\documentclass[12pt,letterpaper,leqno]{amsart}
\usepackage{microtype,amssymb,amsfonts}
\usepackage{geometry}
\usepackage{mathrsfs,tikz,longtable,booktabs,array}
\usetikzlibrary{arrows.meta,calc}
\usepackage{hyperref}
\renewcommand{\Re}{\operatorname{Re}}
\renewcommand{\Im}{\operatorname{Im}}
\newtheorem{theorem}{Theorem}[section]

\newtheorem{proposition}[theorem]{Proposition}
\newtheorem{lemma}[theorem]{Lemma}
\newtheorem{corollary}[theorem]{Corollary}

\numberwithin{equation}{section}
\numberwithin{figure}{section}
\numberwithin{table}{section}
\newcommand{\R}{\mathbb R}
\newcommand{\C}{\mathbb C}
\newcommand{\Sph}{\mathbb S}
\newcommand{\dd}{\,\mathrm d}
\newcommand{\la}{\langle}
\newcommand{\ra}{\rangle}
\newcommand{\eps}{\varepsilon}
\newcommand{\norm}[1]{\left\|#1\right\|}
\newcommand{\ip}[2]{\langle #1,#2\rangle}
\newcommand{\cH}{\mathcal H}
\newcommand{\cR}{\mathcal R}
\newcommand{\cL}{\mathscr L}
\newcommand{\cZ}{\mathcal Z}
\newcommand{\cD}{\mathcal D}
\newcommand{\cA}{\mathcal A}
\newcommand{\cQ}{\mathcal Q}
\newcommand{\cV}{\mathcal V}
\newcommand{\LE}{\operatorname{LE}}
\newcommand{\OO}{\operatorname{OO}}
\newcommand{\TV}{\operatorname{TV}}

\DeclareMathOperator{\conv}{conv}

\DeclareMathOperator{\Ima}{Im}

\newcolumntype{L}[1]{>{\raggedright\arraybackslash}p{#1}}
\space

\space
\hypersetup{pdftitle={Conformal oscillation, a global Kelvin tower, and concentric Type II blow-up for the three-dimensional quintic wave equation},
pdfauthor={Xuwen Chen, Justin Holmer, and Shunlin Shen},
pdfsubject={Nonradial solitary-wave resolution and finite-time concentration of Hopf profiles}}
\title[Conformal oscillation and concentric Hopf scales]{Conformal oscillation, a global Kelvin tower, and concentric Type II blow-up for the three-dimensional quintic wave equation}
\author[X. Chen]{Xuwen Chen}
\address{Department of Mathematics, University of Rochester, Rochester, NY, USA}
\email{xuwenmath@gmail.com}

\author[J. Holmer]{Justin Holmer}
\address{Department of Mathematics, Brown University, Providence, RI, USA}
\email{justin\_holmer@brown.edu}

\author[S. Shen]{Shunlin Shen}
\address{School of Mathematical Sciences, University of Science and Technology of China, Hefei, 230026, China}
\email{slshen@ustc.edu.cn}
\space
\date{}
\begin{document}
\begin{abstract}
We exhibit three forms of nonradial dynamics for the three-dimensional focusing
energy-critical wave equation, arising from the conformal and nodal geometry of
one sign-changing Hopf stationary state. First, a bounded forward-global
solution disproves fixed-profile soliton resolution: after subtraction of a
single free wave, its core repeatedly visits two distinct stationary similarity
classes, and no fixed finite list of traveling profiles in separated frames
resolves it along the full sequence of alternating visits. Second, we construct
a forward-global pure two-core tower with zero asymptotic free radiation, a
common fixed center, one scale tending to one and the other asymptotic to
$ct^{-4}$. Every admissible finite strong traveling-wave decomposition has
precisely these two profiles, up to fixed similarities and permutation. Third,
we construct finite-time concentric two-bubble Type~II blow-up: two nonradial
singular cores share the same fixed spatial center throughout, while both
shrink on strictly separated scales asymptotic to $a(T-t)^4$ and $b(T-t)^{16}$.
Each carries its full stationary energy in the limit, and the solution equals a
fixed regular stationary exterior outside the backward light cone. The two
towers realize dynamics excluded by Shen's radial one-bubble theorem. A regular
nodal point and a quadratic nodal crossing generate dipole and quadrupole
interactions whose powers determine the tower scales. 
\end{abstract}
\maketitle
{\value{tocdepth}=1\relax\tableofcontents}
\section{Introduction}\label{combined:intro}

We consider real solutions of the three-dimensional focusing energy-critical wave equation,
\begin{equation}\label{intro:wave}
 u_{tt}-\Delta u=u^5,\qquad
 \vec u=(u,u_t)\in\mathcal H:=\dot H^1(\mathbb R^3)\times L^2(\mathbb R^3).
\end{equation}
The radial theory determines both the shape and the number of bubbles. Duyckaerts, Kenig and Merle \cite{DKMradial} proved that every radial solution bounded in $\mathcal H$ resolves into a regular part and a sum of signed, rescaled copies of the ground state
$W(x)=(1+\frac{|x|^2}{3})^{-\frac{1}{2}}$. If $J$ denotes the number of bubbles, R.~Shen \cite[Theorem~1.2]{RS} proved the sharper restriction
\[
 J\le1\quad(t\to+\infty),\qquad J=1\quad(t\uparrow T<\infty\text{ in Type~II blow-up}).
\]
At infinite time the regular part is a free wave. At a finite endpoint one can subtract the regular terminal data.

Without radial symmetry, the stationary equation
\begin{equation}\label{main:stationary}
 -\Delta F=F^5,\qquad F\in\dot H^1(\mathbb R^3),
\end{equation}
has the Aubin--Talenti profile $W$ as a sharp Sobolev optimizer \cite{Aubin,Talenti}. Caffarelli, Gidas and Spruck \cite{CGS} classify its positive entire solutions up to translation and dilation. The equation also admits sign-changing solutions \cite{Ding}. We use one such state to separate two questions. Can the stationary shape of a single core be fixed once and for all, modulo the spatial symmetries? Can two stationary profiles persist at the same fixed center on scales whose ratio tends to zero? Our first solution alternates between two inequivalent shapes. Our other two solutions have fixed profiles and a common center. One is a pure two-bubble at infinite time, and the other forms a concentric two-bubble Type~II singularity in finite time.

\subsection{Main results}\label{intro:results}
All three constructions start from $q=h\circ\pi$, where $\pi:\mathbb S^3\to\mathbb S^2$ is the Hopf map and $h$ is positive on the right lune and odd across its sides. The nodal set of $q$ consists of two transversely intersecting tori. A conformal direction at $q$ changes its stereographic profile even modulo translations, positive dilations and orthogonal transformations. We denote this spatial similarity group by $\mathrm{Sim}$. For the two-bubble solutions we use two stereographic representatives of $q$, taken at a regular nodal point and at a crossing. Their leading harmonics have degrees one and two, respectively.

For the first result, $Q_o$ denotes the regular-point representative. A traveling profile is the initial data of a Lorentz transform of a nonzero solution of \eqref{main:stationary}. We use the norm $\|(f,g)\|_{\mathcal H}^2=\|\nabla f\|_2^2+\|g\|_2^2$. Two sequences of frames, with centers $x_{i,n},x_{j,n}$ and scales $s_{i,n},s_{j,n}$, are separated if
\begin{equation}\label{main:profile-separation}
 \frac{s_{i,n}}{s_{j,n}}+\frac{s_{j,n}}{s_{i,n}}+
 \frac{|x_{i,n}-x_{j,n}|^2}{s_{i,n}s_{j,n}}\longrightarrow\infty .
\end{equation}

\begin{theorem}[One core with two limiting shapes]\label{main:counterexample}
There exist a smooth local chart of stationary states $q_\xi$, $\xi\in\mathbb R^4$, with $q_0=Q_o$, and a unit vector $e\in\mathbb R^4$ with the following property. For each sufficiently small $d>0$ there are a bounded forward-global solution $u$ of \eqref{intro:wave}, a free wave $v_+$ and a $C^1$ path $\xi$ such that
\begin{align}
 \|\vec u(0)-(Q_o,0)\|_{\mathcal H}+\|\vec v_+(0)\|_{\mathcal H}
        &\lesssim d^3,\qquad \|\xi\|_{L^\infty}\lesssim d,\label{main:counterexample-bounds}\\
 \|\vec u(t)-\vec v_+(t)-(q_{\xi(t)},0)\|_{\mathcal H}
        &\longrightarrow0\qquad(t\to\infty).\label{main:scattering}
\end{align}
The solution is odd under $(x_1,x_2,x_3)\mapsto(-x_1,-x_2,x_3)$. For $T_j=4^j d^{-5}$,
\begin{equation}\label{main:two-limits}
 \xi(T_{2k})=0,\qquad \xi(T_{2k+1})=de\quad(k\ge1),\qquad
 \operatorname{dist}_{\dot H^1}(q_{de},\mathrm{Sim}\cdot Q_o)\gtrsim d .
\end{equation}
On the entire sequence $(T_j)$, no free wave and no fixed finite list of traveling profiles with fixed velocities $|\ell_j|<1$ approximate $\vec u(T_j)$ in $\mathcal H$, even when the centers, scales and orthogonal transformations of the pairwise separated frames may vary with time.
\end{theorem}

In Conjecture~5.5 of Duyckaerts, Kenig and Merle \cite{DKM}, the list of profiles, their velocities and a free wave are fixed before one chooses a sequence of times. The centers and scales may vary along that sequence. For the solution of Theorem~\ref{main:counterexample}, the even nodes force a single profile in the similarity orbit of $Q_o$, and the odd nodes force one in the orbit of $q_{de}$. By \eqref{main:two-limits} these orbits are distinct. This is an obstruction on the full sequence $(T_j)$. Each of its two subsequences has a resolution, as allowed by the sequential result of Duyckaerts, Jia, Kenig and Merle \cite{DJKM}. In a companion manuscript \cite{CHS3DR6}, we construct a different one-core solution using a balanced single-node isoparametric profile in the family of Fern\'andez and Petean \cite{FPsingle}, with analytic kernel comparisons using identities of Neves and Pistoia \cite{NPsingle}. Its nodal set on $\mathbb S^3$ is one smooth Clifford torus. The two towers in the present paper use the quadratic crossing of the right-lune profile. In the same three-dimensional equation, Donninger and Krieger \cite{DK} constructed radial solutions with a varying scale, and Pillai \cite{PillaiQuintic} constructed radial solutions with free radiation whose ground-state scale can tend to zero and infinity along different subsequences. At finite time, Donninger, Huang, Krieger and Schlag \cite{DHKS} constructed ground-state scales with oscillation between power laws. In all these examples the core belongs to the similarity orbit of $W$. Here the two limits belong to different similarity orbits.

For the two towers we use the energy scaling, its infinitesimal generator and the Kelvin transform,
\begin{equation*}
 F_s(x)=s^{-\frac{1}{2}}F\!\left(\frac{x}{s}\right),\qquad
 \Lambda F=\frac{1}{2}F+x\cdot\nabla F,\qquad
 (\mathcal K F)(x)=|x|^{-1}F\!\left(\frac{x}{|x|^2}\right).
\end{equation*}
The conserved energy $E$ and its density $e$ are
\begin{equation}\label{p8:eq:energy-density}
 E(f,g)=\int_{\mathbb R^3}e(f,g)\,dx,\qquad
 e(f,g)=\frac{1}{2}\big(|\nabla f|^2+|g|^2\big)-\frac{1}{6}|f|^6.
\end{equation}
For the global tower let $Q_{\rm cross}$ be the representative at the quadratic crossing and set $P=\mathcal KQ_{\rm cross}=Q_{\rm cross}(-\cdot)$. If $H_2$ is the quadratic Taylor term of $Q_{\rm cross}$ at zero, the interaction coefficient and the dilation mass are
\begin{equation*}
 C=5\int_{\mathbb S^2}H_2(\omega)^2\,d\omega>0,\qquad
 M=\|\Lambda Q_{\rm cross}\|_2^2=\|\Lambda P\|_2^2 .
\end{equation*}

\begin{theorem}[A forward-global two-core tower]\label{p8:th:tower}
There exist $T_0>0$, a bounded strong solution $u$ of \eqref{intro:wave} on $[T_0,\infty)$, and positive $C^2$ scale functions $\lambda,\mu$ such that, with $S(t)=P_{\lambda(t)}+(Q_{\rm cross})_{\mu(t)}$,
\begin{gather}
 \lim_{t\to\infty}\mu(t)=1,\qquad
 \lim_{t\to\infty}t^4\lambda(t)=64\left(\frac{M}{C}\right)^2,\label{intro:global-tower}\\
 \|\vec u(t)-(S,S_t)(t)\|_{\mathcal H}\lesssim t^{-9},\qquad
 \|\vec u(t)-(S(t),0)\|_{\mathcal H}\lesssim t^{-5}.\label{p8:eq:tower}
\end{gather}
The conserved energy is $2E(Q_{\rm cross},0)$. The solution is nonradial about every point for all sufficiently large times. Along every sequence $t_n\to\infty$, any finite strong traveling-wave decomposition in pairwise separated frames has zero free radiation and exactly two profiles. Up to constant changes of spatial frame and permutation, they are $(P,0)$ and $(Q_{\rm cross},0)$. In the displayed representation both centers are zero.
\end{theorem}

The fixed profiles $P(x)=Q_{\rm cross}(-x)$ and $Q_{\rm cross}(x)$ have zero Lorentz velocities. Their centers agree, and $\frac{\lambda}{\mu}\to0$ separates the two frames in \eqref{main:profile-separation}. Proposition~\ref{p8:pr:refined-scales} records the next terms in the two scales.

For the finite-time tower we keep $P$ and take $B=\mathcal KQ_{\rm wall}$, where $Q_{\rm wall}$ is the regular-wall representative in Section~\ref{common:signed}. With $g=\partial_3 Q_{\rm wall}(0)$ and the orientation in \eqref{typeii:eq:two-forces}, the interaction coefficients and scale constants are
\begin{equation*}
 \begin{aligned}
 C_1&=4\pi g^2>0,&\qquad C_2&=2g\int_{\mathbb R^3}x_1x_2P(x)^5\,dx>0,\\
 a&=\left(\frac{C_1}{8\|\Lambda Q_{\rm wall}\|_2^2}\right)^2,&
 b&=\left(\frac{96\|\Lambda P\|_2^2}{C_2}\right)^2 a^5 .
 \end{aligned}
\end{equation*}

The exterior couples to the dipole of $Q_{\rm wall}$. For the inner bubble, $\int_{\mathbb R^3}x P(x)^5\,dx=0$ by symmetry. The quadratic term of $Q_{\rm wall}$ is therefore the first one that couples to $P$. The two leading interaction energies are $-C_1\mu^{\frac{3}{2}}$ and $-C_2(\frac{\lambda}{\mu})^{\frac{5}{2}}$.

\begin{theorem}[Two concentric singular scales]\label{typeii:seq-main}
For every $T\in\mathbb R$ there exist $\tau_0>0$, a strong solution $u$ of \eqref{intro:wave} on $[T-\tau_0,T)$, and positive scales $\mu,\lambda$ with the following properties. With $\tau=T-t$,
\begin{equation}\label{typeii:eq:finite-asymptotics}
 \begin{gathered}
 \big\|\vec u(t)-(B+(Q_{\rm wall})_{\mu(t)}+P_{\lambda(t)},0)\big\|_{\mathcal H}
     \lesssim \tau^3,\\
 \mu(t)=a\tau^4\big(1+O(\tau^3)\big),\qquad
 \lambda(t)=b\tau^{16}\big(1+O(\tau)\big).
 \end{gathered}
\end{equation}
The solution satisfies $\vec u(t)=(B,0)$ on $|x|>\tau$ and has Type~II blow-up at the single point $(T,0)$. For $z=u-B$, $r_1=\sqrt{\lambda\mu}$ and $r_2=\sqrt{\mu\tau}$, its singular energy satisfies
\begin{equation}\label{intro:localized-energies}
 \lim_{t\uparrow T}\int_{|x|<r_1} e(z,z_t)\,dx=E(P,0)>0,\qquad
 \lim_{t\uparrow T}\int_{r_1<|x|<r_2}e(z,z_t)\,dx=E(Q_{\rm wall},0)>0 .
\end{equation}
\end{theorem}

The two frames are centered at the origin throughout the construction. Their radii satisfy
$\lambda\ll r_1\ll\mu\ll r_2\ll T-t$, so \eqref{intro:localized-energies} separates the energy of each bubble from that of the other. The stationary exterior $B$ is regular at $T$. In the free-wave formulation of \cite[Conjecture~5.5]{DKM}, one can take the free wave with data $(B,0)$ at $T$. By energy continuity, \eqref{typeii:eq:finite-asymptotics} then exhibits the same two fixed profiles along every sequence $t_n\uparrow T$. Theorems~\ref{p8:th:tower} and~\ref{typeii:seq-main} thus realize $J=2$ in the two time regimes in which Shen excludes more than one radial bubble.

\subsection{Relation to previous work}
Kadar \cite[Theorem~1.1]{FTKadar} constructed solutions of \eqref{intro:wave} with several Lorentz-transformed ground states approaching one blow-up point. With $\tau=T-t$, their centers follow distinct timelike rays $x=z_j\tau$ and their radii are comparable to $\tau^\nu$ for a common $\nu>1$. Thus their center separation is of order $\tau$, much larger than either radius. In Theorem~\ref{typeii:seq-main} the separation takes place entirely in scale at one fixed center. For the cubic wave equation in four space dimensions, Collot has announced ongoing joint work with Li on a singularity formed by the collision of collapsing solitons \cite{CollotAnnouncement}. Martel and Rapha\"el \cite{MartelRaphael} previously constructed finite-time concentration of several bubbles at one point for the two-dimensional mass-critical nonlinear Schr\"odinger equation. For the energy-critical wave maps equation from $\mathbb R^{1+2}$ to $\mathbb S^2$, Jendrej and Krieger \cite{FTJendrejKrieger} constructed finite-time concentric two-bubble solutions in the $2$-equivariant class. Krieger and Palacios \cite{FTKriegerPalacios} constructed such bubble trees of arbitrary finite length in the same class. Jendrej, Yin and Zhao \cite{JYZ} showed that finite-time blow-up for $k$-equivariant wave maps to $\mathbb S^2$, $k\ge2$, requires nonsmooth radiation. The exterior $B$ in Theorem~\ref{typeii:seq-main} is smooth for the scalar wave equation.

For the six-dimensional scalar wave equation, Jendrej \cite{J} constructed a pure two-bubble with an exponentially concentrating inner scale. The same paper proves algebraic concentration for equivariant wave maps and predicts $|t|^{-4/(N-6)}$ for scalar ground-state towers in dimension $N>6$ \cite[Remark~1.4]{J}. In dimension seven that prediction gives the exponent $-4$. Here the power $5/2$ in the interaction energy comes from a quadratic nodal jet in dimension three. Jendrej and Martel \cite{JendrejMartel} constructed five-dimensional multi-bubbles of radii $c_jt^{-2}$ at distinct fixed limiting points. Jendrej, Zhang and Zhao \cite{JZZ} proved rigidity of this rate for spatially separated bubbles with comparable scales. For the energy-critical heat equation in dimensions $N\ge7$, del Pino, Musso and Wei \cite{HeatTower} constructed towers of arbitrary finite length concentrating at one point in infinite time.

For the three-dimensional radial equation, see Kenig and Merle \cite{KM} for scattering below the ground-state threshold, and Krieger, Schlag and Tataru \cite{KST} and Krieger and Schlag \cite{KSFull} for finite-time one-bubble solutions. Jendrej \cite{JOpposite} excluded pure radial two-bubbles of opposite signs in dimensions $N\ge3$, an earlier restriction on the configurations covered in dimension three by Shen \cite{RS}. Shen also proved that the regular terminal data of a three-dimensional radial Type~II blow-up solution are nonzero \cite{ShenPure}. For radial resolution in higher dimensions, see \cite{DKModd,JL}. In the three-dimensional nonradial equation, Kadar \cite{KadarGlobal} constructed global multi-solitons with separating centers and no outgoing radiation. For bounded forward-global solutions with compactly supported data in dimensions three through five, Collot, Duyckaerts, Kenig and Merle \cite[Theorems~1.2--1.4]{CDKM} obtain continuous-time resolution when the sequential decomposition has one traveling ground state or two with distinct velocities. In the two-soliton conclusion the velocity parameters may vary with time. For larger numbers of ground states with distinct velocities they obtain a dichotomy involving collisions.

Ding \cite{Ding} proved the existence of sign-changing stationary solutions of \eqref{main:stationary}. Del Pino, Musso, Pacard and Pistoia \cite{DPMPP} constructed large-energy families, and Musso and Wei \cite{MW} studied conformal nondegeneracy for their polygonal family. Hopf reductions for elliptic problems appear in \cite{CP,FPP}. Duyckaerts, Kenig and Merle \cite{DKMcompact} use conformal nondegeneracy to classify solutions with the compactness property. In Theorem~\ref{main:counterexample} the radiation acts on a conformal direction transverse to the similarity orbit. In five dimensions this mechanism produces a continuum of limiting similarity classes \cite{CHS5D}. For harmonic map heat flow, see \cite{Top97,Top04} for winding, \cite{GNT} for oscillation of the scale, and \cite{JLS} for continuous-time bubble decomposition. Jendrej, Lawrie and Schlag explicitly distinguish convergence to the set of bubble configurations from convergence with profiles fixed across time sequences. Adapting Topping's winding mechanism, Engelstein and Mendelson \cite{EM} constructed wave maps with a continuum of different weak bubble limits along time sequences. Jendrej \cite{J5} obtained a quartic concentration radius for a radial five-dimensional wave with prescribed regular part positive at the origin. In Theorem~\ref{typeii:seq-main} the analogous quartic balance couples the first derivatives of the wall profile and its stationary exterior. For prescribed radiation determining wave-map concentration, see Jendrej, Lawrie and Rodriguez \cite{JLR}.

\subsection{Outline of the proofs}\label{intro:outline}
\begin{figure}[htbp]
\centering
\begin{tikzpicture}[x=1cm,y=1cm,font=\small,
 box/.style={draw,align=center,inner xsep=4pt,inner ysep=4pt,text width=4.25cm},
 arr/.style={-{Stealth[length=1.7mm]}}]
 \node[box,text width=7.0cm] (q) at (0,0) {Stationary Hopf state and signed kernels\\
 Section~\ref{sec:common-hopf} and Appendix~\ref{cert:start}};
 \node[box] (wall) at (-4.9,-1.65) {Radiation impulses at $T_j$\\ along a conformal direction};
 \node[box,text width=6.7cm] (mom) at (2.45,-1.65)
 {Harmonic moments at the two poles\\
 $p=k+\frac{1}{2}$ for a harmonic of degree $k$};
 \node[box] (one) at (-4.9,-3.6) {Two limiting shapes\\ Theorem~\ref{main:counterexample}};
 \node[box] (two) at (0,-3.6) {Quadrupole interaction\\
 $\lambda\sim c t^{-4}$\\ Theorem~\ref{p8:th:tower}};
 \node[box] (three) at (4.9,-3.6) {Linear and quadratic jets\\
 $\mu\sim a\tau^4$, $\lambda\sim b\tau^{16}$\\ Theorem~\ref{typeii:seq-main}};
 \draw[arr] (q.south west) -- (wall.north);
 \draw[arr] (q.south east) -- (mom.north);
 \draw[arr] (wall.south) -- (one.north);
 \draw[arr] (mom.south west) -- (two.north);
 \draw[arr] (mom.south east) -- (three.north);
\end{tikzpicture}
\caption{The stationary directions and harmonic moments used in the three constructions.}
\label{intro:proof-map}
\end{figure}
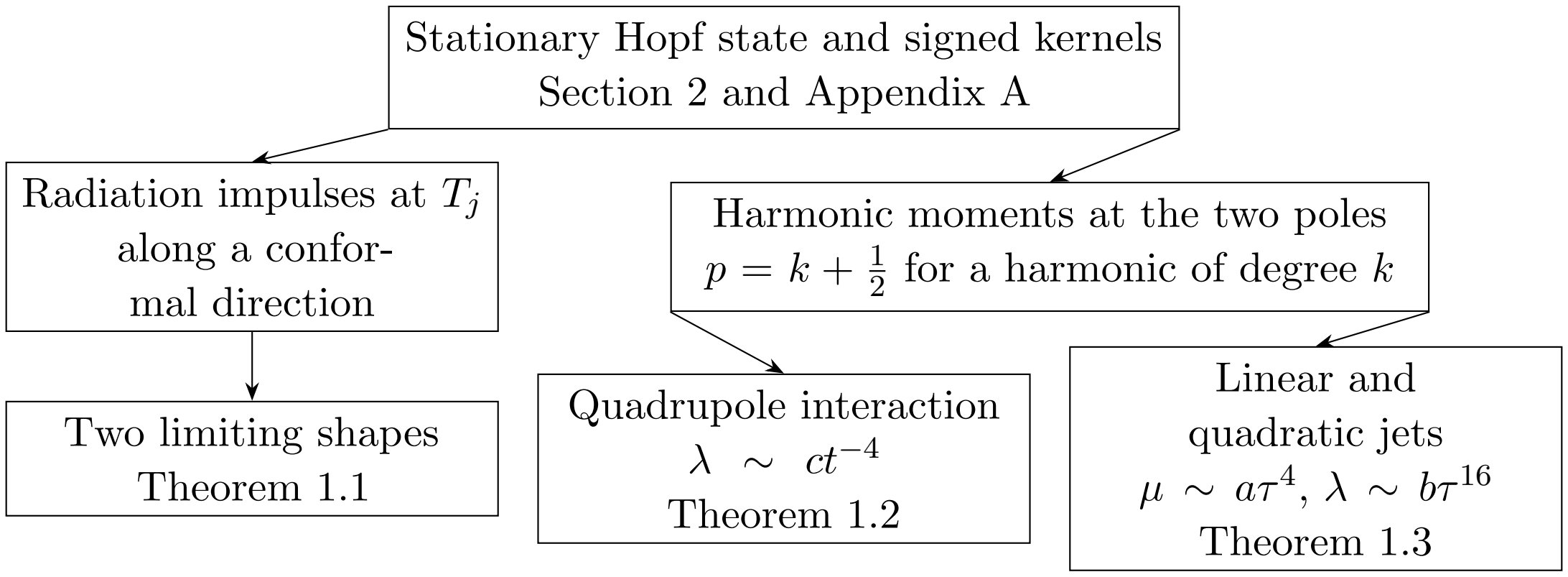

Section~\ref{sec:common-hopf} identifies four odd conformal directions at the regular wall and only the centered dilation in each tower symmetry class. Appendix~\ref{cert:start} constructs the Hopf state by minimization, then uses the charge decomposition to exclude all other zero modes. Figure~\ref{intro:proof-map} indicates how these elliptic conclusions enter the three proofs.

\paragraph{The oscillating core.}
In Sections~\ref{geo:shape-section}--\ref{end:start} we choose a shape direction $e$ transverse to $\mathrm{Sim}\cdot Q_o$ and drive it by a sequence of radiation packets. In the coordinate $\Gamma$ of \eqref{end:ep:coordinate} the desired displacement on $I_j$ is
\[
 \int_{I_j}\Gamma'(t)\,dt=(-1)^j\gamma(de).
\]
The difference of the two adjacent mean velocities prescribes the impulse at $T_j$. Testing against a piecewise affine function supported on the adjacent intervals isolates this impulse. Sections~\ref{cat:sec:reference} and~\ref{force:sec:comparison} estimate the impulse errors. The fixed point in \eqref{end:ep:map} matches the velocities across all nodes. The uniqueness of the free radiation and the comparison of spatial frames in Section~\ref{fp:energy-shape} then exclude a fixed list of traveling profiles.

\paragraph{The two towers.}
If $H_k$ is a homogeneous harmonic polynomial of degree $k$ and $H_kF^5\in L^1$, we change variables $x=sy$ to obtain
\[
 \int_{\mathbb R^3}H_k(x)F_s(x)^5\,dx
   =s^{k+\frac{1}{2}}\int_{\mathbb R^3}H_k(y)F(y)^5\,dy .
\]
The first nonvanishing moment determines the power of the scale ratio in the interaction. For the global Kelvin pair both profiles have a quadratic crossing. In Sections~\ref{sec:tower-1}--\ref{sec:tower-3} we obtain the leading potential $-C(\frac{\lambda}{\mu})^{\frac{5}{2}}$. Since $\mu\to1$, the balance $\frac{M}{2}\dot\lambda^2\sim C\lambda^{\frac{5}{2}}$ on the contracting branch determines the rate $t^{-4}$. Section~\ref{merge:global-wave} controls the wave orthogonal to the dilation modes and chooses the two scales so that both remaining forces vanish.

At a finite endpoint, suppose the larger scale is proportional to $\tau^\alpha$, the smaller one to $\tau^\kappa$, and the potential to their ratio raised to the power $p$. The scale equation \eqref{p1:eq:power-balance} requires
\[
 (2-p)\kappa=2-p\alpha,\qquad
 (p,\alpha)=\left(\frac{3}{2},0\right)\ \Longrightarrow\ \kappa=4,\qquad
 (p,\alpha)=\left(\frac{5}{2},4\right)\ \Longrightarrow\ \kappa=16 .
\]
The first balance is between the exterior and the middle dipole. The second uses the quadrupole of the inner bubble. Section~\ref{typeii:first} constructs the two profile corrections. In Section~\ref{typeii:sec:matchedcone} we match the middle dipole at the light cone. The slope of the corrected middle wave at the origin enters the inner equation through \eqref{typeii:eq:average-slope}. Sections~\ref{typeii:sec:mixedinverse} and~\ref{typeii:merge:closing} select the waves and scales simultaneously and pass to the singular time. The energy limits follow from \eqref{typeii:eq:singular-energy-regions}.
\space
\section{The stationary Hopf state and its signed charts}\label{sec:common-hopf}

On the unit sphere $\mathbb S^2$, the great circles $X=0$ and $Y=0$
bound the right lune $D=\{X>0,\,Y>0\}$. We compute the vanishing order
on its sides and at its vertices, then determine which conformal
directions survive the symmetries of the wave equation.

\subsection{The positive stationary state on the right lune}\label{common:sec:lune}
Let $h$ be the odd reflection of the positive variational solution constructed in Appendix~\ref{ap:profile}.
Its restriction to the lune satisfies
\begin{equation}\label{common:lune-pde}
 (-4\Delta_{\mathbb S^2}+\frac{3}{4})h=h^5\quad\text{in }D,\qquad
 h|_{\partial D}=0.
\end{equation}
The reflection and the symmetries of the minimizing class imply
\begin{equation}\label{common:lune-parities}
 h(-X,Y,Z)=h(X,-Y,Z)=-h(X,Y,Z),\qquad
 h(Y,X,Z)=h(X,Y,-Z)=h(X,Y,Z).
\end{equation}
The strong maximum principle yields $h>0$ in $D$. Hence the analytic extension satisfies
\begin{equation}\label{common:sign}
                 \operatorname{sgn}h=\operatorname{sgn}(XY).
\end{equation}
The peak satisfies $\|h\|_\infty^4<240$ by \eqref{ap:peak}.

\subsection{Regular walls and quadratic crossings}
\label{common:nodal}
We first locate the two types of zeros on the Hopf lift. For $z=(z_1,z_2)\in\mathbb S^3$, we set
\begin{equation}\label{common:hopf-map}
 \pi(z)=(2\Re(z_1\bar z_2),2\Im(z_1\bar z_2),|z_1|^2-|z_2|^2),
 \qquad q=h\circ\pi.
\end{equation}
By odd reflection and the Hopf metric identity, the base and its lift satisfy
\begin{align}
 \left(-4\Delta_{\mathbb S^2}+\frac{3}{4}\right)h&=h^5\quad\text{on }\mathbb S^2,\label{geo:eq:base}\\
 \left(-\Delta_{\mathbb S^3}+\frac{3}{4}\right)q&=q^5,\qquad
 \Delta_{\mathbb S^3}q=4(\Delta_{\mathbb S^2}h)\circ\pi .\label{geo:eq:spherical}
\end{align}
By \eqref{common:sign}, the two nodal surfaces are $\pi^{-1}\{X=0\}$ and
$\pi^{-1}\{Y=0\}$. Each is a Clifford torus, since rotations of the base
lift to isometries and the equatorial preimage is
$\{|z_1|=|z_2|=2^{-\frac{1}{2}}\}$. As $\pi$ is a submersion, the two tori
meet transversely on $\pi^{-1}(0,0,\pm1)$. At every other nodal point the
boundary point lemma on $D$ implies $\nabla q\ne0$.
At the north pole $N$, the coordinates $X,Y$ have Euclidean metric to
leading order. Let $H_m$ be the first nonzero homogeneous term in the
Taylor series of $h$. Equation~\eqref{common:lune-pde} implies
$\Delta_{\mathbb R^2}H_m=0$. In polar coordinates the two reflections in
\eqref{common:lune-parities} require
\begin{equation*}
 H_m(r,\phi)=c r^{2k}\sin(2k\phi),\qquad m=2k,\quad c\ne0 .
\end{equation*}
If $k\ge2$, $H_m$ takes both signs on $0<\phi<\frac{\pi}{2}$,
contradicting \eqref{common:sign} along a ray. Hence $k=1$ and
\begin{equation}\label{common:jet}
 h(X,Y)=\gamma XY+O((X^2+Y^2)^2),\qquad \gamma>0 .
\end{equation}
Oddness in both variables eliminates all cubic terms. The south pole
has the same expansion by the reflection in $Z$.

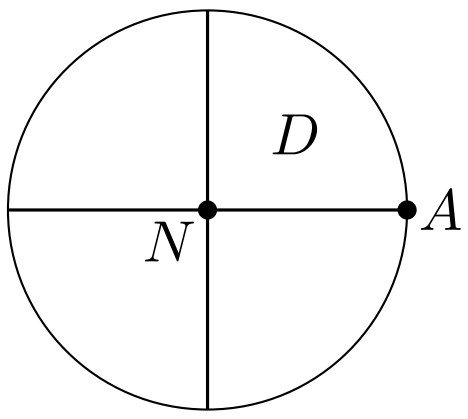
Figure~\ref{combined:fig:marks} records the two choices of pole.
\begin{figure}[htbp]
\centering
\begin{tikzpicture}[x=1cm,y=1cm,font=\small]
 \def\rad{1.32}
 \draw (0,0) circle[radius=\rad];
 \draw[semithick] (-\rad,0)--(\rad,0) (0,-\rad)--(0,\rad);
 \node at (.58,.50) {$D$};
 \fill (0,0) circle[radius=1.8pt];
 \node[below left,inner sep=2pt] at (0,0) {$N$};
 \fill (\rad,0) circle[radius=1.8pt];
 \node[anchor=west,inner sep=2pt] at (\rad,0) {$A$};
 \node[anchor=west] at (2.18,.85) {\textbf{Regular wall} $A=(1,0,0)$};
 \node[anchor=west] at (2.18,.39)
 {$\pi^{-1}(A)=\{2^{-\frac{1}{2}}e^{i\theta}(1,1)\}$};
 \node[anchor=west] at (2.18,-.29) {\textbf{Quadratic crossing} $N=(0,0,1)$};
 \node[anchor=west] at (2.18,-.75)
 {$\pi^{-1}(N)=\{(e^{i\theta},0)\}$};
\end{tikzpicture}
\caption{The northern hemisphere of the Hopf base, projected onto the $(X,Y)$ disk. The two diameters are the nodal walls and $D$ lies in the first quadrant. The circle boundary is the equator. The chosen poles are the points with $\theta=0$ in the indicated fibers.}
\label{combined:fig:marks}
\end{figure}

The oscillation construction uses $p_o=2^{-\frac{1}{2}}(1,1)$ as its
stereographic pole. Its image $(1,0,0)$ lies on an open wall.
The tower uses $p_t=(1,0)$, for which both $p_t$ and $-p_t$ project
to $N$. To fix the normalization, its inverse stereographic map and
conformal factor are
\begin{equation*}
 \sigma_t(x)=
 \left(\frac{|x|^2-1+2i x_3}{1+|x|^2},
       \frac{2(x_1+i x_2)}{1+|x|^2}\right),
 \qquad \Omega=\frac{2}{1+|x|^2}.
\end{equation*}
By the formula for $\sigma_t$ and antipodal invariance of $q$, we have
\begin{equation*}
 Q=\Omega^{\frac{1}{2}}q\circ\sigma_t,\qquad
 \sigma_t\!\left(\frac{x}{|x|^2}\right)=-\sigma_t(-x),\qquad
 P(x)=|x|^{-1}Q\!\left(\frac{x}{|x|^2}\right)=Q(-x).
\end{equation*}
By conformal covariance, both profiles solve $-\Delta F=F^5$. In the crossing chart, with $r=|x|$, the Hopf map reads
\begin{equation*}
 X+iY=\frac{4(r^2-1+2ix_3)(x_1-i x_2)}{(1+r^2)^2}.
\end{equation*}
We substitute in \eqref{common:jet} and use $\Omega^{\frac{1}{2}}=\sqrt{2}+O(r^2)$. With $b=16\sqrt{2}\,\gamma$, we obtain
\begin{align}
 Q(x)&=H_2(x)+H_3(x)+O(r^4),\qquad H_2(x)=-b x_1x_2,\label{common:quad}\\
 H_3(x)&=2b x_3(x_2^2-x_1^2).\label{p2:eq:jet}
\end{align}
Duyckaerts, Kenig and Merle \cite[Remark~3.2, (3.4)]{DKMcompact} describe the leading harmonic tail of a stationary solution through its Kelvin transform. Yuan \cite[Lemma~2.2]{YuanFour} used inversion at a regular nodal point to obtain a dipolar tail for a four-dimensional excited multi-soliton. See also his five-dimensional construction \cite{YuanFive}. Here the crossing yields a quadrupole. Both $P$ and $Q$ have leading term $H_2$ at zero and the leading
term $|x|^{-3}H_2\bigl(\frac{x}{|x|}\bigr)$ at infinity. In particular the tower
interaction coefficient $5\int_{\mathbb S^2}H_2(\omega)^2\,d\omega$
is strictly positive. The change of pole giving these two Euclidean
representatives is a conformal change on stationary data. The time
dependence asserted later comes from solving the wave equation.

\subsection{The conformal stabilizer and the three signed kernels}
\label{common:signed}
We first identify the kernel of the infinitesimal conformal action.
Averaging the remaining generators over each symmetry group then
determines the signed conformal directions.
\begin{proposition}[Connected stabilizer]\label{geo:prop:stabilizer}
The connected conformal stabilizer of $q$ is the Hopf circle
$z\mapsto e^{i\alpha}z$. Its conformal orbit has dimension nine.
\end{proposition}
\begin{proof}
A conformal stabilizer preserves the singular nodal set, since the
action on a stationary state multiplies its pullback by a smooth
positive weight. Its identity component preserves each singular circle
separately. These circles are $\mathbb S^3\cap V_1$ and
$\mathbb S^3\cap V_2$ for the orthogonal complex coordinate planes
$V_1,V_2\subset\mathbb R^4$. Every conformal vector field on the round sphere has the form
$X(\omega)=A\omega+b-(b\cdot\omega)\omega$, where $A^*=-A$ and $b\in\mathbb R^4$.
Tangency at $\omega$ and $-\omega$ on the $j$th circle implies
\begin{equation*}
 0=\Pi_{V_j^\perp}(\pm A\omega+b)
 \quad\Longrightarrow\quad b\in V_j,\qquad A V_j\subset V_j .
\end{equation*}
Since $V_1\cap V_2=0$, only the two plane rotations remain. We write their action
and the induced action on the base as
\begin{equation*}
 (z_1,z_2)\longmapsto(e^{i\alpha}z_1,e^{i\beta}z_2),\qquad
 X+iY\longmapsto e^{i(\alpha-\beta)}(X+iY).
\end{equation*}
By \eqref{common:sign}, preservation of the two meridians requires
$\alpha-\beta\in\frac{\pi}{2}\mathbb Z$. Continuity from the identity
implies $e^{i\alpha}=e^{i\beta}$. Conversely, every common phase fixes
$\pi$ and hence $q$. Thus the conformal orbit has dimension $10-1=9$.
\end{proof}

Let $S(z_1,z_2)=(\bar z_1,\bar z_2)$ and $R(z_1,z_2)=(z_1,iz_2)$.
Both reverse the sign of $q$. We set $L_q=-\Delta_{\mathbb S^3}+\frac{3}{4}-5q^4$
and denote by $J_q$ the infinitesimal conformal action on $q$.
By Proposition~\ref{geo:prop:stabilizer}, the image of $J_q$ has
dimension nine, and Theorem~\ref{thm:stationary-certificate} identifies
it with $\ker L_q$. For a finite group $G$ with $U_g q=\chi(g)q$,
where $U_g$ is pullback by $g^{-1}$, we have
\begin{equation}\label{common:equivariant}
 U_gJ_q(A)=\chi(g)J_q(\operatorname{Ad}_g A).
\end{equation}
If $J_q(A)$ has character $\chi$, averaging \eqref{common:equivariant} shows that $|G|^{-1}\sum_g\operatorname{Ad}_g A$ has the same image. The conformal algebra is $\mathfrak{so}(4)\oplus\mathbb R^4$, and the kernel of $J_q$ is spanned by the Hopf generator. Conjugation $S$ reverses this generator. Thus for every group used below,
\begin{equation}\label{common:fixed-count}
 \dim(\ker L_q)_\chi=\dim\mathfrak{so}(4)^G+\dim(\mathbb R^4)^G .
\end{equation}
For commuting involutions with joint eigenspace dimensions $m_\eta$, a fixed skew map preserves each eigenspace. The two terms on the right are $\sum_\eta\binom{m_\eta}{2}$ and $m_{\bf1}$. For $G=\langle S\rangle$ the real and imaginary planes both have dimension two, so the dimension in \eqref{common:fixed-count} is $1+1+2=4$. The two tower groups also contain $S$ with $\chi(S)=-1$.
The same count applies in the Euclidean energy space at either pole.
The conformal energy identity is
\begin{equation*}
 \int_{\mathbb R^3}|\nabla(\Omega^{\frac{1}{2}}k\circ\sigma)|^2
   =\int_{\mathbb S^3}(|\nabla k|^2+\frac{3}{4}|k|^2).
\end{equation*}
It identifies $H^1(\mathbb S^3)$ with $\dot H^1(\mathbb R^3)$ and
carries the spherical weak kernel equation to its Euclidean counterpart,
since a point has zero $H^1$ capacity in dimension three. All spherical kernel fields are
smooth by elliptic regularity.

We write $\Lambda=\frac{1}{2}+x\cdot\nabla$. For the finite-time tower
we use the same two Euclidean half-turns in both charts,
\begin{equation*}
 S_x=\operatorname{diag}(1,-1,-1),\qquad
 D_x=\operatorname{diag}(-1,1,-1),\qquad \chi(S_x)=\chi(D_x)=-1 .
\end{equation*}
At the regular pole $p=p_o$ put $v=\frac{(1,-1)}{\sqrt{2}}$. The coordinates for
this construction are
\begin{equation*}
 \sigma_{\rm wall}(x)=
 \frac{(|x|^2-1)p+2(x_1v+x_2ip+x_3iv)}{1+|x|^2},\qquad
 Q_{\rm wall}=\Omega^{\frac{1}{2}}q\circ\sigma_{\rm wall}.
\end{equation*}
At $A=(1,0,0)$ the parities \eqref{common:lune-parities} imply
$h=h_Y(A)Y+O((Y^2+Z^2)^{\frac{3}{2}})$. In these wall coordinates,
\begin{equation}\label{common:wall-jet}
 Y=\frac{4(r^2-1)x_3-8x_1x_2}{(1+r^2)^2},\qquad
 Q_{\rm wall}(x)=g(x_3+2x_1x_2)+O(r^3),
\end{equation}
where $g=-4\sqrt{2}\,h_Y(A)<0$ by the boundary point lemma.
In the ordered real basis $(p,v,ip,iv)$ the two half-turns lift to
$S$ and $(z_1,z_2)\mapsto(z_2,z_1)$. The interchange sends
$(X,Y,Z)$ to $(X,-Y,-Z)$ and thus reverses $h$. At the crossing
pole, the corresponding basis is $((1,0),(0,1),(0,i),(i,0))$
from $\sigma_t$, and the lifts are $S$ and $R^2S$.
In either basis, the four joint characters on the pole and tangent axes
are $(+,+),(+,-),(-,+),(-,-)$. Each joint eigenspace has dimension one. Thus \eqref{common:fixed-count} leaves only the fixed gradient along $p$, and the three tangent axes have no common fixed vector.
We put
\begin{equation*}
 Q_{\rm cross}=\Omega^{\frac{1}{2}}q\circ\sigma_t,\qquad
 (\mathcal KF)(x)=|x|^{-1}F\!\left(\frac{x}{|x|^2}\right).
\end{equation*}
The fixed gradient generates centered dilation. By \eqref{common:equivariant},
we obtain
\begin{equation}\label{common:finite-kernel}
 \ker_{\dot H^1_\chi}(-\Delta-5F^4)=\mathbb R\Lambda F,\qquad
 F=Q_{\rm wall},\ \mathcal KQ_{\rm wall},\ Q_{\rm cross},\
 \mathcal KQ_{\rm cross}.
\end{equation}
Indeed $\mathcal K$ commutes with both half-turns and conjugates the
dilation action to its inverse. The regular wall and its antipodal
point have nonzero normal derivative, so the corresponding profiles
and their dilations have dipole decay $O(|x|^{-2})$ at infinity.

The global tower imposes the quarter-turn $R$ in addition to $S$.
Its signed sector is contained in the crossing half-turn sector just
computed, since $\chi(R^2S)=-1$. Conversely the centered dilation
inherits both signs from its profile. With $F=P,Q$ in the crossing
chart, we obtain
\begin{equation}\label{common:restricted}
 \dim\ker(L_q|_{S\text{-odd}})=4,\qquad
 \ker_{\dot H^1_\chi}(-\Delta-5F^4)=\mathbb R\Lambda F .
\end{equation}
These dilations are nonzero by \eqref{common:quad} and lie in $L^2$
by their $O(|x|^{-3})$ decay.
\space
\subsection{The oscillating chart and shape separation}\label{geo:start}
We construct an odd stationary chart $q_\xi$ through $Q$ and choose a
direction $e$ transverse to the spatial similarities \eqref{main:eq:similarity}.
We use the regular nodal pole of the Hopf lift from Section~\ref{sec:common-hopf}.
Appendix~\ref{cert:start} identifies the full odd energy kernel needed
for the linear estimates. For the conformal action we refer to
\cite{DKMcompact}, and for the Hopf reduction to \cite{CP,FPP}.

\subsubsection{The regular-wall chart for the oscillation}\label{geo:hopf}
We take the pole $p=p_o$, the involution $\tau z=\bar z$ and the lift $q$ from \eqref{geo:eq:spherical}. Since $p\in\operatorname{Fix}\tau$ and $q\circ\tau=-q$, we have $q(p)=0$.
In the coordinates of Section~\ref{common:signed} conjugation fixes the first tangent axis.
For the oscillation we permute axes by $C(x_1,x_2,x_3)=(x_3,x_1,x_2)$ and set
$\sigma=\sigma_{\rm wall}\circ C$, $Q_o=Q_{\rm wall}\circ C$. Then
\begin{equation}\label{geo:eq:euclidean}
 Q(x)=\Omega(x)^{\frac{1}{2}}q(\sigma(x)),\qquad -\Delta Q=Q^5.
\end{equation}
In this order the Euclidean action of $\tau$ is
$\mathscr R(x_1,x_2,x_3)=(-x_1,-x_2,x_3)$. Thus $Q\circ\mathscr R=-Q$.
Since $q(p)=0$ and $\sigma(x)-p=O(|x|^{-1})$, Taylor's formula implies
\begin{equation}\label{geo:eq:decay}
 |\partial_x^\beta Q(x)|\lesssim_\beta\langle x\rangle^{-2-|\beta|}.
\end{equation}
In particular $Q\in H^1(\mathbb R^3)$ is real, nonzero and sign-changing.

\subsubsection{The odd stationary chart}\label{geo:chart-section}
We now choose chart coordinates that are linear under projection onto
the tangent space. Let $A:z\mapsto iz$ be the generator of the connected stabilizer in
Proposition~\ref{geo:prop:stabilizer}. The weighted conformal action on
stationary fields is an isometry of $\dot H^1$ and preserves the
Euclidean stationary equation.

\begin{proposition}[Odd stationary chart]\label{geo:prop:chart}
There is a smooth four-dimensional stationary chart
\begin{equation}\label{geo:eq:chart}
 q_\xi=Q+K\xi+R(\xi),\qquad
 K^*K=I_4,\quad K^*R(\xi)=0,\quad R(0)=DR(0)=0.
\end{equation}
Here $K^*f=(\langle K_i,f\rangle)_{i=1}^4$ on $L^2$ and
all chart elements are odd under $\mathscr R$. Every fixed finite order
of parameter differentiation obeys \eqref{geo:eq:decay}, uniformly on a
smaller chart. In particular, with $D=\sqrt{-\Delta}$ and
$Y_s=L^2(\langle x\rangle^{2s}\dd x)$,
\begin{equation}\label{geo:eq:density}
 D^2\partial_\xi^a q_\xi=O(\langle x\rangle^{-10})\in Y_4.
\end{equation}
\end{proposition}
\begin{proof}
We exponentiate the four $\tau$-fixed generators in Section~\ref{common:signed}.
Their derivatives are independent since $A$ is $\tau$-anti-invariant.
Each transformation commutes with $\tau$, so every field on this
four-dimensional local orbit is odd and vanishes at $p$. Taylor's
formula at $p$ proves \eqref{geo:eq:decay} for every fixed parameter derivative.

We choose local orbit coordinates $q(\zeta)$ with $q(0)=Q$ and
$L^2$-orthonormal tangent vectors $Dq(0)=K$. Then
\begin{equation*}
 \left.D_\zeta\bigl[K^*(q(\zeta)-Q)\bigr]\right|_{\zeta=0}=K^*K=I_4.
\end{equation*}
The inverse function theorem permits the choice $\xi=K^*(q-Q)$,
which is precisely \eqref{geo:eq:chart}. Differentiating the stationary
equation yields $D^2\partial_\xi^a q_\xi=\partial_\xi^a(q_\xi^5)$.
Each product on the right contains five factors bounded by
$O(\langle x\rangle^{-2})$, proving \eqref{geo:eq:density}. In particular,
$(-\Delta-5Q^4)K_i=0$ for $1\le i\le4$.
\end{proof}

\subsubsection{A shape direction modulo all spatial similarities}\label{geo:shape-section}
We choose a chart direction transverse to the similarity orbit defined
in Section~\ref{main:start}. The comparison uses the full spatial group,
including similarities that do not preserve oddness.

\begin{proposition}[Quantitative shape separation]\label{geo:prop:shape}
For some unit $e\in\mathbb R^4$ and $c,d_0>0$,
\begin{equation}\label{geo:eq:shape-separation}
 \|q_{de}-Q\|_{\dot H^1}\lesssim |d|,
 \qquad
 \operatorname{dist}_{\dot H^1}(q_{de},\mathrm{Sim}\cdot Q)
       \ge c|d|\quad(0<|d|<d_0).
\end{equation}
\end{proposition}
\begin{proof}
Let $V$ denote the tangent space of $\mathrm{Sim}\cdot Q$ at $Q$. Similarities fix the pole $p$. By Proposition~\ref{geo:prop:stabilizer}, a $\tau$-fixed conformal generator $X$ induces a tangent in $V$ precisely when $X(p)\in\mathbb R ip$. But $X(p)$ belongs to the real line $T_p\mathbb S^3\cap\mathbb R^2$, and the infinitesimal rotation $B(z_1,z_2)=(-z_2,z_1)$ spans this line at $p$. Thus evaluation at $p$ has rank one on the four chart generators, and
\begin{equation*}
 \dim(V\cap\operatorname{ran}K)=3 .
\end{equation*}
The odd projection preserves $V$ by conjugation of similarities under $\mathscr R$. Its image on $V$ is $V\cap\operatorname{ran}K$ by \eqref{common:restricted}. We choose a unit $e$ for which $Ke$ is energy-orthogonal to these three tangents, hence to $V$.

By \eqref{fp:dislocation}, any sequence $S_n Q\to Q$ has a convergent subsequence of parameters. Modulo its stabilizer the orbit is therefore a local $C^2$ graph $Q+v+O_{\dot H^1}(\|v\|^2)$ over $V$. Only $\|v\|_{\dot H^1}\lesssim |d|$ can improve on the competitor $Q$. Projecting $q_{de}=Q+dKe+O_{\dot H^1}(d^2)$ onto $V^\perp$ yields the matching lower bound,
\begin{equation*}
 \operatorname{dist}_{\dot H^1}(q_{de},\mathrm{Sim}\cdot Q)
 = |d|\,\|Ke\|_{\dot H^1}+O(d^2).
\end{equation*}
\end{proof}

\subsubsection{The stationary state and its kernel}\label{geo:validation-section}\label{geo:sec:kernel}

\begin{proposition}[The stationary state and its kernel]\label{geo:ass:stationary}
There is a smooth nonzero real solution $h$ of \eqref{common:lune-pde} with the parities \eqref{common:lune-parities},
specified by the minimization in \eqref{ap:min} and the odd reflection across the sides of the lune.
We construct $Q$, its normalized odd chart and a unit shape vector $e$ by
\eqref{geo:eq:euclidean}, Propositions~\ref{geo:prop:chart} and
\ref{geo:prop:shape}. Then
\begin{equation}\label{geo:eq:odd-kernel}
 \ker_{\dot H^1,\,\mathrm{odd}}(-\Delta-5Q^4)
       =\operatorname{span}\{K_1,K_2,K_3,K_4\}.
\end{equation}
The full spherical linearization has kernel dimension nine, consisting
exactly of the infinitesimal conformal orbit.
\end{proposition}
\begin{proof}
The variational construction and the spherical kernel identity are proved in Theorem~\ref{cert:thm:executed}. The four independent tangents of Proposition~\ref{geo:prop:chart} belong to $H^2$ by \eqref{geo:eq:decay}. They exhaust the odd kernel by \eqref{common:restricted}.
\end{proof}
\space
\section{Recurrent conformal motion and linear wave estimates}\label{main:start}
In Sections~\ref{main:start}--\ref{end:start} we use radiation packets to drive the half-turn-odd core between the two stationary shapes in Theorem~\ref{main:counterexample}. The radiation-driven shape mechanism and the identification of the free wave and surviving profile have a five-dimensional counterpart in \cite[Section~10]{CHS5D}. We begin the proof of Theorem~\ref{main:counterexample}. Throughout Sections~\ref{main:start}--\ref{end:start} we write $Q=Q_o$ for the regular-wall representative. We first specify the competing traveling profiles and identify the impulse that reverses the motion at each node. Sections~\ref{lin:start} and \ref{lin:old-start} then estimate the wave at the current time and after the forcing has ended.

\subsection{Spatial frames and the fixed-profile assertion}\label{main:fixed-profiles}
We denote the spatial similarity group, including reflections, by
$\mathrm{Sim}$. For $a\in\R^3$, $\lambda>0$ and $O\in O(3)$, we set
\begin{equation}\label{main:eq:similarity}
 S_{a,\lambda,O}(f,g)(x)=
 \left(\lambda^{-\frac{1}{2}}f\!\left(\frac{O^{-1}(x-a)}{\lambda}\right),
       \lambda^{-\frac{3}{2}}g\!\left(\frac{O^{-1}(x-a)}{\lambda}\right)\right).
\end{equation}
This action preserves both the conserved energy and the norm of $\mathcal H$.

For a nonzero solution
$P\in\dot H^1$ of \eqref{main:stationary} and $|\ell|<1$, we write
\begin{align*}
 P_\ell(x)&=P\!\left(x+
   \left(\frac{1}{\sqrt{1-|\ell|^2}}-1\right)
   \frac{\ell\cdot x}{|\ell|^2}\ell\right)\quad(\ell\ne0),\qquad P_0=P,\\
 \vec P_\ell&=(P_\ell,-\ell\cdot\nabla P_\ell).
\end{align*}
The function $P_\ell(x-t\ell)$ solves \eqref{intro:wave}. We allow rotations,
enlarging the class in \cite{DKM}. The forward-global assertion has the
quantifier order
\begin{equation}\label{main:eq:quantifiers}
 \begin{gathered}
 \exists\,J\in\mathbb N\cup\{0\},\quad
 \{Q_j,\ell_j\}_{j=1}^J,\quad v_L,
 \qquad -\Delta Q_j=Q_j^5,\quad Q_j\in\dot H^1\setminus\{0\},\quad |\ell_j|<1,\\
 \forall\,t_n\to\infty,\quad
 \exists\,\{a_{j,n},\lambda_{j,n},O_{j,n}\}_{j,n}.
 \end{gathered}
\end{equation}
Here $v_L$ is a free wave. We write
$S_{j,n}=S_{a_{j,n},\lambda_{j,n},O_{j,n}}$. The asserted decomposition is
\begin{equation}\label{main:FP}
 \vec u(t_n)=\vec v_L(t_n)+\sum_{j=1}^J S_{j,n}\vec Q_{j,\ell_j}
               +o_{\mathcal H}(1),\qquad n\to\infty,
\end{equation}
where distinct frames satisfy \eqref{main:profile-separation}.
In \eqref{main:eq:quantifiers} the same $J$, $Q_j$, $\ell_j$ and $v_L$
must work for every sequence. Theorem~\ref{main:counterexample} rules out
\eqref{main:FP} already on the single sequence $t_n=T_n$.

\subsection{The impulse for alternating positions}\label{main:outline}
Set $s_0=0$, $s_j=T_j=4^jd^{-5}$ for $j\ge1$, $I_j=[s_j,s_{j+1}]$ and $\ell_j=s_{j+1}-s_j$. In the coordinate $\Gamma=\gamma(\xi)$ of \eqref{end:ep:coordinate}, the two endpoints are $0$ and $d r=\gamma(de)$. Their displacement on $I_j$ is $(-1)^j d r$. For $j\ge1$ let $\varphi_j$ be the continuous piecewise affine function supported on $I_{j-1}\cup I_j$ with $\varphi_j(s_j)=1$. Integration by parts yields
\begin{equation}\label{main:eq:impulse-block}
 \int\varphi_j\Gamma''=-\int\varphi_j'\Gamma'
 =(-1)^j d r\left(\frac{1}{\ell_{j-1}}+\frac{1}{\ell_j}\right).
\end{equation}
The identity holds by duality when the velocity has jumps. We place one packet group at each $s_j$, with ideal force $\eta_j^2a_j\delta_{s_j}$ and $T_j\eta_j^2=d$. Since $\varphi_j(s_k)=\delta_{jk}$, its normalized impulse equals the slope jump in \eqref{main:eq:impulse-block} divided by $\eta_j^2$. With $\ell_0=4T_0$ and $\ell_j=3T_j$ for $j\ge1$, we obtain
\begin{equation}\label{main:eq:ideal-controls}
 a_1^*=-\frac{4}{3}r,\qquad a_j^*=(-1)^j\frac{5}{3}r\quad(j\ge2).
\end{equation}
Section~\ref{cat:sec:reference} estimates the error when we test the wave packets against $\varphi_j$. In Section~\ref{end:map-section} the same test enforces continuity of the exact velocity.

All spatial norms are over $\R^3$ unless specified otherwise. The implicit
constant in $A\lesssim B$ is independent of $d$ and the terminal time.
The stationary state, packet family and path-domain constant are fixed
before $d$. Parameter derivatives and time derivatives are distinguished
by their subscripts, and we sum over repeated chart indices.
\space
\subsection{The projected Green operator}\label{lin:start}
We prove the estimate \eqref{lin:eq:green} uniformly in the terminal
time. Localized tests remove the zero modes, while one terminal
condition for each growing negative mode selects the initial data.
We compare the projected solution with a free wave, estimate two
weighted observations, and recover the remaining norms from the free
equations.

Krieger and Schlag select the unstable coordinate from future forcing
and separate the resonant part of the propagator \cite{KS}.
Beceanu and Goldberg \cite[Theorem~3]{BecGold} prove reversed Strichartz estimates when zero is neither an eigenvalue nor a resonance. Beceanu \cite{Bec} separates an explicit rank-one resonance term at the ground state. Krieger, Nakanishi and Schlag construct a
centre-stable manifold in the energy space without symmetry assumptions
\cite{KNS}. For a defocusing wave equation with a trapping potential, Jia, Liu, Schlag and Xu \cite{JLSX} construct finite-codimensional center-stable manifolds of unstable excited states. For resolvent methods at the threshold we refer to
\cite{JK,JN}. Here we must retain the Newton observation of the velocity
and remove a multidimensional odd kernel.

\subsubsection{Setting and the Green estimate}\label{lin:green-section}
The stationary input of Proposition~\ref{geo:ass:stationary} satisfies
the hypotheses below with $V=5Q^4$ and $m=4$ by
\eqref{geo:eq:decay} and \eqref{geo:eq:odd-kernel}.

We write $\mathscr R(x_1,x_2,x_3)=(-x_1,-x_2,x_3)$ and work throughout with real
half-turn odd functions, $f(\mathscr R x)=-f(x)$, or their complexification. We set $\rho=\langle x\rangle$ and write
\begin{equation*}
 D=\sqrt{-\Delta},\quad G_0=D^{-2},\quad
 X_s=L^2(\la x\ra^{-2s}\dd x),\quad
 Y_s=L^2(\la x\ra^{2s}\dd x),\quad
 \cH=\dot H^1\times L^2.
\end{equation*}
We assume $V$ is real, smooth, $\mathscr R$-invariant, and
\begin{equation}\label{lin:eq:potential}
 |\partial_x^\alpha V(x)|\lesssim_\alpha\la x\ra^{-8-|\alpha|}.
\end{equation}
Let $L=D^2-V$. We assume its full odd energy kernel has a finite
$L^2$-orthonormal basis $K=(K_1,\ldots,K_m)$ satisfying
\begin{equation*}
 |\partial_x^\alpha K_i|\lesssim_\alpha\la x\ra^{-2-|\alpha|},
 \qquad \nu=D^2K=VK=O(\la x\ra^{-10}).
\end{equation*}
Let $\Psi_l$ be an orthonormal basis of the negative eigenspace,
$L\Psi_l=-\kappa_l^2\Psi_l$, $\kappa_l>0$. It is finite, and the
eigenfunctions decay exponentially. Let $\eta_R$ be a radial cutoff tending to one. Since $K^*\eta_R K\to I_m$, the compact odd tests $\Phi_0=(K^*\eta_R K)^{-1}K^*\eta_R$ are dual to $K$ for large $R$. Fix $c>1+\|V\|_\infty$. Choose a real even $\chi\in C_c^\infty(\R)$ equal to one near $\sqrt c$, with $\chi(\sqrt{c-\kappa_l^2})=0$ for all $l$ and $\chi(\sqrt{c+\lambda})=0$ for $\lambda\ge1$. We set
\begin{equation}\label{lin:eq:spectralgauge}
 \Phi=\Phi_0\chi(\sqrt{L+c}),\qquad
 \Phi K=I_m,\quad \Phi\Psi_l=0,\quad
 P=I-K\Phi,\quad P_0=I-KK^* .
\end{equation}
By spectral calculus, $\Phi\mathbf1_{[1,\infty)}(L)=0$ and
\eqref{lin:eq:spectralgauge} holds. The new tests are odd Schwartz
functions. Indeed, for a component $f$ of $\Phi_0^*$,
we write $\chi(\sqrt{L+c})f$ as the Fourier cosine integral.
If $f$ is supported in $B_R$, finite propagation confines
$\cos(t\sqrt{L+c})f$ to $B_{R+|t|}$ and its $L^2$ norm is at most
$\|f\|_2$. Thus the exterior $L^2$ norm at radius $R+r$ is
$O_N((1+r)^{-N})$, by the rapid decay of $\widehat\chi$. This holds
also after applying $(L+c)^k$, by replacing $f$ with $(L+c)^kf$.
Applying interior elliptic estimates on unit balls and Sobolev embedding,
we obtain the same decay for every pointwise derivative.
The map $P$ is bounded on energy positions. We use $P_0$ only on $L^2$
variables. A pair is called complementary if both entries are annihilated
by $\Phi$ and by every $\Psi_l^*$. The latter pairings are defined on
$\dot H^1$ because the eigenfunctions decay exponentially.

We fix $H>0$. We extend the sources by zero after $H$, and consider
\begin{equation}\label{lin:eq:system}
 x_t=y-G_0\beta,\qquad
 y_t=-PLx+P(F_2+F_1),\qquad \Phi x=\Phi y=0,
\end{equation}
with the compatibility conditions
\begin{equation}\label{lin:eq:compat}
 \Phi G_0\beta=0,\qquad K^*F_2=0.
\end{equation}
We prescribe complementary initial data $h_c$ and set
\begin{equation*}
 a_l=\frac{\Psi_l^*y+\kappa_l\Psi_l^*x}{2\kappa_l},\qquad
 b_l=\frac{\kappa_l\Psi_l^*x-\Psi_l^*y}{2\kappa_l}.
\end{equation*}
The negative-mode conditions are $b_l(0)=0$ and $a_l(H)=0$. The remaining
initial component is $h_c$. We define
\begin{align}
 \|(x,y)\|_{\cR(I)}={}&\| (\nabla x,y)\|_{L_t^\infty L_x^2(I)}
 +\|x\|_{L_t^5L_x^{10}(I)}+\|x\|_{L_t^2X_2(I)}\notag
\\
 &+\|\nabla x\|_{L_t^2X_1(I)}+\|y\|_{L_t^2X_1(I)}
 +\|G_0y\|_{L_t^2X_4(I)}.
\label{lin:eq:green-norm}
\end{align}

We use the observation norms
\begin{equation}\label{lin:eq:observations}
 \begin{aligned}
 E(f)&=\|\nabla f\|_{L_t^\infty L_x^2},&
 S(f)&=\|f\|_{L_t^5L_x^{10}},\\
 \LE(f)&=\|f\|_{L_t^2X_2},&
 \OO(g)&=\|G_0g\|_{L_t^2X_4}.
 \end{aligned}
\end{equation}

\begin{theorem}[Uniform frozen Green estimate]\label{lin:thm:green}
For odd sources satisfying \eqref{lin:eq:compat}, there is a unique energy
solution of \eqref{lin:eq:system} with the stated mixed conditions, and
\begin{equation}\label{lin:eq:green}
 \|(x,y)\|_{\cR([0,\infty))}\lesssim
 \|h_c\|_{\cH}+\|\beta\|_{L_t^2Y_2}
 +\|F_2\|_{L_t^2Y_1}+\|F_1\|_{L_t^1L_x^2}.
\end{equation}
The map is linear, and the constant is independent of $H$. It may depend
on $V$, its full kernel, the localized tests, and the negative spectrum.
\end{theorem}

The future forcing determines the unstable initial amplitudes in
\eqref{lin:eq:negative}. With these amplitudes fixed we solve forward
and extend the solution by zero to $t<0$ for the Laplace transform.

Sections~\ref{lin:free-section}--\ref{lin:obs-section} prove
Theorem~\ref{lin:thm:green}. The observation multiplier will also be used
for the decay estimates in Section~\ref{lin:old-start}.

\subsubsection{Free estimates and kernels}\label{lin:free-section}

We put $C_0(t)=\cos(tD)$ and $S_0(t)=D^{-1}\sin(tD)$. For the free system
\begin{equation*}
 u_t=v-G_0\beta,\qquad v_t=\Delta u+F_2+F_1,
\end{equation*}
we first prove
\begin{equation}\label{lin:eq:freebound}
 \|(u,v)\|_{\cR([0,\infty))}\lesssim
 \|(u_0,v_0)\|_{\cH}+\|\beta\|_{L_t^2Y_2}
 +\|F_2\|_{L_t^2Y_1}+\|F_1\|_{L_t^1L_x^2}.
\end{equation}
The estimate holds for arbitrary odd data and sources in these spaces.

\paragraph{Odd Rellich and the Newton input.}
For smooth compact odd $v$,
\begin{equation}\label{lin:eq:rellich}
 \|r^{-2}v\|_2\le\frac{4}{5}\|\Delta v\|_2.
\end{equation}
We test $-\Delta v$ against $r^{-2}v$ and complete the square in $r$.
Since $v(0)=0$, the boundary term at zero vanishes. In polar coordinates,
\begin{align*}
 \Re\langle-\Delta v,r^{-2}v\rangle-\frac{5}{4}\|r^{-2}v\|_2^2
 &=\int_0^\infty\!\int_{\mathbb S^2}
    \left|\partial_rv-\frac{v}{2r}\right|^2\,d\omega\,dr\\
 &\quad+\int_0^\infty\frac{1}{r^2}\!\int_{\mathbb S^2}
       (|\nabla_\omega v|^2-2|v|^2)\,d\omega\,dr\ge0 .
\end{align*}
Each spherical restriction has mean zero by oddness. The spherical
Poincar\'e inequality therefore bounds the second integral from below by zero.
Cauchy--Schwarz now proves \eqref{lin:eq:rellich}.
We apply \eqref{lin:eq:rellich} by radial cutoff to Newton potentials
of odd Schwartz sources and pass to the closure by density. Duality
and integration by parts yield
\begin{equation}\label{lin:eq:weightedNewton}
 \|G_0b\|_2+\|\la x\ra\nabla G_0b\|_2\lesssim\|b\|_{Y_2}.
\end{equation}
For the weighted gradient we use, with $v=G_0b$,
\begin{equation*}
 \|r\nabla v\|_2^2=\ip{r^2b}{v}+3\|v\|_2^2,
 \qquad\|\nabla v\|_2^2=\ip bv.
\end{equation*}
At each fixed time, duality in \eqref{lin:eq:rellich} also yields
$G_0:L^2_{\rm odd}\to X_2$. For the spacetime observation we use $X_4$
in \eqref{lin:eq:homNewton}.

\paragraph{Energy, local energy and Strichartz.}
We set $f(r)=\arctan r$, $q(r)=\frac{f(r)}{r}$ and $M=f\partial_r+q$. The
energy and Morawetz identities are
\begin{equation}\label{lin:eq:energy-morawetz}
\begin{aligned}
 \frac{d}{dt}\frac{\|\nabla u\|_2^2+\|v\|_2^2}{2}
 &=-\ip u\beta+\ip v{F_2+F_1},\\
 \frac{d}{dt}\ip v{Mu}
 &=-\mathcal B+\ip{F_2+F_1}{Mu}-\ip v{MG_0\beta},
\end{aligned}
\end{equation}
where
\begin{equation*}
 \mathcal B=\int\frac{f'}{2}(v^2+u_r^2)
 +\left(\frac{f}{r}-\frac{f'}{2}\right)|\nabla_\perp u|^2
 -\frac{\Delta q}{2}u^2\dd x.
\end{equation*}
Here $f'=\la r\ra^{-2}$, $-\frac{\Delta q}{2}=\la r\ra^{-4}$, and
$\frac{f}{r}-\frac{f'}{2}\gtrsim\la r\ra^{-2}$. Hardy's inequality controls the
time endpoints, and $\|MG_0\beta\|_{Y_1}\lesssim\|\beta\|_{Y_2}$ by
\eqref{lin:eq:weightedNewton}. We denote the initial energy norm by $E_0$, its supremum on the interval by $E$,
and the sum of the three local-energy norms by $J$. Adding the two integrated bounds from \eqref{lin:eq:energy-morawetz}
with a sufficiently small fixed coefficient on the Morawetz bound, we obtain
\begin{equation*}
 E^2+J^2\lesssim E_0^2+
 J\big(\|\beta\|_{L_t^2Y_2}+\|F_2\|_{L_t^2Y_1}\big)
 +E\|F_1\|_{L_t^1L_x^2}.
\end{equation*}
We absorb both products by Young's inequality, uniformly in the interval length.

The free energy Strichartz estimate \cite{KT}, the adjoints of the
homogeneous position and velocity local-energy maps, and the
Christ--Kiselev truncation \cite{CK} yield the $L_t^5L_x^{10}$ bound.
The relevant time exponents are $2<5$. Explicitly, the full-time
convolutions $-C_0G_0*\beta$ and $S_0*F_2$ factor through bounded maps
into free energy data. We apply retarded truncation to these two
convolutions and Minkowski's inequality to the $F_1$ term.

\paragraph{Odd free resolvent and the Newton observation.}
Let
\begin{equation*}
 G_z=(D^2-z^2)^{-1},\qquad
 G_z(r)=\frac{e^{izr}}{4\pi r},\qquad
 H_z=\frac{G_z-G_0}{z^2},\quad H_0(r)=-\frac{r}{8\pi}.
\end{equation*}
All these kernels act in the odd sector. The scalar term $\frac{i}{4\pi z}$
in $H_z$ cancels.

The odd kernel estimates needed below are
\begin{equation}\label{lin:eq:kernelbounds}
 \|H_z\|_{Y_2\to X_2}
 +\|\partial_z(zH_z)\|_{Y_2\to X_2}
 +\|\partial_zH_z\|_{Y_4\to X_4}\lesssim1.
\end{equation}
The radial derivative of the exact kernel has the integral form
\begin{equation*}
 4\pi\,\partial_rH_z(r)=g(zr),\qquad
 g(w)=-\int_0^1 s e^{isw}\dd s,\qquad g(w)+wg'(w)=-e^{iw}-g(w).
\end{equation*}
The last identity follows by integration by parts. For $\Ima z\ge0$,
the radial derivatives of $H_z$ and $\partial_z(zH_z)$ are therefore
bounded, while $|\partial_r\partial_zH_z(r)|\lesssim r$. For
$r=|x-y|$, $r'=|\mathscr R x-y|$, one has
\begin{equation*}
 |r-r'|\le2\min(|x|,|y|),\qquad |r^2-r'^2|\le4|x||y|.
\end{equation*}
Thus the first two odd kernels in \eqref{lin:eq:kernelbounds} are
bounded by $C\min(|x|,|y|)$, and the third by $C|x||y|$. The third bound is Hilbert--Schmidt from $Y_4$ to $X_4$.

For the mixed weights we write $4\pi zH_z(r)=k(zr)$, where $k(w)=i\int_0^1e^{isw}\dd s$. Since $|k(w)|+|wk'(w)|\lesssim1$, its odd difference is bounded by $C\min(1,\frac{|x|}{|y|})$, using the derivative when $|y|>2|x|$. We put $a=|x|$, $b=|y|$. After angular Cauchy--Schwarz and radial conjugation to $L^2(\dd r)$, the kernel bounds for the two weight pairs satisfy
\begin{equation*}
 \frac{ab\min(a,b)}{\langle a\rangle^2\langle b\rangle^2}
 +\frac{ab\min(1,\frac{a}{b})}{\langle a\rangle^2\langle b\rangle}
 \lesssim\frac{1}{\max(a,b)}.
\end{equation*}
Since $\int_0^\infty\frac{b^{-\frac{1}{2}}}{\max(a,b)}\dd b=4a^{-\frac{1}{2}}$, Schur's test proves \eqref{lin:eq:kernelbounds} and
\begin{equation}\label{lin:eq:mixed-newton}
 \|zH_z\|_{Y_1\to X_2}\lesssim1\qquad(\Ima z\ge0).
\end{equation}
The derivative bounds in \eqref{lin:eq:kernelbounds} also control the differences of $zH_z:Y_2\to X_2$ and $H_z:Y_4\to X_4$ by $C|z-z'|$.

We prove the homogeneous Newton estimate directly on frequency spheres.
For odd $g$ the integral of $\widehat g(\lambda\omega)$ over $\mathbb S^2$ vanishes at almost every $\lambda>0$. Since $|e^{is}-1|\le |s|$, we have
\begin{equation*}
 \left|\int_{\mathbb S^2}e^{i\lambda x\cdot\omega}\widehat g(\lambda\omega)\,d\omega\right|
 \lesssim \lambda |x|\,\|\widehat g(\lambda\,\cdot)\|_{L^2(\mathbb S^2)} .
\end{equation*}
In the polar Fourier formula for $G_0e^{-itD}g$, the multiplier $\lambda^{-2}$ cancels the Jacobian $\lambda^2$. Plancherel in $t$ and then in $x$ yields, for every $s>\frac{5}{2}$,
\begin{equation}\label{lin:eq:homNewton}
 \|G_0e^{-itD}g\|_{L_t^2X_s}^2
 \lesssim \||x|\langle x\rangle^{-s}\|_2^2
 \int_0^\infty\lambda^2\|\widehat g(\lambda\,\cdot)\|_{L^2(\mathbb S^2)}^2\,d\lambda
 \lesssim_s\|g\|_2^2 .
\end{equation}
We first take smooth annular Fourier data and then use density. With $s=4$ this controls the homogeneous velocity and, by Minkowski's inequality, the $F_1$ term. For the retarded terms, the Laplace multipliers of $S_0G_0$ and $C_0G_0$ are $H_z$ and $-izH_z$. We apply time Plancherel with an exponential damping and let the damping decrease to zero. The bounds \eqref{lin:eq:kernelbounds} and \eqref{lin:eq:mixed-newton}, with $X_2\subset X_4$, imply
\begin{equation*}
 \|S_0G_0*\beta\|_{L_t^2X_4}+\|C_0G_0*F_2\|_{L_t^2X_4}
 \lesssim\|\beta\|_{L_t^2Y_2}+\|F_2\|_{L_t^2Y_1}.
\end{equation*}
This proves \eqref{lin:eq:freebound}. For separated time intervals we will use the same retarded kernels in physical space. Kirchhoff's formula yields, for $\tau>0$,
\begin{equation}\label{lin:eq:newtonkernels}
 (S_0(\tau)G_0)(r)=\frac{\min(\tau,r)}{4\pi r},\qquad
 (C_0(\tau)G_0)(r)=\frac{\mathbf1_{\tau<r}}{4\pi r}.
\end{equation}

\subsubsection{Construction and the choice of comparison wave}\label{lin:construction-section}

We choose the negative modes before comparing with a free wave.
We set $g_l=\ip{\Psi_l}{G_0\beta}$ and $f_l=\ip{\Psi_l}{F_2+F_1}$.
Solving their scalar equations with $a_l(H)=b_l(0)=0$, we obtain
\begin{align}\label{lin:eq:negative}
 a_l(t)&=-\int_t^H e^{-\kappa_l(s-t)}
             \frac{f_l-\kappa_lg_l}{2\kappa_l}(s)\dd s,\\
 b_l(t)&=-\int_0^t e^{-\kappa_l(t-s)}
             \frac{f_l+\kappa_lg_l}{2\kappa_l}(s)\dd s.\notag
\end{align}
Their initial unstable amplitudes are bounded by the right side of
\eqref{lin:eq:green}, uniformly in $H$. We select
\begin{equation}\label{lin:eq:selecteddata}
 (x_0,y_0)=h_c+\sum_l a_l(0)(\Psi_l,\kappa_l\Psi_l).
\end{equation}
The map $Vx+K\Phi Lx:\dot H^1\to L^2$ is bounded, so the projected
generator is a bounded perturbation of the free energy generator.
Duhamel's formula defines the unique solution with data \eqref{lin:eq:selecteddata}
on each finite interval.

We first justify the Laplace transform on the dense class with
$x_0\in H^1$. Since $x_t=y-G_0\beta\in L^1([0,H];L^2)$,
$x(H)\in H^1$. After the source stops, let $u$ solve
$u_{tt}+Lu=0$ with data $(x,y)(H)$. Then $(x,y)=(Pu,Pu_t)$
by $LP=L$ and uniqueness. Its negative modes decay because $a_l(H)=0$.
On the nonnegative spectrum, $|\cos(t\sqrt\lambda)|\le1$ and
$\bigl|\frac{\sin(t\sqrt{\lambda})}{\sqrt{\lambda}}\bigr|\le t$, including $\lambda=0$.
By the spectral theorem,
\begin{equation*}
 \|u_t\|_2+\||L|^{\frac{1}{2}}u\|_2\lesssim_H 1,\qquad
 \|u\|_2\lesssim_H 1+t,\qquad
 \|\nabla u\|_2^2=\langle Lu,u\rangle+\int V|u|^2 .
\end{equation*}
Since $V\in L^\infty$ and $P$ is bounded on both $\dot H^1$ and $L^2$, we obtain
\begin{equation}\label{lin:eq:coercive}
 \|(x,y)(t)\|_{\cH}\lesssim_H 1+t\qquad(t\ge H).
\end{equation}
We use the preliminary bound \eqref{lin:eq:coercive}
only to apply damped Plancherel.

We choose the comparison free wave $(u_f,v_f)$ with data $(x_0,P_0y_0)$,
first source $-G_0\beta$, and second source $F_2+P_0F_1$. Its energy
and observations satisfy \eqref{lin:eq:freebound} since $P_0$ is
orthogonal on $L^2$. Since $PK=0$, the discarded
initial velocity and force take values in the kernel and can both
be absorbed into the same localized-gauge multiplier.
We set $z_f=G_0v_f$. We test the free velocity equation against $K$ to obtain the causal identity
\begin{equation}\label{lin:eq:charge}
 \partial_t(\nu^*z_f)=-\nu^*u_f .
\end{equation}
Indeed $\nu^*z_f=K^*v_f$, while $K^*(F_2+P_0F_1)=0$ and
$K^*v_f(0)=0$. All position tests pair with the decaying $\nu$.

\subsubsection{The observation multiplier and norm recovery}\label{lin:obs-section}

We recover the two observations $(x,G_0y)$ from $(u_f,G_0v_f)$
through a bounded time Fourier multiplier. We use
$\widehat f(z)=\int_0^\infty e^{izt}f(t)\dd t$ for $\Ima z>0$ and
include the initial atoms in the causal distributions. We set $\cL_z=(H_z-H_0)\nu$.

\paragraph{The transformed gauge equations.}
We write $X=\widehat x$, $Z=G_0\widehat y$. In the difference between the
causal physical and free velocity equations the forcing is $VX+K\gamma$,
where
\begin{equation*}
 \gamma=\Phi LX-\Phi\widehat{(F_2+F_1)}
                  +K^*\widehat F_1+K^*y_0 .
\end{equation*}
The last term is the transform of the initial velocity difference
$K^*y_0\delta_0$ in kernel coordinates.
The two position equations have the same first source and initial
position. Thus
\begin{align}\label{lin:eq:transformed}
 X&=\widehat u_f+G_zVX+H_z\nu\gamma,\\
 Z&=\widehat z_f-izH_zVX-\frac{i}{z}\cL_z\gamma.\notag
\end{align}
Here
\begin{equation}\label{lin:eq:newtonidentity}
 G_zK=H_z\nu,\qquad
 G_0H_z\nu=\frac{(H_z-H_0)\nu}{z^2}.
\end{equation}
The first identity follows from $G_0\nu=K$ and the resolvent identity. For an odd Schwartz test $f$, \eqref{lin:eq:weightedNewton} implies $G_0f\in L^2$. Hence
\begin{equation*}
 G_{i\sigma}f=\frac{D^2}{D^2+\sigma^2}G_0f\to G_0f\quad\text{in }L^2,
 \qquad \langle H_{i\sigma}\nu,f\rangle=\langle K,G_{i\sigma}f\rangle .
\end{equation*}
Letting $\sigma\downarrow0$ in \eqref{lin:eq:kernelbounds}, we obtain $H_0\nu=G_0K$ by duality. A second application of the resolvent identity proves the second equality in \eqref{lin:eq:newtonidentity}. Taking $f=\nu$ by density, we also obtain $\nu^*H_0\nu=K^*K=I$.

\paragraph{The position and charge equations at zero frequency.}
Let $A_z=I-G_zV$ on $X_2$ and $b=\gamma/z$. Since $G_0\nu=K$ and $K^*V=\nu^*$, we have
\begin{equation}\label{lin:eq:threshold-zeros}
 A_zK=-z^2H_z\nu,\qquad \nu^*A_z=-z^2\nu^*H_zV .
\end{equation}
The free equations, including the initial traces, imply
\begin{equation}\label{lin:eq:free-charge-gauge}
 \nu^*\widehat u_f=iz\nu^*\widehat z_f,\qquad
 -iz\Phi\widehat u_f=(\Phi D^2)\widehat z_f .
\end{equation}
The first identity is \eqref{lin:eq:charge}. The second uses $\Phi x_0=\Phi G_0\beta=0$. We add $K\Phi X=0$ to the position equation in \eqref{lin:eq:transformed}, then test the original equation against $-\nu^*/z$. By \eqref{lin:eq:threshold-zeros}--\eqref{lin:eq:free-charge-gauge},
\begin{equation}\label{lin:eq:threshold-system}
 \begin{pmatrix}
 A_z+K\Phi & -zH_z\nu\\
 z\nu^*H_zV & \nu^*H_z\nu
 \end{pmatrix}\binom{X}{b}=
 \binom{\widehat u_f}{-i\nu^*\widehat z_f}.
\end{equation}
At $z=0$ this matrix is $\operatorname{diag}(A_0+K\Phi,I)$ by \eqref{lin:eq:newtonidentity}. If $(A_0+K\Phi)f=0$, testing by $\nu^*$ yields $(\nu^*K)\Phi f=0$. The matrix $\nu^*K$ is the positive Gram matrix of $\nabla K$. Hence $f=G_0Vf\in\dot H^1$ and $\Phi f=0$. The odd kernel is spanned by $K$, so $f=0$. Since $G_0V$ is Hilbert--Schmidt on $X_2$, the Fredholm alternative proves invertibility. The entries of \eqref{lin:eq:threshold-system} are Lipschitz near zero by \eqref{lin:eq:kernelbounds}, using $\nu\in Y_4$ in the last entry. Its inverse is therefore uniformly bounded and Lipschitz on a small upper half-disk.

\paragraph{Nonzero frequencies and recovery of the velocity.}
For real $z\ne0$, $A_z$ is Fredholm on $X_2$. If $A_z f=0$, then $Vf\in Y_6$, $|f|\le G_0|Vf|$ and $\|f\|_4\lesssim\|Vf\|_{\frac{12}{11}}$ by Hardy--Littlewood--Sobolev. Outgoing injectivity \cite[Lemma~9]{GS} for $z>0$ and conjugation for $z<0$ imply $f=0$. Away from a neighborhood of zero the inverse is uniformly bounded in a thin upper strip, by continuity on compact sets and the following high-frequency bounds. The first is Agmon's weighted free-resolvent estimate \cite{Agmon}, displayed in \cite[(1.5)]{Syl}.
\begin{equation}\label{lin:eq:highG}
 \|G_z\|_{Y_2\to X_2}\lesssim |z|^{-1},\qquad
 \|\partial_zG_z\|_{Y_2\to X_2}\lesssim1\qquad(|z|\ge1).
\end{equation}
The second bound follows directly from the kernel. By \eqref{lin:eq:threshold-zeros}, $A_z(X+Kb/z)=\widehat u_f$. Applying $P$ and $z\Phi$, respectively, and using \eqref{lin:eq:free-charge-gauge}, we obtain
\begin{equation}\label{lin:eq:outgoing-reduction}
 X=P A_z^{-1}\widehat u_f,\qquad
 b=i(\Phi D^2)\widehat z_f+z\Phi G_zV A_z^{-1}\widehat u_f .
\end{equation}
At every frequency the second row of \eqref{lin:eq:transformed} reads
\begin{equation}\label{lin:eq:lowmatrix}
 Z=\widehat z_f-izH_zVX-i(H_z-H_0)\nu b .
\end{equation}
At high frequency $z\Phi G_z=z^{-1}[(\Phi D^2)G_z-\Phi]$ and $zH_z=(G_z-G_0)/z$. Thus \eqref{lin:eq:kernelbounds} and \eqref{lin:eq:highG} control all coefficients and their first differences.
\begin{lemma}[Multiplier for the two wave observations]
\label{lin:lem:observation-multiplier}
We set
\begin{equation*}
 \mathbb X=\mathbb Y=X_2\oplus X_4,\qquad \cZ_f=(u_f,z_f).
\end{equation*}
There is $\epsilon_0>0$ and a family
$\mathbb M_\epsilon(\lambda)\in\mathcal L(\mathbb X,\mathbb Y)$,
$0\le\epsilon\le\epsilon_0$, such that
\begin{equation*}
 \sup_{\epsilon,\lambda}\|\mathbb M_\epsilon(\lambda)\|_{\mathcal L(\mathbb X,\mathbb Y)}
 <\infty,\qquad
 \|\mathbb M_0(\lambda)-\mathbb M_0(\mu)\|_{\mathcal L(\mathbb X,\mathbb Y)}
 \lesssim|\lambda-\mu|.
\end{equation*}
For every real $\lambda$, $\mathbb M_\epsilon(\lambda)$ converges in
operator norm to $\mathbb M_0(\lambda)$ as $\epsilon\downarrow0$.
On the actual transformed observations,
\begin{equation*}
 (X,Z)=\mathbb M_\epsilon(\lambda)\widehat{\cZ_f}
       (\lambda+i\epsilon).
\end{equation*}
The multiplier $\mathbb M_0(D_t)$ is bounded from
$L^2(\R;\mathbb X)$ to $L^2(\R;\mathbb Y)$. For
$\chi\in C_c^\infty(\R)$ and $L\ge1$, it also satisfies
\begin{equation}\label{lin:eq:commutator}
 \|[\chi\bigl(\frac{t}{L}\bigr),\mathbb M_0(D_t)]\mathcal Z\|_{L_t^2\mathbb Y}
 \lesssim L^{-1}\|\mathcal Z\|_{L_t^2\mathbb X}.
\end{equation}
The operator bounds hold for arbitrary, independent coordinates in
$\mathbb X$, without compatibility conditions on $\mathcal Z$.
\end{lemma}

\begin{proof}
We join the maps for $(X,b)$ in \eqref{lin:eq:threshold-system} and \eqref{lin:eq:outgoing-reduction} with a fixed smooth cutoff in $\Re z$, supported in the half-disk for sufficiently small $\Im z$. On the overlap both maps recover the physical pair. For independent inputs the cutoff defines a bounded Lipschitz extension. We then apply \eqref{lin:eq:lowmatrix}. The $L^2$ multiplier bound follows by Plancherel in time.

Lipschitz continuity bounds the Fourier kernel of the commutator
in operator norm by $C|\lambda-\mu|L|\widehat\chi(L(\lambda-\mu))|$.
Its $L^1$ norm in the difference variable is $O(L^{-1})$.
Young's inequality proves \eqref{lin:eq:commutator} for all
$\mathcal Z\in L^2(\R;\mathbb X)$, including observations multiplied
by time cutoffs.
\end{proof}

\paragraph{Plancherel, norm recovery, and the cutoff commutator.}
The preliminary energy construction makes the damped observations
square integrable. Applying the multiplier identity and Plancherel, then
Fatou as $\epsilon\downarrow0$, proves
\begin{equation}\label{lin:eq:obsidentity}
 (x,G_0y)=\mathbb M_0(D_t)\cZ_f,
 \qquad\|x\|_{L_t^2X_2}+\|G_0y\|_{L_t^2X_4}\lesssim\mathfrak D,
\end{equation}
where $\mathfrak D$ is the right side of \eqref{lin:eq:green}.
The global $L^2$ bounds permit passage to the limit in the identity
by dominated convergence on the free inputs and uniform boundedness
of the multipliers.

We set $v=P_0y$ and $c=K^*y=\nu^*G_0y$. The exact free form is
\begin{equation}\label{lin:eq:recovery}
 x_t=v-G_0(\beta-\nu c),\qquad
 v_t=\Delta x+Vx+F_2+P_0F_1.
\end{equation}
By \eqref{lin:eq:obsidentity},
\begin{equation*}
 \|\nu c\|_{L_t^2Y_4}+\|Vx\|_{L_t^2Y_1}
 \lesssim \|G_0y\|_{L_t^2X_4}+\|x\|_{L_t^2X_2}\lesssim\mathfrak D .
\end{equation*}
We apply \eqref{lin:eq:freebound} to \eqref{lin:eq:recovery} and recover
the physical velocity from $y=Pv$.
This proves Theorem~\ref{lin:thm:green}, including the derivative local-energy
norm needed later. We preserve compatibility in the density argument
by applying $P^*$ to each compact force and subtracting a compact smooth Newton-density
correction $D^2\eta\,\Phi G_0\beta$ with $\Phi\eta=I$.

The same Green operator acts on the half-line. For compactly time-supported
sources the future integral selecting each growing negative mode does
not change when the terminal time is increased beyond the support.
Thus the resulting global response is independent of that endpoint.
For general half-line sources, we truncate in time and keep the data fixed. By \eqref{lin:eq:green}, the truncated solutions are Cauchy in $\cR$ and in energy. The limit solves the equations
and has the same growing-mode selection, given by its convergent
future integral to infinity. We may therefore apply
Theorem~\ref{lin:thm:old} on the half-line to any source stopped at
a finite time.

\subsection{Decay after the forcing has ended}\label{lin:old-start}
We estimate a response after its forcing has stopped. The resulting decay
separates earlier radiation from the control on each later interval.
We first prove free estimates and then transfer them to the projected
flow using the observation multiplier from Section~\ref{lin:obs-section}.
We retain the Green system, projections and norms of Section~\ref{lin:start}.

Let $T\ge1$, $H\ge \frac{T}{2}$, and assume that all sources
vanish after $\frac{T}{2}$. Suppose
\begin{equation}\label{lin:eq:oldclasses}
 \begin{gathered}
 \beta\in L_t^2Y_4,\qquad
 F_2\in L_t^2Y_{\frac{7}{4}},\qquad \nabla F_2\in L_t^2L_x^1,\\
 F_1=F_{\rm p}+G_0\Gamma,\qquad
 \Gamma\in L_t^1Y_4,\\
 F_{\rm p}\in L_t^1L_x^2\cap L_t^{\frac{4}{3}}W_x^{1,1},\\
 h_c=(G_0B_c,G_0\Gamma_c),\qquad B_c,\Gamma_c\in Y_4.
 \end{gathered}
\end{equation}
All densities, forces and data are odd. We retain \eqref{lin:eq:compat} and the
complementary condition on $h_c$. Let $A$ be the sum of the
source and density norms in \eqref{lin:eq:oldclasses}. The Sobolev norm is inhomogeneous. The same conclusion holds with $F_2\in L_t^2W_x^{1,\frac{10}{9}}$ in place of the gradient assumption, with the corresponding change in $A$.

\begin{theorem}[Decay of earlier forcing]\label{lin:thm:old}
Under \eqref{lin:eq:oldclasses},
\begin{align}\label{lin:eq:oldobs}
 \|x\|_{L_t^2([bT,\infty);X_2)}
 +\|G_0y\|_{L_t^2([bT,\infty);X_4)}
 &\lesssim_b AT^{-\frac{3}{4}},\qquad b>\frac{1}{2},\\
 \|x\|_{L_t^5([T,\infty);L_x^{10})}&\lesssim AT^{-\frac{1}{10}}.
 \label{lin:eq:oldS}
\end{align}
Its global $\cR$ norm is $O(A)$. The constants are independent of $T,H$.
Increasing $H$ beyond the source support leaves the selected solution
unchanged.
\end{theorem}

\subsubsection{Free separated observations}\label{lin:old-section}

We take sources supported in $[0,\frac{T}{2}]$ and observe on $J_T=[T,4T]$, with $T\ge1$. For a velocity source $F$ let $\mathcal K F=(S_0*F,C_0G_0*F)$, with retarded convolutions. We claim
\begin{equation}\label{lin:eq:tailtable}
 \|\mathcal K F\|_{L^2(J_T;X_2\oplus X_4)}
 \lesssim T^{1-s}\|F\|_{L_t^2Y_s},\qquad 1\le s\le2 .
\end{equation}
The endpoint $s=1$ is \eqref{lin:eq:freebound}. For $s=2$ the time multiplier $(G_z,-izH_z)$ is Lipschitz from $Y_2$ to $X_2\oplus X_4$. Indeed $4\pi|\partial_zG_z(r)|\le1$, $\langle x\rangle^{-2}\in L^2$, and \eqref{lin:eq:kernelbounds} controls the second entry. Choose $\chi(t/T)$ equal to one on $J_T$ and zero on the source interval. Then $\chi(t/T)\mathcal K F=[\chi(t/T),\mathcal K]F$, and the Fourier estimate in the proof of \eqref{lin:eq:commutator} supplies the factor $T^{-1}$. For $1<s<2$ we apply the $s=2$ bound on $|x|\le T$ and the $s=1$ bound on $|x|>T$, obtaining \eqref{lin:eq:tailtable}.

For the Newton data we estimate a source at a single time. For $\ell\ge1$ and odd $b\in Y_4$, we have the entrywise bounds
\begin{equation}\label{lin:eq:newton-time-atoms}
 \begin{pmatrix}
 \|C_0(\ell)G_0 b\|_{X_2}&\|S_0(\ell)G_0 b\|_{X_2}\\
 \|S_0(\ell)G_0 b\|_{X_4}&\|C_0(\ell)G_0^2 b\|_{X_4}
 \end{pmatrix}
 \lesssim \|b\|_{Y_4}
 \begin{pmatrix}\ell^{-2}&\ell^{-\frac{3}{2}}\\
 \ell^{-\frac{7}{2}}&\ell^{-\frac{5}{2}}\end{pmatrix}.
\end{equation}
For the first entry we use $\|C_0(\ell)G_0 b\|_2\lesssim\|b\|_{Y_2}$ on $|x|\ge\frac{\ell}{2}$. Inside this ball, \eqref{lin:eq:newtonkernels} bounds the convolution by $C\ell^{-1}\int_{|y|>\ell/2}|b(y)|\,dy\lesssim\ell^{-\frac{7}{2}}\|b\|_{Y_4}$.

The other three entries use one cancellation. If $K$ is radial, constant on $B_\ell$ and Lipschitz with constant $L$, then $\int b=0$ and $(K*b)(0)=0$. Subtract $K(x)$ in the exterior integral and $K(-y)$ in the interior integral. By Cauchy--Schwarz,
\begin{equation*}
 |K*b(x)|\lesssim L\|b\|_{Y_4}
 \begin{cases}
  1,& |x|\ge\frac{\ell}{2},\\
  |x|\ell^{-\frac{5}{2}},& |x|<\frac{\ell}{2}.
 \end{cases}
\end{equation*}
Here $\||y|b\|_1\lesssim\|b\|_{Y_4}$, and the interior difference vanishes for $|y|\le\frac{\ell}{2}$. Integration against the weights yields $\|K*b\|_{X_p}\lesssim L\ell^{\frac{3}{2}-p}\|b\|_{Y_4}$ for $p=2,4$. For $S_0(\ell)G_0$ we have $L\lesssim\ell^{-1}$ by \eqref{lin:eq:newtonkernels}. For $C_0(\ell)G_0^2$ we have $L\lesssim1$, since integration of $\partial_\ell(C_0(\ell)G_0^2)=-S_0(\ell)G_0$ with initial kernel $H_0(r)=-\frac{r}{8\pi}$ yields
\begin{equation*}
 -\frac{\ell}{4\pi}-\frac{(r-\ell)_+^2}{8\pi r}.
\end{equation*}
This proves \eqref{lin:eq:newton-time-atoms}. For a source at time $v$ we have $\ell=t-v\asymp T$ on $J_T$. Minkowski's inequality costs $T^{\frac{1}{2}}$ for the observation norm and, for an $L_t^2Y_4$ density, at most another $T^{\frac{1}{2}}$ for its time integral. The two columns thus cover initial data, $Y_4$-valued time measures and the first Newton forcing. Both arguments apply on every later dyadic interval with fixed separated endpoints.

\paragraph{The pure nonlinear source channel.}
For $f\in C_c^\infty$ put $m_x(r)=\int_{|x-y|=r}|f(y)|\,dS_y$. By polar integration $\int_0^\infty m_x(r)\,dr=\|f\|_1$. Integrating $-m_x'$ from $r$ to infinity and discarding the nonpositive contribution from differentiating the area factor, we find $m_x(r)\le\|\nabla f\|_1$. Kirchhoff's formula therefore implies, uniformly in $x$ and $L>0$,
\begin{equation}\label{lin:eq:spherical-observation}
 4\pi t|S_0(t)f(x)|\le m_x(t),\qquad
 \|S_0(\cdot)f(x)\|_{L^2(L,\infty)}
 \lesssim L^{-1}\|f\|_1^{\frac{1}{2}}\|\nabla f\|_1^{\frac{1}{2}} .
\end{equation}
For the Newton velocity observation we use oddness directly in \eqref{lin:eq:newtonkernels}. For $a,b>0$,
\begin{equation}\label{lin:eq:step-distance}
 \int_0^\infty\left|\frac{\mathbf1_{t<a}}{a}
          -\frac{\mathbf1_{t<b}}{b}\right|^2dt
   =\frac{1}{a}+\frac{1}{b}-\frac{2\min(a,b)}{ab}
   =\frac{|a-b|}{ab}.
\end{equation}
For $t>L$ we may replace the radii $|x-y|,|\mathscr Rx-y|$ by their maxima with $L$. Their difference is at most $2|x|$. Subtracting the kernels at $x$ and $\mathscr Rx$, then applying Minkowski and \eqref{lin:eq:step-distance}, we obtain
\begin{equation*}
 \|C_0(\cdot)G_0 f(x)\|_{L^2(L,\infty)}
 \lesssim L^{-1}|x|^{\frac{1}{2}}\|f\|_1 .
\end{equation*}
For $L\ge1$ and $|x|\le\frac{L}{2}$, only $|y|>\frac{L}{2}$ contributes. On $|x|>\frac{L}{2}$ we have $\|\mathbf1_{|x|>\frac{L}{2}}|x|^{\frac{1}{2}}\|_{X_4}\lesssim L^{-2}$. Thus for every odd $f\in L^1$,
\begin{equation}\label{lin:eq:odd-tail}
 \|C_0(\cdot)G_0f\|_{L^2((L,\infty);X_4)}
 \lesssim L^{-1}\|\mathbf1_{|y|>\frac{L}{2}}f\|_1+L^{-3}\|f\|_1 .
\end{equation}
If $F_{\rm p}$ is supported in $[0,\frac{T}{2}]$ with $T\ge2$, we take $L=\frac{T}{2}$ in \eqref{lin:eq:spherical-observation} and \eqref{lin:eq:odd-tail}. Since $1\in X_2$, Minkowski and H\"older in the source time imply, for $0\le\theta\le2$,
\begin{equation}\label{lin:eq:ppobs}
 \begin{aligned}
 \|S_0*F_{\rm p}\|_{L^2([T,\infty);X_2)}
 &\lesssim T^{-\frac{3}{4}}\|F_{\rm p}\|_{L_t^{\frac{4}{3}}W_x^{1,1}},\\
 \|C_0G_0*F_{\rm p}\|_{L^2([T,\infty);X_4)}
 &\lesssim T^{-\frac{3}{4}-\theta}\|\langle x\rangle^\theta F_{\rm p}\|_{L_t^{\frac{4}{3}}L_x^1}.
 \end{aligned}
\end{equation}
We use $\theta=0$ for the pure quintic.
\subsubsection{Separated Strichartz bounds}
We estimate the position in $L_t^5L_x^{10}$ after the source stops.
From \eqref{lin:eq:spherical-observation}, Plancherel and dual Sobolev embedding, for $\ell>0$,
\begin{equation*}
 \|S_0(\ell)F\|_\infty\lesssim\ell^{-1}\|\nabla F\|_1,\qquad
 \|S_0(\ell)F\|_2\lesssim
     \min\{\|F\|_{\frac{6}{5}},\,\ell\|F\|_2\}.
\end{equation*}
We interpolate the output by $\|u\|_{10}\le\|u\|_\infty^{\frac{4}{5}}\|u\|_2^{\frac{1}{5}}$. With $\ell=t-s\asymp T$ on the separated intervals, Minkowski and H\"older imply
\begin{equation}\label{lin:eq:ordinary-separated}
\begin{aligned}
 \|S_0*F_2\|_{L_t^5(J_T;L^{10})}
 &\lesssim T^{-\frac{1}{10}}\|\nabla F_2\|_{L_t^2L^1}^{\frac{4}{5}}
                  \|F_2\|_{L_t^2L^{\frac{6}{5}}}^{\frac{1}{5}},\\
 \|S_0*F_{\rm p}\|_{L_t^5(J_T;L^{10})}
 &\lesssim T^{-\frac{1}{5}}\|\nabla F_{\rm p}\|_{L_t^{\frac{4}{3}}L^1}^{\frac{4}{5}}
                  \|F_{\rm p}\|_{L_t^1L^2}^{\frac{1}{5}}.
\end{aligned}
\end{equation}
The source-time integrals cost $T^{\frac{1}{2}}$ and $T^{\frac{1}{5}}$, respectively, and the output time norm costs $T^{\frac{1}{5}}$. We used $Y_{\frac{7}{4}}\hookrightarrow L^{\frac{6}{5}}$ in the first row. For the alternative $W^{1,\frac{10}{9}}$ hypothesis, the frequency-$N$ dispersive and energy bounds interpolate to $\|S_0(\ell)P_NF\|_{10}\lesssim\ell^{-\frac{4}{5}}N^{\frac{3}{5}}\|P_NF\|_{\frac{10}{9}}$. We sum $N^{\frac{3}{5}}$ below one and $N^{-\frac{2}{5}}$ above one after placing a derivative on $F$.

For the two odd Newton inputs we claim
\begin{equation}\label{lin:eq:separatedS}
 \|C_0(\ell)G_0B\|_{10}\lesssim\ell^{-\frac{4}{5}}\|B\|_{Y_4},\qquad
 \|S_0(\ell)G_0\Gamma\|_{10}\lesssim\ell^{-\frac{7}{10}}\|\Gamma\|_{Y_4}
 \quad(\ell\ge1).
\end{equation}
The two kernels in \eqref{lin:eq:newtonkernels} have different behavior on the light cone. For the cosine we use its pointwise bound $(4\pi\ell)^{-1}$ and the unitarity bound $\|C_0(\ell)G_0B\|_2\le\|G_0B\|_2\lesssim\|B\|_{Y_4}$ from \eqref{lin:eq:weightedNewton}. Interpolation proves the first estimate. For the sine write $k_\ell(x)=\frac{1}{4\pi}\min(1,\frac{\ell}{|x|})$. Since $\int\Gamma=0$, we have
\begin{equation*}
 \|k_\ell*\Gamma\|_{10}
 \le\int\!\|k_\ell(\cdot-y)-k_\ell\|_{10}|\Gamma(y)|\,dy
 \le\|\nabla k_\ell\|_{10}\,\||x|\Gamma\|_1
 \lesssim \ell^{-\frac{7}{10}}\|\Gamma\|_{Y_4}.
\end{equation*}
Here $\nabla k_\ell=0$ inside the cone and $|\nabla k_\ell|=\frac{\ell}{4\pi r^2}$ outside. On $J_T$ the output time norm costs $T^{\frac{1}{5}}$. For the first Newton source the input time integral costs at most $T^{\frac{1}{2}}$ more. Thus \eqref{lin:eq:separatedS} supplies rates $T^{-\frac{1}{10}}$ for the first source, $T^{-\frac{1}{2}}$ for the velocity measure and $T^{-\frac{3}{5}}$ for the initial position.

\subsubsection{Transfer to the projected flow}

We first transfer the two observation bounds by the cutoff commutator
\eqref{lin:eq:commutator}, then recover the Strichartz bound from
\eqref{lin:eq:recovery}. The selected negative-mode initial data in \eqref{lin:eq:selecteddata} still
have Newton densities in $Y_4$ with norms $O(A)$, because the
eigenfunctions and their derivatives are exponentially localized.
Formula \eqref{lin:eq:negative} shows that these data are independent of $H$
once $H$ exceeds the source support. Projection of the free
velocity inputs adds only the Newton terms
\begin{equation*}
 P_0F_1-F_1=-G_0\nu K^*F_1,\qquad
 P_0y_0-y_0=-G_0\nu K^*y_0 .
\end{equation*}
These are an early $L^1Y_4$ density and an initial $Y_4$ atom,
of combined norm $O(A)$. By \eqref{lin:eq:freebound},
\eqref{lin:eq:tailtable} at $s=\frac{7}{4}$, \eqref{lin:eq:newton-time-atoms}, and \eqref{lin:eq:ppobs},
\begin{equation}\label{lin:eq:free-late-block}
 \|\cZ_f\|_{L_t^2\mathbb X}\lesssim A,\qquad
 \|\cZ_f\|_{L_t^2([cL,CL];\mathbb X)}\lesssim_{c,C} AL^{-\frac{3}{4}}
 \quad (L\ge T,\ \frac{1}{2}<c<C).
\end{equation}
We fix $b>\frac{1}{2}$ and take $\chi\in C_c^\infty((\frac{1}{2},\infty))$ with
$\chi=1$ on $[b,2b]$. We write $\chi_L(t)=\chi\bigl(\frac{t}{L}\bigr)$. By
\eqref{lin:eq:obsidentity}, \eqref{lin:eq:commutator} and
\eqref{lin:eq:free-late-block},
\begin{align*}
 \|\chi_L(x,G_0y)\|_{L_t^2\mathbb Y}
 &\lesssim \|\chi_L\cZ_f\|_{L_t^2\mathbb X}
          +L^{-1}\|\cZ_f\|_{L_t^2\mathbb X}\\
 &\lesssim_b A(L^{-\frac{3}{4}}+L^{-1}).
\end{align*}
We sum the squares over $L=2^jT$, $j\ge0$, to obtain \eqref{lin:eq:oldobs}.

For the Strichartz tail we use the exact reconstruction
\eqref{lin:eq:recovery}. The global Green estimate, including its weighted
gradient observation, implies
\begin{equation}\label{lin:eq:feedbacksources}
 \|Vx\|_{L_t^2L^{\frac{6}{5}}}+\|\nabla(Vx)\|_{L_t^2L^1}
                 +\|\nu c\|_{L_t^2Y_4}\lesssim A.
\end{equation}
Indeed, for $\rho=\langle x\rangle$,
\begin{align*}
 \|Vx\|_{\frac{6}{5}}&\le\|\rho^2 V\|_3\|x\|_{X_2},\\
 \|\nabla(Vx)\|_1&\le\|\rho^2\nabla V\|_2\|x\|_{X_2}
                      +\|\rho V\|_2\|\nabla x\|_{X_1}.
\end{align*}
The coefficients are finite by \eqref{lin:eq:potential}, and
$c=\nu^*G_0y$ bounds the last term in \eqref{lin:eq:feedbacksources}.

We split the feedback smoothly between $\frac{3T}{4}$ and $\frac{7T}{8}$. For its early
part, we apply \eqref{lin:eq:ordinary-separated} and \eqref{lin:eq:separatedS} and
\eqref{lin:eq:feedbacksources} to obtain $AT^{-\frac{1}{10}}$. Its late part
has $L_t^2Y_1\times L_t^2Y_2$ norm $O(AT^{-\frac{3}{4}})$ by
\eqref{lin:eq:oldobs}. The global free estimate then applies. The external sources and
projected initial data obey the same bound,
including the Newton correction to $F_1$ and $y_0$ just estimated.
We repeat on later dyadic scales and sum fifth powers
to obtain \eqref{lin:eq:oldS}. This proves Theorem~\ref{lin:thm:old}.
\space
\section{Radiation controls and nonlinear response}\label{merge:controls}
\subsection{Radiation packets and quadratic controls}\label{cat:start}
We place one four-dimensional radiation impulse at each interior node and estimate its integral against the piecewise affine tests of Section~\ref{main:outline}.
\subsubsection{The radiation family and reference moments}\label{cat:sec:estimates}
We construct potential waves whose quadratic impulses control the
velocity and position of the stationary core. We use the profile and
chart of Proposition~\ref{geo:ass:stationary}, the spectrally localized coordinate tests
from Section~\ref{lin:start}, and the norms
\eqref{lin:eq:green-norm}--\eqref{lin:eq:observations}.
The finite packet family is fixed before $d$.

We define the fixed-background force densities and impulses by
\begin{equation}\label{cat:eq:force-density}
 w_i=\partial_{\xi_i}(5q_\xi^4)|_0=20Q^3K_i,\qquad
 F_i(v(t))=\frac{1}{2}\int w_i(x)v(t,x)^2\dd x,\qquad
 J_i(v)=\int_\R F_i(v(t))\dd t.
\end{equation}
In particular $w_i=O(\rho^{-8})$.

\begin{theorem}[Compatible finite catalogue]\label{cat:thm:catalogue}
There exist at most eight real odd solutions $v_\nu$ of
$v_{\nu,tt}+Lv_\nu=0$ with pairwise disjoint compact symmetric temporal
bands, all separated from zero, with the following properties.
\begin{enumerate}
\renewcommand{\labelenumi}{(\roman{enumi})}
\item Each wave scatters in energy in both directions, belongs to
$C_b(\R;H^2)$, and is orthogonal, in both components, to every zero
and negative mode. Moreover
\begin{equation*}
 \sup_t\bigl(\|v_\nu(t)\|_2+\|D^{-1}v_{\nu,t}(t)\|_2\bigr)<\infty.
\end{equation*}
\item Its global $\cR$ norm is finite, and for each fixed $k$,
\begin{equation*}
 \la t\ra\partial_t^kv_\nu\in L^2X_2,\quad
 \la t\ra\nabla v_\nu\in L^2X_2,\quad
 \la t\ra G_0v_{\nu,t}\in L^2X_4.
\end{equation*}
For $T\ge1$,
\begin{equation*}
 \LE_{|t|>T}(v_\nu)+\OO_{|t|>T}(v_{\nu,t})\lesssim T^{-1},
 \qquad S_{|t|>T}(v_\nu)\lesssim T^{-\frac{3}{5}}.
\end{equation*}
These bounds also hold for every time translate after recentering the tail.
\item Zero lies in the interior of $\conv\{J(v_\nu)\}$. On a
neighborhood of $\mathcal C=[-3,3]^4$ there are smooth positive $c_\nu(z)$
such that, with $\cV_z=\sum_\nu c_\nu(z)v_\nu$,
\begin{equation}\label{cat:eq:groupimpulse}
 J(\cV_z)=z,\qquad
 \sup_{z\in\mathcal C}\int_\R |t|\,|F(\cV_z(t))|\dd t<\infty.
\end{equation}
\item $\Phi v_\nu(t)=0$ for all $t$. Thus $(v_\nu,v_{\nu,t})$ itself
is the homogeneous gauged Green solution, including on the half-line.
\end{enumerate}
\end{theorem}

For $\eps=d^3$, $T_0=d^{-5}$, we set
\begin{equation}\label{main:scales}
 T_j=4^jT_0,\quad \eta_j=\eps2^{-j},\quad T_j\eta_j^2=d,
 \qquad I_0=[0,4T_0],\quad I_j=[T_j,4T_j]\ (j\ge1).
\end{equation}
For $a_j\in\mathcal C$, $j\ge1$, the radiation train is
\begin{equation}\label{cat:eq:train}
 v_a(t)=\sum_{j\ge1}\eta_j\mathcal V_{a_j}(t-T_j).
\end{equation}
The series converges in uniform energy and in the observation norms since $\sum_{j\ge1}\eta_j=\eps$. For the tests in Section~\ref{main:outline} we set
\begin{equation}\label{cat:eq:moment-map}
 \mathcal M(a)_j=\eta_j^{-2}\int_0^\infty\varphi_j(t)F(v_a(t))\dd t .
\end{equation}
\begin{theorem}[Uniform reference impulses]\label{cat:thm:reference}
Uniformly in the controls in $\mathcal C$,
\begin{equation}\label{cat:eq:reference-moments}
 \sup_{j\ge1}|\mathcal M(a)_j-a_j|\lesssim T_0^{-1}.
\end{equation}
\end{theorem}

\begin{figure}[htbp]
\centering
\begin{tikzpicture}[>=Stealth,font=\small,x=1.8cm,y=1cm]
 \draw[->] (-.4,0)--(5.55,0) node[right] {$t$};
 \draw[thick] (0,0)--(1,1.1)--(5,0);
 \draw[dashed] (1,.04)--(1,1.1);
 \foreach \x/\lab in {0/{s_{j-1}},1/{s_j},5/{s_{j+1}}}{
   \draw (\x,.05)--(\x,-.05) node[below] {$\lab$};}
 \node[above] at (1,1.1) {$\varphi_j(s_j)=1$};
\end{tikzpicture}
\caption{The test for the impulse at $s_j$ vanishes at both neighboring arrivals. The lengths shown have ratio $1:4$, as on $j\ge2$.}
\label{cat:impulse-diagram}
\end{figure}
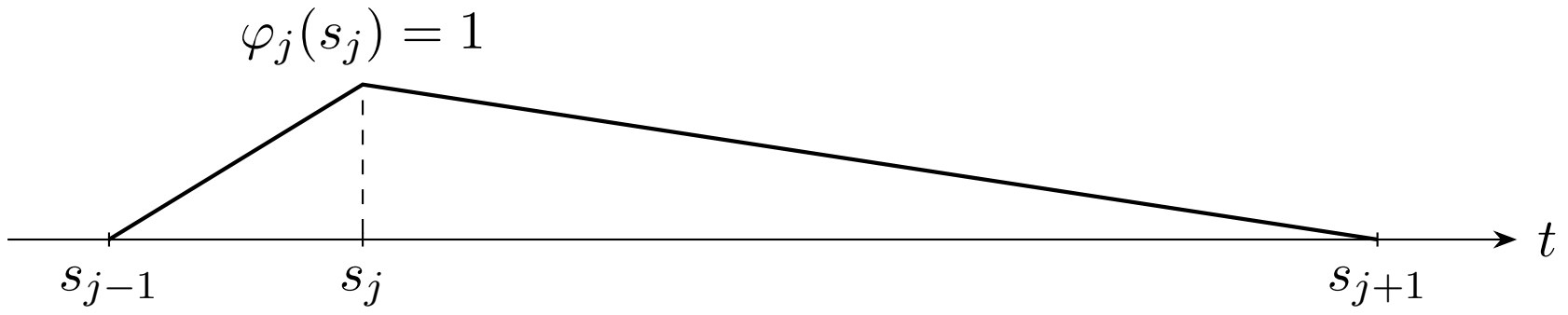
We prove Theorem~\ref{cat:thm:catalogue} in
Sections~\ref{cat:sec:ray}--\ref{cat:sec:positive} and
Theorem~\ref{cat:thm:reference} in Section~\ref{cat:sec:reference}.

\subsubsection{The ray cone and its wave realization}\label{cat:sec:ray}

We first choose lines whose weighted ray integrals surround the origin.
We then approximate these integrals by the impulses \eqref{cat:eq:force-density} of odd potential waves.
For an oriented line $\ell=(b,\omega)$, $b\perp\omega\in\mathbb S^2$,
we write $Xf(\ell)=\int_\R f(b+s\omega)\dd s$. With our normalization of line measure, Fuglede's integral formula \cite{Fuglede} takes the following form. Expanding the square and
integrating in polar coordinates, we obtain
\begin{equation}\label{cat:eq:ray-energy}
 \mathcal B(V):=\int_{\mathbb S^2}\int_{\omega^\perp}(XV)^2\dd b\dd\omega
 =2\iint_{\R^3\times\R^3}\frac{V(x)V(y)}{|x-y|^2}\dd x\dd y.
\end{equation}
The integrals in \eqref{cat:eq:ray-energy} and their first variations
converge absolutely. We verify invariance first under similarities, by
translation, rotation and scaling. For inversion $I(x)=\frac{x}{|x|^2}$,
substitution in the second integral of \eqref{cat:eq:ray-energy} uses
\begin{equation*}
 V^I(x)=|x|^{-4}V(Ix),\qquad
 |Iu-Iv|=\frac{|u-v|}{|u||v|},\qquad d(Iu)=|u|^{-6}du .
\end{equation*}
The powers at each variable cancel since $4+2-6=0$. These transformations
generate the conformal group. Differentiating along $V_\xi=5q_\xi^4$ at
$\xi=0$, we obtain
\begin{equation}\label{cat:eq:balance}
 \int (XV)(Xw_i)\dd\ell=0
\end{equation}
for every chart direction. Suppose $c\ne0$ and $Xw_c\ge0$, where
$w_c=20Q^3Kc$. The continuous integrand $(XV)(Xw_c)$ in
\eqref{cat:eq:balance} is nonnegative, so it vanishes on every line.
Since $V=5Q^4$, we also have
\begin{equation*}
 XV(\ell)=0\quad\Longrightarrow\quad Q|_\ell=0
 \quad\Longrightarrow\quad Xw_c(\ell)=0 .
\end{equation*}
Thus $Xw_c=0$ on all lines. By the Fourier slice identity for the line transform \cite{Helgason}, $w_c=0$.
Thus $Kc=0$ on the dense set $\{Q\ne0\}$, hence everywhere
by continuity. As $K^*K=I$, this contradicts $c\ne0$.
Thus no closed half-space through zero contains all ray columns. By the classical Steinitz theorem \cite[Proposition~1.1]{IN}, we may choose at most eight incident lines whose four-dimensional impulse vectors satisfy
\begin{equation}\label{cat:eq:raycone}
 j_\nu=\frac{1}{4}(Xw_i(\ell_\nu))_{i=1}^4,\qquad
 \delta_0\overline B_{\R^4}\subset\conv\{j_1,\ldots,j_M\},\qquad M\le8 .
\end{equation}
The margin persists under small changes of these finitely many lines.
In particular, we choose directions for which
$\mathscr R\omega\ne\omega$ and $\mathscr R\omega\ne-\omega$.

\paragraph{The ray limit and its normalization.}
We choose a real Schwartz envelope $a_\delta$ with smooth compact Fourier
support and $\|a_\delta\|_2=1$, concentrating at zero as $\delta\downarrow0$.
For a selected line, we define
\begin{equation*}
 u_h^0=\operatorname{Re}\left(h e^{-itD}
                      [e^{i\frac{\omega\cdot x}{h}}a_\delta(x-b)]\right),
 \qquad u_h=\frac{u_h^0-u_h^0\circ\mathscr R}{\sqrt2}.
\end{equation*}
This real odd free wave has energy $O(1)$ and temporal Fourier support in fixed-width
bands $B_h$ about $\pm h^{-1}$. With the envelope and line fixed,
\begin{equation}\label{cat:eq:freeweighted}
 \|\la t\ra u_h\|_{L^2X_2}\lesssim h.
\end{equation}
For unrestricted $g\in L^2$, take the free wave $U=iD^{-1}e^{-itD}g$, so that $\partial_t U=e^{-itD}g$ and its energy norm is $\sqrt{2}\|g\|_2$. Applying the homogeneous Morawetz identity \eqref{lin:eq:energy-morawetz} in both time directions, to real and imaginary parts, gives $\|e^{-itD}g\|_{L_t^2X_1}\lesssim\|g\|_2$. The homogeneous identity holds in every angular sector. If $A_j=\frac{D_j}{D}$, then
\begin{equation*}
 t e^{-itD}g=\sum_j A_j
          [x_je^{-itD}g-e^{-itD}(x_jg)].
\end{equation*}
On the packet supports, we insert an annular cutoff $\chi(hD)$ in $A_j$.
Its convolution kernel has uniformly bounded weighted second $L^1$
moment, hence acts boundedly on $X_2$. We use
$|x|\rho^{-2}\le\rho^{-1}$ and apply the free observation to $g$ and
$x_jg$. This proves \eqref{cat:eq:freeweighted}.

On each bounded time interval, we expand the phase and obtain
\begin{align*}
 \left|\frac{\omega}{h}+\zeta\right|
 &=h^{-1}+\omega\cdot\zeta+O_\delta(h),\\
 h^{-1}u_h^0(t,x)
 &=a_\delta(x-b-t\omega)\cos\!\left(\frac{\omega\cdot x-t}{h}\right)+o_{L_x^2}(1).
\end{align*}
For every continuous $w=O(\rho^{-4})$, averaging $\cos^2$ and then
using the tail in \eqref{cat:eq:freeweighted} yields
\begin{equation}\label{cat:eq:raylimit}
 h^{-2}J_w(u_h^0)\longrightarrow
 \frac{1}{4}\int |a_\delta(y)|^2Xw(b+y,\omega)\dd y
 \longrightarrow_{\delta\downarrow0}\frac{1}{4}Xw(b,\omega),
\end{equation}
where $J_w(v)=\frac{1}{2}\iint wv^2$.
The omitted time tail is $O(h^2T^{-2})$ before letting $T\to\infty$.
For the actual force densities $w_i$, $w_i\circ\mathscr R=w_i$.
The two diagonal contributions of $u_h$ agree, and the cross terms
vanish in the normalized limit because their spatial phases have nonzero
gradients $\omega\pm\mathscr R\omega$ on bounded time intervals, and
the tails are controlled by \eqref{cat:eq:freeweighted}. Thus
\eqref{cat:eq:raylimit} holds for $u_h$.

The factor $\frac{1}{4}$ in \eqref{cat:eq:raylimit} is the product of
the force normalization in \eqref{cat:eq:force-density} and the mean of $\cos^2$.

\paragraph{The retarded packet.}
On functions with temporal Fourier support in $B_h$, we solve the incoming equation
\begin{equation}\label{cat:eq:duhamel}
 v_h=u_h+\mathcal Dv_h,\qquad
 (\mathcal Df)(t)=\int_{-\infty}^t S_0(t-s)Vf(s)\dd s .
\end{equation}
With $\widehat f(\tau)=\int e^{it\tau}f(t)\dd t$, the retarded multiplier is $G_{\tau+i0}V$. Indeed, exponential damping of the retarded kernel replaces $i0$ by $i\delta$, and the weighted kernels converge in $Y_2\to X_2$ on the compact bands. On $B_h$, \eqref{lin:eq:highG} implies
\begin{equation*}
 \|G_{\tau+i0}V\|_{X_2\to X_2}\lesssim h,\qquad
 \|(\partial_\tau G_{\tau+i0})V\|_{X_2\to X_2}\lesssim1 .
\end{equation*}
For small $h$ we therefore set $\widehat v_h=(I-G_{\tau+i0}V)^{-1}\widehat u_h$ on $B_h$ and zero off $B_h$. The inverse preserves the odd symmetry and the conjugation relation at opposite frequencies. Differentiating the multiplier equation, we obtain
\begin{equation*}
 (I-G_{\tau+i0}V)\partial_\tau\widehat v_h
 =\partial_\tau\widehat u_h+(\partial_\tau G_{\tau+i0})V\widehat v_h .
\end{equation*}
By Plancherel and \eqref{cat:eq:freeweighted},
\begin{equation}\label{cat:eq:packet-perturbation}
 \|v_h-u_h\|_{L^2X_2}\lesssim h^2,\qquad
 \|\la t\ra v_h\|_{L^2X_2}\lesssim h.
\end{equation}
In particular $Vv_h\in L^1L^2$. The damped integrals converge to the retarded energy integral in \eqref{cat:eq:duhamel}, so this field has the prescribed incoming wave.
The decay of $V$ and Cauchy--Schwarz imply
\begin{equation}\label{cat:eq:forcing}
 \|\la t\ra Vv_h\|_{L^2Y_2}+\|Vv_h\|_{L^1L^2}\lesssim h.
\end{equation}
The energy and Strichartz estimates in \eqref{cat:eq:duhamel} now control the packet globally. Its energy distance from $u_h$ is $O(h)$ and tends to zero at $-\infty$. The future scattering state exists by integrability of $Vv_h$. For the impulse, $|w_i|\lesssim\rho^{-4}$ and the difference of squares imply
\begin{equation}\label{cat:eq:forceerrorh}
 \int(1+|t|)|F_i(v_h(t))-F_i(u_h(t))|\dd t
 \lesssim\|v_h-u_h\|_{L^2X_2}\|\la t\ra(v_h+u_h)\|_{L^2X_2}\lesssim h^3 .
\end{equation}

A temporal primitive supplies the extra regularity and the low-frequency
velocity bound at once. We choose a smooth cutoff $\theta_h$ equal to one
near $B_h$, supported in slightly larger bands of fixed width. We define
$T_h$ by the temporal Fourier multiplier $\frac{\theta_h(\tau)}{-i\tau}$.
Its convolution kernel has $L^1$ norm $O(h)$, by translation of the
two bands and integration by parts. On the packet $U=T_hv_h$ one has
\begin{equation*}
 U_t=v_h,\qquad -LU=v_{h,t},\qquad
 \|\nabla U\|_{C_bL^2}+\|v_h\|_{C_bL^2}\lesssim h .
\end{equation*}
The first estimate applies $T_h$ to $\nabla v_h$ and the second uses
$v_h=T_hv_{h,t}$. These identities hold first in spacetime
distributions by the band support. Since $v_{h,t}\in C_bL^2$ is band-limited,
$v_{h,tt}\in C_bL^2$. Thus $D^2v_h=Vv_h-v_{h,tt}\in C_bL^2$ and $v_h\in C_bH^2$. Using $v_{h,t}=-D^2U+VU$ and Sobolev, we obtain
\begin{equation*}
 \|D^{-1}v_{h,t}\|_2\lesssim
 \|\nabla U\|_2+\|V\|_{\frac{3}{2}}\|U\|_6\lesssim h .
\end{equation*}
This bound is for the unprojected packet.

\subsubsection{Time moments and positive controls}\label{cat:sec:positive}

We establish the time moments and spectral orthogonality of each wave.
We then combine waves in disjoint bands so that their impulses add exactly.
Multiplication by a smooth cutoff times $\tau^k$ preserves
$H^1_\tau(X_2)$, hence $\la t\ra\partial_t^kv_h\in L^2X_2$.
We integrate $D^2v=Vv-v_{tt}$ against $\rho^{-4}v$ to obtain
\begin{equation*}
 \|\rho^{-2}\nabla v\|_2^2
 \lesssim\|v\|_{X_2}(\|v_{tt}\|_{X_2}+\|v\|_{X_2}).
\end{equation*}
Multiplying by $\la t\ra^2$ and integrating proves
$\la t\ra\nabla v\in L^2X_2$. For the Newton observation we use the same temporal primitive:
\begin{equation}\label{cat:eq:newtonband}
 G_0v_t=-(I-G_0V)T_h v,\qquad I-G_0V:X_2\longrightarrow X_4 .
\end{equation}
Indeed $G_0v_t=-U+G_0VU$ in the energy realization since
$v_t=-D^2U+VU$, $U\in\dot H^1$ and $VU\in L^{\frac{6}{5}}$.
The mapping property follows from the weighted Newton bound and the decay of $V$. The smooth band multiplier of $T_h$ preserves $H^1_\tau(X_2)$,
hence $\la t\ra G_0v_t\in L^2X_4$. The first time moments imply the $T^{-1}$ tails for both observations.

We prove pointwise decay for the full packet. Set $g(y)=|V(y)|\sup_s\langle s\rangle |v_h(s,y)|$. By the one-dimensional Sobolev inequality in $s$ and the weighted bounds on $v_h,\partial_s v_h$,
\begin{equation*}
 \|g\|_1+\|g\|_2\lesssim
 (\|\rho^2V\|_2+\|\rho^2V\|_\infty)
 \|\langle s\rangle v_h\|_{H_s^1 X_2}\lesssim_h1 .
\end{equation*}
In the retarded Kirchhoff formula \eqref{cat:eq:duhamel}, write $r=|x-y|$ and use $\langle t\rangle\le\langle t-r\rangle+r$. Then
\begin{equation}\label{cat:eq:packet-decay}
 \langle t\rangle |v_h(t,x)-u_h(t,x)|
 \le \frac{1}{4\pi}\int_{\mathbb R^3}(1+r^{-1})g(y)\,dy
 \lesssim\|g\|_1+\|g\|_2 .
\end{equation}
The last bound uses $|x|^{-1}\in L^2(B_1)$. The free packet has Schwartz data on a fixed annulus and satisfies $\|u_h(t)\|_\infty\lesssim_h\langle t\rangle^{-1}$ by the frequency-localized wave estimate \cite[p.~972]{KT}. Interpolating \eqref{cat:eq:packet-decay} and the free bound with $\sup_t\|v_h(t)\|_2\lesssim h$, we obtain
\begin{equation*}
 \|v_h(t)\|_{10}\lesssim_h\langle t\rangle^{-\frac{4}{5}},\qquad
 S_{|t|>T}(v_h)\lesssim_h T^{-\frac{3}{5}} .
\end{equation*}

For spectral orthogonality, we choose the packet bands in $\{|\tau|>1\}$ and
put $E=\mathbf1_{(-\infty,1]}(L)$. This projection contains all gauge
tests and all nonpositive modes. On $\operatorname{ran}E$ the spectrum
lies in $[-\|V\|_\infty,1]$. For $\tau$ in the packet bands,
\begin{equation*}
 (L-\tau^2)\widehat{Ev}(\tau)=0,\qquad
 \bigl\|(L|_{\operatorname{ran}E}-\tau^2)^{-1}\bigr\|\le\frac{1}{\tau^2-1}.
\end{equation*}
The inverse is smooth in $\tau$ on both bands, so $Ev=0$ as a temporal
distribution. Hence $\Phi v=\Phi v_t=0$ and $Pv=v$ at every time.
The packet and each of its translates have zero negative amplitudes.
By Green uniqueness they are the homogeneous selected solutions.

\paragraph{Disjoint bands and positive controls.}\label{cat:sec:disjoint}

We choose the finitely many envelopes first, close enough to their rays
for the cone margin in \eqref{cat:eq:raycone}. We choose the frequencies
$h_\nu^{-1}$ successively large so that the enlarged symmetric
bands $B_\nu$ are pairwise disjoint and contained in $\{|\tau|>1\}$. We normalize the potential waves by
\begin{equation*}
 v_\nu=h_\nu^{-1}v_{\nu,h_\nu}.
\end{equation*}
By \eqref{cat:eq:raylimit} and \eqref{cat:eq:forceerrorh}, the columns
$J_\nu=J(v_\nu)$ still satisfy \eqref{cat:eq:raycone} with a smaller positive margin. Their energy constants may be large, of order $h_\nu^{-1}$,
but are fixed before $d$ and independent of all later time translations.

Multiplication by $\rho^4 w_i$ is bounded on $L_x^2$ and preserves the
temporal support of $\rho^{-2}v_\nu\in L_t^2L_x^2$. Translation and
multiplication by $t$ also preserve temporal Fourier support. By Plancherel,
for $\nu\ne\mu$, $k=0,1$, and all $s,s'\in\R$,
\begin{equation}\label{cat:eq:band-orthogonality}
 \iint t^k w_i(x)v_\nu(t-s,x)v_\mu(t-s',x)\dd x\dd t=0.
\end{equation}

Choose positive barycentric weights $\beta_\nu$ for zero in the interior of $\conv\{J_\nu\}$. Since the columns span $\mathbb R^4$, their covariance matrix $G=\sum_\nu\beta_\nu J_\nu J_\nu^T$ is positive definite. For a sufficiently large constant $L$ we take
\begin{equation*}
 c_\nu(z)=\sqrt{\beta_\nu(L+J_\nu^T G^{-1}z)},\qquad
 J(\cV_z)=\sum_\nu c_\nu(z)^2J_\nu=z .
\end{equation*}
The radicands are positive near $\mathcal C$. The equality follows from \eqref{cat:eq:band-orthogonality} and $\sum_\nu\beta_\nu J_\nu=0$.
For the first moment, part~(ii) and the coefficient bounds imply
\begin{equation*}
 \int_\R |t|\,|F(\cV_z(t))|\dd t
 \lesssim \|t\cV_z\|_{L^2X_2}\|\cV_z\|_{L^2X_2}\lesssim1
 \qquad(z\in\mathcal C).
\end{equation*}
The same bounds hold after one control derivative by bilinearity and
the coefficient bounds. This proves Theorem~\ref{cat:thm:catalogue}.

\subsubsection{Testing the radiation at the nodes}\label{cat:sec:reference}
Write $v_k=\eta_k\mathcal V_{a_k}(t-T_k)$ for $k\ge1$ and $v_0=0$. On $I_j$ the only arrivals are at its endpoints. Every other packet is at distance $\gtrsim T_j$ for $k<j$ and $\gtrsim T_k$ for $k>j+1$. By Theorem~\ref{cat:thm:catalogue}(ii),
\begin{equation}\label{cat:eq:off}
 \LE_{I_j}(v_a-v_j-v_{j+1})
 \lesssim T_j^{-1}\sum_{1\le k<j}\eta_k+\sum_{k>j+1}T_k^{-1}\eta_k
 \lesssim T_0^{-1}\eta_j .
\end{equation}
On each half of $I_j$ the more distant endpoint packet also has $\LE$ norm $O(T_j^{-1}\eta_j)$. Expanding the square in \eqref{cat:eq:force-density} and applying Cauchy--Schwarz on each half, we find
\begin{equation}\label{cat:eq:crossbound}
 \|F(v_a)-F(v_j)-F(v_{j+1})\|_{L^1(I_j)}\lesssim T_0^{-1}\eta_j^2 .
\end{equation}
The zero extension of $\varphi_j$ has Lipschitz constant $\lesssim T_j^{-1}$. On $I_{j-1}\cup I_j$ the diagonal terms in \eqref{cat:eq:crossbound} have indices $k=j-1,j,j+1$. For the central packet we use $|\varphi_j(t)-1|\lesssim |t-T_j|/T_j$ on the entire line. For either neighbor use $|\varphi_j(t)|\lesssim |t-T_k|/T_j$ on its interval. The first absolute moment in \eqref{cat:eq:groupimpulse} then bounds the total error by
\begin{equation*}
 C T_0^{-1}(\eta_{j-1}^2+\eta_j^2)
 +\frac{C}{T_j}\sum_{k=j-1}^{j+1}\eta_k^2\lesssim T_0^{-1}\eta_j^2 .
\end{equation*}
This proves \eqref{cat:eq:reference-moments}, including $j=1$.

\paragraph{Reference observations and finite energy.}
The global bounds follow from $\sum_{j\ge1}\eta_j=\eps$. The two endpoint packets have block norms $O(\eta_j)$. For all others we use \eqref{cat:eq:off} for the first-moment observations and the same sum with $T^{-\frac{3}{5}}$ for the Strichartz norm.
Consequently the reference field $r=v_a$ satisfies
\begin{equation}\label{cat:eq:blockobservations}
 \LE_{I_j}(r)+S_{I_j}(r)+\OO_{I_j}(v_{a,t})\lesssim\eta_j.
\end{equation}
Its uniform energy bound is $C\eps$.

Summing the geometric block bounds, we may use \eqref{cat:eq:blockobservations}
also on $J_j=[\frac{T_j}{3},\infty)$ for $j\ge1$, since $J_j$ is covered by
the blocks $I_k$, $k\ge j-1$, and $\sum_{k\ge j-1}\eta_k=4\eta_j$.
The energy bound on each such tail is $C\eps$.

\subsection{The prescribed-path response and force estimates}\label{non:start}
Theorem~\ref{non:thm:main}
constructs the field for a prescribed path, and
Theorem~\ref{cat:thm:blocks} bounds the force error in each nodal test. Section~\ref{end:start}
uses these tests to choose the path.

\subsubsection{The prescribed-path problem}\label{non:sec:field}
For fixed $\xi$ and $v_a$, we solve a projected equation and identify the
force that will determine $\xi$. We retain $D,G_0,L,V,K,\Phi,P,P_0$ and the norms
from Section~\ref{lin:start}, and set $\cD=P\Phi^*$.
For the product estimates we use $M(f)=\|f\|_{L_t^\infty L_x^2}$.

\paragraph{The stationary chart.}
We use the chart \eqref{geo:eq:chart}, with $m=4$. For small $\xi$,
its derivatives satisfy the decay bounds \eqref{geo:eq:decay} and
\begin{align*}
 R_i=\partial_{\xi_i}R&=O(|\xi|\rho^{-2}),&
 R_{ij}=\partial_{\xi_j}R_i&=O(\rho^{-2}),\\
 \|R_i\|_{H^1}+\|D^2R_i\|_{Y_4}&\lesssim |\xi|,&
 \|R_{ij}\|_{H^1}+\|D^2R_{ij}\|_{Y_4}&\lesssim1.
\end{align*}

The density bounds follow from \eqref{geo:eq:density}.
We use $K^*R=0$ for the stationary chart and the localized gauge $\Phi$
for energy-class radiation.

\paragraph{The reference radiation.}
The reference field $v_a$ solves $v_{a,tt}+Lv_a=0$ and is orthogonal to
all $K_i$ and $\Psi_l$, in both its position and velocity components.
Since its packets lie in $\{L>1\}$, we have $Pv_a=v_a$ and $Pv_{a,t}=v_{a,t}$
by \eqref{lin:eq:spectralgauge}. Summing their bounds in Theorem~\ref{cat:thm:catalogue}, we obtain
\begin{equation}\label{non:eq:radminus}
 \|(v_a,v_{a,t})\|_{\cR}+\|v_a(0)\|_2+
 \|D^{-1}v_{a,t}(0)\|_2\lesssim\eps .
\end{equation}
All fixed derivatives converge on compact time intervals since the
catalogue is finite and the train coefficients are summable.

Let $0<d\ll1$, $\eps=d^3$, $0<\mathsf T\le\infty$, and prescribe
$\xi\in H^1_{\rm loc}([0,\mathsf T);\R^4)$ with
\begin{equation}\label{non:eq:path}
 \xi(0)=0,\qquad \|\xi\|_\infty\lesssim d,\qquad
 \|\xi'\|_2\lesssim\sqrt d\,\eps .
\end{equation}
The block bounds $\|\xi'\|_{L^2(I_j)}\lesssim\sqrt d\,\eta_j$ imply
the global derivative bound since $\eta_j=\eps2^{-j}$.

\begin{theorem}[Prescribed-path nonlinear closure]\label{non:thm:main}
For the chart \eqref{geo:eq:chart}, the train \eqref{cat:eq:train} and a path
satisfying \eqref{non:eq:path}, the system
\eqref{non:eq:field}--\eqref{non:eq:initial} has a unique small correction
in its contraction class on $[0,\mathsf T]$. Uniformly in $\mathsf T$,
\begin{gather}
 \|(w,\bar v)\|_{\cR}\lesssim\eps,\qquad
 \|(w-v_a,\bar v-v_{a,t})\|_{\cR}\lesssim d\eps,\qquad |z|\lesssim d\eps,
 \label{non:eq:mainbasic}\\
 \|H\|_{L^1}\lesssim\eps^2.
 \label{non:eq:maincharges}
\end{gather}
The full field satisfies $\|w\|_{L_t^{\frac{20}{3}}L_x^5}+\|w\|_{L_t^\infty\dot H_x^{\frac{3}{4}}}\lesssim\eps$.
There is an explicit two-component dipole response $h$, defined in
\eqref{non:eq:halo}, whose driving measure has total variation
$\delta\lesssim\eps^2$, such that
\begin{equation}\label{non:eq:mainmass}
 G=w-h,\qquad M(G)\lesssim\eps.
\end{equation}
The exact Green sources admit the enhanced decomposition
\begin{align}
 F_2&=P^*[(V_\xi-V)w+B_{q_\xi}(w)],\notag
\\
 F_1&=F_p+G_0\Gamma_{\rm c},\qquad
 F_p=w^5,\quad
 \Gamma_{\rm c}=D^2\Phi^*(H-K^*w^5),
 \label{non:eq:enhanced}
\end{align}
Here $V_\xi=5q_\xi^4$. We separate the pure quintic by writing
\begin{equation*}
 B_q(w)=(q+w)^5-q^5-5q^4w-w^5=10q^3w^2+10q^2w^3+5qw^4 .
\end{equation*}
The source terms satisfy
\begin{align}
 \|\beta\|_{L^2Y_4}&\lesssim d\eps,\notag
\\
 \|F_2\|_{L^2Y_2}+\|\nabla F_2\|_{L^2L^1}
 &\lesssim d\eps+\eps^2,\notag
\\
 \|F_p\|_{L^1L^2}+\|F_p\|_{L^{\frac{4}{3}}W^{1,1}}
 &\lesssim\eps^5,\notag
\\
 \|\Gamma_{\rm c}\|_{L^1Y_4}&\lesssim\eps^2.
 \label{non:eq:mainenhanced}
\end{align}
We also have $\|F_2\|_{L^2W^{1,\frac{10}{9}}}\lesssim d\eps+\eps^2$ and $\|F_p\|_{L^{\frac{5}{4}}W^{1,\frac{10}{9}}}\lesssim\eps^5$.
The stronger $Y_2$ bound implies the $Y_{\frac{7}{4}}$ bound required in
Theorem~\ref{lin:thm:old}. All sources retain the exact full-kernel compatibilities.
\end{theorem}

\paragraph{The field equations.}

We write $\cQ_i=\partial_{\xi_i}q_\xi=K_i+R_i$ and set
\begin{equation*}
 M_\xi=I+\Phi DR(\xi),\quad
 p=-M_\xi^{-*}DR(\xi)^*\bar v,\quad
 \widetilde v=\bar v+\cD p.
\end{equation*}
The matrix $M_\xi$ is invertible for small $d$. We define
\begin{equation}\label{non:eq:force-definition}
\begin{aligned}
 N_q(w)&=(q+w)^5-q^5-5q^4w=B_q(w)+w^5,\\
 \cA&=(V_\xi-V)w+N_{q_\xi}(w),\\
 H_i&=\xi_j'\la R_{ij}(\xi),\widetilde v\ra
       +\la\cQ_i(\xi),N_{q_\xi}(w)\ra.
\end{aligned}
\end{equation}
Repeated chart indices are summed. The prescribed-path equations are
\begin{equation}\label{non:eq:field}
 \begin{split}
 w_t&=\bar v+\cD p-PDR(\xi)\xi',\\
 \bar v_t&=-PLw+PP^*\cA+\cD H,\qquad \Phi w=\Phi\bar v=0.
 \end{split}
\end{equation}
By the stationary equation,
\begin{equation}\label{non:eq:chartids}
 LR_i=(V_\xi-V)\cQ_i,\qquad K^*R_i=K^*R_{ij}=0.
\end{equation}
Consequently,
\begin{equation}\label{non:eq:algebra}
 DR^*\cD=M_\xi^*-I,\qquad p=-DR^*\widetilde v,
 \qquad PP^*=P-\cD K^*.
\end{equation}
These identities use the $L^2$ chart normalization and the localized gauge
in different roles.

The exact frozen Green sources are
\begin{equation}\label{non:eq:basicsources}
 \begin{split}
 \beta&=D^2[PDR(\xi)\xi'-\cD p],\\
 F_2&=P^*[(V_\xi-V)w+B_{q_\xi}(w)],\\
 F_1&=P^*w^5+\Phi^*H.
 \end{split}
\end{equation}
Indeed \eqref{non:eq:field} becomes
$w_t=\bar v-G_0\beta$,
$\bar v_t=-PLw+P(F_2+F_1)$, and
\begin{equation}\label{non:eq:compat}
 \Phi G_0\beta=0,\qquad K^*F_2=0,\qquad K^*F_1=H.
\end{equation}

Since $DR(0)=0$, every input to the prescribed-path construction obeys
\begin{equation}\label{non:eq:p0}p(0)=0.\end{equation}
The mixed Green operator selects the negative amplitudes $z_l$, giving
\begin{equation}\label{non:eq:initial}
 (w,\bar v)(0)=(v_a,v_{a,t})(0)
       +\sum_l z_l(\Psi_l,\kappa_l\Psi_l).
\end{equation}

The contraction controls the field and charges. Removing the two
horizontal first moments from its Newton density will close the regular
position mass and the full quintic estimate together.

\subsubsection{Contraction and exact identities}\label{non:sec:contraction}

We estimate the Green sources in $\cR$, contract around $v_a$, then
differentiate the chart constraint. We set $\mathcal N=\|(w,\bar v)\|_{\cR}$.
Since $R_i=G_0D^2R_i$, Newton duality implies
\begin{equation*}
 |\langle R_i,y\rangle|=|\langle D^2R_i,G_0y\rangle|
 \le\|D^2R_i\|_{Y_4}\|G_0y\|_{X_4}\lesssim d\|G_0y\|_{X_4}.
\end{equation*}
We apply the same identity to $R_{ij}$ and use the formula for $p$ to obtain
\begin{equation}\label{non:eq:pbound}
 |\langle R_{ij},y\rangle|\lesssim\|G_0y\|_{X_4},\qquad
 \|p\|_2+\|p\|_\infty\lesssim d\mathcal N .
\end{equation}
\begin{lemma}[Weighted polynomial estimates]\label{non:lem:poly}
For $p=2,3,4$,
\begin{equation*}
 \|q^{5-p}w^p\|_{L^2Y_2}
 \lesssim \LE(w)^{\frac{4-p}{3}}
 S(w)^{\frac{5(p-1)}{6}}E(w)^{\frac{p-1}{2}},
 \qquad \|w^5\|_{L^1L^2}=S(w)^5.
\end{equation*}
The polynomial force satisfies
\begin{equation}\label{non:eq:forces}
 \|\langle\cQ,N_q(w)\rangle\|_{L^1}
 \lesssim\sum_{p=2}^5
 \LE(w)^{\frac{8-p}{3}}S(w)^{\frac{5(p-2)}{6}}E(w)^{\frac{p-2}{2}}.
\end{equation}
Each summand bounds the corresponding absolute spacetime integral
$\iint|\cQ_iq^{5-p}w^p|$.
\end{lemma}
\begin{proof}
By H\"older and Sobolev,
$\||w|^4\|_{L^2_{t,x}}\lesssim S(w)^{\frac{5}{2}}E(w)^{\frac{3}{2}}$.
Interpolating this bound with $\|\rho^{-2}w\|_{L^2_{t,x}}\le\LE(w)$, we obtain, for $1\le k\le4$,
\begin{equation}\label{non:eq:weighted-interpolation}
 \|\rho^{-\frac{2(4-k)}{3}}|w|^k\|_{L^2_{t,x}}
 \lesssim\LE(w)^{\frac{4-k}{3}}S(w)^{\frac{5(k-1)}{6}}E(w)^{\frac{k-1}{2}} .
\end{equation}
Since $|q|+|\cQ_i|\lesssim\rho^{-2}$, \eqref{non:eq:weighted-interpolation} applies to $\rho^2q^{5-p}w^p$ with $k=p$ for $2\le p\le4$. For the force we take $k=p-1$ and use Cauchy--Schwarz,
\begin{equation*}
 J_p:=\iint|\cQ_iq^{5-p}w^p|
 \lesssim\|\rho^{-2}w\|_{L^2_{t,x}}
 \|\rho^{-\frac{2(5-p)}{3}}|w|^{p-1}\|_{L^2_{t,x}},\qquad 2\le p\le5 .
\end{equation*}
This proves \eqref{non:eq:forces}. On any time interval, H\"older applied directly to the force integrands also yields
\begin{equation}\label{non:eq:force-interpolation}
 J_p\le J_2^{\frac{5-p}{3}}J_5^{\frac{p-2}{3}},\qquad p=3,4 .
\end{equation}
\end{proof}

Substituting the polynomial estimates and \eqref{non:eq:pbound} in
\eqref{non:eq:basicsources}, we obtain
\begin{equation}\label{non:eq:basic-source-bounds}
\begin{aligned}
 \|\beta\|_{L^2Y_4}&\lesssim d(\|\xi'\|_2+\mathcal N),\\
 \|F_2\|_{L^2Y_2}&\lesssim d\mathcal N+
                    \mathcal N^2+\mathcal N^3+\mathcal N^4,\\
 \|H\|_1+\|F_1\|_{L^1L^2}
 &\lesssim\|\xi'\|_2\mathcal N+
                  \mathcal N^2+\mathcal N^3+\mathcal N^4+\mathcal N^5.
\end{aligned}
\end{equation}
The geometric term of $H$ uses \eqref{non:eq:pbound} and Cauchy--Schwarz in
time. The radiation $(v_a,v_{a,t})$ is a homogeneous gauged solution by
its spectral support. We apply the mixed Green map with zero
complementary initial correction on the closed ball
\begin{equation*}
 \|(w-v_a,\bar v-v_{a,t})\|_{\cR}\le C_1d\eps
\end{equation*}
in the gauged energy-continuous class.
On this ball, \eqref{non:eq:basic-source-bounds} bounds the inhomogeneous
correction by $O(d\eps+\eps^2)$. We choose $C_1$ large, then $d$ small. By
polarization, the Lipschitz constant is at most
\begin{equation*}
 C\left(d+\eps+\|\xi'\|_2+
           \sum_{j=1}^4(\mathcal N_1^j+\mathcal N_2^j)\right)<\frac{1}{2}.
\end{equation*}
This proves the contraction and \eqref{non:eq:mainbasic}. By the mixed Green formulas we also have $|z|\lesssim d\eps$. Equations
\eqref{non:eq:forces} and \eqref{non:eq:path} now prove
\eqref{non:eq:maincharges}. Replacing $\bar v$ by $\widetilde v$ preserves
the stated field bounds, by \eqref{non:eq:pbound}.

The $H^{-1}$--$H^1$ product rule applies to $M_\xi^*p=-DR^*\bar v$, since $R_i(\xi)\in H^1_{t,\rm loc}H_x^1$ and $\bar v_t\in L^1_{t,\rm loc}H_x^{-1}$. Thus $p\in W^{1,1}_{\rm loc}$. Differentiating the stationary equation, we obtain $(\Delta+V_\xi)\cQ_i=0$. Consequently
\begin{equation*}
 \cQ^*(-Lw+\cA)=\cQ^*N_{q_\xi}(w),\qquad
 \widetilde v_t=P(-Lw+\cA)+\cD(p'+H-K^*\cA).
\end{equation*}
The second identity follows from \eqref{non:eq:field}. Since $DR^*P=DR^*$ and $K^*L=0$, the definition of $H$ implies
\begin{equation*}
 0=(p+DR^*\widetilde v)'=(I+DR^*\cD)(p'+H-K^*\cA).
\end{equation*}
The first factor is $M_\xi^*$, so
\begin{equation}\label{non:eq:constraint}
 p'+H=K^*\cA,\qquad w_t=\widetilde v-PDR\xi',\qquad
 \widetilde v_t=P(-Lw+\cA).
\end{equation}
Using $DR^*\widetilde v=-p$, $K^*\cD=I-\Phi\Phi^*$ and $\cQ^*P=\cQ^*-\Phi$, we obtain
\begin{equation}\label{non:eq:tangent-observation}
 \cQ^*\widetilde v=K^*\bar v-\Phi\Phi^*p,\qquad
 (\cQ^*\widetilde v)'=H-\Phi(-Lw+\cA).
\end{equation}

\subsubsection{Odd radiation below the energy regularity}\label{non:sec:dipole}

We estimate the full Newton forcing by a free wave estimate at regularity $\frac{3}{4}$. The two dipole moments enter only in identifying the bounded $L^2$ remainder. We write
\begin{equation*}
 X(f)=\|f\|_{L_t^{\frac{20}{3}}L_x^5},\qquad X(w)^5=\|w^5\|_{L_t^{\frac{4}{3}}L_x^1}.
\end{equation*}

To identify the Newton source, we set
$z_v=P_0\bar v$, $c=K^*\bar v=\nu^*G_0\bar v$ and
$B=\beta-\nu c$. The exact free system is
\begin{equation}\label{non:eq:masssystem}
 \begin{split}
 w_t&=z_v-G_0B,\\
 (z_v)_t&=\Delta w+F_{\rm wt}+w^5+\Phi^*a-G_0\nu H,
 \end{split}
\end{equation}
where
\begin{equation*}
 F_{\rm wt}=Vw+P^*[(V_\xi-V)w+B_{q_\xi}(w)],\qquad a=H-K^*w^5 .
\end{equation*}
By Lemma~\ref{non:lem:poly}, the basic estimates and
$\|w^5\|_{L^1L^2}=S(w)^5$, we have
\begin{equation*}
 \|B\|_{L^2Y_4}+\|F_{\rm wt}\|_{L^2Y_2}\lesssim\eps,\qquad
 \|a\|_{L^1}\lesssim\eps^2+\eps^5 .
\end{equation*}
By spectral orthogonality the entire initial velocity is ordinary,
\begin{equation}\label{non:eq:z0}
 z_v(0)=z_{\rm ord}=v_{a,t}(0)+\sum_l\kappa_l z_l\Psi_l.
\end{equation}
The odd Newton density in \eqref{non:eq:masssystem} is therefore
\begin{equation}\label{non:eq:measure}
 \mathrm d\gamma(s)=-\nu H(s)\,\mathrm ds,\qquad
 \|\gamma\|_{\TV(Y_4)}\lesssim\eps^2 .
\end{equation}
\paragraph{The spacetime estimate and the regular position mass.}\label{non:sec:mass}
Write $h_K=-S_0*(KH)$ for the entire Newton term, with retarded convolution in time. Since $|K_i|\lesssim\langle x\rangle^{-2}$, the dual Sobolev inequality implies
\begin{equation*}
 \|K_i\|_{\dot H^{-\frac{1}{4}}}\lesssim\|K_i\|_{\frac{12}{7}}<\infty .
\end{equation*}
The homogeneous wave Strichartz estimate \cite[Corollary~1.3]{KT} with $(q,r)=(\frac{20}{3},5)$ has position regularity $\frac{3}{4}$. By Minkowski's inequality, this estimate and the energy-class Strichartz estimate with $(q,r)=(5,10)$ imply
\begin{equation}\label{non:eq:charge-strichartz}
 X(h_K)+S(h_K)+\sup_t\|h_K(t)\|_{\dot H^{\frac{3}{4}}}
 \lesssim\|H\|_{L^1}\sum_i\big(\|K_i\|_{\dot H^{-\frac{1}{4}}}+\|K_i\|_2\big)
 \lesssim\eps^2 .
\end{equation}
For the ordinary forcing, duality in the homogeneous free estimates of Section~\ref{lin:free-section} implies
\begin{equation}\label{non:eq:dualities}
 M(C_0G_0*B)\lesssim\|B\|_{L^2Y_4},\qquad
 M(S_0*F)\lesssim\|F\|_{L^2Y_2+L^{\frac{5}{4}}L^{\frac{10}{9}}}.
\end{equation}
The initial data are $O(\eps)$ in $L^2\times\dot H^{-1}$ by \eqref{non:eq:radminus} and \eqref{non:eq:z0}, and $\|\Phi^*a\|_{L^1\dot H^{-1}}\lesssim\eps^2$. For $g=w-h_K$, H\"older and \eqref{non:eq:masssystem} imply
\begin{equation}\label{non:eq:mass-absorption}
 M(g)\lesssim\eps+\|w^5\|_{L^{\frac{5}{4}}L^{\frac{10}{9}}}
 \le C\eps+C X(w)^4 S(w),\qquad X(g)^4\le M(g)S(g)^3 .
\end{equation}
On bounded intervals $M(g)<\infty$ by the initial data and local time integration in $L^2$. Since $S(g)+S(w)\lesssim\eps$, we absorb $C\eps^4 M(g)$ in \eqref{non:eq:mass-absorption} and conclude $M(g)+X(w)\lesssim\eps$ uniformly in the interval. Interpolation with the energy bound and \eqref{non:eq:charge-strichartz} also yields $\sup_t\|w(t)\|_{\dot H^{\frac{3}{4}}}\lesssim\eps$.

To identify the two dipoles, take radial $\psi\in C_c^\infty$ of integral one, $\rho_j=-\partial_j\psi$ and $b_j(f)=\int x_j f$ for $j=1,2$. Define
\begin{equation}\label{non:eq:halo}
 h(t)=\sum_{j=1}^2\int_{[0,t]} S_0(t-s)G_0\rho_j\,b_j(\mathrm d\gamma(s)).
\end{equation}
Differentiating \eqref{lin:eq:newtonkernels} in space yields
\begin{equation}\label{non:eq:dipole-kernel}
 S_0(t)G_0\rho_j=\left(\frac{t x_j}{4\pi|x|^3}\mathbf1_{\{|x|>t\}}\right)*\psi .
\end{equation}
For odd $f\in Y_4$, $f^\circ=f-\sum_j b_j(f)\rho_j$ has zero mean and first moments. For $|\zeta|\le1$, $|\widehat{f^\circ}(\zeta)|\lesssim |\zeta|^2\||x|^2 f^\circ\|_1$. By Plancherel, $\sup_t\|S_0(t)G_0 f^\circ\|_2\le\|D^{-3}f^\circ\|_2\lesssim\|f\|_{Y_4}$. Thus $M(h_K-h)\lesssim\|\gamma\|_{\TV(Y_4)}$, and \eqref{non:eq:mainmass} follows from $M(g)\lesssim\eps$.

\subsubsection{Enhanced sources and projected evolution}\label{non:sec:enhanced}
The decomposition \eqref{non:eq:enhanced} separates three estimates.
The localized polynomial belongs to $L^2Y_2$, the pure quintic has
the integrable time norms of \eqref{non:eq:mainenhanced}, and the
localized correction is a Newton density in $L^1Y_4$.
These are precisely the source classes in
Theorem~\ref{lin:thm:old}. The identity $P^*=I-\Phi^*K^*$ also keeps
$K^*F_2=0$ exact.

\paragraph{The localized polynomial source.}
The $Y_2$ bound for $B_q(w)$ is Lemma~\ref{non:lem:poly}.
For $p=2,3,4$ take $\frac{1}{s_p}=\frac{4-p}{6}$. By \eqref{geo:eq:decay},
\begin{equation*}
 \|\rho q^{5-p}\|_{s_p}+\|\rho^2\nabla(q^{5-p})\|_{s_p}\lesssim1,
 \qquad \frac{1}{s_p}+\frac{p-1}{6}+\frac{1}{2}=1 .
\end{equation*}
When the derivative falls on $w$ we use $\rho^{-1}\nabla w$, and when it falls on $q^{5-p}$ we use $\rho^{-2}w$. Thus one product estimate covers the three powers,
\begin{equation*}
 \|\nabla(q^{5-p}w^p)\|_{L_t^2L^1}
 \lesssim E(w)^{p-1}\big(\|\nabla w\|_{L_t^2X_1}+\|w\|_{L_t^2X_2}\big).
\end{equation*}
The bounds $|V_\xi-V|\lesssim d\rho^{-8}$ and $|\nabla(V_\xi-V)|\lesssim d\rho^{-9}$ imply the same gradient estimate for $(V_\xi-V)w$ with coefficient $d$ in place of $E(w)^{p-1}$, and its $L^2Y_2$ bound is $O(d\LE(w))$. The projection subtracts the Schwartz function $\Phi^*K^*F$, with $\|K^*F\|_{L_t^2}\lesssim\|F\|_{L^2Y_2}$. This proves the $F_2$ bound in \eqref{non:eq:mainenhanced}.

\paragraph{The entire quintic source.}
Since $\|w\|_{L^{\frac{16}{3}}L^8}\le X(w)^{\frac{1}{4}}S(w)^{\frac{3}{4}}$, the product rule yields
\begin{equation}\label{non:eq:strongpure}
 \|w^5\|_{L^1L^2}=S(w)^5,\qquad
 \|w^5\|_{L^{\frac{4}{3}}W^{1,1}}
 \lesssim X(w)^5+E(w)X(w)S(w)^3\lesssim\eps^5 .
\end{equation}

For the additional Sobolev bounds in Theorem~\ref{non:thm:main}, we use $Y_2\hookrightarrow L^{\frac{10}{9}}$ on undifferentiated local terms. The coefficient derivatives satisfy the same weighted envelope as the coefficients in \eqref{non:eq:weighted-interpolation}. With one derivative on $w$, spatial H\"older uses $\rho q^{5-p}\in L^{\frac{30}{17-5p}}$ for $p=2,3$. For $p=4$ and the pure quintic the remaining bounds are
\begin{align*}
 \|qw^3\nabla w\|_{L^2L^{\frac{10}{9}}}
 &\le\|q\|_{15}\|w\|_{L^6L^9}^3E(w)
 \lesssim E(w)^{\frac{3}{2}}S(w)^{\frac{5}{2}},\\
 \|w^5\|_{L^{\frac{5}{4}}W^{1,\frac{10}{9}}}
 &\lesssim X(w)^4S(w)+E(w)S(w)^4 .
\end{align*}
Here we interpolated $L^6L^9$ between $L^\infty L^6$ and $L^5L^{10}$. The decaying potential difference and the Schwartz correction obey the same bound by H\"older.

\paragraph{Exact compatibility and enhanced initial data.}
The decomposition \eqref{non:eq:enhanced} is an exact rewriting of
\eqref{non:eq:basicsources}, since
\begin{equation*}
 P^*w^5+\Phi^*H=w^5+\Phi^*(H-K^*w^5)
                =w^5+G_0\Gamma_{\rm c}.
\end{equation*}
Using $G_0D^2\Phi^*=\Phi^*$ in the canonical decaying realization, we obtain
\begin{equation*}
 K^*F_2=0,\qquad K^*(F_p+G_0\Gamma_{\rm c})=H,
 \qquad \|\Gamma_{\rm c}\|_{L^1Y_4}\lesssim\eps^2+\eps^5.
\end{equation*}
The density $\Gamma_{\rm c}$ has vanishing mean and first moments by
integration by parts. In the free velocity equation the Newton density
is $-\nu H$ as in \eqref{non:eq:measure}. By \eqref{non:eq:initial} the
correction $(w-v_a,\bar v-v_{a,t})$ has zero complementary initial data
and zero initial kernel charge. Its enhanced source norms sum to
\begin{equation*}
 A_*\lesssim d\eps+\eps^2+\eps^5\lesssim d\eps.
\end{equation*}
Multiplication of these sources by $\mathbf1_{[0,\frac{T}{2}]}$ preserves
their norms and both compatibilities in \eqref{non:eq:compat}. Applying
the mixed Green operator separately to the earlier and later sources
selects their respective negative amplitudes by linearity.
By Theorem~\ref{lin:thm:old}, the observations of the earlier part decay
at rate $T^{-\frac{3}{4}}$ and its Strichartz norm at rate $T^{-\frac{1}{10}}$.

The contraction, the regular mass estimate and the enhanced source
bounds prove Theorem~\ref{non:thm:main}.

\paragraph{The projected evolution.}\label{non:sec:projected}

We set
\begin{equation*}
 W=w+PR(\xi),\qquad \Gamma_{\mathrm{coord}}=\xi+\Phi R(\xi),\qquad
 U=Q+W+K\Gamma_{\mathrm{coord}}=q_\xi+w.
\end{equation*}
Differentiating $W$ in \eqref{non:eq:constraint} and using
$\Delta U+U^5=-Lw+\cA$, we obtain
\begin{equation}\label{non:eq:projected}
 W_t=\widetilde v,\qquad
 W_{tt}=P[\Delta U+U^5],\qquad
 P[U_{tt}-\Delta U-U^5]=0.
\end{equation}
Thus an arbitrary prescribed path solves the projected nonlinear
problem and leaves a kernel-valued residual. Section~\ref{end:start}
chooses the path so that this residual vanishes. Since the projected
equation contains $\Gamma_{\rm coord}$ with no time derivative,
supplying the localized coordinate acceleration recovers the full wave
equation at the endpoint fixed point.

\subsubsection{The force and the current response}\label{force:start}

We compare $H$ with the quadratic force $F(v_a)$ in $L^1(I_j)$ and then test the difference against $\varphi_j$. We estimate the radiation created on recent blocks by a
recurrence, and use time decay for the radiation created earlier.

\begin{theorem}[Current response and nonlinear force]\label{cat:thm:blocks}
Let $\xi$ be an admissible prescribed path for Theorem~\ref{non:thm:main}, with
\begin{equation*}
 \|\xi\|_\infty\lesssim d,\quad \xi(0)=0,\quad
 \|\xi'\|_{L^2(I_j)}\lesssim\sqrt d\,\eta_j.
\end{equation*}
Let $(w,\bar v)$ be its prescribed field driven by \eqref{cat:eq:train},
with $p,\widetilde v,H$ as in \eqref{non:eq:force-definition}--\eqref{non:eq:field}. For the splitting
in \eqref{cat:eq:source-split}, we have
\begin{equation*}
 \|\boldsymbol c_j\|_{\cR}\lesssim d\eta_j,\qquad
 \LE_{I_j}(w-v_a)+\OO_{I_j}(\bar v-v_{a,t})\lesssim d\eta_j.
\end{equation*}
Its force satisfies
\begin{equation}\label{cat:eq:forcecomparison}
 H=F(v_a)+R_{\rm ext},\qquad
 \|R_{\rm ext}\|_{L^1(I_j)}\lesssim\sqrt d\,\eta_j^2,
\end{equation}
and its tested impulses satisfy
\begin{equation}\label{cat:eq:firsterror}
 \sup_{j\ge1}\left|\eta_j^{-2}\int\varphi_j H\,\dd t-a_j\right|
 \lesssim T_0^{-1}+\sqrt d .
\end{equation}
\end{theorem}

\paragraph{The exact force.}\label{cat:sec:nonlinear}
By the spectral separation in Theorem~\ref{cat:thm:catalogue},
\begin{equation}\label{cat:eq:gaugeforce}
 \Phi v_a=0,\qquad r=v_a .
\end{equation}
Subtracting the quadratic force \eqref{cat:eq:force-density} from
\eqref{non:eq:force-definition}, we obtain
\begin{equation}\label{cat:eq:Rext}
 (R_{\rm ext})_i=
 \la\cQ_i,N_{q_\xi}(w)\ra-F_i(v_a)
 +\xi_k'\la R_{ik}(\xi),\widetilde v\ra .
\end{equation}
By Theorem~\ref{non:thm:main}, the response
$\boldsymbol e=(w-v_a,\bar v-v_{a,t})$ has
$\|\boldsymbol e\|_{\cR}\lesssim d\eps$ and $\|H\|_1\lesssim\eps^2$.
Its enhanced sources in \eqref{non:eq:mainenhanced}, including the
pure term $w^5$, have total size
$A\lesssim d\eps+\eps^2+\eps^5\lesssim d\eps=d^4$.
We first estimate the current response.

\paragraph{Earlier forcing and the current response.}\label{cat:sec:current}

For $j\ge1$, we choose a smooth cutoff $\chi_j$ equal to one on $t\le \frac{T_j}{3}$
and zero on $t\ge \frac{T_j}{2}$. We split the enhanced sources with multipliers
$\chi_j$ and $1-\chi_j$. We apply the mixed Green operator to each part
with zero complementary data, selecting its negative amplitudes separately.
By linearity the two responses satisfy
\begin{equation}\label{cat:eq:source-split}
 \boldsymbol e=\boldsymbol o_j+\boldsymbol c_j,\qquad
 \boldsymbol o_j=(o_j,y_{o_j}),\quad \boldsymbol c_j=(c_j,y_{c_j}),\quad
 \boldsymbol o_0=0,\quad \boldsymbol c_0=\boldsymbol e.
\end{equation}
The source cutoff preserves the compatibilities in \eqref{non:eq:compat}.
Theorem~\ref{lin:thm:old} applies to $\boldsymbol o_j$, and the global
Green bound applies to both parts:
\begin{equation}\label{cat:eq:oldbounds}
 \begin{split}
 \LE_{I_j}(o_j)+\OO_{I_j}(y_{o_j})&\lesssim AT_j^{-\frac{3}{4}},\\
 S_{I_j}(o_j)&\lesssim AT_j^{-\frac{1}{10}},\qquad E(o_j)\lesssim A,\\
 \|\boldsymbol o_j\|_{\cR}+\|\boldsymbol c_j\|_{\cR}&\lesssim A.
 \end{split}
\end{equation}
We set $J_j=[\frac{T_j}{3},\infty)$ for $j\ge1$ and
$C_j=\|\boldsymbol c_j\|_{\cR}$. To estimate the sources of
$\boldsymbol c_j$, we use the \emph{preceding} split
$w=r+o_{j-1}+c_{j-1}$ on the whole tail $J_j$. The old sources for this
split have ended by $\frac{T_{j-1}}{2}=\frac{T_j}{8}$. We may therefore apply Theorem~\ref{lin:thm:old} to obtain
\begin{equation*}
 \LE_{J_j}(o_{j-1})+\OO_{J_j}(y_{o_{j-1}})\lesssim AT_j^{-\frac{3}{4}},
 \qquad S_{J_j}(o_{j-1})\lesssim AT_j^{-\frac{1}{10}}.
\end{equation*}
For $j=1$ this old field is zero. By the catalogue tail bounds and the
sum of the squared block bounds on $\xi'$, we have
$\|\xi'\|_{L^2(J_j)}\lesssim\eta_j$ and
\begin{equation*}
 \LE_{J_j}(w)+\OO_{J_j}(\bar v)
       \lesssim\eta_j+AT_j^{-\frac{3}{4}}+C_{j-1}.
\end{equation*}

\paragraph{A two-power bound for the whole block nonlinearity.}
We split $N_q(w)$ on $|w|\le|q|$ and its complement. The two odd pieces
are bounded by $C|q|^3|w|^2$ and $C|w|^5$, respectively. We apply $P^*$
to each to preserve the kernel compatibility. We use this split only in
the basic Green source norm $L^2Y_1+L^1L^2$, after the original smooth
enhanced sources have supplied the old response. On $J_j$, the two
endpoint powers of $w=r+o_{j-1}+c_{j-1}$ are bounded by the sum of the
same powers of its three components. The estimates are
\begin{equation}\label{cat:eq:thresholdbounds}
\begin{array}{c|cc}
 f & \||q|^3 |f|^2\|_{L^2Y_1} & \||f|^5\|_{L^1L^2} \\ \hline
 r       & O(\eps^{\frac{1}{2}}\eta_j^{\frac{3}{2}}) & O(\eta_j^5)\\
 o_{j-1} & O(A^2 T_j^{-\frac{7}{12}})          & O(A^5 T_j^{-\frac{1}{2}})\\
 c_{j-1} & O(A C_{j-1})               & O(A^4 C_{j-1}).
\end{array}
\end{equation}
We used Lemma~\ref{non:lem:poly} in the quadratic column and
$\|f^5\|_{L_t^1L_x^2}=S(f)^5$ in the quintic column. In the last row all factors have
the uniform bound $O(A)$ and one is bounded by $C_{j-1}$.

\paragraph{The potential and the charge.}
By the decay of $V_\xi-V$ and Newton duality,
\begin{equation*}
 \|(V_\xi-V)w\|_{L^2(J_j;Y_2)}\lesssim d\LE_{J_j}(w),\qquad
 \|p\|_{L^2(J_j)}\lesssim d\OO_{J_j}(\bar v).
\end{equation*}
For $\beta=D^2[PDR(\xi)\xi'-\cD p]$ we thus obtain
\begin{equation*}
 \|\beta\|_{L^2(J_j;Y_4)}
 \lesssim d\big(\eta_j+AT_j^{-\frac{3}{4}}+C_{j-1}\big).
\end{equation*}
The geometric term in \eqref{non:eq:force-definition} is bounded by
$C\|\xi'\|_{L^2(J_j)}\OO_{J_j}(\widetilde v)$. In \eqref{non:eq:forces},
the exponent $\frac{8-p}{3}$ is at least one for $2\le p\le5$. We bound one
$\LE$ factor on $J_j$ and the remaining factors on the whole half-line to obtain
\begin{equation}\label{cat:eq:Hblock}
 \|H\|_{L^1(J_j)}\lesssim\eps(\eta_j+AT_j^{-\frac{3}{4}}+C_{j-1}).
\end{equation}

\paragraph{Closing the current-response estimate.}
The sources of $\boldsymbol c_j$ are supported in $J_j$ and its
complementary initial data vanish. Applying Theorem~\ref{lin:thm:green}
with \eqref{cat:eq:thresholdbounds} and \eqref{cat:eq:Hblock}, we obtain
\begin{equation}\label{cat:eq:blocksource}
 C_j\lesssim(d+\eps)(\eta_j+AT_j^{-\frac{3}{4}}+C_{j-1})
        +A^2T_j^{-\frac{7}{12}}+A^5T_j^{-\frac{1}{2}}.
\end{equation}
As $\eta_j=\sqrt{\frac{d}{T_j}}$, $A\lesssim d^4$ and $T_j\ge1$, the three
explicit tail errors are $O(d\eta_j)$. Dividing \eqref{cat:eq:blocksource} by $\eta_j=\frac{\eta_{j-1}}{2}$, we obtain,
with $K$ independent of $d,j$,
\begin{equation}\label{cat:eq:muclose}
 \frac{C_j}{\eta_j}\le Kd+2K(d+\eps)\frac{C_{j-1}}{\eta_{j-1}},
 \qquad \frac{C_0}{\eta_0}\lesssim d .
\end{equation}
We take $2K(d+\eps)\le\frac{1}{2}$. Iterating \eqref{cat:eq:muclose}, we obtain
\begin{equation*}
 \frac{C_j}{\eta_j}
 \le Kd\sum_{k=0}^{j-1}2^{-k}+2^{-j}\frac{C_0}{\eta_0}\lesssim d .
\end{equation*}
Substitution in \eqref{cat:eq:oldbounds} yields
\begin{equation}\label{cat:eq:currentbound}
 \|\boldsymbol c_j\|_{\cR}\lesssim d\eta_j,\qquad
 \LE_{I_j}(w-r)+\OO_{I_j}(\bar v-v_{a,t})\lesssim d\eta_j .
\end{equation}
The recurrence starts from the global correction at $j=0$. Every
response in it uses the negative-mode integral on the entire future half-line.

\subsubsection{The normalized force error}\label{force:sec:comparison}

We estimate the quadratic, higher-order and path-derivative terms in
\eqref{cat:eq:Rext} separately. We set $n_j=r+c_j$, so $w=n_j+o_j$. By \eqref{cat:eq:blockobservations}, \eqref{cat:eq:oldbounds} and
\eqref{cat:eq:currentbound}, we have
\begin{equation*}
 \begin{gathered}
 \LE_{I_j}(n_j)+S_{I_j}(n_j)\lesssim\eta_j,\qquad E(n_j)\lesssim\eps,\\
 \LE_{I_j}(o_j)\lesssim AT_j^{-\frac{3}{4}},\quad
 S_{I_j}(o_j)\lesssim AT_j^{-\frac{1}{10}},\quad E(o_j)\lesssim A.
 \end{gathered}
\end{equation*}

For the quadratic term in \eqref{cat:eq:Rext}, we use the difference of squares to obtain
\begin{equation}\label{cat:eq:quadratic-force-error}
 \int_{I_j}\!\int\rho^{-8}|w^2-r^2|\dd x\dd t
 \lesssim\LE_{I_j}(w-r)[\LE_{I_j}(w)+\LE_{I_j}(r)]
 \lesssim d\eta_j^2.
\end{equation}
For the coefficient difference, \eqref{geo:eq:density} yields
\begin{equation*}
 |\cQ_i(\xi)q_\xi^3-K_iQ^3|\lesssim d\rho^{-8},\qquad
 \iint_{I_j\times\mathbb R^3}d\rho^{-8}|w|^2\dd x\dd t\lesssim d\eta_j^2 .
\end{equation*}

We apply \eqref{non:eq:force-interpolation} on $I_j\times\mathbb R^3$ after dividing each integral by $\eta_j^2$. The normalized quadratic term is bounded by $C\eta_j^{-2}\LE_{I_j}(w)^2\lesssim1$. It remains to estimate only the fifth power. We use $|w|^5\lesssim|n_j|^5+|o_j|^5$ and the endpoint $p=5$ of \eqref{non:eq:forces}. Since $\frac{3}{4}+\frac{5}{2}\cdot\frac{1}{10}=1$ and $T_j\eta_j^2=d$, we obtain
\begin{equation*}
 \frac{J_5}{\eta_j^2}\lesssim\eta_j^{-2}\left(\eta_j^{\frac{7}{2}}\eps^{\frac{3}{2}}+A^5 T_j^{-1}\right)
 \lesssim \eps^3+\frac{A^5}{d}\lesssim\eps^3 .
\end{equation*}
Here $\eta_j\le\eps=d^3$ and $A\lesssim d^4$. Hence
\begin{equation}\label{cat:eq:oldhighernormalized}
 \eta_j^{-2}\sum_{p=3}^5\iint_{I_j\times\mathbb R^3}|\cQ_i q^{5-p}w^p|\,\dd x\,\dd t
 \lesssim\sum_{p=3}^5\eps^{p-2}\lesssim d^3 .
\end{equation}
For the geometric term, $D^2R_{ik}(\xi)$ is uniformly bounded in
$Y_4$. By Newton duality and \eqref{cat:eq:currentbound}, we have
$\OO_{I_j}(\bar v)\lesssim\eta_j$ and
$\|p\|_{L^2(I_j)}\lesssim d\eta_j$, hence
$\OO_{I_j}(\widetilde v)\lesssim\eta_j$. Cauchy--Schwarz yields
\begin{equation}\label{cat:eq:geometric}
 \|\xi'_k\langle R_{ik}(\xi),\widetilde v\rangle\|_{L^1(I_j)}
 \lesssim\|\xi'\|_{L^2(I_j)}\OO_{I_j}(\widetilde v)
 \lesssim\sqrt d\,\eta_j^2 .
\end{equation}

Adding \eqref{cat:eq:quadratic-force-error}, the chart-coefficient error,
\eqref{cat:eq:oldhighernormalized} and \eqref{cat:eq:geometric}, we obtain
\eqref{cat:eq:forcecomparison}. Since $0\le\varphi_j\le1$ and $\eta_{j-1}=2\eta_j$, the error tested over the two neighboring blocks is $O(\sqrt d\,\eta_j^2)$. We apply \eqref{cat:eq:reference-moments} to the reference force.
This proves \eqref{cat:eq:firsterror} and completes Theorem~\ref{cat:thm:blocks}.
\space
\section{Selection of a recurrent solution on the half-line}\label{end:start}
We prescribe both endpoints of the chart coordinate on every block and solve its acceleration equation there. The packet controls match the resulting velocities at the common endpoints. We then recover the wave equation and identify the free wave along both sequences.
\subsection{The endpoint variables and admissible domain}\label{end:domain-section}
We prescribe the endpoint positions and bound the two terms that
enter the endpoint map, the force $H$ and the velocity observation $Z$.
We use the scales \eqref{main:scales}, the train \eqref{cat:eq:train}
and the field of Theorem~\ref{non:thm:main}. The velocity observation
$Z=\cQ^*\widetilde v$ is the forcing in the chart coordinate
\begin{equation}\label{end:ep:coordinate}
 \gamma(\xi)=\Phi(q_\xi-Q)=\xi+\Phi R(\xi),\qquad
 D\gamma(\xi)=M_\xi.
\end{equation}
On a fixed small convex chart ball, $\|M_\xi-I\|\le\frac{1}{2}$.
Thus $\gamma$ is a smooth bi-Lipschitz coordinate, its inverse has
Lipschitz constant at most two, and its image contains a neighborhood
of zero. We use the maximum norm for the control cube and the Euclidean norm for the block velocities.

\paragraph{The domain and coordinate forcing.}
The coordinate in \eqref{end:ep:coordinate} satisfies $\gamma(0)=0$ and $D\gamma(0)=I$. For the alternating endpoints of Section~\ref{main:outline} we therefore have
\begin{equation}\label{end:eq:control-cube}
 r=\frac{\gamma(de)}{d}=e+O(d),\qquad \mathcal A=\prod_{j\ge1}[-2,2]^4 .
\end{equation}
We write $c_j=\frac{\ell_j}{T_j}$, so $c_0=4$ and $c_j=3$ for $j\ge1$.
We take the chart velocity $g=(\gamma(\xi))'$ as the unknown. For $C_2\ge1$ let
\begin{equation}\label{end:ep:path-domain}
 \mathcal Y=\left\{g\in L^2(0,\infty;\mathbb R^4):
  \int_{I_j}g\,\dd t=(-1)^j d r,\quad
  \|g\|_{L^2(I_j)}\le C_2\sqrt d\,\eta_j\quad(j\ge0)\right\}.
\end{equation}
We recover the path by
\begin{equation}\label{end:eq:path-primitive}
 \Gamma(t)=\int_0^t g(s)\dd s,\qquad
 \xi(t)=\gamma^{-1}(\Gamma(t)).
\end{equation}
The block integrals cancel in pairs, so $\Gamma(s_{2k})=0$ and $\Gamma(s_{2k+1})=dr$. On $I_j$ we have $|\Gamma(t)-\Gamma(s_j)|\le\sqrt{\ell_j}\,\|g\|_{L^2(I_j)}\le2C_2d$. For fixed $C_2$ and small $d$ the primitive thus lies in the coordinate neighborhood. By the chain rule,
\begin{equation*}
 \|\xi\|_\infty\lesssim(1+C_2)d,\qquad
 \|\xi'\|_{L^2(I_j)}\le 2C_2\sqrt d\,\eta_j .
\end{equation*}
The block constants $g=\frac{(-1)^j d r}{\ell_j}$ belong to $\mathcal Y$ for $C_2$ sufficiently large. We place $\mathcal A$ in the weighted sequence space with squared norm $\sum_j2^{-j}|a_j|^2$. Its norm topology on the cube is coordinatewise convergence. In the product of this space and $L^2(0,\infty;\mathbb R^4)$, the set $\mathcal A\times\mathcal Y$ is nonempty, closed, bounded and convex.

\paragraph{The prescribed wave on the half-line.}
\begin{lemma}\label{end:cont-field}
The prescribed field on $[0,\infty)$ depends continuously in $\cR$ on $(a,g)\in\mathcal A\times\mathcal Y$.
\end{lemma}
\begin{proof}
Let $(a_n,g_n)\to(a,g)$. The packet tail in $\cR$ is bounded by $C\sum_{j\ge J}\eta_j$, so convergence of the finite sums implies $v_{a_n}\to v_a$ in $\cR$. By \eqref{end:eq:path-primitive}, $\Gamma_n\to\Gamma$ uniformly on compact intervals. The chain rule and dominated convergence imply
\begin{equation*}
 \|\xi_n'-\xi'\|_2\le 2\|g_n-g\|_2+
 \|[D\gamma^{-1}(\Gamma_n)-D\gamma^{-1}(\Gamma)]g\|_2\longrightarrow0 .
\end{equation*}
Let $E_n,E_*$ be the fixed corrections of $\mathcal T_n,\mathcal T$. Since their contraction constants are at most $\frac{1}{2}$,
\begin{equation*}
 \|E_n-E_*\|_{\cR}\le 2\|\mathcal T_nE_*-\mathcal TE_*\|_{\cR}.
\end{equation*}
In the right side we may replace $v_{a_n}$ and $\xi_n'$ by $v_a$ and $\xi'$ at cost $o(1)$ by polarization of \eqref{non:eq:basic-source-bounds}. Only the chart coefficients then vary. They converge pointwise in time with the same spatial weights. On this fixed input, Lemma~\ref{non:lem:poly} dominates the polynomial sources in their respective $L^1_t$ or $L^2_t$ norms. The terms containing $p$ or $DR$ are dominated there by the Newton observation, $|\xi'|$, or their product, by \eqref{non:eq:pbound}. Dominated convergence and \eqref{lin:eq:green} now imply $E_n\to E_*$.
\end{proof}
By \eqref{cat:eq:force-density} and \eqref{cat:eq:blockobservations},
\begin{equation*}
 \|F(v_a)\|_{L^1(I_j)}\lesssim\LE_{I_j}(v_a)^2\lesssim\eta_j^2 .
\end{equation*}
Thus \eqref{cat:eq:forcecomparison} and \eqref{cat:eq:firsterror} imply
\begin{equation}\label{end:eq:block-force-moments}
 \|H\|_{L^1(I_j)}\le C_F\eta_j^2,\qquad
 \left|\eta_j^{-2}\int\varphi_j H\,\dd t-a_j\right|
 \lesssim T_0^{-1}+C_{\rm force}\sqrt d\quad(j\ge1).
\end{equation}
The first estimate also holds for $j=0$. Once $C_{\rm force}\sqrt d\le1$, its constant depends only on the reference family.
\begin{lemma}[Small coordinate forcing]\label{end:ep:k}
The coordinate forcing satisfies
\begin{equation}\label{end:ep:k-bounds}
 \|Z\|_{L^2(I_j)}\le C_Z d\eta_j.
\end{equation}
The constant may depend on the fixed path bound $C_2$.
\end{lemma}
\begin{proof}
By \eqref{non:eq:tangent-observation}, $K^*v_{a,t}=0$ and $K^*=\nu^*G_0$ on velocities,
\begin{equation*}
 \|Z\|_{L^2(I_j)}\lesssim
 \OO_{I_j}(\bar v-v_{a,t})+d\,\OO_{I_j}(\bar v)\lesssim d\eta_j .
\end{equation*}
We used \eqref{non:eq:pbound} on $I_j$ and the block response estimate of Theorem~\ref{cat:thm:blocks}.

\end{proof}

\subsection{The endpoint map}\label{end:map-section}
For $U=q_\xi+w$ and $B=\Phi(\Delta U+U^5)$, the second identity in \eqref{non:eq:tangent-observation} reads
\begin{equation}\label{end:cont-primitive-identity}
 Z'=H-B,\qquad B=\Phi(-Lw+\cA).
\end{equation}
On $I_j=[s_j,s_{j+1}]$ let $\langle\cdot\rangle_j$ denote the interval mean. We define the output velocity by its derivative and its prescribed mean, and update each control by the jump at its node,
\begin{equation}\label{end:ep:map}
 \widehat g'=B\ \text{on }I_j^\circ,\qquad
 \langle\widehat g\rangle_j=\frac{(-1)^jdr}{\ell_j},\qquad
 \widehat a_j=a_j+\frac{\widehat g(s_j^+)-\widehat g(s_j^-)}{\eta_j^2}\quad(j\ge1).
\end{equation}
The velocity conditions include $I_0$. On each open block $(\widehat g+Z)'=H$. Testing its distributional derivative against the hat $\varphi_j$ and using the two prescribed means, we obtain
\begin{equation}\label{end:eq:moment-ends}
 \widehat g(s_j^+)-\widehat g(s_j^-)
 =\eta_j^2 a_j^*-\int\varphi_j H\,\dd t-\int\varphi_j'Z\,\dd t .
\end{equation}
Here $Z$ is continuous by \eqref{end:cont-primitive-identity}, so only $\widehat g$ contributes to the jump.
\begin{proposition}[Preservation of the endpoint domain]\label{end:ep:preservation}
For $C_2$ sufficiently large and $d$ sufficiently small, \eqref{end:ep:map} preserves $\mathcal A\times\mathcal Y$ and has relatively compact image.
\end{proposition}
\begin{proof}
By Cauchy--Schwarz and \eqref{end:ep:k-bounds},
\begin{equation*}
 \left|\int\varphi_j' Z\,\dd t\right|
 \le\sum_{k=j-1}^{j}\ell_k^{-\frac{1}{2}}\|Z\|_{L^2(I_k)}
 \lesssim C_Z\sqrt d\,\eta_j^2 .
\end{equation*}
Thus \eqref{end:eq:block-force-moments} and \eqref{end:eq:moment-ends} imply, for every input,
\begin{equation}\label{end:ep:control-invariance}
 \sup_{j\ge1}|\widehat a_j-a_j^*|
 \lesssim T_0^{-1}+(C_{\rm force}+C_Z)\sqrt d .
\end{equation}
Since $\sup_j |a_j^*|\le\frac{5}{3}(1+O(d))$ by \eqref{main:eq:ideal-controls}, all output controls belong to the cube.

On $I_j$, the mean-zero part of $\widehat g+Z$ is the primitive of $H$ with kernel $\mathbf1_{s<t}-\frac{s_{j+1}-s}{\ell_j}$. Its squared $L^2_t(I_j)$ norm is $\frac{(s-s_j)(s_{j+1}-s)}{\ell_j}\le\frac{\ell_j}{4}$. Subtracting the mean-zero part of $Z$ and using $\sqrt{\ell_j}\eta_j=\sqrt{c_jd}$, we obtain
\begin{equation}\label{end:ep:path-output}
 \frac{\|\widehat g\|_{L^2(I_j)}}{\sqrt d\,\eta_j}
 \le\frac{|r|}{\sqrt{c_j}}+\frac{\sqrt{c_j}}{2}C_F+C_Z\sqrt d .
\end{equation}
We first take $C_2>2C_F+3$, using the reference bound for $C_F$, then take $d$ small enough for both estimates. The Schwartz tests satisfy $|B|\lesssim\|\nabla U\|_2+\|U\|_6^5\lesssim1$. Thus \eqref{end:ep:map} bounds the outputs in $H^1$ on each fixed block. Rellich compactness there and the uniform tail $\sum_{j\ge J}C_2^2d\eta_j^2\to0$ imply relative compactness in $L^2(0,\infty)$. The control cube is compact in its weighted norm.
\end{proof}

\subsection{The fixed point and the wave equation}\label{end:cont-section}
For convergent inputs, Lemma~\ref{end:cont-field} implies local uniform convergence of $U$ in $\dot H^1$, hence of $B=\Phi(\Delta U+U^5)$. The derivative and mean in \eqref{end:ep:map} determine $\widehat g$ uniformly on each closed block, including its one-sided endpoints. The uniform tail bound \eqref{end:ep:path-output} proves continuity in the domain topology. Schauder's theorem now supplies a fixed point.

At a fixed point all jumps in \eqref{end:ep:map} vanish, so $g'=B$ on the half-line. Thus $\pi=g+Z$ belongs to $L^2$ and $\pi'=H\in L^1$ by \eqref{end:cont-primitive-identity}. Its limit at infinity exists and is zero. By \eqref{end:eq:block-force-moments},

\begin{equation}\label{end:eq:future-momentum}
 \pi=g+Z,\qquad \|\pi\|_{L^\infty(I_j)}\lesssim\eta_j^2,\qquad
 \pi(t)=-\int_t^\infty H(s)\,\dd s .
\end{equation}

\paragraph{Recovery of the equation and initial data.}\label{end:cont-recovery}
Since $B$ is continuous, $\Gamma=\gamma(\xi)$ and $\xi$ are $C^2$, and
\begin{equation}\label{end:eq:physical-coordinate}
 \Gamma'=M_\xi\xi'=\pi-Z,\qquad \Gamma''=\Phi(\Delta u+u^5),\qquad u=q_\xi+w .
\end{equation}

The field identities imply
\begin{equation}\label{end:cont-physical}
 u_t=\widetilde v+KM_\xi\xi',\qquad
 \cQ^*u_t=\pi,\qquad
 \Delta u+u^5=-Lw+\cA.
\end{equation}
In $H^{-1}$ we differentiate $u=Q+W+K\Gamma$ and use \eqref{non:eq:projected} and \eqref{end:eq:physical-coordinate} to obtain
\begin{equation*}
 u_{tt}=P(\Delta u+u^5)+K\Gamma''=(P+K\Phi)(\Delta u+u^5)=\Delta u+u^5 .
\end{equation*}
On compact time intervals $q_\xi+w\in L^5L^{10}$, hence $u^5\in L^1L^2$. Uniqueness for the free wave with this source and these energy traces identifies $u$ with its Duhamel solution.

At $t=0$, the packets and negative modes are orthogonal to $K$. By \eqref{non:eq:p0} and \eqref{non:eq:initial}, $Z(0)=K^*\bar v(0)=0$. Since $M_0=I$, \eqref{end:eq:future-momentum} implies
\begin{equation*}
 \xi'(0)=\pi(0),\qquad |\pi(0)|\lesssim\eta_0^2=\varepsilon^2 .
\end{equation*}
The initial data are
\begin{equation}\label{end:cont-initial}
 \begin{aligned}
 u(0)&=Q+v_a(0)+\sum_l z_l\Psi_l,\\
 u_t(0)&=v_{a,t}(0)+\sum_l\kappa_l z_l\Psi_l+K\pi(0),\qquad
 |z|\lesssim d\eps .
 \end{aligned}
\end{equation}
Since the columns of $K$ belong to $L^2$, \eqref{non:eq:radminus}
implies $\|\vec u(0)-(Q,0)\|_{\cH}\lesssim\eps$.
\subsection{Scattering and the alternating limits}\label{end:global-section}
We identify the radiation shared by both endpoint subsequences. The half-line
fixed point $u=q_\xi+w$ satisfies
\begin{equation*}
 \|\xi\|_\infty\lesssim d,\qquad
 \|\xi'\|_{L^2(0,\infty)}\lesssim\sqrt d\,\eps,\qquad
 \|(w,\bar v)\|_{\cR([0,\infty))}\lesssim\eps.
\end{equation*}

\begin{proposition}[One free radiation]\label{end:global-scattering}
There is a unique free energy wave $v_+$ with data of size $O(\eps)$ such
that
\begin{equation}\label{end:global-radiation}
 \norm{(u,u_t)-(q_\xi,0)-(v_+,v_{+,t})}_{\mathcal H}\longrightarrow0.
\end{equation}
\end{proposition}
\begin{proof}
Set $w=u-q_\xi$ and let $U_0(t)$ be the free wave group. In the interaction representation the physical equation reads
\begin{equation}\label{end:global-physical-sources}
 \partial_t\!\left[U_0(-t)\binom{w}{u_t}\right]
 =U_0(-t)\binom{-\cQ(\xi)\xi'}{V_\xi w+B_{q_\xi}(w)+w^5}.
\end{equation}
Since $D^2\cQ_i=5q_\xi^4\cQ_i$ is bounded in $Y_2$, the free estimate \eqref{lin:eq:freebound} on $I=[a,b]$ implies
\begin{equation}\label{end:global-duality}
 \left\|\left[U_0(-t)\binom{w}{u_t}\right]_a^b\right\|_{\mathcal H}
 \lesssim \|\xi'\|_{L^2(I)}+\|V_\xi w+B_{q_\xi}(w)\|_{L^2(I;Y_1)}
          +\|w\|_{L^5(I;L^{10})}^5 .
\end{equation}
By Lemma~\ref{non:lem:poly} all three norms are $O(\eps)$ on the half-line and tend to zero on its tails. The energy limit of $U_0(-t)(w,u_t)(t)$ therefore exists. By \eqref{end:cont-initial} its norm is $O(\eps)$. Unitarity identifies it as the unique datum in \eqref{end:global-radiation}.
\end{proof}

By \eqref{end:ep:path-domain} and \eqref{end:eq:path-primitive}, $\xi(s_{2k})=0$ and $\xi(s_{2k+1})=de$. Since $s_j=T_j$ for $j\ge1$, \eqref{end:global-radiation} implies
\begin{equation}\label{end:eq:two-core-limits}
 (\vec u-\vec v_+)(T_{2k})\longrightarrow(Q,0),\qquad
 (\vec u-\vec v_+)(T_{2k+1})\longrightarrow(q_{de},0)
 \quad\text{in }\mathcal H.
\end{equation}

\subsection{Rigidity of a strong core limit}\label{fp:start}\label{fp:radiation-energy}

We use weak orthogonality of free wave frames \cite{BG,DKMprofiles}. For the frame $S_{c,s,O}\mathcal W(-b/s)$ the parameters are scale $s$, center $c$, rotation $O$ and time $b$. Two frames are orthogonal if
\begin{equation*}
 \frac{s}{r}+\frac{r}{s}+\frac{|c-d|+|b-b'|}{s}\longrightarrow\infty .
\end{equation*}
Their relative action then tends weakly to zero. Otherwise, on a subsequence the relative parameters converge, including the rotations, and the action converges strongly on each datum. In particular, for spatial frames,
\begin{equation}\label{fp:dislocation}
 S_n\to S\ \Longrightarrow\ S_nf\to Sf,\qquad
 S_n\text{ escaping}\ \Longrightarrow\ S_nf\rightharpoonup0
 \quad(f\in\mathcal H).
\end{equation}
A free wave $\vec v(t_n)$ at $t_n\to\infty$ has frame $(1,0,-t_n)$, orthogonal to every sequence of spatial frames, since either the scales separate or the normalized time translations diverge. Hence
\begin{equation}\label{fp:free-separation}
 \sup_{S\in\mathrm{Sim}}|\langle\vec v(t),Sf\rangle_{\mathcal H}|\longrightarrow0
 \qquad(t\to\infty).
\end{equation}

\begin{proposition}[Matching separated profiles]\label{fp:comparison}
Let $t_n\to\infty$ and let $v$ be a free energy wave. Suppose two finite families of nonzero energy pairs have internally separated spatial frames as in \eqref{main:profile-separation}, and
\begin{equation}\label{fp:core-representation}
 \vec v(t_n)+\sum_j S_{j,n}\mathbf f_j-\sum_k T_{k,n}\mathbf g_k\to0
 \quad\text{in }\mathcal H .
\end{equation}
Then $v=0$. After extraction the two families match bijectively, with
$S_{j,n}^{-1}T_{k,n}\to U_j$ and $\mathbf f_j=U_j\mathbf g_k$ for each matched pair. In particular a nonzero strong limit has exactly one profile.
\end{proposition}
\begin{proof}
After extraction, each pair of spatial frames has either convergent or escaping relative parameters. Let $\mathcal P$ consist of the pairs $(j,k)$ with $S_{j,n}^{-1}T_{k,n}\to U_{jk}$. Internal separation permits at most one partner for each index. By \eqref{fp:dislocation} and \eqref{fp:free-separation}, the squared norm in \eqref{fp:core-representation} tends to
\begin{equation*}
 \|\vec v(0)\|_{\mathcal H}^2+
 \sum_{(j,k)\in\mathcal P}\|\mathbf f_j-U_{jk}\mathbf g_k\|_{\mathcal H}^2+
 \sum_{j\ {\rm unmatched}}\|\mathbf f_j\|_{\mathcal H}^2+
 \sum_{k\ {\rm unmatched}}\|\mathbf g_k\|_{\mathcal H}^2=0 .
\end{equation*}
Every profile is nonzero. Hence every index has a partner, the matched profiles agree in their common frame, and $v=0$.
\end{proof}

\subsection{Exclusion of fixed profiles}\label{fp:energy-shape}
\begin{proof}[Proof of Theorem~\ref{main:counterexample}]
By \eqref{end:cont-initial} and Proposition~\ref{end:global-scattering}, $u$ satisfies \eqref{main:counterexample-bounds} and \eqref{main:scattering}. The alternating endpoints yield the two strong limits in \eqref{end:eq:two-core-limits}, and \eqref{geo:eq:shape-separation} separates their similarity orbits. Assume that \eqref{main:FP} and
\eqref{main:profile-separation} hold on $T_n$. We apply Proposition~\ref{fp:comparison}
on $T_{2k}$ with $v=v_L-v_+$ and $\mathbf f_j=\vec Q_{j,\ell_j}$. Then
\begin{equation*}
 v_L=v_+,\qquad J=1,\qquad (Q,0)=S_*\mathbf f_1.
\end{equation*}
Applying the same proposition on $T_{2k+1}$, we obtain another similarity $S'$ with
\begin{equation*}
 (q_{de},0)=S'\mathbf f_1=S'S_*^{-1}(Q,0).
\end{equation*}
Its position component contradicts
$\operatorname{dist}_{\dot H^1}(q_{de},\mathrm{Sim}\cdot Q)\gtrsim d>0$.
A full-time fixed-profile resolution would hold on $T_n$, so it is
excluded as well.
\end{proof}
\space
\section{The global Kelvin tower and its approximation}\label{merge:global-profile}

\begingroup

\def\R{\mathbb R}
\def\HH{\mathcal H}
\def\D{\mathcal D}
\def\W{\mathcal W}
\def\St{\mathcal S}
\def\dd{\,\mathrm d}
\def\vct#1{\boldsymbol{#1}}

Sections~\ref{sec:tower-1} and~\ref{sec:tower-2} compute the first interaction of two Kelvin-related profiles and determine their scales. Section~\ref{sec:tower-3} removes the leading error in the wave equation.

\subsection{The Kelvin pair and its scale expansion}
Here and in Section~\ref{merge:global-wave} we write $Q=Q_{\rm cross}$ and $P=\mathcal KQ$. We denote the free wave group on $\mathcal H$ by $\mathcal W(t)$. We use the energy density \eqref{p8:eq:energy-density} for the two profiles of Theorem~\ref{p8:th:tower}.
In the stereographic coordinates of Section~\ref{sec:tower-1}, the pair satisfies
\begin{equation}\label{p8:eq:kelvin-pair}
 P(x)=|x|^{-1}Q\bigl(\frac{x}{|x|^2}\bigr)=Q(-x),\qquad F_s(x)=s^{-\frac{1}{2}}F\bigl(\frac{x}{s}\bigr).
\end{equation}
Both solve $-\Delta F=F^5$ and have the same quadratic leading term $H_2(x)=-b x_1x_2$ at the origin, with $b=16\sqrt{2}\,\gamma>0$ in the chosen orientation. We define $Y_2=H_2|_{\mathbb S^2}$, $\Lambda=\frac{1}{2}+x\cdot\nabla$, and
\begin{equation}\label{p8:eq:coefficients}\begin{aligned}
 M&=\|\Lambda P\|_2^2=\|\Lambda Q\|_2^2,\qquad
 C=5\int_{\mathbb S^2}Y_2^2\dd\omega,\qquad a=64\bigl(\frac{M}{C}\bigr)^2,\\
 B_0&=-\frac{5}{4}\int_{\R^3}|x|^{-3}Y_2\bigl(\frac{x}{|x|}\bigr)Q(x)\dd x .
\end{aligned}\end{equation}
Here $M,C>0$ and $B_0<0$ by Sections~\ref{sec:tower-1}--\ref{sec:tower-2}.

\begin{proposition}[Refined scale asymptotics]\label{p8:pr:refined-scales}
The solution in Theorem~\ref{p8:th:tower} may be chosen so that, with $L=at^{-4}$ and $a$ as in \eqref{p8:eq:coefficients},
\begin{equation}\label{p8:eq:scale}
 \lambda=L+\frac{32}{15}L^2+O(t^{-12}),\qquad
 \mu=1-\frac{5}{18}L^2-\frac{2B_0}{5M}L^{\frac{5}{2}}+O(t^{-12}) .
\end{equation}
\end{proposition}

We prove the proposition together with Theorem~\ref{p8:th:tower} in Section~\ref{sec:tower-8}. The time origin matters in \eqref{p8:eq:scale}. If $t=T_0+\tau$ places the initial data at $\tau=0$, then $L=a(T_0+\tau)^{-4}$ throughout the expansion.
\endgroup
\space
\subsection{The signed Kelvin moment and the leading scales}\label{sec:tower-1}

We identify the first nonzero harmonic moment of the Kelvin pair
\eqref{p8:eq:kelvin-pair} and determine the scale exponent it selects.
In the stereographic chart $\sigma_t$ of Section~\ref{common:nodal}, the
pole $(1,0)$ lies on the critical Hopf fiber and $Q=\Omega^{\frac{1}{2}}q\circ\sigma_t$.
Here $P$ and $Q$ are finite-energy stationary states. In the
orientation of \eqref{common:quad}, we put $b=16\sqrt 2\,\gamma>0$ and
$H_2(x)=-b x_1x_2$. Both profiles start with $H_2(x)+O(|x|^3)$ at
zero. By Kelvin inversion, both $F=P,Q$ satisfy
\begin{equation}\label{p1:eq:kelvin-tail}
 F(r\omega)=r^{-3}H_2(\omega)+O(r^{-4}),\qquad
 \partial_r F(r\omega)=-3r^{-4}H_2(\omega)+O(r^{-5}).
\end{equation}
The differentiated remainders follow from Taylor expansion in the
smooth chart at the other end. The same expansion can be differentiated
any fixed number of times with respect to $r\partial_r$ and angular
variables.

\subsubsection{The moment entering the scale equation}
Using the ordinary area measure on $\mathbb S^2$, we compute
\begin{equation}\label{p1:eq:C}
 C=5\int_{\mathbb S^2}H_2(\omega)^2\,d\omega
 =5b^2\frac{4\pi}{15}=\frac{4\pi}{3}b^2=\frac{2048\pi}{3}\gamma^2>0 .
\end{equation}
The same flux computes the moments for both tower constructions.
Suppose a smooth stationary profile $F$ has the differentiated expansion
\begin{equation*}
 F(r\omega)=\sum_{j=0}^k r^{-j-1}Y_j(\omega)+O(r^{-k-2}),
 \qquad -\Delta_{\mathbb S^2}Y_j=j(j+1)Y_j .
\end{equation*}
For any homogeneous harmonic polynomial $H_k$ of degree $k$ with
$H_k F^5\in L^1$, we apply Green's identity to obtain
\begin{equation}\label{p1:eq:harmonic-flux}
 \int_{\mathbb R^3}H_k F^5
  =\lim_{R\to\infty}\int_{\partial B_R}
       (F\partial_r H_k-H_k\partial_r F)\,dS
  =(2k+1)\int_{\mathbb S^2}H_k(\omega)Y_k(\omega)\,d\omega .
\end{equation}
Spherical orthogonality removes the lower degrees. At degree $k$ the
two radial derivatives contribute $k$ and $k+1$. In particular
$\int H_2P^5=C>0$ for the Kelvin pair.

We also need the pairing with the dilation. If $H$ is homogeneous of
degree $k$ and the products below are integrable, we change scale and then
differentiate at $s=1$ to obtain
\begin{equation}\label{typeii:eq:force-moment}
 \int F_s^5H=s^{k+\frac{1}{2}}\int F^5H,\qquad
 5\int F^4\Lambda F\,H=-\left(k+\frac{1}{2}\right)\int F^5H .
\end{equation}
We use this identity with the quadratic jet and its Kelvin transform in
\eqref{p2:eq:signcheck}, and with the incident dipole in \eqref{typeii:eq:linear-responses}.
\paragraph{From a harmonic moment to the scale and its linearization.}
If $j$ is the degree of the first interacting harmonic, we set $p=j+\frac{1}{2}$.
For child and parent scales $s,m$, the potential $-C\bigl(\frac{s}{m}\bigr)^p$ produces the force
\begin{equation}\label{p1:eq:leading-force}
 Ms''=-\partial_s[-C\bigl(\frac{s}{m}\bigr)^p]=pCs^{p-1}m^{-p}.
\end{equation}
With $s=B\tau^k$ and $m=A\tau^a$, matching powers and coefficients in
\eqref{p1:eq:leading-force} yields
\begin{equation}\label{p1:eq:power-balance}
 (2-p)k=2-pa,\qquad
 B^{2-p}=\frac{pC}{MA^p k(k-1)} .
\end{equation}
The second time derivative is the same for $\tau=t$ and $\tau=T-t$.
To linearize about $s_*=B\tau^k$ at fixed $m_*=A\tau^a$, we use
\eqref{p1:eq:power-balance} once more and find
\begin{equation*}
 \frac{1}{M}\,\partial_s\big(pCs^{p-1}m^{-p}\big)\big|_{(s_*,m_*)}
 =(p-1)k(k-1)\tau^{-2}.
\end{equation*}
Thus $\delta s=s_*\tau^{-\gamma}$ solves the homogeneous linearized
equation precisely when
\begin{equation}\label{p1:eq:indicial}
 (k-\gamma)(k-\gamma-1)=(p-1)k(k-1).
\end{equation}
The three choices in \eqref{p1:eq:power-balance} are
\begin{equation*}
 \begin{array}{c|ccc}
 &\text{global inner scale}&\text{finite-time middle scale}&\text{finite-time inner scale}\\[2pt]
 (p,a)&(\frac{5}{2},0)&(\frac{3}{2},0)&(\frac{5}{2},4)\\[2pt]
 k&-4&4&16\\
 (p-1)k(k-1)&30&6&360
 \end{array}
\end{equation*}
The coefficient identity determines $B=64\bigl(\frac{M}{C}\bigr)^2$
in the first column and the two constants of Theorem~\ref{typeii:seq-main}
in the other columns.
For the true scales we estimate the wave and force errors in Sections~\ref{sec:tower-5} and~\ref{sec:tower-6}.

\subsubsection{Center and dilation modes}\label{p1:sec:kernel}
In the crossing chart, the signed rotations $R(x)=(-x_2,x_1,x_3)$ and $S(x)=(x_1,-x_2,-x_3)$ have no common fixed direction. The quarter-turn excludes spherical harmonics of degrees zero and one. By \eqref{common:restricted},
\begin{equation}\label{p1:eq:kernel}
 \ker_{\dot H^1_\chi}(-\Delta-5F^4)=\mathbb R\Lambda F,\qquad
 F=P,Q,\qquad M=\|\Lambda P\|_2^2=\|\Lambda Q\|_2^2\in(0,\infty).
\end{equation}
The equality of the masses follows by reflection, and their finiteness by \eqref{p1:eq:kelvin-tail}.
\space
\subsection{Interaction energy and the reference path}\label{sec:tower-2}

\begingroup

\def\R{\mathbb R}
\def\Sp{\mathbb S}
\def\dd{\mathrm d}
\def\dx{\,\dd x}
\def\dy{\,\dd y}
\def\dw{\,\dd\omega}
\def\ip#1#2{\langle #1,#2\rangle}
\def\Os{O_{\mathrm{sym}}}
\def\one{\frac{1}{2}}

We compute the static energy and kinetic form on the pair of scales. Their Euler--Lagrange equations determine the first corrections to the law $\lambda\sim at^{-4}$.

\subsubsection{The two-scale result}

Let $S=Q_\mu+P_\lambda$ and $\rho=\frac{\lambda}{\mu}$, with the Kelvin pair
and energy density of \eqref{p8:eq:kelvin-pair} and \eqref{p8:eq:energy-density}.
We use the real $L^2(\R^3)$ pairing, the constants in \eqref{p8:eq:coefficients},
and the signed dilation kernels in \eqref{p1:eq:kernel}.

\begin{theorem}\label{p2:th:result}
The excess static energy $V=E(S,0)-2E(Q,0)$ and the kinetic cross coefficient $B=\ip{\partial_\lambda S}{\partial_\mu S}$ depend only on $\rho$. As $\rho\downarrow0$,
\begin{align}
 V(\rho)&=-C\rho^{\frac{5}{2}}\big(1+\frac{16}{5}\rho+\Os(\rho^2)\big),
 \label{p2:eq:potential}\\
 B(\rho)&=B_0\rho^{\frac{3}{2}}+\Os(\rho^{\frac{5}{2}}),\qquad
 B_0=-\frac{5}{4}\int_{\R^3}(\,|x|^{-3}H_2\bigl(\frac{x}{|x|}\bigr)\,)Q(x)\dx<0 .
 \label{p2:eq:mass}
\end{align}
The symbol $\Os(\rho^a)$ means an $O(\rho^a)$ bound after any fixed number of applications of $\rho\partial_\rho$. The bounds are uniform in the common scale.

For positive $C^2$ scales, the two conditions
$\ip{S_{tt}-\Delta S-S^5}{\partial_s S}=0$, $s=\lambda,\mu$, are exactly
\begin{align}
 M\ddot\lambda+B\ddot\mu
   &=\mu^{-1}(-V'+\rho B'\dot\mu^2), \label{p2:eq:inner}\\
 B\ddot\lambda+M\ddot\mu
   &=\mu^{-1}(\rho V'-B'\dot\lambda^2). \label{p2:eq:outer}
\end{align}
Primes on $V,B$ denote $\rho$-derivatives. We put $a=64\bigl(\frac{M}{C}\bigr)^2$ and $L=at^{-4}$. We define the reference path $q_*=(\lambda_*,\mu_*)$ by
\begin{equation}\label{p2:eq:reference}
 \lambda_*=L+\frac{32}{15}L^2,\qquad
 \mu_*=1-\frac{5}{18}L^2-\frac{2B_0}{5M}L^{\frac{5}{2}}.
\end{equation}
For this path both left sides of \eqref{p2:eq:inner}--\eqref{p2:eq:outer} differ from their right sides by $O(t^{-14})$.
\end{theorem}

The finite polynomial \eqref{p2:eq:reference} fixes the time origin and the limiting outer scale for the wave construction.

\subsubsection{Static interaction and differentiation}

\paragraph{The two harmonic moments.}
The crossing terms $H_2,H_3$ in \eqref{common:quad}--\eqref{p2:eq:jet} are homogeneous harmonic polynomials of degrees two and three.
Set $Y_j=H_j|_{\Sp^2}$ and $S_j(x)=r^{-j-1}Y_j(\frac{x}{r})$ for $j=2,3$. Kelvin inversion implies $P=S_2+S_3+O(r^{-5})$ at infinity, also after any fixed number of applications of $D=x\cdot\nabla$. In both charts $|D^kF|\lesssim_k\min(r^2,r^{-3})$ for $F=Q,P$. By \eqref{p1:eq:harmonic-flux},
\begin{equation}\label{p2:eq:moments}
 C_j:=\int H_jP^5\dx=(2j+1)\|Y_j\|_{L^2(\Sp^2)}^2,\qquad j=2,3 .
\end{equation}

\paragraph{Testing the stationary equations against the opposite profile}
By simultaneous dilation we may take $\mu=1$ and $\lambda=\rho$.
We test $-\Delta Q=Q^5$ against $P_\rho$ and
$-\Delta P_\rho=P_\rho^5$ against $Q$. Since both profiles belong to
$\dot H^1\cap L^6$, we obtain
\begin{equation}\label{p2:eq:mutual-test}
 I:=\int\nabla Q\cdot\nabla P_\rho\dx
 =\int Q^5P_\rho\dx=\int QP_\rho^5\dx
 =\rho^{\frac{1}{2}}\int P(y)^5Q(\rho y)\dy .
\end{equation}
Expanding the sextic potential and using \eqref{p2:eq:mutual-test} twice,
we obtain the identity
\begin{equation}\label{p2:eq:poly}
 -V-I=\frac{5}{6}\int Q^2P_\rho^2
 \big[Q^2+P_\rho^2+2(Q+P_\rho)^2\big]\dx .
\end{equation}
The Kelvin involution at radius $\sqrt\rho$ exchanges the profiles:
\begin{equation}\label{p2:eq:kelvin-exchange}
 (\mathcal K_{\sqrt\rho} f)(x)=\frac{\sqrt\rho}{|x|}
           f\!\left(\frac{\rho x}{|x|^2}\right),\qquad
 \mathcal K_{\sqrt\rho} Q=P_\rho,\quad \mathcal K_{\sqrt\rho}P_\rho=Q .
\end{equation}
Its conformal factor to the sixth power is the Jacobian. Thus
$\int Q^4P_\rho^2=\int Q^2P_\rho^4$, and \eqref{p2:eq:poly} requires only one overlap:
\begin{equation}\label{p2:eq:neck}
 0\le -V-I\lesssim\int Q^4P_\rho^2\dx
 \lesssim\rho^5\int |x|^{-6}Q^4\dx\lesssim\rho^5 .
\end{equation}
For the symbol bound, take $\rho\partial_\rho$ under the integrals. The differentiated inner factors obey the same envelope as $P_\rho$. The positive envelope $m_0(r)=\min(r^2,r^{-3})$ is itself Kelvin invariant at unit radius. Hence \eqref{p2:eq:kelvin-exchange} exchanges its two scaled copies as well, and the same bound applies to every fixed derivative.

We also record the localized energy identity for inversion in a sphere of radius $R$. Substitution of $r=R^2/s$ preserves the angular and sextic integrals. The radial gradient difference is $\frac{1}{2}\int_0^R\partial_s(sf^2)\,ds$ on each ray. Thus
\begin{equation}\label{p2:eq:kelvin-flux}
 \int_{|x|>R}e(\mathcal K_R f,0)\dx-
 \int_{|x|<R}e(f,0)\dx=\frac{R}{2}\int_{\mathbb S^2} f(R\omega)^2\,d\omega .
\end{equation}
The identity extends from smooth data to $\dot H^1$ by the trace bound $\sqrt R\,\|f(R\,\cdot)\|_2\le\|\nabla f\|_2$, obtained by radial integration.

Taylor's formula and $\int |P(y)|^5|y|^4\dy<\infty$ bound the remainder in \eqref{p2:eq:mutual-test}, including its dilation derivatives. Thus
\begin{equation}\label{p2:eq:energy-moments}
 I=C_2\rho^{\frac{5}{2}}+C_3\rho^{\frac{7}{2}}+\Os(\rho^{\frac{9}{2}}),
 \qquad V=-I+\Os(\rho^5).
\end{equation}
The spherical means of $\omega_1^2\omega_2^2$ and $\omega_3^2(\omega_2^2-\omega_1^2)^2$ are $\frac{1}{15}$ and $\frac{4}{105}$. Equations \eqref{p2:eq:jet} and \eqref{p2:eq:moments} imply $C_2=C$ and $\frac{C_3}{C}=\frac{16}{5}$, which proves \eqref{p2:eq:potential}. In particular, for any fixed $i,j\ge0$,
\begin{equation}\label{p2:eq:partial}
 \left|\partial_\lambda^i\partial_\mu^j
   (V+C\rho^{\frac{5}{2}}+C_3\rho^{\frac{7}{2}})\right|
   \lesssim_{i,j}\frac{\rho^{\frac{9}{2}}}{\lambda^i\mu^j}.
\end{equation}
Here $\lambda\partial_\lambda=\rho\partial_\rho$ and
$\mu\partial_\mu=-\rho\partial_\rho$. Consequently
\begin{equation}\label{p2:eq:forces}
 -\partial_\lambda V=\frac{5C}{2\mu}\rho^{\frac{3}{2}}
   \left(1+\frac{112}{25}\rho+O(\rho^2)\right),\qquad
 -\partial_\mu V=-\rho(-\partial_\lambda V).
\end{equation}

\subsubsection{The kinetic form and the two projections}

For either profile, the tangent and its mass satisfy
\begin{equation*}
 \partial_s F_s(x)=-s^{-\frac{3}{2}}(\Lambda F)\!\left(\frac{x}{s}\right),\qquad
 \|\partial_s F_s\|_2^2=\|\Lambda F\|_2^2=M .
\end{equation*}
Pairing the two scale velocities, we obtain
\begin{equation}\label{p2:eq:kinetic}
 T=\frac{M}{2}(\dot\lambda^2+\dot\mu^2)+B(\rho)\dot\lambda\dot\mu,\qquad
 B(\rho)=\rho^{-\frac{3}{2}}\int(\Lambda Q)(x)(\Lambda P)\bigl(\frac{x}{\rho}\bigr)\dx .
\end{equation}
We compute $B$ before differentiating either profile in space. Write $I_0(\rho)=\langle P_\rho,Q\rangle$, so $\langle P_\lambda,Q_\mu\rangle=\mu^2 I_0(\rho)$. Since $P=S_2+O(r^{-4})$ at infinity, with the same bound after Euler differentiation of the remainder, and $\int |x|^{-4}|Q(x)|\dx<\infty$,
\begin{equation}\label{p2:eq:kinetic-overlap}
 I_0(\rho)=\rho^{\frac{5}{2}}\int S_2 Q\dx+\Os(\rho^{\frac{7}{2}}),\qquad
 B=\partial_\mu\partial_\lambda[\mu^2 I_0(\rho)]=I_0'-\rho I_0'' .
\end{equation}
On $\rho^a$ the operator in \eqref{p2:eq:kinetic-overlap} acts by $a(2-a)\rho^{a-1}$. At $a=\frac{5}{2}$ the factor is $-\frac{5}{4}$, which proves \eqref{p2:eq:mass} except for its sign, established in Section~\ref{p2:sec:sign}. In particular $M^2-B^2>0$ for small $\rho$.

We put
\begin{equation*}
 q=(\lambda,\mu),\qquad G=\begin{pmatrix}M&B\\ B&M\end{pmatrix},\qquad
 p=Gq'=\partial_{q'}T,\qquad
 \nabla_q\rho=\frac{1}{\mu}(1,-\rho).
\end{equation*}
Since $T-V$ depends on $q$ only through $\rho$, both projected equations
take the form
\begin{equation}\label{p2:eq:momentum}
 \binom{\ip{\mathcal R(S)}{\partial_\lambda S}}{\ip{\mathcal R(S)}{\partial_\mu S}}
 =p'-\frac{B'(\rho)\dot\lambda\dot\mu-V'(\rho)}{\mu}
                  \binom{1}{-\rho},\qquad
 \mathcal R(S)=S_{tt}-\Delta S-S^5 .
\end{equation}
Indeed $\partial_t\ip{S_t}{\partial_s S}-\partial_s T
=\ip{S_{tt}}{\partial_s S}$, and $\partial_s V$ is the static pairing.
Expanding $p'$ in \eqref{p2:eq:momentum} recovers \eqref{p2:eq:inner}--\eqref{p2:eq:outer}. In
particular the projected dynamics conserve $T+V$ and obey
$p_o'=-\rho p_i'$. For use with arbitrary scale velocities the differentiated kinetic error is
\begin{equation}\label{p2:eq:bprime}
 B'=\frac{3}{2} B_0\rho^{\frac{1}{2}}+\Os(\rho^{\frac{3}{2}}).
\end{equation}
The errors in \eqref{p2:eq:forces} and \eqref{p2:eq:bprime} are uniform in the two velocities.

In the physical residual the leading mixed terms are $-5QP_\rho^4$ on
the inner scale and $-5Q^4P_\rho$ on the outer scale. Kelvin inversion
implies $\int Q^5S_2=\int P^5H_2=C$. We apply \eqref{typeii:eq:force-moment}
with degrees $2$ and $-3$, respectively, to obtain
\begin{equation}\label{p2:eq:signcheck}
 5\int P^4\Lambda P\,H_2\dx=-\frac{5}{2} C,\qquad
 5\int Q^4\Lambda Q\,S_2\dx=\frac{5}{2} C .
\end{equation}
The crossing estimates \eqref{p1:eq:kelvin-tail} justify differentiation
under both integrals.
In both charts the other mixed terms have smaller pairings by the differentiated estimates in \eqref{p2:eq:neck}. The factors converting these pairings to $\partial_\lambda V$ and $\partial_\mu V$ are $\frac{\rho^{\frac{3}{2}}}{\mu}$ and $\frac{\rho^{\frac{5}{2}}}{\mu}$, respectively.

\subsubsection{The reference path and its force defect}\label{p2:sec:branch}

Write $r_i=\ip{\mathcal R(S)}{\partial_\lambda S}$ and $r_o=\ip{\mathcal R(S)}{\partial_\mu S}$. The energy identity and the common dilation in \eqref{p2:eq:momentum} read
\begin{equation}\label{p2:eq:energy-dilation}
 (T+V)'=\lambda' r_i+\mu' r_o,\qquad r_o+\rho r_i=p_o'+\rho p_i'.
\end{equation}
We take $\lambda=L+dL^2$, $\mu=1+eL^2+fL^{\frac{5}{2}}$. Since $\frac{M}{2}(L')^2=CL^{\frac{5}{2}}$ and $(L^m)''=\frac{4m(4m+1)C}{8M}L^{m+\frac{1}{2}}$, equations \eqref{p2:eq:potential}, \eqref{p2:eq:mass} and \eqref{p2:eq:kinetic} imply
\begin{align*}
 T+V&=C\left(\frac{3d}{2}-\frac{16}{5}\right)L^{\frac{7}{2}}+O(L^{\frac{9}{2}}),\\
 p_o'+\rho p_i'&=C\left(9e+\frac{5}{2}\right)L^{\frac{5}{2}}
       +C\left(\frac{55f}{4}+\frac{11B_0}{2M}\right)L^3+O(L^{\frac{7}{2}}).
\end{align*}
In the second line we used $p_o=M\mu'+B\lambda'$ and $\rho p_i'=MLL''+O(L^{\frac{7}{2}})$. Vanishing of the three displayed coefficients determines precisely \eqref{p2:eq:reference}. Along this explicit path the symbol bounds differentiate the energy remainder to $O(t^{-19})$. Since $\lambda'-\rho\mu'\sim-4at^{-5}$, \eqref{p2:eq:energy-dilation} implies
\begin{equation}\label{p2:eq:defect}
 r_i=\frac{(T+V)'-\mu'(p_o'+\rho p_i')}{\lambda'-\rho\mu'}=O(t^{-14}),
 \qquad r_o=p_o'+\rho p_i'-\rho r_i=O(t^{-14}).
\end{equation}
Differentiating \eqref{p2:eq:reference}, we obtain
\begin{align}
 \dot\lambda_*&=-4at^{-5}(1+O(t^{-4})),&
 \ddot\lambda_*&=20at^{-6}(1+O(t^{-4})),\nonumber\\
 \dot\mu_*&=\frac{20}{9}a^2t^{-9}(1+O(t^{-2})),&
 \ddot\mu_*&=-20a^2t^{-10}(1+O(t^{-2})).
 \label{p2:eq:rates}
\end{align}

\subsubsection{The sign of the kinetic cross coefficient}\label{p2:sec:sign}

\paragraph{Positivity on a spherical lune}
By Section~\ref{sec:common-hopf}, $h>0$ on the open lune $\mathcal D=\{X>0,\ Y>0\}$ of the base sphere. Odd reflection across its edges implies $h(\theta,\phi)\sin(2\phi)>0$ off the meridians when $0<\theta<\pi$. The cosine moment vanishes by oddness under $\phi\mapsto-\phi$. Hence
\begin{equation}\label{p2:eq:angularsign}
 \int_0^{2\pi} h(\theta,\phi)\sin(2\phi)\,\dd\phi>0,\qquad
 \int_0^{2\pi} h(\theta,\phi)\cos(2\phi)\,\dd\phi=0 .
\end{equation}

\paragraph{The Kelvin integral}
We write $z_1=c e^{i\alpha},\,z_2=s e^{i\delta}$ with $c=\cos\chi,\,s=\sin\chi$, $0<\chi<\frac{\pi}{2}$. Then $\theta=2\chi$ and $\phi=\alpha-\delta$ on the base. We set $u=\Re z_1$, $w_1=\Re z_2$, $w_2=\Im z_2$. The inverse stereographic chart satisfies $x_j=\frac{w_j}{1-u}$ for $j=1,2$, $r^2=\frac{1+u}{1-u}$ and $\dx=(1-u)^{-3}\dd\sigma$. Thus
\begin{equation*}
 \int QS_2\dx=-b\int_{\Sp^3}w_1w_2\,h(2\chi,\alpha-\delta)\,
       W(c\cos\alpha)\,\dd\sigma,\qquad
 W(u)=(1-u)^{-2}(1+u)^{-\frac{5}{2}}.
\end{equation*}
The integral is absolutely convergent by the origin and far-field bounds for $Q$. The measure is $cs\,\dd\chi\,\dd\alpha\,\dd\delta$, and $w_1w_2=\frac{s^2}{2}\sin(2\delta)$. By \eqref{p2:eq:angularsign} the angular integral factors as
\begin{equation*}
 \int_0^{2\pi}\sin(2\delta)h(2\chi,\alpha-\delta)\,\dd\delta
 =-\cos(2\alpha)\int_0^{2\pi} h(2\chi,\phi)\sin(2\phi)\,\dd\phi .
\end{equation*}
For $0<c<1$, the even sum $W_e(u)=W(u)+W(-u)$ increases strictly on $[0,c]$, since $(\log W)''=\frac{2}{(1-u)^2}+\frac{5}{2(1+u)^2}>0$ and hence $W$ is strictly convex. Pairing the quadrants, then exchanging $\alpha$ and $\frac{\pi}{2}-\alpha$, we find
\begin{equation}\label{p2:eq:convex-moment}
 \int_0^{2\pi}\! W(c\cos\alpha)\cos(2\alpha)\,\dd\alpha
 =2\int_0^{\frac{\pi}{4}}\!\cos(2\alpha)
 [W_e(c\cos\alpha)-W_e(c\sin\alpha)]\,\dd\alpha>0 .
\end{equation}
Together with \eqref{p2:eq:angularsign}, this proves $\int QS_2\dx>0$ and hence $B_0<0$.
\endgroup
\space
\subsection{The inner correction}\label{sec:tower-3}

\begingroup
\def\R{\mathbb R}
\def\RR{\mathcal R}
\def\ip#1#2{\langle #1,#2\rangle}
\def\Tr{\mathcal T_\rho}

The reference scales cancel the two dilation forces, but the full
residual still has an inner term of order $\rho^{\frac{3}{2}}$ in $L^2$.
We remove it by solving one elliptic equation on the inner profile.

We set $\RR(u)=u_{tt}-\Delta u-u^5$, $S=P_\lambda+Q_\mu$ and
$\rho=\frac{\lambda}{\mu}$. For $z=\frac{x}{\mu}$, $y=\frac{x}{\lambda}$, we set
$\widehat f(z)=\mu^{\frac{5}{2}} f(\mu z)$ and $\Tr g(z)=\rho^{-\frac{5}{2}}g\bigl(\frac{z}{\rho}\bigr)$.
Thus $\|f\|_2=\mu^{-1}\|\widehat f\|_2$ and
$\|\Tr g\|_2=\rho^{-1}\|g\|_2$. All scales in this section have $\mu\asymp1$.
We write $\Lambda=\frac{1}{2}+x\cdot\nabla$, $K=(1+\Lambda)\Lambda$,
$L_F=-\Delta-5F^4$ and set
\begin{equation}\label{p3:eq:j}
 p=\frac{5}{2},\qquad \alpha=\frac{pC}{M},\qquad \beta=\frac{2C}{M},\qquad
 J=-\alpha\Lambda P+\beta KP-5P^4H_2 .
\end{equation}
Let $\chi$ be a smooth radial cutoff equal to one for $r\le1$ and zero
for $r\ge2$. We shall use the field
\begin{equation}\label{p3:eq:ansatz}
 U=S+w,\qquad w(t,x)=\mu^{-\frac{1}{2}}\rho^2\chi(|z|)W\bigl(\frac{z}{\rho}\bigr),
 \qquad L_P W=-J .
\end{equation}
We write $e_i=\partial_\lambda P_\lambda$ and $e_o=\partial_\mu Q_\mu$.

\begin{theorem}\label{p3:th:main}
On the reference path \eqref{p2:eq:reference}, the correction in
\eqref{p3:eq:ansatz} can be chosen in the signed sector so that
\begin{align}
 \|\RR(U)(t)\|_2&\lesssim t^{-10},\label{p3:eq:mainbound}\\
 |\ip{\RR(U)(t)}{e_i}|+|\ip{\RR(U)(t)}{e_o}|&\lesssim t^{-14},
     \label{p3:eq:proj}\\
 \|\nabla w(t)\|_2&\lesssim t^{-10},\qquad
 \|w_t(t)\|_2\lesssim t^{-13}.\label{p3:eq:energybounds}
\end{align}
\end{theorem}

\subsubsection{The inner source and its elliptic inverse}

The same two exponents control all mixed products in the tower.
In outer coordinates, let $I_\rho(z)=\rho^{-\frac{1}{2}}I\bigl(\frac{z}{\rho}\bigr)$ and $O(z)$
be profiles bounded at unit scale by $m_0(r)=\min(r^2,r^{-3})$.
For $a,b>0$, $1\le p<\infty$ and a radial weight $r^{-\sigma}$, suppose $2(a+b)-\sigma>-\frac{3}{p}$ and $3(a+b)+\sigma>\frac{3}{p}$. Set $u=\frac{3}{p}-\sigma-\frac{a}{2}+2b$ and $v=\frac{5a}{2}$. Then
\begin{equation}\label{p3:eq:overlap}
 \big\|r^{-\sigma}|I_\rho|^a|O|^b\big\|_p
 \lesssim \rho^{\min(u,v)}
 \begin{cases}1,&u\ne v,\\ (1+|\log\rho|)^{\frac{1}{p}},&u=v.\end{cases}
\end{equation}
Using the two pieces of $m_0$ on $r<\rho$, $\rho<r<1$ and $r>1$,
we obtain
\begin{equation*}
 \big\|r^{-\sigma}|I_\rho|^a|O|^b\big\|_p^p
 \lesssim \rho^{pu}+\rho^{pv}\!\int_\rho^1 r^{p(u-v)}\frac{\mathrm{d}r}{r}
             +\rho^{pv}.
\end{equation*}
We take $\sigma=0$ unless stated otherwise. Only $u=v$ produces a logarithm in \eqref{p3:eq:overlap}.
Different factors may be different profiles with this envelope.
In particular the same rule applies after logarithmic scale
differentiations of products of whole powers, by the product rule.
A physical $L^2$-normalized test at scale $s$ carries one additional
factor $s^{-1}$ relative to an energy-scaled profile.

With $p,\alpha,\beta$ as in \eqref{p3:eq:j}, the reference path
\eqref{p2:eq:reference} satisfies
\begin{equation}\label{p3:eq:accelerations}
 \lambda\ddot\lambda=\alpha\rho^p+O(\rho^{p+1}),\quad
 \dot\lambda^2=\beta\rho^p+O(\rho^{p+1}),\quad
 |\mu\ddot\mu|+\dot\mu^2\lesssim\rho^p .
\end{equation}
We used \eqref{p2:eq:rates}. In outer units the inner acceleration is $\Tr(-\lambda\ddot\lambda\Lambda P+\dot\lambda^2 KP)$. In the term $-5P_\rho^4Q$ we use
\begin{equation*}
 Q(\rho y)=\rho^2H_2(y)+O(\rho^3|y|^3),\qquad
 \|P_\rho^4(Q-H_2)\|_{L^2_z}
 \lesssim \rho^{\frac{5}{2}}\||y|^3P^4\|_2 .
\end{equation*}
For the remaining powers, \eqref{p3:eq:overlap}
with $b=5-a$ reads
\begin{equation*}
 \|P_\rho^3 Q^2\|_2\lesssim\rho^4,\qquad
 \|P_\rho^2 Q^3\|_2\lesssim\rho^5,\qquad
 \|P_\rho Q^4\|_2\lesssim\rho^{\frac{5}{2}} .
\end{equation*}
Combining these three bounds with \eqref{p3:eq:accelerations}, we obtain
\begin{equation}\label{p3:eq:raw}
 \|\widehat{\RR(S)}-\rho^{\frac{5}{2}}\Tr J\|_{L^2_z}\lesssim\rho^{\frac{5}{2}}.
\end{equation}
The leading term has norm $\rho^{\frac{3}{2}}\|J\|_2$. At infinity we have
\begin{equation*}
 J(y)=\frac{55C}{4M}|y|^{-3}H_2\!\left(\frac{y}{|y|}\right)+O(|y|^{-4}).
\end{equation*}
\begin{lemma}\label{p3:lem:ell}
There is a signed energy solution $W$ of $L_P W=-J$. Its Euler
derivatives of every fixed order are bounded by $C\min(|y|^2,|y|^{-1})$.
\end{lemma}
\begin{proof}
We separate the homogeneous tail before solving the elliptic equation. This also treats the dipole in Section~\ref{typeii:sec:stationarycells}. Let $F$ be either signed profile, with leading tail $r^{-k-1}Y_k(\omega)$, and let $H_k$ be the incident harmonic polynomial of degree $k$. Here $k=2$ for the crossing and $k=1$ for the wall. Write $p_k=k+\frac{1}{2}$ and $S_F=\frac{1}{4}|y|^2\Lambda F$. From $L_F\Lambda F=0$ and $[-\Delta,|y|^2]=-6-4y\cdot\nabla$ we obtain
\begin{equation}\label{p3:eq:homogeneous-tail}
 L_FS_F=-KF,\qquad
 \left(\Lambda+\frac{K}{p_k-1}\right)F
 =\frac{(\Lambda+p_k)\Lambda F}{p_k-1}=O(r^{-k-2}).
\end{equation}
The last Euler factor kills the leading harmonic. To solve $L_F W=\alpha\Lambda F-\beta KF+5F^4H_k$, we set
\begin{equation}\label{p3:eq:shared-response}
 W=\left(\beta+\frac{\alpha}{p_k-1}\right)S_F+T,\qquad
 L_FT=\alpha\left(\Lambda+\frac{K}{p_k-1}\right)F+5F^4H_k .
\end{equation}
The new source belongs to $L^{\frac{6}{5}}$. Since $1+\Lambda$ is skew-adjoint in $L^2$, we have $\langle KF,\Lambda F\rangle=0$. The other pairing vanishes precisely when $\alpha\|\Lambda F\|_2^2=p_k\int F^5H_k$, by integration by parts. Under this condition $(-\Delta)^{-1}L_F$ is a compact perturbation of the identity on the signed energy space, with kernel $\mathbb R\Lambda F$ by \eqref{common:finite-kernel}. By the Fredholm alternative an energy solution $T$ exists and is unique modulo this kernel.
We estimate $T$ by its Newton potential. Write $g=-\Delta T$. Local elliptic estimates on unit balls, $T\in L^6$ and the bounded coefficients imply $T\in L^\infty$. Thus $g=O(\langle y\rangle^{-k-2})$ by \eqref{p3:eq:shared-response}. Since $g\in L^{\frac{6}{5}}$, its Newton potential is the unique $\dot H^1$ solution. In the signed sector the moments of $g$ of spherical degree less than $k$ vanish on every sphere. For $r=|x|\ge2$ subtract the Taylor polynomial of $|x-y|^{-1}$ through degree $k-1$ on $|y|<\frac{r}{2}$. The constant, and for $k=2$ the linear term, integrate to zero. We obtain
\begin{equation}\label{p3:eq:newton-moments}
 |T(x)|\lesssim r^{-k-1}\int_{|y|<\frac{r}{2}}|y|^k|g(y)|\,dy
       +\int_{|y|\ge\frac{r}{2}}\frac{|g(y)|}{|x-y|}\,dy
 \lesssim r^{-k}.
\end{equation}
The first integral before its prefactor is $O(r)$, and dilation by $r$ bounds the second. Rescaled elliptic estimates give the corresponding derivative bounds.
For $F=P$ the moment condition is \eqref{p2:eq:signcheck}. The quarter-turn excludes degrees zero and one, so \eqref{p3:eq:newton-moments} with $k=2$ implies $T=O(r^{-2})$. Thus
\begin{equation}\label{p3:eq:crossing-response}
 V_\times=W-\beta S_P=\frac{2}{3}\alpha S_P+O(r^{-2}),\qquad
 L_PV_\times=\alpha\Lambda P+5P^4H_2 .
\end{equation}
The term $S_P$ belongs to $\dot H^1$, with tail $O(r^{-1})$. We may fix the remaining multiple of $\Lambda P$ by taking $V_\times$ gradient-orthogonal to it. At zero the signed symmetries remove the constant and linear Taylor terms. The differentiated stationary expansions and rescaled elliptic estimates prove the asserted Euler bounds.
\end{proof}

\subsubsection{The residual and its two projections}

We first collect the spatial residual on an arbitrary pair of nearby scales. Its two leading terms will also control differences of paths in Section~\ref{p5:sec:tube}. Let $\mathcal D_i=\partial_{\log\lambda}$, $\mathcal D_o=\partial_{\log\mu}$
at fixed physical $x$. With $m(s)=\min(s^2,s^{-1})$ and $r=|z|$,
Lemma~\ref{p3:lem:ell} bounds the correction, also after one spatial Euler derivative, by
\begin{equation}\label{p3:eq:env}
 |\mathcal D^\nu w(t,\mu z)|\lesssim_\nu
    \rho^2 m\bigl(\frac{r}{\rho}\bigr)\mathbf 1_{r\le2}.
\end{equation}
In inner units the gradient and Laplacian of $W$ have square-integrable
tails. The $L^2$ norm of $W$ up to radius $\frac{2}{\rho}$ is $O(\rho^{-\frac{1}{2}})$.
After rescaling and estimating the cutoff derivatives on $r\asymp1$, we obtain
\begin{equation}\label{p3:eq:symbols}
 \|\mathcal D^\nu w\|_6+\|\nabla\mathcal D^\nu w\|_2\lesssim_\nu\rho^{\frac{5}{2}},
 \quad \|\mathcal D^\nu w\|_2\lesssim_\nu\rho^3,\quad
 \|\Delta\mathcal D^\nu w\|_2\lesssim_\nu\frac{\rho^{\frac{5}{2}}}{\lambda} .
\end{equation}
On the reference path the first two logarithmic time derivatives of each
scale are $O(t^{-1})$ and $O(t^{-2})$. Thus $\|w_t\|_2\lesssim t^{-13}$
and $\|w_{tt}\|_2\lesssim t^{-14}$. The same bounds hold on paths with
the corresponding size and derivative bounds.

Before applying the cutoff we have
\begin{equation*}
 L_{P_\rho}\!\left[\rho^2W\!\left(\frac{z}{\rho}\right)\right]
    =-\rho^{\frac{5}{2}}\Tr J .
\end{equation*}
On $r\ge1$ the omitted source is $O(\rho^3 r^{-3})$. The two cutoff
commutators are $O(\rho^3)$ in $L_z^2$ by the annular derivative bounds.
We put $A=P_\rho+Q$, $v=\mu^{\frac{1}{2}}w(t,\mu z)$ and $N(A,v)=(A+v)^5-A^5-5A^4v$. Then
\begin{align*}
 \|(A^4-P_\rho^4)v\|_2
 &\lesssim\|P_\rho^3Qv\|_2+\|Q^4v\|_2\lesssim\rho^4+\rho^3,\\
 \|N(A,v)\|_2
 &\lesssim\||A|^3v^2\|_2+\|v^5\|_2\lesssim\rho^4 .
\end{align*}
In the first line we use the inner bound on $Q$ and
$|v|\lesssim\min(\rho^2,\frac{\rho^3}{r})$. For the second line we use
$\|A\|_6\lesssim1$, $\|v\|_\infty\lesssim\rho^2$ and
\eqref{p3:eq:symbols}. By the stationary equations, the leading term of $-\Delta S-S^5$ in outer units is $-5\rho^{\frac{5}{2}}\Tr(P^4H_2)$. The elliptic correction contributes $-\rho^{\frac{5}{2}}\Tr J$. The terms in $P^4H_2$ cancel by \eqref{p3:eq:j}. Write $g^{[s]}(x)=s^{-\frac{3}{2}}g(\frac{x}{s})$. In physical units, for $\mu\asymp1$,
\begin{equation}\label{p3:eq:spatial-remainder}
 -\Delta U-U^5=\frac{\rho^{\frac{5}{2}}}{\lambda}
       (\alpha\Lambda P-\beta KP)^{[\lambda]}+R(\lambda,\mu),\qquad
 \|\mathcal D^\nu R\|_2\lesssim_\nu\rho^{\frac{5}{2}}\quad(|\nu|\le1).
\end{equation}
The bound follows from the mixed products above, the Taylor remainder of $Q$ and the cutoff commutators. Each logarithmic derivative preserves these bounds: it differentiates a smooth shape by an Euler derivative or a cutoff on $|z|\asymp1$. Adding the accelerations of the two cores, we obtain
\begin{equation}\label{p3:eq:residual-decomposition}
\begin{aligned}
 \RR(U)={}&\left(\lambda''-\frac{\alpha\rho^{\frac{5}{2}}}{\lambda}\right)e_i
       +\frac{(\lambda')^2-\beta\rho^{\frac{5}{2}}}{\lambda}(KP)^{[\lambda]}\\
 &+\mu'' e_o+\frac{(\mu')^2}{\mu}(KQ)^{[\mu]}+w_{tt}+R(\lambda,\mu).
\end{aligned}
\end{equation}
On the reference path, \eqref{p3:eq:accelerations} bounds the first two coefficients by $O(\rho^{\frac{5}{2}})$. Hence \eqref{p3:eq:residual-decomposition} proves \eqref{p3:eq:mainbound}.
For $F_i=P,F_o=Q$ and $s_i=\lambda,s_o=\mu$ we use
$L_{(F_j)_{s_j}}e_j=0$. On every path, integration by parts yields
\begin{equation}\label{p3:eq:projid}
 \ip{\RR(U)-\RR(S)}{e_j}
 =\ip{w_{tt}}{e_j}
 -5\ip{w}{(S^4-(F_j)_{s_j}^4)e_j}-\ip{N(S,w)}{e_j}.
\end{equation}
By \eqref{p3:eq:overlap} the mixed quartics are $O(\rho^{\frac{5}{2}})$
in $L^{\frac{3}{2}}$. The opposite pure well times a test uses $(a,b)=(1,4)$
for the inner row and $(4,1)$ for the outer row, at $p=\frac{6}{5}$.
After normalizing the tests, these two products are
$O(\rho^{\frac{3}{2}})$ and $O(\rho^{\frac{5}{2}})$, respectively.
Since $\|e_j\|_6\lesssim s_j^{-1}$, the full coefficient in
\eqref{p3:eq:projid} satisfies, also after one logarithmic scale derivative,
\begin{equation*}
 \|(S^4-(F_j)_{s_j}^4)e_j\|_{\frac{6}{5}}
 \lesssim \begin{cases}t^{-6},&j=i,\\ t^{-10},&j=o.\end{cases}
\end{equation*}
By \eqref{p3:eq:symbols}, $\|w\|_6\lesssim t^{-10}$. Thus the middle
pairing is $O(t^{-16})$ in both rows. For the last pairing we use
\begin{equation*}
 |\ip{N(S,w)}{e_j}|\lesssim s_j^{-1}(\|w\|_6^2+\|w\|_6^5)
 \lesssim t^{-16}.
\end{equation*}
Finally $\|e_j\|_2^2=M$, so the inertial pairing is $O(t^{-14})$.
The bare forces on the reference path are the two defects computed in
\eqref{p2:eq:defect}. This proves \eqref{p3:eq:proj} and completes
the proof of Theorem~\ref{p3:th:main}.

The same path has the energy accuracy
\begin{equation}\label{p3:eq:energydefect}
 E(U,U_t)=2E(Q,0)+O(t^{-18}).
\end{equation}
Indeed $U_t=\lambda' e_i+\mu' e_o+w_t$. Testing the residual against $U_t$ and applying \eqref{p3:eq:mainbound}--\eqref{p3:eq:energybounds} yields
\begin{equation}\label{p3:eq:energy-work}
 \left|\frac{d}{dt}E(U,U_t)\right|
 =\left|\lambda'\ip{\RR(U)}{e_i}+\mu'\ip{\RR(U)}{e_o}+\ip{\RR(U)}{w_t}\right|
 \lesssim t^{-5}t^{-14}+t^{-10}t^{-13}\lesssim t^{-19}.
\end{equation}
By \eqref{p2:eq:potential}, $E(S,0)\to2E(Q,0)$. Since $\|\nabla w\|_2+\|U_t\|_2\to0$, $E(U,U_t)$ has the same limit. Integration of \eqref{p3:eq:energy-work} from $t$ to infinity proves \eqref{p3:eq:energydefect}.
\endgroup
\space
\section{Wave correction, scale selection and identification of the tower}\label{merge:global-wave}
We first solve for the wave orthogonal to the two dilation modes
in Sections~\ref{sec:tower-4} and~\ref{sec:tower-5}. Section~\ref{sec:tower-6}
then varies the scales to cancel the two remaining forces.
Section~\ref{sec:tower-8} identifies the energy and both profiles of the resulting solution.

\subsection{Energy estimates for the two wells}\label{sec:tower-4}

\begingroup

\def\R{\mathbb R}
\def\Sp{\mathbb S}
\def\dd{\mathrm d}
\def\ip#1#2{\langle #1,#2\rangle}
\def\RR{\mathcal R}
\def\En{\dot H^1\times L^2}
\def\N{\mathcal N}
\def\Trho{\mathcal T_\rho}
\def\dil{\Lambda}
\def\dx{\,\dd x}

The wave equation on the complement of the two dilations has a moving
potential. We first prove uniform coercivity and a backward energy
estimate, then apply both to the residual of \eqref{p3:eq:ansatz}.

\subsubsection{The first wave response}

Let $S=P_\lambda+Q_\mu$ and $U=S+w$, with the inner correction
\eqref{p3:eq:ansatz}. We write $s_i=\lambda$, $s_o=\mu$, $F_i=P$, $F_o=Q$,
$e_j=\partial_{s_j}(F_j)_{s_j}$, $V=5U^4$ and $L(t)=-\Delta-V$. First
take the explicit reference path \eqref{p2:eq:reference}. All fields have
the two signed symmetries of Section~\ref{sec:tower-1}. We set
$\RR(u)=u_{tt}-\Delta u-u^5$.

\begin{theorem}\label{p4:th:main}
On $[T_0,\infty)$ for $T_0$ sufficiently large there is a selected signed solution of
\begin{align}
 h_{tt}+L(t)h&=-\RR(U)+c_i e_i+c_o e_o,\qquad
 \ip{h}{e_i}=\ip{h}{e_o}=0 . \label{p4:eq:wave}
\end{align}
It has zero future free scattering state and obeys
\begin{equation}\label{p4:eq:rates}
 X(t):=\|(h,h_t)(t)\|_{\En}\lesssim t^{-9},\qquad
 |c_i(t)|\lesssim t^{-10},\qquad |c_o(t)|\lesssim t^{-14}.
\end{equation}
Only the dilation forces appear in \eqref{p4:eq:wave}. All negative eigenmodes are retained in the evolution, with the endpoint selection \eqref{p5:eq:selection}.
The first response satisfies
\begin{align}
 |E(U+h,U_t+h_t)-E(U,U_t)|&\lesssim t^{-18},\label{p4:eq:Egoal}\\
 |\ip{\RR(U+h)}{\partial_\lambda U}|&\lesssim t^{-10},\qquad
 |\ip{\RR(U+h)}{\partial_\mu U}|\lesssim t^{-14}. \label{p4:eq:projgoal}
\end{align}
The same projection bounds hold for $e_i,e_o$. All constants are independent of a late terminal cutoff.
\end{theorem}

On the reference path this field is precisely the linear response
$h_{q_*}^{(1)}$ in Corollary~\ref{p5:cor:nonlin}, equivalently the
first iterate of \eqref{p5:eq:compare} starting from zero. The inverse uses the energy estimate \eqref{p4:eq:Hprime}.

\subsubsection{Inputs and estimates for the moving potential}\label{p4:sec:inputs}

For the energy estimate we allow any $C^2$ path with $\lambda'\le0$ and
\begin{equation}\label{p4:eq:path}
 \lambda\asymp t^{-4},\quad\mu\asymp1,\quad
 |\lambda^{(m)}|\lesssim t^{-4-m},\quad
 |\mu^{(m)}|\lesssim t^{-8-m}\qquad(m=1,2).
\end{equation}
In particular $\rho=\frac{\lambda}{\mu}\asymp t^{-4}$. The constants in Lemma~\ref{p4:lem:pot}
depend only on the constants in \eqref{p4:eq:path} and on the fixed
profiles. We set $V_j=5(F_j)_{s_j}^4$. Each signed one-core energy kernel
is spanned by dilation. The stationary profiles and their dilation
derivatives are bounded by a constant times $\min(r^2,r^{-3})$.

We use the $L^2$ isometry
\begin{equation*}
 \psi^{[s]}(x)=s^{-\frac{3}{2}}\psi\!\left(\frac{x}{s}\right),\qquad
 e_j=(-\dil F_j)^{[s_j]} .
\end{equation*} In the following estimates $\psi_j$ denotes either $e_j$ or the $L^2$ rescaling of a normalized one-core negative eigenfunction. The negative eigenfunctions are smooth and exponentially decreasing with their first two derivatives. The signed symmetries make their value and gradient vanish at the origin. Thus all the tests and their first two dilation derivatives have the same two-sided envelope as the stationary profiles before scaling.

\begin{lemma}\label{p4:lem:pot}
With $D=x\cdot\nabla$, set $V_m=V-V_i-V_o$. Uniformly for large $t$,
\begin{align}
 \|V_m\|_{\frac{3}{2}}&\lesssim\rho^{\frac{5}{2}}, &
 \|V_m\|_3&\lesssim\rho^{\frac{3}{2}},\label{p4:eq:Vmix}\\
 \|(V-V_i)\psi_i\|_{\frac{6}{5}}&\lesssim t^{-6},&
 \|(V-V_o)\psi_o\|_{\frac{6}{5}}&\lesssim t^{-10}.\label{p4:eq:almostmode}
\end{align}
Moreover $|\ip{\psi_i}{\psi_o}|\lesssim\rho^{\frac{3}{2}}$ and
\begin{equation}\label{p4:eq:td}
 \|\dot\psi_j\|_2\lesssim \frac{|s_j'|}{s_j},\quad
 \|\dot\psi_j\|_{\frac{6}{5}}\lesssim |s_j'|,\quad
 \|\ddot\psi_j\|_{\frac{6}{5}}\lesssim |s_j''|+\frac{|s_j'|^2}{s_j} .
\end{equation}
\end{lemma}
\begin{proof}
In outer coordinates $\mu\asymp1$. For $p=\frac{3}{2},3$, \eqref{p3:eq:overlap} bounds every mixed monomial in $(Q+P_\rho)^4$ by interpolating between $P_\rho^3Q$ and $P_\rho Q^3$. In particular,
\begin{equation*}
 \|P_\rho^3 Q\|_3\lesssim\rho^{\frac{3}{2}}
    \||y|^2P^3\|_3,\qquad
 \|P_\rho Q^3\|_3\lesssim\rho^{\frac{5}{2}}\||x|^{-3}Q^3\|_3 .
\end{equation*}
The corresponding $L^{\frac{3}{2}}$ bounds are $O(\rho^{\frac{5}{2}})$. By \eqref{p3:eq:env}, $\|w\|_6\lesssim\rho^{\frac{5}{2}}$ and $\|w\|_{12}\lesssim\rho^{\frac{9}{4}}$. For the quartic terms containing $w$, H\"older's inequality with three further $L^6$ or $L^{12}$ factors proves \eqref{p4:eq:Vmix}, since $\|P_\rho\|_6\lesssim1$ and $\|P_\rho\|_{12}\lesssim\rho^{-\frac{1}{4}}$.

For \eqref{p4:eq:almostmode}, H\"older's inequality and the two bounds used in \eqref{p3:eq:projid} imply
\begin{align*}
 \|V_m\psi_j\|_{\frac{6}{5}}&\le\|V_m\|_{\frac{3}{2}}\|\psi_j\|_6
           \lesssim \frac{\rho^{\frac{5}{2}}}{s_j},\\
 \|V_o\psi_i\|_{\frac{6}{5}}&\lesssim\rho^{\frac{3}{2}},\qquad
 \|V_i\psi_o\|_{\frac{6}{5}}\lesssim\rho^{\frac{5}{2}} .
\end{align*}
We use $\mu\asymp1$ and $\lambda\asymp\rho$ in the first line. For $\ip{\psi_i}{\psi_o}$ we apply \eqref{p3:eq:overlap} with $a=b=p=1$, including the extra factor $O(\rho^{-1})$ in the inner $L^2$ rescaling.
These bounds also hold for a logarithmic
derivative of the spatial products, as all test profiles and their
Euler derivatives obey the same envelope. Finally \eqref{p4:eq:td} follows by differentiating the $L^2$ rescaling, whose $L^{\frac{6}{5}}$ norm carries one factor of $s_j$.
\end{proof}

\subsubsection{Static coercivity and the two constraints}\label{p4:sec:coer}

We first remove the negative modes and the dilation of each separate well. Let $Y_{j,k}$ be a real orthonormal basis of negative eigenfunctions of $-\Delta-5F_j^4$ in the signed $L^2$ sector, with eigenvalues $-\nu_{j,k}^2<0$, counted with multiplicity. We use $Y_{j,k}^{[s_j]}$ at scale $s_j$. The family is finite. Indeed multiplication by $5F_j^4$ defines a compact form on $\dot H^1$ by $L^{\frac{3}{2}}$ approximation and local compactness, so the negative index is finite. The decaying bounded potential has essential spectrum $[0,\infty)$ in the full space and no negative essential spectrum in the sector. Elliptic regularity at negative spectral parameter implies the decay used in Lemma~\ref{p4:lem:pot}.

Let $\mathcal B_t(h)=\int(|\nabla h|^2-Vh^2)\dx$. There are constants $\delta>0$, $K_0<\infty$ such that
\begin{equation}\label{p4:eq:coer}
 \mathcal B_t(h)\ge\delta\|\nabla h\|_2^2
 -K_0\sum_{j=i,o}s_j^{-2}
  \left(|\ip{h}{e_j}|^2+\sum_k|\ip{h}{Y_{j,k}^{[s_j]}}|^2\right).
\end{equation}
The same compactness argument treats one or two cores and any continuous unit-scale displacement test $\ell_j$ with $\ell_j(\Lambda F_j)\ne0$. On the complement of the negative eigenfunctions the single-well form is nonnegative. This follows from the spectral theorem on $H^1$, then by approximation on $\dot H^1$, subtracting the finitely many negative projections. A vector of zero form in this complement pairs to zero in the form with the whole complement by polarization, and with each negative eigenfunction by orthogonality. Hence
\begin{equation}\label{p4:eq:one-core-gap}
 g\perp\{Y_{j,k}\}_k,\quad \langle L_{F_j}g,g\rangle=0,\quad \ell_j(g)=0
 \quad\Longrightarrow\quad g\in\mathbb R\Lambda F_j,\quad g=0 .
\end{equation}

If the penalized estimate \eqref{p4:eq:coer} fails uniformly for small $\rho$, we can choose $\rho_n\to0$ and $\|\nabla h_n\|_2=1$ so that the form plus $n$ times the sum of the squared unit-scale observations is at most $1/n$. The forms are uniformly bounded below by Sobolev, so every observation tends to zero. After extraction $s_{j,n}^{\frac{1}{2}}h_n(s_{j,n}\,\cdot)\rightharpoonup g_j$ in $\dot H^1$. Thus $\ell_j(g_j)=0$ and $g_j\perp Y_{j,k}$. By frame orthogonality \eqref{fp:dislocation}, compactness of each potential form, and \eqref{p4:eq:Vmix},
\begin{equation}\label{p4:eq:frame-gap}
 \lim_n\mathcal B_{t_n}(h_n)
 =\left(1-\sum_j\|\nabla g_j\|_2^2\right)
     +\sum_j\langle L_{F_j}g_j,g_j\rangle\le0 .
\end{equation}
The parenthesis is nonnegative by expansion of $\|h_n-\sum_j(g_j)_{s_{j,n}}\|_{\dot H^1}^2$. The remaining terms are nonnegative as well. By \eqref{p4:eq:one-core-gap} each $g_j=0$, whereas vanishing of the parenthesis requires their squared norms to sum to one. This proves \eqref{p4:eq:coer}. For a single well the same argument uses one weak limit.

On a finite interval replace $-\RR(U)$ in \eqref{p4:eq:wave} by a given continuous signed $L^2$ source $f$, and put $v=h_t$, $X=\|(\nabla h,v)\|_2$. Two derivatives of $\ip{h}{e_j}=0$ determine the reactions,
\begin{equation}\label{p4:eq:reaction}
 \sum_{\ell=i,o}G_{j\ell}c_\ell
 =\ip{h}{L(t)e_j}-\ip{f}{e_j}
    -2\ip{v}{\dot e_j}-\ip{h}{\ddot e_j},
 \qquad G_{j\ell}=\ip{e_j}{e_\ell}.
\end{equation}
Here $G_{ii}=G_{oo}=M$ and $G_{io}=O(t^{-6})$. By Lemma~\ref{p4:lem:pot},
\begin{equation}\label{p4:eq:gencreact}
 |c_i|+|c_o|\lesssim\|f\|_2+t^{-1}X,\qquad
 |\ip{v}{e_j}|=|\ip{h}{\dot e_j}|\lesssim |s_j'|X .
\end{equation}
We have used $L(t)e_j=-(V-V_j)e_j$. Equation \eqref{p4:eq:reaction} defines a bounded finite-rank feedback on each compact time interval. Since $V(t)\in L^3$, the free wave Duhamel formula yields energy well-posedness in both time directions. Imposing the position constraint and its first derivative at one endpoint preserves it throughout by \eqref{p4:eq:reaction}.

\subsubsection{A backward energy estimate for the concentrating well}\label{p4:sec:virial}

In the wave equation we put the mixed potential on the right. By \eqref{p4:eq:Vmix} and Sobolev,
\begin{equation*}
 h_{tt}=\Delta h+(V_i+V_o)h+f+c_i e_i+c_o e_o+V_m h,\qquad
 \|V_m h\|_2\lesssim t^{-6}X .
\end{equation*}
We use Jendrej's localized dilation multiplier for two bubbles \cite[Section~3]{J}, which adapts the energy--virial method of Rapha\"el and Szeftel \cite{RaphaelSzeftel}. Choose a smooth $0\le\chi\le1$ equal to one on $(-\infty,\frac{1}{3}]$ and zero on $[\frac{2}{3},\infty)$. Set
\begin{equation}\label{p4:eq:neck-cutoffs}
 u=\frac{\log(r/\lambda)}{\log(\mu/\lambda)},\qquad
 c=\frac{\lambda'}{\lambda}\chi(u),\qquad
 r_-=\lambda\rho^{-\frac{1}{3}},\quad r_+=\mu\rho^{\frac{1}{3}} .
\end{equation}
The transition lies in $[r_-,r_+]$. With $\Lambda=\frac{1}{2}+D$, we define
\begin{equation}\label{p4:eq:H}
 \mathcal H=\frac{1}{2}\left(\|v\|_2^2+\|\nabla h\|_2^2
       -\int(V_i+V_o)h^2\dx\right)+\ip{v}{c\Lambda h}.
\end{equation}
Hardy's inequality on $r\le r_+$ implies $\|c\Lambda h\|_2\lesssim t^{-1}r_+X$, and \eqref{p4:eq:Vmix} bounds the change in the potential form by $O(\rho^{\frac{5}{2}}X^2)$. Thus \eqref{p4:eq:coer} implies, for sufficiently late $T_0$,
\begin{equation}\label{p4:eq:Hcoer}
 X^2\le C_0\mathcal H+C_1\sum_{j,k'}s_j^{-2}
           |\ip{h}{Y_{j,k'}^{[s_j]}}|^2,\qquad
 |\mathcal H|\lesssim X^2 .
\end{equation}
The constants are uniform for all sufficiently late $T_0$ in a fixed path class.

\begin{lemma}\label{p4:lem:vir}
For all sufficiently large $T_0$, on every $[T_0,T]$ we have
\begin{equation}\label{p4:eq:Hprime}
 \mathcal H'(t)\ge-\frac{C}{t\log t}X(t)^2-C\|f(t)\|_2X(t).
\end{equation}
\end{lemma}
\begin{proof}
For any smooth radial $c$ of compact spatial support and $\gamma_j=s_j'/s_j$, the identity $\partial_t V_j=-\gamma_j(2+D)V_j$ and spatial integration by parts imply
\begin{equation}\label{p4:eq:virial-cancellation}
\begin{aligned}
 -\frac{1}{2}\int(\partial_t V_j)h^2\dx+\ip{V_j h}{c\Lambda h}
 &=\frac{1}{2}\int[(\gamma_j-c)(2+D)V_j-(Dc)V_j]h^2\dx,\\
 \ip{v}{c\Lambda v}&=-\int(c+\frac{1}{2} Dc)v^2\dx,\\
 \ip{\Delta h}{c\Lambda h}&=\frac{1}{2}\int(Dc)|\nabla h|^2\dx
             -\int(\nabla c\cdot\nabla h)\Lambda h\dx .
\end{aligned}
\end{equation}
Since $c\le0$, the term $-\int cv^2\dx$ is nonnegative. On $r<r_-$ the inner potential cancels exactly. On $r>r_+$ we have $c=0$. The bounds $F_j=O(r^2)$ at zero and $F_j=O(r^{-3})$ at infinity hold also for $DF_j$. Hence the $L^{\frac{3}{2}}$ norms of $(2+D)V_i$ on $r>r_-$ and $(2+D)V_o$ on $r<r_+$ are $O(\rho^{\frac{10}{3}})$. At fixed $r$ in the transition interval,
\begin{equation*}
 \partial_t u=-\frac{(1-u)\gamma_i+u\gamma_o}{\log(\mu/\lambda)},\qquad
 |Dc|\lesssim\frac{1}{t|\log\rho|},\qquad
 |\partial_t c|\lesssim t^{-2}\mathbf1_{r\le r_+}.
\end{equation*}
Thus
\begin{equation*}
 \sum_{j=i,o}\|(\gamma_j-c)(2+D)V_j\|_{\frac{3}{2}}
 \lesssim t^{-1}\rho^{\frac{10}{3}}+|\gamma_o|,\qquad
 \|(\partial_t c)\Lambda h\|_2\lesssim r_+t^{-2}X .
\end{equation*}
Hardy's inequality bounds the spatial cutoff terms in \eqref{p4:eq:virial-cancellation} by $CX^2/(t|\log\rho|)$. By \eqref{p4:eq:gencreact} we have
\begin{equation*}
 |\ip{e_j}{v+c\Lambda h}|\lesssim (|s_j'|+t^{-1}r_+)X .
\end{equation*}
We use $r_+\lesssim t^{-\frac{4}{3}}$, $\rho\asymp t^{-4}$ and $\|V_m h\|_2\lesssim t^{-6}X$ in the remaining source pairings. This proves \eqref{p4:eq:Hprime}. Approximation on compact intervals extends the integrated inequality to energy data.

\end{proof}

\subsubsection{The response on the explicit path}

\begin{proof}[Proof of Theorem~\ref{p4:th:main}]
For $f=-\RR(U)$ on \eqref{p2:eq:reference}, \eqref{p3:eq:mainbound} and \eqref{p3:eq:proj} imply
\begin{equation*}
 \|f(t)\|_2\lesssim t^{-10},\qquad \sum_j|\ip{f(t)}{e_j}|\lesssim t^{-14}.
\end{equation*}
By Theorem~\ref{p5:th:inverse} with $b=9$, $X(t)\lesssim t^{-9}$, since $\int_s^{2s}\|f(t)\|_2\,dt\lesssim s^{-9}$.
In \eqref{p4:eq:reaction}, the differentiated tests contribute $O(t^{-10})$ in the inner row and $O(t^{-18})$ in the outer row. The off-diagonal entry of the Gram inverse is $O(t^{-6})$, so
\begin{equation*}
 |c_i|\lesssim t^{-10}+t^{-14},\qquad
 |c_o|\lesssim t^{-14}+t^{-18}+t^{-6}t^{-10}.
\end{equation*}
This proves \eqref{p4:eq:rates}, and $X(t)\to0$ proves scattering to zero in energy.

For $\N(U,h)=(U+h)^5-U^5-5U^4h$, we have
\begin{equation*}
 \RR(U+h)=\sum_j c_j e_j-\N(U,h),\qquad
 \|\N(U,h)\|_{\frac{6}{5}}\lesssim X^2+X^5 .
\end{equation*}
By \eqref{p3:eq:symbols} and scaling, the tests satisfy
\begin{equation*}
 \|e_j\|_6\lesssim s_j^{-1},\qquad
 \|\partial_{s_j}w\|_6\lesssim s_j^{-1}\rho^{\frac{5}{2}},\qquad
 \|\partial_{s_j}w\|_2\lesssim s_j^{-1}\rho^3 .
\end{equation*}
Pairing the residual with $e_j+\partial_{s_j}w$ proves \eqref{p4:eq:projgoal}.

The static force of $S$ in $L^{\frac{6}{5}}$ is $O(\rho^{\frac{5}{2}})$ by
\eqref{p3:eq:overlap}. Adding $w$ changes its $\dot H^{-1}$ norm by at
most $C\|\nabla w\|_2$, hence
\begin{equation}\label{p4:eq:dual}
 \|-\Delta U-U^5\|_{\dot H^{-1}}\lesssim\rho^{\frac{5}{2}}.
\end{equation}
For the kinetic first variation we differentiate $\ip{h}{e_j}=0$ and obtain
\begin{equation*}
 \ip{S_t}{h_t}=-\sum_j s_j'\ip{h}{\dot e_j}
     =O\!\left(X\sum_j|s_j'|^2\right)=O(t^{-19}).
\end{equation*}
Here we used \eqref{p4:eq:td}. The $w_t$ pairing is $O(t^{-22})$,
the static first variation is $O(t^{-19})$, and the Taylor remainder
is $O(X^2+X^6)$. Together with \eqref{p4:eq:dual}, these estimates prove
\eqref{p4:eq:Egoal} and preserve the accuracy of
\eqref{p3:eq:energydefect}.
\end{proof}
\endgroup
\space
\subsection{The projected wave inverse}\label{sec:tower-5}

\begingroup

\def\R{\mathbb R}
\def\dd{\mathrm d}
\def\ip#1#2{\langle #1,#2\rangle}
\def\RR{\mathcal R}
\def\N{\mathcal N}
\def\GG{\mathcal G}
\def\XX{\mathcal X}
\def\YY{\mathcal Y}
\def\St{\mathsf S}
\def\M{\mathsf M}
\def\dx{\,\dd x}

We construct the inverse for the projected linear equation and estimate
its response to nonlinear forcing. For the scale comparison we shall use
this inverse at the single reference path \eqref{p2:eq:reference}.

\subsubsection{The linear estimate and its application}

We work in the signed sector of Section~\ref{p1:sec:kernel}. For
the scaled profiles in \eqref{p8:eq:kelvin-pair}, we write
$q=(s_i,s_o)=(\lambda,\mu)$ and $\RR(u)=u_{tt}-\Delta u-u^5$.
On each path put
\begin{equation*}
  U_q=P_\lambda+Q_\mu+w(\lambda,\mu),\qquad L_q=-\Delta-5U_q^4,
\end{equation*}
where $w$ is the inner correction \eqref{p3:eq:ansatz}, cut off at the outer scale. The dilation tests are $e_i=\partial_\lambda P_\lambda$, $e_o=\partial_\mu Q_\mu$. For $c\in\R^2$ write $E_qc=\sum_jc_j e_j$, $E_q^*h=(\ip{h}{e_j})_j$, and $G_q=E_q^*E_q$.

We prescribe positive $C^2$ paths on $[T_0,T]$, $T\le\infty$, with $\lambda'\le0$ and the bounds \eqref{p4:eq:path}. We bound all comparison constants by a fixed $K>1$ large enough to include \eqref{p2:eq:reference}. For $T_0\le S<T$ let $I_S=[S,\min(2S,T)]$, set $\St(I)=L^5_t L^{10}_x(I\times\R^3)$, and write $H_h(t)=\|(\nabla h,h_t)(t)\|_{L^2\times L^2}$. We define for $b>0$
\begin{align}
 \|f\|_{\YY_b}&=\sup_S S^b\int_{I_S}\|f(t)\|_2\,\dd t,\nonumber\\
 \|h\|_{\XX_b}&=\sup_t t^b H_h(t)
                 +\sup_S S^{b-1}\|h\|_{\St(I_S)} . \label{p5:eq:norms}
\end{align}
The solution space in \eqref{p5:eq:norms} consists of energy-continuous pairs $(h,h_t)$ with the stated spacetime norm. We use the distributional time derivative in this definition.

\begin{theorem}\label{p5:th:inverse}
For fixed $K,b$ and sufficiently large $T_0$, the selected inverse $h=\GG_q f$ solves
\begin{equation}\label{p5:eq:wave}
 h_{tt}+L_q h=f+E_qc,\qquad E_q^*h=0,\qquad
 \|h\|_{\XX_b}\lesssim_{K,b}\|\Pi_q f\|_{\YY_b},\quad
 \Pi_q=1-E_qG_q^{-1}E_q^*.
\end{equation}
The estimate is uniform in the path, the terminal time and the scale ratio. Sources are signed and locally $L^1_tL^2_x$ with $\Pi_q f\in\YY_b$. The two constraints hold at every time and all one-core negative modes evolve. At finite $T$ we impose the negative-mode conditions and the terminal data restrictions in \eqref{p5:eq:selection}. On the half-line there is a unique solution with $\sup t^b H_h(t)<\infty$ and the same initial negative-mode conditions.
\end{theorem}

We use this inverse only on the reference path $q_*$ of \eqref{p2:eq:reference}.
We write $E_*=E_{q_*}$, $\Pi_*=\Pi_{q_*}$, and $\GG_*=\GG_{q_*}$.
For every trial $q$ impose the same dilation constraints and endpoint conditions,
all evaluated on $q_*$. We shall solve
\begin{equation}\label{p5:eq:complement}
 \RR(U_q+h_q)=E_*c_q,\qquad E_*^*h_q=0,\qquad
 H_{h_q}(t)\lesssim t^{-9},\quad \|h_q\|_{\St(I_S)}\lesssim S^{-8}.
\end{equation}
The neighborhood and its difference estimate are in Section~\ref{p5:sec:tube}. Section~\ref{sec:tower-6} solves the scalar equations $c_q=0$.

\subsubsection{Coercivity and energy with time-integrated forcing}\label{p5:sec:coercivity}

We use the dilation tests, negative eigenfunctions $Y_{j,k}$, Gram matrix
and localized form $\mathcal H$ from Section~\ref{sec:tower-4}.
We write $\rho=\frac{\lambda}{\mu}$, $V=5U_q^4$, $V_j=5(F_j)_{s_j}^4$,
$V_m=V-V_i-V_o$, $D=x\cdot\nabla$ and $v=h_t$.
We apply the coercivity estimate \eqref{p4:eq:coer}. System \eqref{p4:eq:reaction} determines
the reactions, and \eqref{p4:eq:gencreact} bounds them.
Adding $E_q a$ to the source translates the reactions by $-a$ and leaves the wave fixed. We may therefore take $f=\Pi_q f$. Approximation on compact
intervals extends the energy Cauchy problem and the integrated inequality
of Lemma~\ref{p4:lem:vir} to locally $L^1_tL^2_x$ sources.

For each $\psi=Y_{j,k}^{[s_j]}$ set
\begin{equation*}
  u_\psi=s_j^{-1}\ip{h}{\psi},\quad v_\psi=\ip{v}{\psi},\quad
  z_\psi^\pm=v_\psi\pm\nu_{j,k}u_\psi .
\end{equation*}
Their equations are
\begin{equation}
 (z_\psi^\pm)'=\pm\frac{\nu_{j,k}}{s_j} z_\psi^\pm+
                         O(\|f\|_2+t^{-1}H_h).\label{p5:eq:negative}
\end{equation}
Let $\mathcal T_j(t)$ be the span of $e_j(t)$ and all $Y_{j,k}^{[s_j(t)]}$ at that core. The conditions on $z_\psi^\pm$ below hold for every negative mode. We select the solution by imposing
\begin{equation}\label{p5:eq:selection}
\begin{aligned}
 &z_\psi^-(T_0)=0 &&\text{on every interval},\\
 &z_\psi^+(T)=0 &&(T<\infty),\\
 &h(T)\in\sum_{j=i,o}s_j(T)\mathcal T_j(T),\qquad
 h_t(T)\in\sum_{j=i,o}\mathcal T_j(T) &&(T<\infty),\\
 &E_q(T)^* h(T)=0,\qquad
 E_q(T)^* h_t(T)=-\dot E_q(T)^*h(T) &&(T<\infty).
\end{aligned}
\end{equation}
The last line imposes the dilation constraints and their first time derivatives. For given negative moments it fixes the dilation moments of the terminal data. The same-core moment matrices are diagonal with positive entries. Off-core displacement entries are $O(\rho^{\frac{1}{2}})$ even after the factor $\frac{s_o}{s_i}$. The velocity matrix has $O(\rho^{\frac{3}{2}})$ off-core entries. Consequently $H_h(T)\lesssim\sum_\psi |z_\psi^-(T)|$. Equation \eqref{p5:eq:negative} follows by differentiating the moments and applying \eqref{p4:eq:almostmode} and \eqref{p4:eq:gencreact}. Its inner exponential uses $\int \frac{\dd t}{\lambda}\asymp t^5$ over a late dyadic interval, and its outer exponential uses $\int \frac{\dd t}{\mu}\asymp t$.

\paragraph{The energy estimate and the terminal trace.}
For the selected problem we take $p=\|\Pi_q f\|_2$ and $Z=0$.
The same argument also allows arbitrary terminal data, keeping
$E_q^*h=0$ and $z^-(T_0)=0$, if $Z=T^b H_h(T)$.
We write $B=\|\Pi_q f\|_{\YY_b}+Z$ and $A_b=\sup t^bH_h(t)$. In either case
\begin{equation}\label{p5:eq:defect-bound}
                         A_b\lesssim_{K,b} B .
\end{equation}
Both rates $\frac{\nu_{j,k}}{s_j}$ in \eqref{p5:eq:negative} are bounded below. By the dyadic bound in \eqref{p5:eq:norms},
\begin{equation*}
 \int_{T_0}^{T}e^{-c|t-s|}p(s)\,ds\lesssim_b Bt^{-b},\qquad
 \int_{T_0}^{T}e^{-c|t-s|}\frac{H_h(s)}{s}\,ds\lesssim_b A_b t^{-b-1}.
\end{equation*}
For $s\ge \frac{t}{2}$ we sum the dyads starting at $\frac{t}{2}$ (clipped at $T_0$). For $s<\frac{t}{2}$ we use exponential decay. The endpoint conditions and the terminal span estimate now imply
\begin{equation*}
 t^b\bigl(|z^+(t)|+|z^-(t)|\bigr)
 \lesssim_b B+T_0^{-1}A_b .
\end{equation*}
At finite $T$, $T^bH_h(T)$ has the same bound, by the terminal span estimate in the selected case and by the definition of $Z$ otherwise. Integrating Lemma~\ref{p4:lem:vir} backwards, we use $\int_t^T pH_h\lesssim_b BA_b t^{-2b}$ and \eqref{p4:eq:Hcoer} to obtain
\begin{equation*}
 A_b^2\le C_b(B+T_0^{-1}A_b)^2+\frac{C_b}{\log T_0}A_b^2+C_bBA_b .
\end{equation*}
Young's inequality and a sufficiently large $T_0$ prove \eqref{p5:eq:defect-bound}.

For existence, we prescribe the free vector $z^-(T)$ in the given terminal span and solve backward. The map to $z^-(T_0)$ is affine, with equal domain and target dimensions. Its homogeneous kernel is zero by the estimate at $B=0$, which proves solvability for \eqref{p5:eq:selection}. We take a weakly convergent subsequence of data at $T_0$. The bounded linear propagator on each compact interval preserves weak convergence, so the limit solves the half-line equation with the same estimate. For a homogeneous half-line difference in the stated decay class,
the future mode and energy boundary terms vanish upon integration to
infinity at each fixed $t$. With $p=0$, the preceding estimate reduces to
\begin{equation*}
 A_b^2\le\left(C_b T_0^{-2}+\frac{C_b}{\log T_0}\right)A_b^2 .
\end{equation*}
Our choice of $T_0$ forces $A_b=0$. The half-line inverse is therefore unique and linear.

\subsubsection{The critical spacetime norm}\label{p5:sec:str}

For $z_t=v+g$, $v_t=\Delta z+f$ on $J$, we apply the homogeneous free-wave estimate \cite[Corollary~1.3]{KT} to each Duhamel datum $(g,f)$. The pair $(4,12)$ is wave-admissible in dimension three and $\frac{1}{4}+\frac{3}{12}=\frac{1}{2}$. By Minkowski's inequality, for $t_*\in J$ we have
\begin{equation}\label{p5:eq:mixed-free-wave}
 \|z\|_{L^4_t L^{12}_x(J)}
 \lesssim \|(\nabla z,v)(t_*)\|_2+
     \int_J(\|\nabla g\|_2+\|f\|_2)\,\dd t .
\end{equation}
The same estimate holds for $(5,10)$. We keep the energy inside the forcing integral. With $H(t)=\|(\nabla z,v)(t)\|_2$, either endpoint $c$ of $I$, and $(q,r)=(4,12),(5,10)$, we claim
\begin{equation}\label{p5:eq:variation-strichartz}
 \|z\|_{L^q(I;L^r)}^q\lesssim H(c)^q+
 \int_I H(s)^{q-1}(\|\nabla g(s)\|_2+\|f(s)\|_2)\,\dd s .
\end{equation}
Let $w_s$ be the free wave with data $(z,v)(s)$ at $s$. Then $w_s(s)=z(s)$. For smooth inputs, $\partial_s w_s$ is the free wave with data $(g,f)(s)$ at $s$. On $I=[a,d]$, H\"older and the homogeneous estimate imply
\begin{equation*}
 \left|\frac{d}{ds}\int_s^d\|w_s(t)\|_r^q\,\dd t+\|z(s)\|_r^q\right|
 \lesssim H(s)^{q-1}(\|\nabla g(s)\|_2+\|f(s)\|_2).
\end{equation*}
Integration proves \eqref{p5:eq:variation-strichartz} for $c=a$. Integrating over $[a,s]$ proves it for $c=d$. The energy and Duhamel estimates allow approximation in the data and in $L^1(I;\dot H^1\times L^2)$ for the forcing.

On $I_S$ we may take $f=\Pi_qf$. By \eqref{p4:eq:gencreact}, scaling and \eqref{p4:eq:Vmix}, the full free source satisfies $\|h_{tt}-\Delta h\|_2\lesssim_K t^4H_h(t)+\|f(t)\|_2$. Taking $(q,r)=(5,10)$, $g=0$ in \eqref{p5:eq:variation-strichartz} and applying Young's inequality, we obtain
\begin{equation}\label{p5:eq:short}
 \|h\|_{\St(I_S)}\lesssim_K S\sup_{I_S}H_h+\int_{I_S}\|f(t)\|_2\,\dd t .
\end{equation}
On the half-line, if $H_h(t)\le M t^{-b}$ with $b>1$, the right endpoint term tends to zero. For $S\ge T_0$,
\begin{equation*}
 \|h\|_{\St([S,\infty))}^5\lesssim_{K,b}
 M^5 S^{5-5b}+M^4\int_S^\infty t^{-4b}\|f(t)\|_2\,\dd t .
\end{equation*}
Equations \eqref{p5:eq:short} and \eqref{p5:eq:defect-bound} prove Theorem~\ref{p5:th:inverse}. The same scaling bounds the background field,
\begin{equation*}
 \|(F_j)_{s_j}\|_{10}=s_j^{-\frac{1}{5}}\|F_j\|_{10},\qquad
 \|U_q\|_{\St(I_S)}^5\lesssim_K\int_{I_S}(\lambda^{-1}+\mu^{-1})\,\dd t
 \lesssim_K S^5 .
\end{equation*}
The correction satisfies the same spatial estimate as the cores. For $\N(U,h)=(U+h)^5-U^5-5U^4h$ we have
\begin{equation*}
 \N(U,h)-\N(U,g)
 =5(h-g)\int_0^1\big[(U+g+\theta(h-g))^4-U^4\big]\,d\theta .
\end{equation*}
Expanding the fourth power and using Young's inequality, we obtain
\begin{equation}\label{p5:eq:quintic-pointwise}
 |\N(U,h)-\N(U,g)|\lesssim |h-g|
 \big[|U|^3(|h|+|g|)+|h|^4+|g|^4\big].
\end{equation}
For $2\le m\le5$, H\"older's inequality on $I_S$ reads
\begin{equation}\label{p5:eq:power}
 \|U^{5-m}h^m\|_{L^1_tL^2_x(I_S)}
 \le\|U\|_{\St(I_S)}^{5-m}\|h\|_{\St(I_S)}^m
 \lesssim_K S^{5-mb}\|h\|_{\XX_b}^{m}.
\end{equation}
If $\|h\|_{\XX_b},\|g\|_{\XX_b}\le R_0$, we place the difference factor in $\St(I_S)$ with weight $\beta-1$. For $b\ge5$ and $\beta>0$, \eqref{p5:eq:quintic-pointwise} then implies
\begin{equation}\label{p5:eq:prodDiff}
 \|\N(U,h)-\N(U,g)\|_{\YY_\beta}
 \lesssim_K \big[R_0T_0^{5-b}+R_0^4T_0^{5-4b}\big]\|h-g\|_{\XX_\beta}.
\end{equation}
At $b=9$ the leading quadratic term belongs to $\YY_{13}$ by \eqref{p5:eq:power}, and the Lipschitz bound in $\YY_9$ gains $T_0^{-4}$.

\subsubsection{Differences between scale paths}\label{p5:sec:difference}

We set $V_q=5U_q^4$ and $A_j=\partial_{\log s_j}$ at fixed physical $x$.
The change of potential obeys
\begin{equation}\label{p5:eq:Vdiff}
 \|A_i V_q\|_3\lesssim_K\lambda^{-1},\qquad
 \|A_o V_q\|_3\lesssim_K1 .
\end{equation}
The individual wells have the stated scaling exponents. Logarithmic derivatives preserve the core envelopes in \eqref{p3:eq:overlap} and the $L^{12}$ bound on the correction in \eqref{p3:eq:env}. Applying the proof of \eqref{p4:eq:Vmix} to each differentiated product bounds the mixed terms by $O_K(\rho^{\frac{3}{2}})$ in $L^3$.
Using $L_{q_*}=-\Delta-V_*$, equation \eqref{p5:eq:complement} becomes
\begin{equation}\label{p5:eq:compare}
 h_q=\GG_*\big[-\Pi_*\RR(U_q)+(V_q-V_*)h_q+\N(U_q,h_q)\big].
\end{equation}
Its boundary spaces and projection are identical for every trial path.
All path differences now occur in the displayed source.

\subsubsection{The actual residual on neighboring paths}\label{p5:sec:tube}

We now estimate the change in the residual when the prescribed scales vary. Let $q_*=(\lambda_*,\mu_*)$ be the explicit reference path \eqref{p2:eq:reference}. For varying the scales we use the supremum over time of the sum of the six weighted differences
\begin{equation*}
\begin{array}{c|ccc}
 &\text{position}&\text{velocity}&\text{acceleration}\\ \hline
 i&t^{12}|\delta\lambda|&t^{13}|\delta\lambda'|&t^{10}|\delta\lambda''|\\
 o&t^{11}|\delta\mu|&t^{12}|\delta\mu'|&t^{13}|\delta\mu''|
\end{array}
\end{equation*}
We denote it by $d(q,\bar q)$ and take the class $d(q,q_*)\le A$, for fixed $A$. For $T_0$ large depending on $A$ the paths are positive, contracting in the inner variable and satisfy \eqref{p4:eq:path}, by $\lambda_*'\sim-4a t^{-5}$. These weights allow the inner acceleration to vary more than its twice-integrated effect. The scalar force identity \eqref{p6:eq:exact} explains these weights.

\begin{proposition}\label{p5:pr:res}
For any two paths in this class, uniformly in $T$,
\begin{align}
 \|V_q-V_{\bar q}\|_3&\lesssim_A d(q,\bar q)t^{-4},\label{p5:eq:tubed}\\
 \|\Pi_*[\RR(U_q)-\RR(U_{\bar q})]\|_2
                &\lesssim_A d(q,\bar q)t^{-14}. \label{p5:eq:tuber}
\end{align}
The unprojected residual satisfies $\|\RR(U_q)\|_2\lesssim_A t^{-10}$ and
$\|\RR(U_q)-\RR(U_{\bar q})\|_2\lesssim_A d(q,\bar q)t^{-10}$.
\end{proposition}
\begin{proof}
We write $d=d(q,\bar q)$, $\delta g=g_q-g_{\bar q}$ and $\ell_j=\log s_j$. Dividing the path bounds by the corresponding scales, we obtain componentwise
\begin{equation}\label{p6:eq:logs}
\begin{aligned}
 (|\delta\ell_i|,|\delta\ell_i'|,|\delta\ell_i''|)
    &\lesssim_A d\,(t^{-8},t^{-9},t^{-6}),\\
 (|\delta\ell_o|,|\delta\ell_o'|,|\delta\ell_o''|)
    &\lesssim_A d\,(t^{-11},t^{-12},t^{-13}).
\end{aligned}
\end{equation}
Integrating \eqref{p5:eq:Vdiff} on the segment in $(\ell_i,\ell_o)$ proves \eqref{p5:eq:tubed}.

We now take differences in \eqref{p3:eq:residual-decomposition}. We denote the two inner coefficients by
\begin{equation*}
 a_\lambda=\lambda''-\frac{\alpha\rho^{\frac{5}{2}}}{\lambda},\qquad
 b_\lambda=\frac{(\lambda')^2-\beta\rho^{\frac{5}{2}}}{\lambda}.
\end{equation*}
By \eqref{p3:eq:accelerations} they are $O(t^{-10})$ on $q_*$. On the whole ball,
\begin{equation}\label{p5:eq:split}
 |\delta a_\lambda|\lesssim_A dt^{-10},\qquad |\delta b_\lambda|\lesssim_A dt^{-14},\qquad
 |a_\lambda|+|b_\lambda|\lesssim_A t^{-10}.
\end{equation}
Indeed $\frac{\rho^{\frac{5}{2}}}{\lambda}\asymp t^{-6}$ and its relative variation is $O_A(dt^{-8})$. For $\frac{(\lambda')^2}{\lambda}$ the velocity variation is bounded by $C_A t^4t^{-5}|\delta\lambda'|\lesssim_A dt^{-14}$. For the outer coefficients $|\delta\mu''|\le dt^{-13}$ and $\big|\delta\big(\frac{(\mu')^2}{\mu}\big)\big|\lesssim_A dt^{-21}$. All four $L^2$-rescaled shapes have bounded log derivatives.

For the correction, the chain rule reads
\begin{equation}\label{p5:eq:shape-acceleration}
 w_{tt}=\sum_j\ell_j'' A_jw+\sum_{j,k}\ell_j'\ell_k'A_j A_kw .
\end{equation}
Since $\|A^\nu w\|_2\lesssim_\nu t^{-12}$, \eqref{p6:eq:logs} implies $\|\delta w_{tt}\|_2\lesssim_A dt^{-18}$. The same bound for $\delta R$ follows directly from \eqref{p3:eq:spatial-remainder}.

It remains to project the two acceleration terms. Since $\Pi_* e_{*,j}=0$, we have $\|\Pi_*e_j\|_2\lesssim_A t^{-8}$ and $\|\delta e_j\|_2\lesssim_A dt^{-8}$. For $(b_i,b_o)=(a_\lambda,\mu'')$ this implies
\begin{equation*}
 \|\Pi_*\delta(b_j e_j)\|_2
 \le |\delta b_j|\,\|\Pi_* e_j\|_2+|\bar b_j|\,\|\delta e_j\|_2
 \lesssim_A dt^{-18}.
\end{equation*}
The remaining terms of \eqref{p3:eq:residual-decomposition} vary by at most $C_A dt^{-14}$ in $L^2$, which proves \eqref{p5:eq:tuber}. Before projection the largest variation is $O_A(dt^{-10})$ by \eqref{p5:eq:split}. The reference bound \eqref{p3:eq:mainbound} now implies the asserted bound on $\RR(U_q)$.
\end{proof}

\begin{corollary}\label{p5:cor:nonlin}
Uniformly on this class of prescribed paths, \eqref{p5:eq:complement} has a selected solution with $\|h_q\|_{\XX_9}\le C_A$, unique in the ball specified in the proof. For the linear response $h_q^{(1)}$ obtained by setting $\N=0$ in
\eqref{p5:eq:compare},
\begin{equation*}
 \|h_q-h_q^{(1)}\|_{\XX_{13}}\lesssim_A1,\qquad
 \|h_q-h_{\bar q}\|_{\XX_{12}}\lesssim_A d(q,\bar q).
\end{equation*}
The assertions hold both at a finite terminal cutoff and on the half-line.
\end{corollary}
\begin{proof}
On $I_S$, Sobolev and \eqref{p5:eq:tubed} imply
\begin{equation*}
 \|(V_q-V_*)h\|_{L^1_tL^2_x(I_S)}
 \lesssim_A S^{-3-\beta}\|h\|_{\XX_\beta},\qquad \beta=9,12,13 .
\end{equation*}
We absorb this term into the linear inverse \eqref{p5:eq:wave} for $T_0$ large. By \eqref{p3:eq:mainbound} and \eqref{p5:eq:tuber}, $\|h_q^{(1)}\|_{\XX_9}\lesssim_A1$. On the $\XX_9$ ball of radius $\max(1,2\sup_{q,T}\|h_q^{(1)}\|_{\XX_9})$, \eqref{p5:eq:power} and \eqref{p5:eq:prodDiff} bound the nonlinear Lipschitz constant by $C_AT_0^{-4}$. Thus \eqref{p5:eq:compare} has a unique fixed point in this ball, and
\begin{equation*}
 \|h_q-h_q^{(1)}\|_{\XX_{13}}\lesssim_A\|\N(U_q,h_q)\|_{\YY_{13}}\lesssim_A1 .
\end{equation*}

We subtract \eqref{p5:eq:compare} for $q$ and $\bar q$. With $d=d(q,\bar q)$, \eqref{p5:eq:tubed}, \eqref{p5:eq:tuber} and differentiation in the log scales imply
\begin{align*}
 \|(V_q-V_{\bar q})h_{\bar q}\|_{\YY_{12}}+
 \|\Pi_*[\RR(U_q)-\RR(U_{\bar q})]\|_{\YY_{13}}&\lesssim_A d,\\
 \|U_q-U_{\bar q}\|_{\St(I_S)}&\lesssim_A dS^{-7}.
\end{align*}
In \eqref{p5:eq:power}, replacing one $U$ by $U_q-U_{\bar q}$ gains $dS^{-8}$. By \eqref{p5:eq:prodDiff}, the divided difference in $h$ has operator norm $O_A(T_0^{-4})$ from $\XX_\beta$ to $\YY_\beta$ for $\beta=9,12$. We solve the difference equation in $\XX_{12}$ and identify its solution in $\XX_9$ by uniqueness. All bounds are uniform in the terminal time.
\end{proof}

\endgroup
\space
\subsection{Closing the two scale equations}\label{sec:tower-6}

\begingroup

\def\R{\mathbb R}
\def\HH{\dot H^1(\R^3)\times L^2(\R^3)}
\def\res{\mathcal R}
\def\N{\mathcal N}
\def\ip#1#2{\langle #1,#2\rangle}
\def\dd{\,\mathrm d}

Section~\ref{sec:tower-5} constructs a wave whose residual lies in the
span of the two reference dilations. We now project on the dilations
at the actual scales and choose those scales so that both coefficients
vanish.

\subsubsection{The nonlinear closure theorem}

We continue to use the signed profiles \eqref{p8:eq:kelvin-pair},
the constants in \eqref{p8:eq:coefficients} and the residual
$\res(u)=u_{tt}-\Delta u-u^5$. Let
$q_*=(\lambda_*,\mu_*)$ be the polynomial path \eqref{p2:eq:reference},
with outer limit one and clock $L=at^{-4}$. In particular,
\begin{equation}\label{p6:eq:branch}
 \lambda_*=a t^{-4}+O(t^{-8}),\quad \lambda_*'=-4at^{-5}+O(t^{-9}),
 \qquad 1-\mu_*\sim\frac{5}{18}a^2t^{-8},\quad \mu_*'=O(t^{-9}).
\end{equation}
For a scale path $q=(s_i,s_o)=(\lambda,\mu)$ denote the bare pair by $S=P_\lambda+Q_\mu$ and the corrected pair by $U=U_q=S+w(q)$, with $w$ from \eqref{p3:eq:ansatz}. For a perturbation $y=(y_i,y_o)$ on $I=[T_0,T]$ with $T>T_0$, or on $I=[T_0,\infty)$, we set
\begin{equation}\label{p6:eq:norm}
 \|y\|_I=\sup_{t\in I}\big(t^{12}|y_i|+t^{13}|y_i'|+t^{10}|y_i''|
               +t^{11}|y_o|+t^{12}|y_o'|+t^{13}|y_o''|\big).
\end{equation}

\begin{theorem}\label{p6:th:main}
There exist $A,T_0$ such that on every one of these intervals there is a path $q=q_*+y$ with $\|y\|_I\le A$ for which both reference dilation multipliers in \eqref{p5:eq:complement} vanish. The field $u=U_q+h_q$ is an exact strong energy solution, and
\begin{equation}\label{p6:eq:wavebound}
 \|(\nabla h_q,h_{q,t})(t)\|_2\lesssim t^{-9},\qquad
 \|h_q\|_{L^5_tL^{10}_x([s,\min(2s,T)])}\lesssim s^{-8}.
\end{equation}
At finite $T$ take $y(T)=y'(T)=0$ and select the wave by all the conditions in \eqref{p5:eq:selection}, including the restrictions on its terminal displacement and velocity. On the half-line impose \eqref{p6:eq:norm} on the scales and select the decaying wave with $z_\psi^-(T_0)=0$ for every negative mode. For the complementary wave selected in the contraction ball of Section~\ref{sec:tower-5}, the scale path is unique in the ball $\|y\|_I\le A$. The constants are independent of $T$.
\end{theorem}

All wave constraints and endpoint tests in \eqref{p5:eq:selection} are
evaluated on $q_*(t)$, for every trial path. In particular $E_*^*h=0$,
and the negative coordinates at the endpoints use the same two
reference scales. The parameter $q$ in $U_q$ remains free.
\subsubsection{Exact forces on a prescribed path}

Let $Ec=\sum c_j e_j$, $G=E^*E$, $S=P_\lambda+Q_\mu$,
$b(\rho)=\ip{e_i}{e_o}$, and $V=E_{\rm stat}(S)-2E_{\rm stat}(Q)$.
The diagonal entries of $G$ are $M$. By \eqref{p2:eq:momentum},
\begin{equation}\label{p6:eq:bare}
 E^*\res(S)=(Gq')'-\mathcal F(q,q'),\qquad
 \mathcal F=\frac{b'(\rho)\lambda'\mu'-V'(\rho)}{\mu}
                    \binom{1}{-\rho}.
\end{equation}
A star denotes evaluation on $q_*$. By \eqref{p2:eq:defect},
\begin{equation*}
 r_*=M^{-1}\big[(G_*q_*')'-\mathcal F_*\big]=O(t^{-14}).
\end{equation*}
We test \eqref{p5:eq:complement} against the actual dilation $e_j$. The product rule in time reads
\begin{equation}\label{p6:eq:moment-derivative}
 \ip{h_{tt}}{e_j}=\partial_t\big(\ip{h_t}{e_j}-\ip{h}{\dot e_j}\big)
                    +\ip{h}{\ddot e_j}.
\end{equation}
We put $v=w+h$, so the full field on every trial path is $S+v$. Its static residual differs from that of $S$ by $L_Sv-\N(S,v)$, where $L_S=-\Delta-5S^4$ and $\N(S,v)=(S+v)^5-S^5-5S^4v$. For $k\ne j$ set
\begin{align}
 2M B_j&=\ip{h}{\dot e_j}-\ip{h_t}{e_j}
                  -b(\rho)s_k'+b(\rho_*)s_{*,k}',\nonumber\\
 MD_j&=\ip{w_{tt}}{e_j}+\ip{v}{L_Se_j}+\ip{h}{\ddot e_j}
                     -\ip{\N(S,v)}{e_j}.\label{p6:eq:D}
\end{align}
Substituting \eqref{p6:eq:moment-derivative} and \eqref{p6:eq:bare} in the tested equation, we obtain
\begin{equation}\label{p6:eq:exact}
 (y'-2B)'=M^{-1}(\mathcal F-\mathcal F_*)-r_*-D+M^{-1}E^*E_*c .
\end{equation}
The matrix $E^*E_*$ is invertible, since $G_*^{-1}$ is uniformly bounded and $\|E-E_*\|_{\mathbb R^2\to L^2}\lesssim_A t^{-8}$. We first obtain \eqref{p6:eq:exact} in distributions. The dual estimate in Lemma~\ref{p6:lem:feedback} and energy continuity imply $B\in C^1$.

\subsubsection{Bounds on scale feedback}

We fix $A<\infty$ and take $T_0$ sufficiently large depending on $A$. In the ball \eqref{p6:eq:norm} about $q_*$, $\lambda\asymp t^{-4}$, $\mu\asymp1$, $\lambda'<0$ and the derivative bounds of Section~\ref{sec:tower-5} hold. For two paths we use $\delta$ for the first quantity minus the second, and $d=\|q-\bar q\|_I$. The scale difference is measured about the same reference path throughout. All bounds are uniform in the right endpoint.

\begin{lemma}\label{p6:lem:bare}
With the clock of \eqref{p6:eq:branch},
\begin{equation*}
 |M^{-1}\delta\mathcal F_i-30t^{-2}\delta\lambda|\lesssim_A d\,t^{-17},
 \qquad |M^{-1}\delta\mathcal F_o|\lesssim_A d\,t^{-18}.
\end{equation*}
\end{lemma}
\begin{proof}
At fixed velocities $\mathcal F_i$ in \eqref{p6:eq:bare} is homogeneous of degree $-1$ in $(\lambda,\mu)$, and $\mathcal F_o=-\rho\mathcal F_i$. For a variation of the path, Euler's identity yields
\begin{equation}\label{p6:eq:force-differential}
 \begin{aligned}
 \dd\mathcal F_i&=\partial_\lambda\mathcal F_i\,(\dd\lambda-\rho\,\dd\mu)
       -\frac{\mathcal F_i}{\mu}\dd\mu
       +\frac{b'(\rho)}{\mu}(\mu'\dd\lambda'+\lambda'\dd\mu'),\\
 \dd\mathcal F_o&=-\rho\,\dd\mathcal F_i-\mathcal F_i\,\dd\rho .
 \end{aligned}
\end{equation}
Only the inner position derivative needs an expansion. By \eqref{p2:eq:forces} and \eqref{p2:eq:bprime},
\begin{equation*}
 M^{-1}\partial_\lambda\mathcal F_i
 =\frac{15C}{4M}\lambda^{\frac{1}{2}}\mu^{-\frac{5}{2}}(1+O(\rho))+O_A(t^{-12})
 =30t^{-2}+O_A(t^{-6}).
\end{equation*}
The $O_A(t^{-12})$ term comes from $b''(\rho)\lambda'\mu'/\mu^2$. Since $|\mathcal F_i|\lesssim t^{-6}$ and $|b'(\rho)|\lesssim t^{-2}$, integration of \eqref{p6:eq:force-differential} along the segment joining the paths implies
\begin{align*}
 |M^{-1}\delta\mathcal F_i-30t^{-2}\delta\lambda|
 &\lesssim_A t^{-6}(|\delta\lambda|+|\delta\mu|)
       +t^{-11}|\delta\lambda'|+t^{-7}|\delta\mu'|\lesssim_A dt^{-17},\\
 |\delta\mathcal F_o|&\lesssim_A t^{-4}dt^{-14}+t^{-6}dt^{-12}\lesssim_A dt^{-18}.
\end{align*}
\end{proof}

\begin{lemma}\label{p6:lem:feedback}
On the reference path, with its own selected wave for the endpoint under consideration, the bounds are
\begin{equation}\label{p6:eq:base}
\begin{array}{c|ccc}
  & |D_j| & |B_j| & |B_j'|\\ \hline
 i & O(t^{-14})&O(t^{-14})&O(t^{-10})\\
 o & O(t^{-14})&O(t^{-18})&O(t^{-18}).
\end{array}
\end{equation}
For two paths in the ball, differences divided by $d$ satisfy
\begin{equation}\label{p6:eq:diff}
\begin{array}{c|ccc}
  & \frac{|\delta D_j|}{d} & \frac{|\delta B_j|}{d} & \frac{|\delta B_j'|}{d}\\ \hline
 i & O_A(t^{-17}) & O_A(t^{-17}) & O_A(t^{-13})\\
 o & O_A(t^{-17}) & O_A(t^{-19}) & O_A(t^{-16}).
\end{array}
\end{equation}
\end{lemma}
\begin{proof}
We use the log-scale differences \eqref{p6:eq:logs}, with $\mathcal D_j=\partial_{\ell_j}$ at fixed $x$. Corollary~\ref{p5:cor:nonlin} bounds the wave in the energy norm $H_g=\|(\nabla g,g_t)\|_2$,
\begin{equation*}
 H_h\lesssim_A t^{-9},\qquad H_{\delta h}\lesssim_A d\,t^{-12}.
\end{equation*}

\paragraph{The force.}
For $v=w+h$ we have $\|v\|_6\lesssim_A t^{-9}$ and $\|\delta v\|_6\lesssim_A dt^{-12}$ by \eqref{p3:eq:symbols} and the wave bound. The spatial part of \eqref{p6:eq:D} uses only the bare well error
\begin{equation}\label{p6:eq:well}
 L_Se_j=-5\big(S^4-(F_j)_{s_j}^4\big)e_j,\qquad
 \|\mathcal D^\alpha L_Se_j\|_{\frac{6}{5}}
 \lesssim \begin{cases}t^{-6},&j=i,\\t^{-10},&j=o,\end{cases}\qquad |\alpha|\le1 .
\end{equation}
This is the estimate following \eqref{p3:eq:projid}. For $(A_0,z)=(S,v)$ or $(U,h)$, we use H\"older and the quadratic vanishing of $\N$ to obtain
\begin{equation}\label{p6:eq:nonlinear-variation}
 \|\N(A_0,z)\|_{\frac{6}{5}}\lesssim_A t^{-18},\qquad
 \|\delta\N(A_0,z)\|_{\frac{6}{5}}
 \lesssim_A t^{-9}\|\delta z\|_6+t^{-18}\|\delta A_0\|_6
 \lesssim_A d t^{-21}.
\end{equation}
Here the implicit constant depends on the path ball, and $\|\delta S\|_6+\|\delta U\|_6\lesssim_A dt^{-8}$. For the inertial term we already have
\begin{equation*}
 \|w_{tt}\|_2\lesssim_A t^{-14},\qquad
 \|\delta w_{tt}\|_2\lesssim_A dt^{-18}
\end{equation*}
from \eqref{p5:eq:shape-acceleration}. For the remaining tests, $e_j$ and its log derivatives have norms $O(1),O(s_j),O(s_j^{-1})$ in $L^2,L^{\frac{6}{5}},L^6$. Using $\dot e_j=\ell_j'\mathcal D_j e_j$ and $\ddot e_j=\ell_j''\mathcal D_j e_j+(\ell_j')^2\mathcal D_j^2e_j$, we find
\begin{equation}\label{p6:eq:test-derivatives}
\begin{array}{c|ccc}
 &\|\dot e_j\|_2&\|\dot e_j\|_{\frac{6}{5}}&\|\ddot e_j\|_{\frac{6}{5}}\\ \hline
 i&O(t^{-1})&O(t^{-5})&O(t^{-6})\\
 o&O(t^{-9})&O(t^{-9})&O(t^{-10})\\[2pt]
 &\frac{\|\delta\dot e_j\|_2}{d}&\frac{\|\delta\dot e_j\|_{\frac{6}{5}}}{d}&
                       \frac{\|\delta\ddot e_j\|_{\frac{6}{5}}}{d}\\ \hline
 i&O_A(t^{-9})&O_A(t^{-13})&O_A(t^{-10})\\
 o&O_A(t^{-12})&O_A(t^{-12})&O_A(t^{-13})
\end{array}
\end{equation}
For instance $s_i|\delta\ell_i''|\lesssim_A dt^{-10}$. All entries follow from the same rescaling and \eqref{p6:eq:logs}. Substituting in \eqref{p6:eq:D}, we obtain, uniformly in both rows,
\begin{align*}
 |D_j|&\lesssim_A t^{-14}+t^{-9}t^{-6}+s_j^{-1}t^{-18}\lesssim_A t^{-14},\\
 |\delta D_j|&\lesssim_A d\big(t^{-18}+t^{-12}t^{-6}+t^{-9}t^{-10}
       +s_j^{-1}(t^{-21}+t^{-18}t^{-8})\big)\lesssim_A d t^{-17}.
\end{align*}
In particular these estimates suffice for the first columns of \eqref{p6:eq:base} and \eqref{p6:eq:diff}.

\paragraph{The momentum.}
Set $\phi_j=e_j-e_{*,j}$ and $C_j=\ip{h}{\dot e_j}-\ip{h_t}{e_j}$. Since $E_*^*h=0$,
\begin{equation}\label{p6:eq:wave-current}
\begin{aligned}
 C_j&=2\ip{h}{\dot e_{*,j}}+\ip{h}{\dot\phi_j}-\ip{h_t}{\phi_j},\\
 C_j'&=2\ip{h_t}{\dot e_{*,j}}+2\ip{h}{\ddot e_{*,j}}
                +\ip{h}{\ddot\phi_j}-\ip{h_{tt}}{\phi_j}.
\end{aligned}
\end{equation}
We eliminate the reaction before estimating the last pairing. Twice differentiating $E_*^*h=0$ in \eqref{p5:eq:complement}, we obtain
\begin{equation*}
 h_{tt}=\Pi_*[-L_qh-\res(U)+\N(U,h)]
      -E_*G_*^{-1}(2\dot E_*^*h_t+\ddot E_*^*h).
\end{equation*}
The projection $\Pi_*$ is uniformly bounded on $H^{-1}$. Indeed it is self-adjoint on $L^2$, and its $\dot H^1$ bound follows from $\|\nabla e_j\|_2\lesssim s_j^{-1}$, $\|e_k\|_{\frac{6}{5}}\lesssim s_k$ and $|(G_*^{-1})_{jk}|s_k/s_j\lesssim1$ at the reference scales. Using \eqref{p5:eq:tuber}, \eqref{p6:eq:nonlinear-variation} and \eqref{p6:eq:test-derivatives}, we find
\begin{equation*}
 \|h_{tt}\|_{H^{-1}}\lesssim_A t^{-9},\qquad
 \|\delta h_{tt}\|_{H^{-1}}\lesssim_A d\,t^{-12}.
\end{equation*}
Here $\|V_q\|_{\frac{3}{2}}\lesssim1$ and $\|\delta V_q\|_{\frac{3}{2}}\lesssim_A dt^{-8}$ by rescaling. The fixed projection preserves time continuity.

On $q_*$ we have $\phi=0$. Away from $q_*$ its derivatives satisfy the difference rows of \eqref{p6:eq:test-derivatives}, with $d$ replaced by $A$. Rescaling also yields
\begin{equation*}
 (\|\phi_i\|_2,\|\phi_o\|_2)\lesssim_A (t^{-8},t^{-11}),\qquad
 (\|\phi_i\|_{H^1},\|\phi_o\|_{H^1})\lesssim_A(t^{-4},t^{-11}),
\end{equation*}
and the same bounds times $d$ for $\delta\phi_j$. Applying these bounds to \eqref{p6:eq:wave-current} yields
\begin{align*}
 (|C_{*,i}|,|C_{*,o}|)&\lesssim(t^{-14},t^{-18}),&
 (|C_{*,i}'|,|C_{*,o}'|)&\lesssim(t^{-10},t^{-18}),\\
 (|\delta C_i|,|\delta C_o|)&\lesssim_A d(t^{-17},t^{-20}),&
 (|\delta C_i'|,|\delta C_o'|)&\lesssim_A d(t^{-13},t^{-20}).
\end{align*}
For example, the inner acceleration term in the last entry is at most $C_A d(t^{-12}+t^{-9})t^{-4}$ by $H^{-1},H^1$ duality.
The remaining term of $2MB_j$ vanishes at $q_*$. By \eqref{p2:eq:mass} and \eqref{p2:eq:bprime},
\begin{equation*}
 (|b|,|\dot b|)\lesssim(t^{-6},t^{-7}),\qquad
 (|\delta b|,|\delta\dot b|)\lesssim_A d(t^{-14},t^{-15}).
\end{equation*}
The product rule implies
\begin{align*}
 (|\delta(b\mu')|,|\delta(b\lambda')|)&\lesssim_A d(t^{-18},t^{-19}),\\
 (|\delta(b\mu')'|,|\delta(b\lambda')'|)&\lesssim_A d(t^{-19},t^{-16}).
\end{align*}
The last power comes from $b\,\delta\lambda''$. These estimates prove \eqref{p6:eq:base} and \eqref{p6:eq:diff}.

\end{proof}

\subsubsection{Integration in the corrected momentum}

For a trial path we set $k_i=5,k_o=0$ and
\begin{equation*}
 f_j=M^{-1}(\mathcal F_j-\mathcal F_{*,j})-r_{*,j}-D_j
                      -k_j(k_j+1)t^{-2}y_j .
\end{equation*}
On the reference input $|f_i|+|f_o|\lesssim t^{-14}$. For both rows,
\begin{equation}\label{p6:eq:fdiff}
 |\delta f_j|\lesssim_A d\,t^{-17}
\end{equation}
by Lemmas~\ref{p6:lem:bare} and~\ref{p6:lem:feedback}.
For each row $k\in\{0,5\}$ we use the corrected momentum in \eqref{p6:eq:exact} and solve directly for the output position $Y$,
\begin{equation}\label{p6:eq:linear}
 (Y'-2B)'=k(k+1)t^{-2}Y+f .
\end{equation}
\begin{lemma}\label{p6:lem:green}
If $|f|\le F t^{-m}$, $|B|\le H t^{-n}$ with $m>k+2$ and $n>k+1$,
the selected solution satisfies
\begin{equation}\label{p6:eq:gin}
 |Y|+t|Y'|\lesssim_{m,n} F t^{2-m}+H t^{1-n}.
\end{equation}
\end{lemma}
\begin{proof}
Take $\phi=t^{-k}$, so $\phi''=k(k+1)t^{-2}\phi$, and set $J=\phi(Y'-2B)-\phi'Y$. The terms containing $Y$ in $J'$ cancel by \eqref{p6:eq:linear}, leaving
\begin{equation}\label{p6:eq:wronskian}
 J'=\phi f-2\phi'B,\qquad
 \left(\frac{Y}{\phi}\right)'=\frac{J}{\phi^2}+\frac{2B}{\phi}.
\end{equation}
At a finite $T$, $Y(T)=Y'(T)=0$ prescribes $J(T)=-2\phi(T)B(T)$ and $(Y/\phi)(T)=0$. At infinity both terminal values vanish in the class \eqref{p6:eq:gin}. Integrating successively from the terminal end, we obtain
\begin{equation*}
 |J(t)|\lesssim Ft^{1-m-k}+Ht^{-n-k},\qquad
 \left|\frac{Y(t)}{\phi(t)}\right|\lesssim Ft^{2-m+k}+Ht^{1-n+k}.
\end{equation*}
Indeed $|J(T)|\le2HT^{-n-k}$, and the exponents in the second integrand are $1-m+k<-1$ and $k-n<-1$. These integrals determine $Y$ uniquely. The identity $Y'=(J+\phi'Y)/\phi+2B$ proves the velocity bound.

Only the acceleration uses a derivative of $B$,
\begin{equation}\label{p6:eq:acceleration}
 Y''=k(k+1)t^{-2}Y+f+2B' .
\end{equation}

\end{proof}

\begin{proof}[Proof of Theorem~\ref{p6:th:main}]
We use \eqref{p6:eq:linear} and Lemma~\ref{p6:lem:green} to define $y\mapsto Y$ on the closed ball in the weighted $C^2$ space \eqref{p6:eq:norm}, a complete metric space also on the half-line. At $y=0$,
\begin{equation*}
 |Y_i|+t|Y_i'|\lesssim t^{-12},\qquad
 |Y_o|+t|Y_o'|\lesssim t^{-12},\qquad
 |Y_i''|\lesssim t^{-10},\quad |Y_o''|\lesssim t^{-14}.
\end{equation*}
The position and velocity use $B_i=O(t^{-14})$ directly. The weaker bound on $B_i'$ enters only the physical acceleration.

For two inputs, \eqref{p6:eq:diff}, \eqref{p6:eq:fdiff} and the Green estimate imply, for $j=i,o$,
\begin{equation*}
 |\delta Y_j|+t|\delta Y'_j|\lesssim_A dt^{-15},\qquad
 |\delta Y_i''|\lesssim_A dt^{-13},\quad |\delta Y_o''|\lesssim_A dt^{-16}.
\end{equation*}
For the last two bounds we used \eqref{p6:eq:acceleration} with $k_i=5$ and $k_o=0$. Let $A$ exceed twice the bound for the image $Y(0)$ of the reference correction $y=0$. By increasing $T_0$ we obtain, in the norm \eqref{p6:eq:norm},
\begin{equation*}
 \|Y(y)-Y(\bar y)\|_I\le C_AT_0^{-3}\|y-\bar y\|_I,\qquad
 \|Y(y)\|_I\le\|Y(0)\|_I+\frac{1}{2} A\le A .
\end{equation*}
At the fixed point, \eqref{p6:eq:linear} and \eqref{p6:eq:exact} imply $E^*E_*c=0$. This matrix is invertible, so $c_i=c_o=0$. Section~\ref{sec:tower-5} supplies energy continuity and the local $L^5_tL^{10}_x$ bound. Positivity, $\lambda'<0$ and the outer normalization follow from \eqref{p6:eq:branch} and \eqref{p6:eq:norm}.

\end{proof}

\subsubsection{The half-line solution}

\begin{corollary}\label{p6:cor:tower}
The half-line solution of Theorem~\ref{p6:th:main} is bounded and obeys
\begin{equation*}
 \|(\nabla(u-S),u_t-S_t)(t)\|_2=O(t^{-9}),\qquad
 \|(u,u_t)(t)-(P_{\lambda(t)}+Q_{\mu(t)},0)\|_{\HH}=O(t^{-5}).
\end{equation*}
The two profiles are centered, with $\mu\to1$ and $\frac{\lambda}{\mu}\to0$. Moreover,
\begin{equation}\label{p6:eq:outer-sharp}
 |\mu-\mu_*|+t|\mu'-\mu_*'|+t^2|\mu''-\mu_*''|\lesssim t^{-12}.
\end{equation}
\end{corollary}
\begin{proof}
Since $u=S+w+h$, \eqref{p3:eq:symbols} and \eqref{p6:eq:wavebound} control the first error by $Ct^{-9}$. The bare velocity satisfies $\|S_t\|_2\le\sqrt M(|\lambda'|+|\mu'|)\lesssim t^{-5}$ by \eqref{p6:eq:branch} and \eqref{p6:eq:norm}, which also imply the scale limits. On the fixed path, \eqref{p6:eq:base}--\eqref{p6:eq:diff} and \eqref{p6:eq:fdiff} imply
\begin{equation*}
 |f_o|\lesssim t^{-14},\qquad |B_o|\lesssim t^{-18},\qquad |B_o'|\lesssim t^{-16}.
\end{equation*}
Taking $k=0$ in \eqref{p6:eq:gin} and \eqref{p6:eq:acceleration} proves \eqref{p6:eq:outer-sharp}.
\end{proof}

\endgroup
\space
\subsection{The exact tower and its two profiles}\label{sec:tower-8}
\begingroup

\def\R{\mathbb R}
\def\HH{\mathcal H}
\def\D{\mathcal D}
\def\W{\mathcal W}
\def\St{\mathcal S}
\def\dd{\,\mathrm d}
\def\vct#1{\boldsymbol{#1}}

We complete Theorem~\ref{p8:th:tower} by reading off the scale
asymptotics, splitting the energy on the two spatial scales, and
identifying all profiles in any finite strong resolution.

\begin{proof}[Construction in Theorem~\ref{p8:th:tower} and proof of Proposition~\ref{p8:pr:refined-scales}]
Corollary~\ref{p6:cor:tower} proves \eqref{p8:eq:tower}. By \eqref{p6:eq:norm} for the inner scale and \eqref{p6:eq:outer-sharp} for the outer scale, both positions differ from \eqref{p2:eq:reference} by $O(t^{-12})$, and both velocities by $O(t^{-13})$. Thus \eqref{p8:eq:scale} holds and
\begin{equation*}
 \lambda'=-4at^{-5}-\frac{256}{15}a^2t^{-9}+O(t^{-13}),\qquad
 \mu'=\frac{20}{9}a^2t^{-9}+\frac{4B_0}{M}a^{\frac{5}{2}}t^{-11}+O(t^{-13}).
\end{equation*}
In particular $\lambda'<0<\mu'$ for large $t$. A nonzero time translation changes $at^{-4}$ at order $t^{-5}$, so \eqref{p8:eq:scale} fixes the time origin. We prove the energy and profile assertions below.
\end{proof}

\subsubsection{Energy on either side of the matching radius}\label{p8:sec:energy}

For a measurable set $D$ write
$N_D(\vec F)=\|(\nabla F,F_1)\|_{L^2(D)}+\|F\|_{L^6(D)}$.
For pairs of energy norm at most $R$, Cauchy--Schwarz and H\"older imply
\begin{equation}\label{p8:eq:density-continuity}
 \|e(\vec F)-e(\vec G)\|_{L^1(D)}\lesssim (R+R^5)N_D(\vec F-\vec G).
\end{equation}
For a stationary $F$, $\|\nabla F\|_2\le R$ and $\|\nabla w\|_2\le1$, we also have
\begin{equation}\label{p8:eq:local-stationary}
 \left|\int_D[e(F+w,k)-e(F,0)]\dd x\right|
 \lesssim_R \|\nabla w\|_2^2+\|k\|_2^2
 +\|\nabla w\|_2\min\{N_D(F,0),N_{D^c}(F,0)\}.
\end{equation}
Indeed $-\Delta F=F^5$, so the integrals of $\nabla F\cdot\nabla w-F^5w$ on $D$ and $D^c$ sum to zero. We bound either integral by H\"older and the remaining terms by Taylor's formula and Sobolev's inequality.

For $\vec z,\vec F,\vec G$ in this ball and $D=B_r$, we compare $e(\vec F_\lambda+\vec G_\mu)$ with the inner density on $D$ and with the outer density on $D^c$. Applying \eqref{p8:eq:density-continuity} also to each discarded density and zero, we obtain
\begin{equation}\label{p8:eq:density-decoupling}
 \begin{aligned}
 &\|e(\vec z)-e(\vec F_\lambda)-e(\vec G_\mu)\|_1\\
 &\qquad\lesssim_R\|\vec z-\vec F_\lambda-\vec G_\mu\|_{\HH}
        +N_{|x|>\frac{r}{\lambda}}(\vec F)+N_{|x|<\frac{r}{\mu}}(\vec G).
 \end{aligned}
\end{equation}
If $\frac{\lambda}{\mu}\to0$, the two tails vanish for $r=\sqrt{\lambda\mu}$. The same estimate holds for the gradient density alone.
\begin{proposition}\label{p8:pr:density}
For every $r(t)$ with $\frac{\lambda(t)}{r(t)}\to0$ and $\frac{r(t)}{\mu(t)}\to0$,
\begin{equation}\label{p8:eq:half}
 \int_{|x|<r(t)}e(\vec u)\dd x\longrightarrow E(P,0),\qquad
 \int_{|x|>r(t)}e(\vec u)\dd x\longrightarrow E(Q,0).
\end{equation}
At the Kelvin radius $R=\sqrt{\lambda\mu}$ we have the signed correction
\begin{equation}\label{p8:eq:kelvin-energy-split}
 \int_{|x|<R}e(\vec u)\dd x-E(P,0)
 = E(Q,0)-\int_{|x|>R}e(\vec u)\dd x
 =\frac{3C}{10}\rho^{\frac{5}{2}}+O(t^{-14}).
\end{equation}
For every bounded continuous spatial test function the signed measures converge as
\begin{equation*}
 e(\vec u(t))\dd x \ \rightharpoonup\
                    e(Q,0)\dd x + E(P,0)\delta_0 .
\end{equation*}
The velocity of the concentrating core is resolved one order more precisely by
\begin{equation}\label{p8:eq:kin}
 \big\|t^5\lambda^{\frac{3}{2}}u_t(t,\lambda\,\cdot)-4a\Lambda P\big\|_2
       =O(t^{-4}) .
\end{equation}
\end{proposition}

\begin{proof}
For the static gradient product we apply \eqref{p3:eq:overlap} with $p=a=b=1$, $\sigma=2$. For the mixed sextic terms Young's inequality and \eqref{p2:eq:kelvin-exchange} reduce the estimate to one integral. In outer units,
\begin{align*}
 \int |\nabla P_\rho|\,|\nabla Q|\dd x&\lesssim\rho^{\frac{5}{2}}(1+|\log\rho|),\\
 \|(P_\rho+Q)^6-P_\rho^6-Q^6\|_1&\lesssim\int|P_\rho|\,|Q|^5\dd x
 \lesssim\rho^{\frac{5}{2}}\int |x|^{-3}|Q|^5\dd x .
\end{align*}
By \eqref{p8:eq:tower}, $\|\nabla(u-S)\|_2\lesssim t^{-9}$ and $\|u_t\|_2^2\lesssim t^{-10}$. We apply \eqref{p8:eq:density-continuity} only to the positions and add the kinetic density $\frac{1}{2}|u_t|^2$. Together with the two overlap bounds this yields
\begin{equation*}
 \|e(\vec u)-e(P_\lambda,0)-e(Q_\mu,0)\|_1\lesssim t^{-9}.
\end{equation*} Integration on either side of $r(t)$ proves \eqref{p8:eq:half}. Testing against $\phi\in C_b(\R^3)$ in the two frames and using dominated convergence proves the measure limit. Energy conservation and $-\Delta Q=Q^5$ imply $E(\vec u)=2E(Q,0)$ and $E(Q,0)=\frac{1}{3}\|\nabla Q\|_2^2>0$.

For \eqref{p8:eq:kelvin-energy-split}, the Kelvin transform in $|x|=R$ fixes $S=P_\lambda+Q_\mu$. By \eqref{p2:eq:kelvin-flux},
\begin{equation*}
 \int_{|x|>R}e(S,0)\dd x=\frac{1}{2}E(S,0)+\frac{R}{4}\int_{\mathbb S^2}S(R\omega)^2\,d\omega .
\end{equation*}
Here $S(R\omega)=2\mu^{-\frac{1}{2}}Q(\sqrt\rho\,\omega)$. The terms $H_2,H_3$ in \eqref{common:quad}--\eqref{p2:eq:jet} are orthogonal on the sphere, so the boundary term is $\frac{C}{5}\rho^{\frac{5}{2}}+O(\rho^{\frac{7}{2}})$. Together with \eqref{p2:eq:potential}, this proves the exterior formula for $(S,0)$.

Write $w_1=u-S$. In the exterior first variation, the bulk term is bounded by $\|\Delta S+S^5\|_{\frac{6}{5}}\|w_1\|_6\lesssim \rho^{\frac{5}{2}}t^{-9}$ by \eqref{p3:eq:overlap}. Differentiating the Kelvin identity at $R$ yields $2R\partial_r S=-S$ there. The trace estimate following \eqref{p2:eq:kelvin-flux} therefore bounds the boundary pairing by $O(t^{-9}\rho^{\frac{5}{4}})$. The Taylor remainder is $O(t^{-18})$, and \eqref{p8:eq:tower} bounds the exterior velocity norm by $O(t^{-9}+|\mu'|+|\lambda'|(\lambda/R)^{\frac{3}{2}})=O(t^{-8})$. This proves the exterior estimate for $\vec u$. The interior estimate follows by conservation.
For $D\ge1$, the Kelvin tails $|Q|+r|\nabla Q|\lesssim r^2$ at zero and $|P|+r|\nabla P|\lesssim r^{-3}$ at infinity also yield
\begin{equation*}
 \|\nabla u\|_{L^2(D\lambda<|x|<\frac{\mu}{D})}
 +\|u\|_{L^6(D\lambda<|x|<\frac{\mu}{D})}
 \lesssim D^{-\frac{5}{2}}+t^{-9}.
\end{equation*}
For \eqref{p8:eq:kin}, the rescaling $g\mapsto\lambda^{\frac{3}{2}}g(\lambda\,\cdot)$ preserves the $L^2$ norm. Thus
\begin{equation*}
 \left\|\lambda^{\frac{3}{2}}u_t(t,\lambda\,\cdot)+\lambda'\Lambda P\right\|_2
 \le \|u_t-S_t\|_2+|\mu'|\sqrt M\lesssim t^{-9}.
\end{equation*}
We insert $\lambda'=-4at^{-5}+O(t^{-9})$ to obtain \eqref{p8:eq:kin}. In particular,
\begin{equation*}
 \|u_t\|_2=4a\sqrt M\,t^{-5}+O(t^{-9}),\qquad
 8a^2M-Ca^{\frac{5}{2}}=0 .
\end{equation*}
The second identity follows from $C\sqrt a=8M$. It cancels the order $t^{-10}$ terms in the kinetic and static energies.
\end{proof}

For the half-line fixed point, the signed symmetries and \eqref{p8:eq:density-decoupling} imply
\begin{equation*}
 u(t,Rx)=-u(t,x),\qquad R(x)=(-x_2,x_1,x_3),\qquad
 \|\nabla u(t)\|_2^2\longrightarrow2\|\nabla Q\|_2^2>0 .
\end{equation*}
If $u(t)$ were radial about $c$, it would also be radial about $Rc$. Composing the two central reflections would give a period $2(Rc-c)$. A nonzero $L^6$ function has no nonzero period, so $Rc=c$. Radiality would then imply $u(t,Rx)=u(t,x)$, contrary to the sign change.
This places the tower outside the radial class in Shen's nonexistence theorem \cite{RS}.

\subsubsection{Rigidity of the two profile frames}\label{p8:sec:profiles}

Fix $t_n\to\infty$. Set $\vct F_1=(P,0)$, $\vct F_2=(Q,0)$ and $(s_{1,n},s_{2,n})=(\lambda(t_n),\mu(t_n))$. In the notation $\D_{s,c}=S_{c,s,I}$ of Section~\ref{fp:start}, suppose there is a finite strong traveling-wave resolution

\begin{equation}\label{p8:eq:djkm}
 \vec u(t_n)=\W(t_n)\vct v+
  \sum_{j=1}^{N}\D_{\sigma_{j,n},c_{j,n}}\,\vct A_j+o_{\HH}(1),
 \qquad
 \vct A_j=(A_j,-\ell_j\cdot\nabla A_j),
\end{equation}
where $\vct v\in\HH$, the nonzero $\vct A_j$ are fixed traveling-wave snapshots, $|\ell_j|<1$, and their spatial frames are pairwise orthogonal. These are the properties needed from the good-time resolution of Duyckaerts--Jia--Kenig--Merle \cite[Thm.~1.1]{DJKM}. By \eqref{p8:eq:tower}, the data also have a two-profile decomposition with strong energy remainder.

\begin{proof}[Identification in Theorem~\ref{p8:th:tower}]
Proposition~\ref{fp:comparison} applied to \eqref{p8:eq:tower} and \eqref{p8:eq:djkm} gives $\vct v=0$ and $N=2$. After passing to a common subsequence there is a permutation $\pi$ of $\{1,2\}$ matching each proposed profile $j$ to the core $k=\pi(j)$. On this subsequence the relative frame parameters satisfy
\begin{equation*}
 \frac{\sigma_{j,n}}{s_{k,n}}\longrightarrow d_j>0,\qquad
 \frac{c_{j,n}}{s_{k,n}}\longrightarrow b_j,\qquad
 \D_{d_j,b_j}(A_j,-\ell_j\cdot\nabla A_j)=\vct F_k .
\end{equation*}
The second component implies $\ell_j\cdot\nabla A_j=0$. If $\ell_j\ne0$, Fubini's theorem and $A_j\in L^6$ force $A_j=0$ along almost every line parallel to $\ell_j$. Thus both velocities vanish. Since the profiles and their velocities in \eqref{p8:eq:djkm} are fixed, this identifies the unordered pair of profile classes for the given decomposition. If varying orthogonal maps are allowed in its frames, compactness of the orthogonal group allows us to take their limits on the same subsequence and absorb them into the profiles and velocities.

\end{proof}

For a complete Bahouri--G\'erard decomposition \cite{BG}, \eqref{p8:eq:tower} already supplies two orthogonal free profiles with an $o_{\HH}(1)$ remainder. By \cite[Lemma~3.3]{DKMprofiles}, every such decomposition has two nonzero profiles. After extraction and permutation the relative frame parameters converge to fixed changes of frame, including constant rescaled time shifts. In the representatives with zero time offset the data are $(P,0)$ and $(Q,0)$, so their nonlinear evolutions are stationary.

The finer scale coefficients in Proposition~\ref{p8:pr:refined-scales} use the two dilation constraints and the selected time origin. In particular the outer $t^{-10}$ coefficient comes from the modulation equation, since the whole-space error in \eqref{p8:eq:tower} is $O(t^{-9})$.

\endgroup
\space
\space

\section{Finite-time concentration and its profile corrections}\label{typeii:first}
We begin the proof of Theorem~\ref{typeii:seq-main}. In Sections~\ref{typeii:first}--\ref{typeii:merge:closing} we write $Q=Q_{\rm wall}$, $P=\mathcal K Q_{\rm cross}$ and $B=\mathcal KQ$. We set $\tau=T-t$, $M_o=\|\Lambda Q\|_2^2$ and $M_i=\|\Lambda P\|_2^2$. For a dilated profile $F_s$ we use the $L^2$-unitary dilation and its derivative,
\begin{equation}\label{typeii:eq:unitary-dilation}
 D_s f=s^{-\frac{3}{2}} f\!\left(\frac{\cdot}{s}\right),\qquad
 -\partial_s F_s=D_s\Lambda F,\qquad
 \partial_s D_s f=-\frac{1}{s}D_s\big(\frac{3}{2}+x\cdot\nabla\big)f .
\end{equation}
We set $K=(1+\Lambda)\Lambda$ and $L_F=-\Delta-5F^4$. Section~\ref{common:signed} identifies the signed kernels. To obtain a solution of the wave equation we must select the two corrections together with those paths.

\paragraph{Geometry and measurements.}
Near the origin, $Q=gx_3+2g x_1x_2+O(r^3)$ and $P=O(r^2)$. At infinity, $Q$ has a dipole of order $r^{-2}$, whereas
$P=d_\times r^{-3}\omega_1\omega_2+O(r^{-4})$ with $d_\times\ne0$. We set $H_2=2g x_1x_2$, $C_1=4\pi g^2$ and $C_2=\int P^5H_2>0$.
We work with functions odd under both $\operatorname{diag}(1,-1,-1)$ and $\operatorname{diag}(-1,1,-1)$, as are $Q,P$ and $B$. These two half-turns fix the origin. By \eqref{common:finite-kernel}, centered dilation spans each stationary kernel in this class. We set $\rho=\frac{\lambda}{\mu}$ and $w=5P^4\Lambda P$. For the correction $h_o$ arriving from scale $\mu$, $H(z)=\mu^{\frac{1}{2}}h_o(t,\mu z)$ has the following quadrupole moment in the inner variables
\begin{equation}\label{typeii:eq:measurement}
 \mathcal M_\rho H=\rho^{-2}\int w(y)H(\rho y)\,dy.
\end{equation}
The signed Taylor expansion and the moment identity \eqref{typeii:eq:quadrupolereadout}
remove the linear jet and isolate $\partial_1\partial_2H(0)$.
The physical factor $\lambda^{-1}\rho^{\frac{5}{2}}\asymp\tau^{14}$
has the order of the inner acceleration. To see which norm this force
requires, we choose a smooth, compactly supported signed $\phi$ with
$\int w\phi\ne0$. For $H_\rho(z)=\rho^2\phi\bigl(\frac{z}{\rho}\bigr)$,
we have
\begin{equation*}
 \|\nabla H_\rho\|_2=\rho^{\frac{5}{2}}\|\nabla\phi\|_2\to0,
 \qquad \mathcal M_\rho H_\rho=\int w\phi\ne0 .
\end{equation*}
\paragraph{Proof structure.}
We prove Theorem~\ref{typeii:seq-main} in three steps.
\begin{enumerate}
\item We match the middle dipole on the light cone in Section~\ref{typeii:sec:matchedcone}. The physical $L^2$ residual is $O(\tau^5)$, so the selected middle correction has energy norm and local rescaled derivatives of order $\tau^6$ (Section~\ref{typeii:sec:nonlinearmiddle}).
\item We insert its slope at the origin in the inner displacement and prescribe the momentum independently. The transfer bound \eqref{typeii:truepoltransfer} requires only one time derivative of this slope. We apply the constrained wave estimate of Section~\ref{typeii:sec:mixedinverse} to both cores and control the remaining wave at order $\tau^{16}$ (Section~\ref{typeii:sec:transport}).
\item At each positive terminal cutoff we select the two waves and solve both scale equations by a fixed point. We pass to the limit in the exact wave equation in Section~\ref{typeii:sec:closing} and compute the two limiting energies in Section~\ref{typeii:sec:terminalcutoffs}.
\end{enumerate}

\subsection{Stationary responses at the two concentrating scales}\label{typeii:sec:stationarycells}
We compute the two forces and cancel the leading residual at each scale. By \eqref{common:wall-jet} and antipodal invariance,
\begin{equation*}
 Q(x)=g x_3+2g x_1x_2+O(|x|^3),\qquad B(x)=Q(-x).
\end{equation*}
For the Kelvin crossing profile $P$ of Section~\ref{sec:tower-1}, the coefficient $d_\times$ in \eqref{common:quad} is negative, so $g d_\times>0$.
We apply the harmonic moment identity \eqref{p1:eq:harmonic-flux} to
the two incident polynomials. The dipole of $Q$ is
$Y_1=-g\omega_3$ by Kelvin inversion, and the incident polynomial of $B$
is $H_B=-g x_3$. For $P$ the lower moments vanish and $Y_2=d_\times\omega_1\omega_2$.
Integrating the two spherical harmonics, we obtain
\begin{equation}\label{typeii:eq:two-forces}
 C_1=3g^2\int_{\mathbb S^2}\omega_3^2=4\pi g^2,\qquad
 C_2=10gd_\times\int_{\mathbb S^2}\omega_1^2\omega_2^2=\frac{8\pi}{3}gd_\times>0.
\end{equation}
By \eqref{common:finite-kernel}, the kernels of $L_Q$, $L_P$ and $L_B$ in the common signed class are spanned by $\Lambda Q$, $\Lambda P$ and $\Lambda B$, respectively. The negative eigenspaces are finite dimensional, with smooth exponentially decaying eigenfunctions.

We set
\begin{equation*}
 \alpha_o=\frac{3C_1}{2M_o},\quad \beta_o=\frac{2C_1}{M_o},\qquad
 \alpha_i=\frac{5C_2}{2M_i}.
\end{equation*}
The regular-wall profile requires a dipole correction $Z$, and the
crossing profile requires a correction $W$ for the incident linear slope. They solve
\begin{equation}\label{typeii:eq:linear-responses}
 L_QZ=\alpha_o\Lambda Q-\beta_o KQ+5Q^4H_B,\qquad
 H_B(z)=\nabla B(0)\cdot z,\qquad L_PW=5P^4y_3 .
\end{equation}
Since $1+\Lambda$ is skew-adjoint on $L^2$, we have $\langle KF,\Lambda F\rangle=0$.
We apply \eqref{typeii:eq:force-moment} at $k=1$ to the right sides of \eqref{typeii:eq:linear-responses} and obtain
\begin{align*}
 \langle \alpha_o\Lambda Q-\beta_o KQ+5Q^4 H_B,\Lambda Q\rangle
 &=\alpha_o M_o-\frac{3}{2}C_1=0,\\
 \langle 5P^4 y_3,\Lambda P\rangle&=-\frac{3}{2}\int P^5 y_3=0 .
\end{align*}
In the last equality we used the quarter-turn symmetry.
We solve for $W$ on the signed energy complement and average over the quarter-turn, which preserves its equation and the two half-turn signs. The Newton source for $W$ is $5P^4(W+y_3)$. Its mean vanishes on every sphere and it decays as $O(r^{-11})$ by local elliptic boundedness. In \eqref{p3:eq:newton-moments} at $k=1$ its first absolute moment is finite, whence $W=O(r^{-2})$.

For the wall, \eqref{p3:eq:shared-response} with $k=1$ reads
\begin{equation}\label{typeii:eq:dipole-remainder}
 Z=(\beta_o+2\alpha_o)S_Q+Z_0,\qquad
 L_QZ_0=\alpha_o(\Lambda+2K)Q+5Q^4H_B .
\end{equation}
The moment condition was just verified. Either signed half-turn removes degree zero, so \eqref{p3:eq:newton-moments} implies $Z_0=O(r^{-1})$ with all fixed Euler derivatives. In particular if
$Q=r^{-2}Y+O(r^{-3})$ at infinity, then
\begin{equation}\label{typeii:eq:explicit-dipole-tail}
 Z=Z_\infty+O(r^{-1}),\qquad
 Z_\infty=-\frac{3}{8}(\beta_o+2\alpha_o)Y=-15\sqrt a\,Y .
\end{equation}

We write $H_\times=-b_\times y_1y_2$ for the crossing polynomial, where $b_\times=-d_\times>0$. By the flux identity and the wall expansion,
\begin{equation*}
 C_\times=\int P^5H_\times=5\int_{\mathbb S^2}H_\times(\omega)^2\,d\omega,\qquad
 H_2=2g y_1y_2=\kappa H_\times,\qquad \kappa=-\frac{2g}{b_\times}>0.
\end{equation*}
Thus $C_2=\kappa C_\times$. Scaling the crossing response from \eqref{p3:eq:crossing-response}, we obtain
\begin{equation*}
 L_P(\kappa V_\times)=\alpha_i\Lambda P+5P^4H_2 .
\end{equation*}
Its tail and fixed Euler derivatives are $O(r^{-1})$ by
Lemma~\ref{p3:lem:ell} and the explicit expression for $S_P$.

By $S_F=\frac{|y|^2\Lambda F}{4}$, the tail of $S_F$ is $O(r^{-1})$ for $P$ and $O(1)$ for $Q$,
and its gradient vanishes at zero for both. Linearity combines the quadratic response and the full velocity-square correction into the identity
\begin{equation}\label{typeii:eq:inner-cancellation}
 L_P\big(\kappa\rho^{\frac{5}{2}}V_\times+p_i^2 S_P\big)
 =\rho^{\frac{5}{2}}(\alpha_i\Lambda P+5P^4H_2)-p_i^2 KP .
\end{equation}
We write $d$ for the linear coefficient of the actual normalized middle
background. The full inner position before its cone cutoff is
\begin{equation*}
 \lambda^{-\frac{1}{2}}\!\left\{P(y)+\chi\bigl(\frac{y}{R}\bigr)
 \left[\rho^{\frac{3}{2}}d W(y)+\kappa\rho^{\frac{5}{2}}V_\times(y)+p_i^2 S_P(y)\right]\right\},
 \qquad R=\rho^{-\frac{1}{2}}.
\end{equation*}
Both terms in \eqref{typeii:eq:inner-cancellation} have the same cutoff $\chi\bigl(\frac{y}{R}\bigr)$, so we can estimate their spatial commutators and parameter derivatives together.
We keep the middle velocity-excess cell
$(p_o^2-\beta_o\mu^{\frac{3}{2}})S_Q$ with its own cutoff when matching the dipole
on the physical cone.

The leading acceleration coefficient is $\frac{\alpha_i\rho^{\frac{5}{2}}}{\lambda}$.
In inner units we set $D=\lambda^{\frac{1}{2}}U_b(t,\lambda y)$ and put
$u_\rho(r)=\rho^{\frac{1}{2}}\min(\rho r,1)$. The deterministic background
satisfies $|D|\lesssim u_\rho$. On $r\le2R$,
\begin{equation*}
 D=\rho^{\frac{3}{2}}d y_3+\rho^{\frac{5}{2}}H_{2,t}(y)+O(\rho^{\frac{7}{2}}r^3),
 \qquad H_{2,t}-H_2=O(\mu^{\frac{3}{2}})
\end{equation*}
in polynomial coefficients. After the elliptic cancellations, the two cutoff errors in normalized
$L^2$ are $O(\rho^{\frac{3}{2}}R^{-\frac{5}{2}})$ and
$O(\rho^{\frac{5}{2}}R^{-\frac{3}{2}})$. For the Taylor and exterior
terms we use $r^3P^4\in L^2$ and $u_\rho\le\rho^{\frac{3}{2}}r$ to obtain
\begin{equation*}
 \|\rho^{\frac{7}{2}}r^3P^4\|_2\lesssim\rho^{\frac{7}{2}},
 \qquad \|P^4D\|_{L^2(r>2R)}
   \lesssim \rho^{\frac{3}{2}}\|rP^4\|_{L^2(r>2R)}
   \lesssim\rho^{\frac{3}{2}}R^{-\frac{19}{2}}.
\end{equation*}
For the mixed products we use $u_\rho\le\rho^{\frac{3}{2}}r$ and $|P|\lesssim r^{-3}$. Scaling $z=\rho y$ in the second bound, we obtain
\begin{equation}\label{typeii:eq:mixedendpoints}
 \|P^3u_\rho^2\|_2\le\rho^3\|r^2P^3\|_2,\qquad
 \|P u_\rho^4\|_2\lesssim\rho^{\frac{7}{2}}\big\||z|^{-3}\min(|z|,1)^4\big\|_2 .
\end{equation}
Both fixed norms are finite. By interpolation, $\|P^2u_\rho^3\|_2\lesssim\rho^{\frac{13}{4}}$.
On $r\le2R$ the cutoff cells are
$O(\rho^{\frac{3}{2}}\langle y\rangle^{-2})$. Inserting a cell in the
mixed terms contributes $O(\rho^3)$ after its linearized term is removed.
As $R=\rho^{-\frac{1}{2}}$, the largest normalized cutoff error is $O(\rho^{\frac{11}{4}})$. Restoring physical units, we obtain the inner residual $O(\lambda^{-1}\rho^{\frac{11}{4}})=O(\tau^{17})$. The cutoff of the innermost tail on $r\asymp\tau$ is smaller there by $P_\lambda=O(\lambda^{\frac{5}{2}}r^{-3})$. For the first-order phase below, only the bare inner dilation enters the prescribed momentum. The time derivative of the leading polar position contributes $O(\frac{\rho^{\frac{3}{2}}}{\tau})$ in $\dot H^1$ for a slowly varying incident slope. We estimate the wave contribution to the slope in Section~\ref{typeii:sec:outertransfer}.

\subsubsection{Dilation probes and the reference residual}\label{typeii:sec:pathbook}
We first estimate the dilation tests on general scale paths, then compute the residual on the reference paths. We set $q_o=\mu$, $q_i=\lambda$, $\rho=\frac{\lambda}{\mu}$, and write
$T_{s,F}=D_s\Lambda F$ for the physical dilation test in \eqref{typeii:eq:unitary-dilation}.
We assume $q_o\asymp\tau^4$, $q_i\asymp\tau^{16}$ and
$|\dot q_j|\lesssim \frac{q_j}{\tau}$. For the residual computation we take $q_o=a\tau^4$, $q_i=b\tau^{16}$.

\paragraph{Gauge.}
Let $j_{q,F}=D_q j_F$ be the fixed annular position probes chosen
in Section~\ref{typeii:sec:mixedinverse}. They detect $\Lambda F$
and annihilate its negative modes. The physical velocity test is
$E_{q,F}=T_{q,F}$ times a smooth cone cutoff equal to one on both
probe supports. The resulting matrix $A_{kl}=\langle j_{q_k,F_k},
E_{q_l,F_l}\rangle$ has fixed positive diagonal entries. Since $\Lambda Q\bigl(\frac{x}{\mu}\bigr)=O(\rho)$ on the inner annulus and $\Lambda P\bigl(\frac{x}{\lambda}\bigr)=O(\rho^3)$ on the outer one,
\begin{equation}\label{typeii:eq:probe-overlap}
 A_{io}=O(\rho^{\frac{5}{2}}),\qquad A_{oi}=O(\rho^{\frac{3}{2}}).
\end{equation}
The Kelvin bounds persist under fixed logarithmic scale derivatives. In the order $(i,o)$, we can write \eqref{typeii:eq:probe-overlap} as
\begin{equation*}
 \left(\frac{q_l A_{kl}}{q_k}\right)_{kl}
  =\begin{pmatrix} A_{ii}&O(\rho^{\frac{3}{2}})\\
               O(\rho^{\frac{5}{2}})&A_{oo}\end{pmatrix},
 \qquad \det A=A_{ii}A_{oo}+O(\rho^4).
\end{equation*}
Both matrices have uniformly bounded inverses. Scaling the fixed probe and its dilation derivative, we obtain
\begin{equation*}
 \|j_q\|_{\frac{6}{5}}+\tau\|\dot j_q\|_{\frac{6}{5}}\lesssim_J q,
 \qquad \|j_q\|_2+\tau\|\dot j_q\|_2\lesssim_J1 .
\end{equation*}

\paragraph{Residual on the reference paths.}
We write $P_o=\mu^{\frac{3}{2}}$, $P_i=\rho^{\frac{5}{2}}$ and $\zeta=\zeta(\frac{|x|}{\tau})$ for the cone cutoff.
The identity
\begin{equation}\label{typeii:eq:path-acceleration}
 \partial_t^2 F_{q_j}=D_{q_j}\big[-q_j''\Lambda F_j+(\frac{(q_j')^2}{q_j})KF_j\big]
\end{equation}
holds on every path. The definitions of $a,b$ imply
$q_j''=\frac{\alpha_jP_j}{q_j}$ and $(q_o')^2=\beta_oP_o$ on the reference paths. Inserting these coefficients into \eqref{typeii:eq:path-acceleration}, we cancel the leading accelerations, the incident Taylor terms and the full inner velocity square $(q_i')^2$ by \eqref{typeii:eq:linear-responses} and \eqref{typeii:eq:inner-cancellation}.
We use the full coefficient of $z_3$ in $\mu^{\frac{1}{2}}U_b(t,\mu z)$ in the first
inner cell. The residual is $F_o+F_i$, where $F_o$ is the corrected
middle residual including its cone commutators and $F_i$ the remaining
inner insertion. They satisfy
\begin{equation}\label{typeii:eq:reference-residual}
 \|F_o\|_2\lesssim\tau^3,\qquad \|F_i\|_2\lesssim\tau^{17},\qquad
 |\mu^{\frac{5}{2}}F_o(t,\mu z)|\lesssim\mu^{\frac{3}{2}}|z|\quad(|z|\le1).
\end{equation}
The middle norm is the reference case of \eqref{typeii:eq:middle-symbol} below. Near the origin the middle residual is smooth and odd, and its first normalized derivatives are $O(\mu^{\frac{3}{2}})$ by \eqref{typeii:eq:path-acceleration} and the elliptic cells. This yields the last bound.
For the inner
profile we combine Section~\ref{typeii:sec:stationarycells} with the
two time derivatives of the corrections. In physical $L^2$ their joint bound is
\begin{equation*}
 \|F_i\|_2\lesssim\frac{1}{\lambda}
 \left[\rho^{\frac{11}{4}}
 +\left(\frac{\lambda}{\tau}\right)^2
  (\rho^{\frac{3}{2}}+\rho^{\frac{5}{2}}R^{\frac{1}{2}})\right]
 \lesssim\tau^{17},\qquad R=\rho^{-\frac{1}{2}} .
\end{equation*}
On $|x|\asymp\tau$ both intermediate cutoffs vanish, the stationary inserted
fields are $O(\tau^4)$ and their first derivatives are $O(\tau^3)$, so the
cone commutators contribute $O(\tau^{\frac{7}{2}})$ in physical $L^2$.
We match the middle field to its dipole on the light cone in Section~\ref{typeii:sec:matchedcone}.

\subsection{The common virial and endpoint estimates}\label{typeii:sec:terminalestimate}
The middle wave and the two-core correction use the same estimate for
a moving potential. We formulate the estimate in physical time. We write $s_o=\mu$,
$s_i=\lambda$, $V_j=5(F_j)_{s_j}^4$ and $L=-\Delta-\sum V_j$. For one core the sum contains only $j=o$. We set $\gamma_j=\frac{\dot s_j}{s_j}$. For the constructed Lipschitz paths we use the rates of the reference powers,
\begin{equation*}
 \gamma_o^0=-\frac{4}{\tau},\qquad \gamma_i^0=-\frac{16}{\tau},\qquad
 \max_j|\gamma_j-\gamma_j^0|\lesssim\tau^{d-1}\quad(d>0).
\end{equation*}
On a smooth path one may also take $\gamma_j^0=\gamma_j$ whenever $\gamma_i\le\gamma_o<0$, $|\gamma_j|\lesssim\tau^{-1}$ and $|\dot\gamma_j|\lesssim\tau^{-2}$.

\paragraph{One localized virial identity.} We combine two cutoffs of the localized wave energy--virial correction of \cite[Section~3]{J}, one for each changing scale.
We work on the position space orthogonal to a compact dual of the dilation
and to the negative eigenfunctions at every active scale. For these constraints we claim
\begin{equation}\label{typeii:eq:coercive-gap}
 \langle LZ,Z\rangle\ge c_J\|\nabla Z\|_2^2 .
\end{equation}
By \eqref{common:finite-kernel}, the compact displacement dual detects the only kernel direction at each scale. We apply \eqref{p4:eq:frame-gap} with these duals and the negative eigenfunctions. Their pairings vanish by the constraints, which proves \eqref{typeii:eq:coercive-gap}.

We choose large radial cutoffs $\chi_j(\frac{x}{s_j})$ equal to one on the
corresponding dual supports, with $|r\chi_j'|\le\delta$ by a slow
transition in $\log r$. We write $c=\gamma_o^0\chi_o$ at one scale and
$c=\gamma_o^0\chi_o+(\gamma_i^0-\gamma_o^0)\chi_i$ at two scales. The
inner cutoff lies inside the outer plateau on a late interval. We set
\begin{equation*}
 E=\frac{1}{2}(\langle LZ,Z\rangle+\|W\|_2^2),\qquad
 E_*=E+\langle W,c\Lambda Z\rangle .
\end{equation*}
Then $|E_*-E|\lesssim_{J,\chi}(\frac{\mu}{\tau})E$. For
$Z_t=W+F_1$ and $W_t=-LZ+F_2$, we differentiate to obtain
\begin{align}
 (E_*)_t &\ge -C_J\left(\frac{\delta}{\tau}+\frac{C_\chi\mu}{\tau^2}
       +\max_j|\gamma_j-\gamma_j^0|\right)E \notag\\
     &\quad{}
     +\langle LZ,F_1\rangle+\langle W,c\Lambda F_1\rangle
     +\langle W+c\Lambda Z,F_2\rangle . \label{typeii:eq:commonvirial}
\end{align}
Indeed, we apply the three identities \eqref{p4:eq:virial-cancellation} with $(h,v)=(Z,W)$ and this nested cutoff $c$. The terms containing $F_1,F_2$ are exactly the three pairings in \eqref{typeii:eq:commonvirial}. The potential term splits according to $\gamma_j-c=(\gamma_j-\gamma_j^0)+(\gamma_j^0-c)$. The first coefficient multiplies $(2+D)V_j$, whose $L^{\frac{3}{2}}$ norm is scale invariant. We first choose the plateau radii large and then take $\tau_0$ small. The $L^{\frac{3}{2}}$ tails of $(2+D)V_j$ and the bound $\frac{|x|}{\mu}\lesssim_\chi\rho$ on the inner cutoff imply
\begin{equation*}
 \|r\nabla c\|_\infty+\sum_j\|(\gamma_j^0-c)(2+D)V_j\|_{\frac{3}{2}}
   \lesssim\frac{\delta}{\tau}.
\end{equation*}
Since $c\le0$, the term $-\int cW^2$ is nonnegative. Hardy's and Sobolev's inequalities bound the remaining spatial errors by $C_J(\frac{\delta}{\tau}+\max_j|\gamma_j-\gamma_j^0|)E$. Time differentiation of $c$ uses the prescribed coefficients and the first derivatives of the cutoff radii. On the support of $c_t$ we have $|x|\lesssim_\chi\mu$, whence
\begin{equation*}
 |\langle W,c_t\Lambda Z\rangle|\lesssim_\chi
    \frac{\mu}{\tau^2}\|W\|_2\|\nabla Z\|_2 .
\end{equation*}
Together with \eqref{typeii:eq:coercive-gap}, these estimates prove \eqref{typeii:eq:commonvirial}.
\paragraph{Weighted terminal selection.}
We integrate the continuous energy in physical time and the negative modes in the time variables of their respective scales. For $p>-1$ we fix $0<\eta<p+1$ and use the clocks $\sigma_j'=\frac{1}{s_j}$ in physical time. Let $n$ be the continuous energy norm and $U_j$ the modal sizes. Write $f_c$ for the size of the applied force in the continuous estimate and $f_j$ for its direct contribution to mode $j$. Suppose that
\begin{equation}\label{typeii:eq:terminal-kernels}
\begin{aligned}
 n(\tau)&\lesssim\int_\epsilon^\tau (\frac{\tau}{v})^\eta
          [f_c(v)+v^{-1}\textstyle\sum_j U_j(v)]\,dv,\\
 U_j(\tau)&\lesssim\int_\epsilon^{\tau_0}
 e^{-c_0|\sigma_j(T-\tau)-\sigma_j(T-v)|}
 [f_j(v)+v^{-1}(n(v)+\textstyle\sum_k U_k(v))]\,dv .
\end{aligned}
\end{equation}
The first inequality uses zero terminal continuous data. In the second we integrate the stable and unstable modes from opposite endpoints. We measure the two components by
\begin{equation*}
 N=\sup_{\epsilon\le\tau\le\tau_0}\frac{n(\tau)}{\tau^{p+1}},\qquad
 D=\sup_{j,\,\epsilon\le\tau\le\tau_0}\frac{U_j(\tau)}{s_j(T-\tau)\tau^p}.
\end{equation*}
For each $j$ set $y(v)=\sigma_j(T-v)$, $w(v)=s_j(T-v)v^p$ and $\nu=\sup_{\epsilon\le v\le\tau_0}\frac{\mu(T-v)}{v}$. We need only
\begin{equation}\label{typeii:eq:clock-mass}
 \int_{[r,2r]\cap[\epsilon,\tau_0]} f_c(v)\,dv\le Mr^{p+1},\qquad
 \sup_{y_*}\int_{\substack{\epsilon\le v\le\tau_0\\ |y(v)-y_*|\le1}}
              \frac{f_j(v)}{w(v)}\,dv\le M .
\end{equation}
These hypotheses include pointwise forcing $f_c,f_j\le Mv^p$ up to a constant, and the following cell bounds. On a dyad $v\asymp r$ let $\ell\asymp s_{\min}(T-r)$ be the nominal cell length, using the smallest active scale and retaining $\ell$ for clipped cells. Suppose each force has integral at most $M\ell r^p$ on each cell. Summation over a dyad proves the first bound. Since $|dv|=s_j(T-v)|dy|$ and $|\partial_y\log v|+|\partial_y\log w|\lesssim_p\nu$, a unit clock interval has physical length $O(s_j(T-v_*))$ for any $v_*$ on it. The interval meets only a bounded number of dyads, and its two end cells have nominal lengths $O(s_{\min}(T-v_*))$. Thus all cells meeting it have total nominal length $O(s_j(T-v_*))$. Dividing their masses by $w(v_*)$ proves the second bound, uniformly in the clipped lengths.

Choose $\tau_0$ so that $|\partial_y\log w|\le\frac{c_0}{2}$. In the continuous kernel of \eqref{typeii:eq:terminal-kernels} sum over $v\asymp 2^{-k}\tau$. In the modal kernel sum over $k\le|y(v)-y(\tau)|<k+1$. By \eqref{typeii:eq:clock-mass},
\begin{align*}
 \int_\epsilon^\tau\left(\frac{\tau}{v}\right)^\eta f_c(v)\,dv
 &\lesssim M\tau^{p+1}\sum_{k\ge0}2^{-k(p+1-\eta)},\\
 \int_\epsilon^{\tau_0}e^{-c_0|y(v)-y(\tau)|} f_j(v)\,dv
 &\lesssim_p M w(\tau)\sum_{k\ge0} e^{-\frac{c_0}{2}k}.
\end{align*}
Both series converge. Since $v^{-1}(n+\sum_jU_j)\lesssim (N+\nu D)v^p$, the same two estimates applied in \eqref{typeii:eq:terminal-kernels} yield
\begin{equation}\label{typeii:eq:weighted-bounds}
 N\lesssim_p M+\nu D,\qquad D\lesssim_p M+N+\nu D .
\end{equation}
Substituting the first bound in the second and using $\nu\lesssim\tau_0^3$, we obtain $N+D\lesssim_p M$, uniformly in the positive cutoff.
\subsection{Interior regularity and the signed quadrupole}\label{typeii:sec:outermoment}
A local propagation estimate suffices for the middle wave used in the
construction. We write $v=\tau$ as increasing time, fix a target $v=\sigma$,
and set $m=\mu(\sigma)\asymp\sigma^4$, $L=\frac{\sigma}{m}$,
$z=\frac{x}{m}$, $b=\frac{v}{m}$, $H(b,z)=m^{\frac{1}{2}}h(v,mz)$. The scales $\mu(v)$ and
$m$ are comparable for $\frac{\sigma}{8}\le v\le\sigma$. Suppose the
cone-supported wave has energy $\le A v^d$ and, on every such slab,
satisfies
\begin{equation*}
 (\partial_b^2-\Delta_z-V_{\rm ext})H=WH+F_{\rm sm}+P .
\end{equation*}
In macroscopic coordinates $(\theta,X)=(b/L,z/L)$, assume the exterior coefficient
$\mathscr V_\sigma(\theta,X)=L^2V_{\rm ext}(L\theta,LX)$ is smooth and,
for each fixed integer $M$, uniformly bounded in $C^M([1/8,1]\times B_3)$
independently of the slab and the trial path. These bounds include mixed
time and space derivatives. A fixed spatial cutoff in $X$ gives a compactly
supported extension agreeing on $B_2$.
Suppose $W$ is supported in $|z|\le \frac{b}{8}$ and, for $0\le j\le3$,
\begin{equation*}
 |\nabla_z^j W|\lesssim_j\langle z\rangle^{-5-j},\qquad
 |\nabla_z^j P|\lesssim_j A\sigma^d\langle z\rangle^{-2-j},
 \qquad \|F_{\rm sm}\|_{H^3_z}\lesssim \frac{A\sigma^d}{L} .
\end{equation*}
The pressure $P$ is supported in the cone. At a positive terminal
cut assume that the data have the same global middle-scaled energy
bounds after up to three spatial commutations. Then
\begin{equation}\label{typeii:movingouterjet}
 \|\mu^{\frac{1}{2}}h(v,\mu z)\|_{C^{2,\delta}(|z|\le2)}
 +\|\mu^{\frac{3}{2}}\partial_v h(v,\mu z)\|_{C^{1,\delta}(|z|\le2)}
 \lesssim_\delta A v^d,\qquad 0<\delta<\frac{1}{2} .
\end{equation}
For the mixed system, the momentum component has the stated local bound at every time. The first equation gives the same bound for $\partial_v h$ almost everywhere on a Lipschitz trial path when its first-row source $G$ satisfies the spatial bounds of Section~\ref{typeii:sec:perturbedmatched}. This argument uses only spatial derivatives of $G$. Sections~\ref{typeii:sec:matchedcone} and~\ref{typeii:sec:perturbedmatched} check the hypotheses at $d=6$.

\paragraph{Compact terminal data.}
We cut each negative eigenpacket off inside $r<\frac{\tau}{32}$ before projecting
it onto the constraint spaces in Section~\ref{typeii:sec:mixedinverse}.
The cutoff changes the packet and its fixed dilation derivatives by
$O_N(\tau^N)$ for every $N$. The Gram corrections are smooth cone-supported columns. In the
single middle frame the terminal displacement and velocity coefficients
are $\frac{p_j}{\mu}$ and $b_j$. By \eqref{typeii:eq:mixed-wave-bound} their absolute values are $\lesssim A\mu\tau^{d-1}\lesssim A\tau^d$.
The fixed profiles and their cone cutoffs have uniform $H^4$ and
$H^3$ norms, respectively. The terminal continuous data vanish. Finite propagation therefore keeps the selected middle wave
inside its cone when the applied sources have that support.

\paragraph{Interior smoothing by propagation on spheres.}
Starting with the energy and local $L^6$ bounds at order one, we raise the spatial derivative order from $n$ to $n+1$ for $1\le n\le3$. The bound is $C_n A v^d$ in local $L^2$ in units of $\mu(v)$, on unit balls centered in $|x_0|\le \frac{v}{4}$. At the induction step we may use a fixed larger ball by a finite covering when needed. Slab comparability transfers the inductive estimates to units of $m$ at any $v\in[\frac{\sigma}{8},\sigma]$. Every ball intersecting the differentiated core potential has its center well inside the previous cone, because the potential is supported in $|x|\le \frac{v}{8}$. Constants from the finite covering are independent of $L$.

We first propagate the initial data. In macroscopic units the retarded kernel of the smooth exterior potential is smooth strictly inside the characteristic cone. On compact subsets of that interior its derivatives in the source and target variables are uniformly bounded. In the finite Hadamard expansion \cite[Chapter~2]{BGP}, each retained transport coefficient and each fixed derivative of it is bounded by finitely many of the stated $C^M$ seminorms of $\mathscr V_\sigma$. Truncating beyond the order of the target derivatives and the one source derivative needed for displacement data leaves a remainder controlled by differentiated energy estimates and Sobolev embedding on this fixed slab. The shrinking coefficient $W$ and the pressure are treated below as sources, using their spatial bounds.
If the slab starts at $v_0=\frac{\sigma}{8}$, its initial data vanish outside $|x|\le v_0$. At $|x_0|\le\frac{\sigma}{4}$ the initial sphere lies in $|x|\ge\frac{5\sigma}{8}$. Every source point in the support of the initial data lies strictly inside the past cone, at a fixed macroscopic distance from its boundary. In macroscopic units we pair the kernel derivatives with the velocity in $L^2$ and the displacement in $L^6$ on a fixed ball. Each pairing is $O(A\sigma^d)$. Returning to middle units gains $L^{-j-\frac{1}{2}}$ for $j$ displacement derivatives. If the slab starts at $v=\epsilon$, the prescribed data already have the three commuted energy bounds. For these data and for $F_{\rm sm}$ we apply the differentiated energy inequality on a slab of length $O(L)$, using
\begin{equation*}
 \|\nabla^j V_{\rm ext}\|_3\lesssim_j L^{-1-j},\qquad
 \int\|F_{\rm sm}\|_{H^3}\,db\lesssim A\sigma^d .
\end{equation*}
Here we extend the exterior potential by a cutoff on scale $L$ and use the Sobolev inequality for the undifferentiated displacement.

It remains to treat $F=WH+P$ and the first-row source $G$. The induction hypothesis and the source bounds imply, on every unit ball,
\begin{equation*}
 \|F\|_{H^{n-1}(B_1(y))}+\|G\|_{H^n(B_1(y))}
 \lesssim A\sigma^d\langle y\rangle^{-2}.
\end{equation*}
Since $\langle y\rangle^{-2}\in L^2(\mathbb R^3)$, the global norms are $O(A\sigma^d)$. Let $(J,K)$ solve the mixed equation with zero data and these two sources,
\begin{equation*}
 J_b=K+G,\qquad K_b=\Delta J+V_{\rm ext}J+F .
\end{equation*}
Commute at most $n-1$ spatial derivatives. By $\|\nabla^j V_{\rm ext}\|_3\lesssim L^{-1-j}$ and Sobolev, the potential contributes at most $C L^{-1}$ times the sum of the commuted energies. Gronwall's inequality on a slab of length $O(L)$ yields
\begin{equation*}
 \sup_b\big(\|J\|_6+\|\nabla J\|_{H^{n-1}}+\|K\|_{H^{n-1}}\big)\lesssim L A\sigma^d .
\end{equation*}
We now use the free Duhamel formula for the same pair. Its extra second-row source is $V_{\rm ext}J$. Since $n\ge1$, the product rule yields
\begin{equation}\label{typeii:eq:slab-potential-energy}
 \|\nabla^n(V_{\rm ext}J)\|_2
 \lesssim L^{-1-n}\|J\|_6+
   \sum_{k=1}^n L^{-2-n+k}\|\nabla^k J\|_2
 \lesssim \frac{A\sigma^d}{L}.
\end{equation}
Here the undifferentiated $J$ uses $\nabla^n V_{\rm ext}\in L^3$, and all other terms use the macroscopic $L^\infty$ coefficient bounds. Integration over the slab costs $O(L)$, so free energy controls this entire contribution at order $n$ by $O(A\sigma^d)$.
For the free response to $(G,F)$ we use Kirchhoff's formula and energy conservation on the Cauchy pair $(G,F)$. At elapsed time $u$, the output on $B_1(z_0)$ uses only the data in $\{\big||z-z_0|-u\big|\le2\}$. After $n$ spatial commutations we cut off $\nabla^n(G,F)$ on a shell enlarged by a fixed width, with bounded cutoff gradient. Both output rows are bounded by the $L^2$ norms there of $\nabla^n F$, $\nabla^{n+1}G$ and $\nabla^nG$. These have local bounds $O(A\sigma^d\langle z\rangle^{-p})$, with $p=5$ for the core product and $p\ge 2+n$ for the other terms. It suffices to use
\begin{equation}\label{typeii:eq:sphere-shell}
 \int_0^\infty\left(\int_{|\,|z-z_0|-u\,|\le5}\langle z\rangle^{-2p}\,dz\right)^{\frac{1}{2}}du
 \lesssim_p1\qquad (p>2).
\end{equation}
For $a=|z_0|$, polar integration about $z_0$ yields
\begin{equation*}
 \int_{|z-z_0|=r}\frac{dS_z}{(1+|z|^2)^2}
 =\frac{4\pi r^2}{[1+(r-a)^2][1+(r+a)^2]}
 \le\frac{4\pi}{1+(r-a)^2}.
\end{equation*}
Since $|z|\ge|r-a|$ on the sphere, the surface integral of $\langle z\rangle^{-2p}$ is $O(\langle r-a\rangle^{2-2p})$ for $p\ge2$. Integrating over $|r-u|\le5$ bounds the integrand in \eqref{typeii:eq:sphere-shell} by $C_p\langle u-a\rangle^{1-p}$. At $a=0$ it is comparable to $u^{1-p}$ for $u\ge10$, so $p>2$ is sharp.

Smooth-data approximation and the uniform energy estimate justify the argument for the selected finite-terminal solutions. At $n=3$ we have $H^4$ for displacement and $H^3$ for velocity
on each unit ball. Since $H^4\hookrightarrow C^{2,\delta}$ and $H^3\hookrightarrow C^{1,\delta}$ for $0<\delta<\frac{1}{2}$, we obtain \eqref{typeii:movingouterjet}.

\paragraph{Application to the projected quadrupole.}
For a signed $H\in\dot H^1\cap C^{2,\delta}(B_1)$ the observation
in \eqref{typeii:eq:measurement} has the explicit expansion
\begin{equation}\label{typeii:eq:quadrupolereadout}
 \mathcal M_\rho H=-\frac{5C_2}{4g}\,\partial_1\partial_2 H(0)
 +O\!\left(\rho^\delta\|H\|_{C^{2,\delta}(B_1)}
                  +\rho^{10}\|\nabla H\|_2\right).
\end{equation}
The two half-turns restrict the Taylor polynomial to two monomials. For $|z|\le1$ we have
\begin{equation*}
 H(z)=\partial_3H(0)z_3+\partial_{12}H(0)z_1z_2
         +O\big(\|H\|_{C^{2,\delta}(B_1)}|z|^{2+\delta}\big).
\end{equation*}
For $w=5P^4\Lambda P$, the quarter-turn symmetry and \eqref{typeii:eq:force-moment} with $k=2$ imply
\begin{equation*}
 \int_{|y|<R}y_3w\,dy=0\quad(R>0),\qquad
 \int y_1y_2w\,dy=-\frac{5}{2}\int y_1y_2P^5\,dy=-\frac{5C_2}{4g}.
\end{equation*}
We insert the expansion at $z=\rho y$ in \eqref{typeii:eq:measurement} on $|y|\le\rho^{-1}$. The remainder is bounded by $C\rho^\delta\|H\|_{C^{2,\delta}(B_1)}$, since $|w(y)|\lesssim\langle y\rangle^{-15}$. On the complement we use
\begin{align*}
 \int_{|y|>\rho^{-1}}|y|^2|w|\,dy&\lesssim\rho^{10},\\
 \rho^{-2}\|w\|_{L^{\frac{6}{5}}(|y|>\rho^{-1})}\|H(\rho\,\cdot)\|_6
 &\lesssim\rho^{-2}\rho^{\frac{25}{2}}\rho^{-\frac{1}{2}}\|\nabla H\|_2
 =\rho^{10}\|\nabla H\|_2 .
\end{align*}
This proves \eqref{typeii:eq:quadrupolereadout}. For $H(z)=\mu^{\frac{1}{2}}h(v,\mu z)$, \eqref{typeii:movingouterjet} bounds the observation by $CAv^d$. Fixed Euler derivatives of $w$ have the same parity and decay.
\subsubsection{Transfer from the middle scale to the inner scale}
\label{typeii:sec:outertransfer}
We use $\rho=\frac{\lambda}{\mu}$ and $R=\rho^{-\frac{1}{2}}$, with
\begin{equation*}
 \mu\asymp\tau^4,\qquad \lambda\asymp\tau^{16},\qquad
 \frac{|\dot\mu|}{\mu}+\frac{|\dot\lambda|}{\lambda}\lesssim\tau^{-1}.
\end{equation*}
By \eqref{typeii:movingouterjet} and the first middle equation, we may assume
\begin{equation*}
 \|\nabla h_o\|_2+
 \|\mu^{\frac{1}{2}}h_o(t,\mu\,\cdot)\|_{C^2(B_1)}+
 \|\mu^{\frac{3}{2}}\partial_t h_o(t,\mu\,\cdot)\|_{C^1(B_1)}
 \le A\le1 .
\end{equation*}
The two half-turn symmetries imply $\partial_1h_o(t,0)=\partial_2h_o(t,0)=0$. Its coefficient in the middle frame is
\begin{equation}\label{typeii:eq:average-slope}
 d_h=\mu^{\frac{3}{2}}\partial_3 h_o(t,0).
\end{equation}
We differentiate at the origin, which is fixed under dilation, and obtain
\begin{equation}\label{typeii:eq:slope-test}
 \mu\dot d_h=\frac{3}{2}\dot\mu\,d_h+\mu^{\frac{5}{2}}\partial_3\partial_t h_o(t,0).
\end{equation}
For Lipschitz paths this identity holds almost everywhere. Indeed, on compact time intervals the local spatial bounds in \eqref{typeii:movingouterjet} allow evaluation of $\partial_3 h_o$ and $\partial_3\partial_t h_o$ at zero in the distributional time derivative. Since $|\dot \mu|\lesssim\frac{\mu}{\tau}$, we obtain
\begin{equation}\label{typeii:eq:slope-rate}
 |d_h|+\mu|\dot d_h|\lesssim A .
\end{equation}

We add the position and momentum pair $(w_h,0)$ to the deterministic
insertion, where $W$ solves \eqref{typeii:eq:linear-responses} and
\begin{equation*}
 w_h=\lambda^{-\frac{1}{2}}\rho^{\frac{3}{2}}d_h\,\chi\bigl(\frac{y}{R}\bigr)W(y),\qquad y=\frac{x}{\lambda} .
\end{equation*}
Let $U_b$ and $U$ be the deterministic middle and two-core positions,
respectively.
The position defect and the spatial terms left after extracting the
terms with at least two explicit factors of $h_o$ satisfy
\begin{equation}\label{typeii:truepoltransfer}
\begin{aligned}
 \|\partial_t w_h\|_{\dot H^1}&\lesssim \frac{A\,\rho^{\frac{3}{2}}}{\mu},\\
 \|\Delta w_h+(U+w_h)^5-U^5+5[(U+w_h)^4-U_b^4]h_o\|_2
      &\lesssim \frac{A\,\rho^{\frac{5}{2}}}{\lambda} .
\end{aligned}
\end{equation}
Since $\rho=\frac{\lambda}{\mu}$, the two powers in \eqref{typeii:truepoltransfer} are equal,
\begin{equation}\label{typeii:eq:transfer-scale}
 \frac{\rho^{\frac{3}{2}}}{\mu}=\frac{\rho^{\frac{5}{2}}}{\lambda}
 \asymp\tau^{14}.
\end{equation}
\begin{proof}[Proof of \eqref{typeii:truepoltransfer}]
The cutoff profile $\chi\bigl(\frac{y}{R}\bigr)W$ and its scale derivatives are bounded
in $\dot H^1_y$ since $|\nabla^j W|\lesssim\langle y\rangle^{-2-j}$.
We differentiate $w_h$ in physical time and use \eqref{typeii:eq:slope-rate} and $\mu\lesssim\tau$ to obtain
\begin{equation*}
 \|\partial_t w_h\|_{\dot H^1}\lesssim
 \rho^{\frac{3}{2}}\left(|\dot d_h|+\frac{|d_h|}{\tau}\right)
 \lesssim \frac{A\rho^{\frac{3}{2}}}{\mu}.
\end{equation*}
In inner units we set $H=\lambda^{\frac{1}{2}}h_o(t,\lambda y)$, $D=\lambda^{\frac{1}{2}}U_b(t,\lambda y)$.
On $|y|\le\rho^{-1}$ the signed Taylor expansion and \eqref{typeii:eq:average-slope} imply
\begin{equation*}
 \left|H-\rho^{\frac{3}{2}}d_hy_3\right|
 \lesssim A\rho^{\frac{5}{2}}|y|^2,
 \qquad
 |H|\lesssim A\rho^{\frac{3}{2}}|y|,\quad |D|\lesssim\rho^{\frac{3}{2}}|y|.
\end{equation*}
By \eqref{typeii:eq:linear-responses}, the elliptic term cancels $5P^4\rho^{\frac{3}{2}}d_hy_3$. On $|y|\le\rho^{-1}$ the Taylor and mixed terms obey
\begin{align*}
 \|P^4(H-\rho^{\frac{3}{2}}d_hy_3)\|_2
 &\lesssim A\rho^{\frac{5}{2}}\||y|^2P^4\|_2,\\
 \big\|\big(|P|^3|D|+|P|^2|D|^2+|P||D|^3\big)H\big\|_2
 &\lesssim A\big(\rho^3+\rho^{\frac{13}{4}}+\rho^{\frac{7}{2}}\big).
\end{align*}
For the second line we use $|H|\lesssim A u_\rho$, $|D|\lesssim u_\rho$ and \eqref{typeii:eq:mixedendpoints}. The cutoff commutator and the omitted source tail are $O(A\rho^{\frac{3}{2}}R^{-\frac{5}{2}})$. Each normalized inner cell, including the added polarization, is bounded by $C\rho^{\frac{3}{2}}\langle y\rangle^{-2}$ on $|y|\le2R$. With $|P|\lesssim\langle y\rangle^{-3}$, any such insertion in the mixed quartic costs $O(A\rho^3)$ in $L^2_y$. The same estimate applies to the nonlinear part of the polarization by
\begin{equation*}
 (U+w_h)^5-U^5=5w_h\int_0^1(U+\theta w_h)^4\,d\theta .
\end{equation*}
For $|y|>\rho^{-1}$ the inner insertion is $O(\langle y\rangle^{-3})$,
$|D|\lesssim\rho^{\frac{1}{2}}$, and $\|H\|_6\lesssim A$. H\"older in
$L^3,L^\infty,L^6$ bounds the exterior potential difference times $H$
by $O(A\rho^{\frac{7}{2}})$. Dividing the normalized $L^2_y$ bound by $\lambda$
proves the momentum estimate.
\end{proof}

\section{The coupled wave inverse and dipole matching}\label{typeii:merge:coupled}
\subsection{Negative-mode columns and the corrected potential}\label{typeii:sec:pairinverse}
For the two-core estimate we set $L=-\Delta-V_i-V_o$ with
$V_i=5P_\lambda^4$ and $V_o=5Q_\mu^4$, omitting $V_i$ for the
single-core evolution. Let $\Pi_{\rm pos}$ be the $L^2$ projection onto
the annular position probes of Section~\ref{typeii:sec:pathbook}.
We choose a fixed smooth cutoff $\chi$ equal to one near zero and supported in $|x|<\frac{1}{32}$. For each unit negative eigenfunction $n_j$ we use the same packets throughout,
\begin{equation*}
 \eta_j=\chi\bigl(\frac{x}{\tau}\bigr)D_{s_j}n_j,\qquad
 \phi_j=(I-\Pi_{\rm pos})\eta_j .
\end{equation*}
Exponential decay, after rescaling, implies for $0\le k\le2$,
\begin{equation*}
 s_j^k\|\nabla^k(D_{s_j}n_j-\eta_j)\|_2+
 s_j^{k-1}\|\nabla^k(D_{s_j}n_j-\eta_j)\|_{\frac{6}{5}}
 \lesssim_k e^{-c\frac{\tau}{s_j}}.
\end{equation*}
The same bound holds after $\tau\partial_t$, since $|\dot s_j|\lesssim\frac{s_j}{\tau}$. On the position complement $Z\perp\phi_j$ implies $Z\perp\eta_j$, whence
\begin{equation*}
 s_j^{-1}|\langle Z,D_{s_j}n_j\rangle|
 \le s_j^{-1}\|Z\|_6\|D_{s_j}n_j-\eta_j\|_{\frac{6}{5}}
 \lesssim e^{-c\frac{\tau}{s_j}}\|\nabla Z\|_2.
\end{equation*}
Since $\frac{\tau}{s_j}\to\infty$, the weak limit at each scale is orthogonal to the full negative eigenspace, as required in the proof of \eqref{typeii:eq:coercive-gap}.

\paragraph{Dressed negative columns.}
In the signed class every smooth eigenfunction vanishes at the origin, and its fixed dilation derivatives obey the same bound $O(|x|)$ there. Its tail and derivatives decay exponentially. An outer negative field has overlap $O(\rho^{\frac{5}{2}})$ with the inner dual, while an inner field has exponentially small overlap with the outer dual up to fixed powers of $\rho^{-1}$. Fields on opposite scales have $L^2$ overlap $O(\rho^{\frac{5}{2}})$. Thus the Gram matrix of the dressed negative fields has a uniformly bounded inverse. For the defect $r_j=(L+\frac{\kappa_j^2}{s_j^2})\phi_j$, where $-\kappa_j^2$ is the isolated eigenvalue, we have
\begin{equation}\label{typeii:eq:negative-column}
 s_j\|r_j\|_2+\|r_j\|_{\frac{6}{5}}\lesssim 1.
\end{equation}
For the outer negative mode, in the coordinates $y=\frac{x}{\lambda}$ we have
\begin{equation*}
 |V_i\eta_o|\lesssim\lambda^{-2}\mu^{-\frac{3}{2}}\rho |y|P(y)^4,
 \qquad
 \mu\|V_i\eta_o\|_2\lesssim\mu^{-1}\rho^{\frac{1}{2}},\quad
 \|V_i\eta_o\|_{\frac{6}{5}}\lesssim\mu^{-1}\rho^{\frac{3}{2}} .
\end{equation*}
Applying $L$ to the inner dressing of size $O(\rho^{\frac{5}{2}})$
introduces a factor $\lambda^{-2}$ in $L^2$, so it obeys the same bounds.
For the inner negative mode we use exponential localization and
$|V_o(x)|\lesssim\mu^{-2}|\frac{x}{\mu}|^4$ near the origin.
Likewise
\begin{equation*}
 \|\partial_t^k\phi_j\|_2\lesssim\tau^{-k},\qquad
 \|\partial_t^k\phi_j\|_{\frac{6}{5}}\lesssim s_j\tau^{-k}
 \quad(0\le k\le1).
\end{equation*}
The $L^{\frac{6}{5}}$ bound also holds after one fixed spatial dilation derivative. The outer dressing is supported at the smaller scale and the opposite dressing has exponentially small coefficient. All these estimates hold for the fixed finite collections of fields.

\paragraph{Corrected potential.}
We separate the stationary exterior from the two concentrating potentials. With $V_i=5P_\lambda^4$ and $V_o=5Q_\mu^4$, we have
\begin{equation}\label{typeii:eq:corrected-potential}
 \|5U^4-V_i-V_o-5B^4\|_3\lesssim\tau^2,\qquad
 \|(5U^4-V_i-V_o)h\|_2\lesssim\|\nabla h\|_2 .
\end{equation}
Put $v=U-P_\lambda-Q_\mu$ and $r=|x|$. The middle correction is $O(\min(r,\mu))$ at the core by its signed zero and bounded tail. The inner cells are $O(\lambda^{-\frac{1}{2}}\rho^{\frac{3}{2}})=O(\tau^{10})$. At $r\asymp\tau$ the cut tails and the matched dipole are $O(\tau^4)$. Thus
\begin{equation*}
 |v|\lesssim\min(r,1)+\tau^{10},\qquad
 |v-B|\lesssim
 \begin{cases}
 \tau^4,& r\le\tau,\\
 \mu^{\frac{3}{2}}r^{-2}+\lambda^{\frac{5}{2}}r^{-3},&r>\tau .
 \end{cases}
\end{equation*}
In particular $\|v-B\|_{12}\lesssim\tau^{\frac{17}{4}}$, and $\|v^4-B^4\|_3\lesssim\tau^{\frac{17}{4}}$ by H\"older. For $F_s=P_\lambda,Q_\mu$, scaling and $F,\,rF^3\in L^3$ imply
\begin{equation*}
 \|F_s^3v\|_3+\|F_sv^3\|_3
 \lesssim s^{\frac{1}{2}}+\tau^{10}s^{-\frac{1}{2}}\lesssim\tau^2 .
\end{equation*}
For the two bare cross terms we use $|Q_\mu|\lesssim\mu^{-\frac{3}{2}}r$ and $|P_\lambda|\lesssim\lambda^{\frac{5}{2}}r^{-3}$. Hence
\begin{equation*}
 \|P_\lambda^3 Q_\mu\|_3\lesssim\lambda^{\frac{1}{2}}\mu^{-\frac{3}{2}}\lesssim\tau^2,\qquad
 \|Q_\mu^3P_\lambda\|_3\lesssim\lambda^{\frac{5}{2}}\mu^{-\frac{7}{2}}\lesssim\tau^{26}.
\end{equation*}
Every mixed quartic is bounded by a sum of these extreme products, by Young's inequality. This proves the first bound in \eqref{typeii:eq:corrected-potential}. The second follows from $B^4\in L^3$ and $\dot H^1\hookrightarrow L^6$. For one core we omit $P_\lambda$ throughout.
\subsection{The nonlinear estimate on inner time cells}\label{typeii:sec:cellnonlinearity}
We bound the quintic remainder in the forcing norm by applying Strichartz on the shortest active time scale. Partition each dyad $\tau\asymp\sigma$ into cells $I$ of nominal
length $\lambda_I\asymp\sigma^{16}$, retaining this length for
clipped end cells. For a pair of applied sources we set
\begin{equation*}
 M_p=\sup_I(\lambda_I\tau_I^p)^{-1}
        \int_I (\|\nabla f_1(t)\|_2+\|f_2(t)\|_2)\,dt .
\end{equation*}
The cell version of \eqref{typeii:eq:terminal-kernels} controls
the energy in units $\tau^q$, $q=p+1$, once the selected inverse is
known. We apply the estimate on each cell to the mixed position--velocity system.

Let $N=\|\nabla h\|_2+\|k\|_2$. We apply the pressure bounds \eqref{typeii:eq:two-pressures} below, obtained by differentiating the two constraints in \eqref{typeii:eq:mixed-wave}. Using $\|\nabla E_j\|_2\lesssim s_j^{-1}$, $\|E_j\|_2\lesssim1$, $\|D_j\|_2\lesssim\tau^{-1}$ and $|I|\lesssim \min_j s_j$, we obtain
\begin{equation*}
 \int_I\sum_j\left(\frac{|\ell_j|}{s_j}+|m_j|+\frac{|\ell_j|}{\tau}\right)dt
 \lesssim_J \sup_I N+\int_I(\|\nabla f_1\|_2+\|f_2\|_2)\,dt .
\end{equation*}
By scaling and the pointwise correction bounds, $\|5U^4\|_{L^1(I;L^3)}=5\|U\|_{L^4(I;L^{12})}^4\lesssim1$. We put $5U^4h$ on the right, using $\|h\|_6\lesssim N$. By \eqref{p5:eq:mixed-free-wave} and interpolation with $L^\infty(I;L^6)$,
\begin{equation}\label{typeii:eq:cell-strichartz}
 \|h\|_{L^4(I;L^{12})}+\|h\|_{L^5(I;L^{10})}
 \lesssim_J \sup_I N+\int_I(\|\nabla f_1\|_2+\|f_2\|_2)\,dt .
\end{equation}
In particular, for an energy field $v$ on the same cell, H\"older and Sobolev imply
\begin{equation}\label{typeii:eq:quartic-multiplier}
 \||h|^4v\|_{L^1(I;L^2)}
 \le\|h\|_{L^4(I;L^{12})}^4\|v\|_{L^\infty(I;L^6)}.
\end{equation}
The same estimates hold on a middle cell for the one-core equation, with $|I|\lesssim\mu$.

Interpolation gives $\|U\|_{L^5(I;L^{10})}^5\le\|U\|_{L^\infty(I;L^6)}\|U\|_{L^4(I;L^{12})}^4\lesssim1$. We write $N_U(h)=(U+h)^5-U^5-5U^4h$. We expand in $h$ and differentiate in $U$ to obtain
\begin{equation}\label{typeii:eq:quintic-taylor}
 \begin{aligned}
 N_U(h)&=20h^2\int_0^1(1-\theta)(U+\theta h)^3\,d\theta,\\
 \partial_U N_U(h)&=60h^2\int_0^1(1-\theta)(U+\theta h)^2\,d\theta .
 \end{aligned}
\end{equation}
We use the second identity in Section~\ref{typeii:sec:transport} to vary $U$. We use \eqref{p5:eq:quintic-pointwise} on each inner cell. For two waves $h,k$ we have
\begin{equation}\label{typeii:eq:quintic-difference}
 |N_U(h)-N_U(k)|\lesssim |h-k|\big[|U|^3(|h|+|k|)+|h|^4+|k|^4\big].
\end{equation}
Assume the energy and $L^4(I;L^{12})$ norms of $h,k$ are bounded by $R\tau_I^q$, $q>15$. In \eqref{typeii:eq:quintic-difference} the bracket has $L^1(I;L^3)$ norm $\lesssim_R\tau_I^q$ by \eqref{typeii:eq:quartic-multiplier} and H\"older in the four coefficient factors. At fixed time its $L^{\frac{3}{2}}$ norm has the same bound by Sobolev. We apply H\"older with $h-k$ in $L^6$ and obtain
\begin{equation}\label{typeii:eq:cell-remainder}
 \begin{aligned}
 \int_I\|N_U(h)-N_U(k)\|_2\,dt
 &\lesssim_R \tau_I^q\sup_I\|\nabla(h-k)\|_2,\\
 |\langle N_U(h)-N_U(k),E_j\rangle|
 &\lesssim_R s_j^{-1}\tau_I^q\|\nabla(h-k)\|_2 .
 \end{aligned}
\end{equation}
The second line uses $\|E_j\|_6\lesssim s_j^{-1}$. With $k=0$ the cell integral is $O_R(\tau_I^{2q})$. Since $s_j\gtrsim\lambda_I\asymp\tau_I^{16}$, the gain relative to both source weights is $\tau_I^{q-15}$.
\subsection{The coupled position and velocity equations}
We choose separate tests for position and velocity, then solve the resulting constrained evolution with mixed terminal conditions on the negative modes.

\subsubsection{Choice of tests}

In a unit frame write $e=\Lambda F$ and let $n_1,\ldots,n_l$ be the orthonormal negative eigenfunctions. We form the column map $D=(e/\|e\|_2,n_1,\ldots,n_l)$. For smooth radial annular cutoffs $0\le\chi_R\le1$ tending to one away from the origin, put
\begin{equation}\label{typeii:eq:compact-duals}
 D_R=\chi_R D(D^*\chi_R D)^{-1},\qquad
 D^*D_R=I,\qquad D_R\longrightarrow D\quad\text{in column }L^2 .
\end{equation}
Indeed, $D^*D=I$ and $\chi_R D\to D$. We take $j_F$ to be the first column of $D_R$ times $\|e\|_2$, choosing $R$ so that $\|j_F-e\|_2\le\varepsilon_1$. This annular test detects $e$ and annihilates every negative mode.

The displacement and velocity require different tests at the middle scale. In middle units let $e=\Lambda Q$, $L=\frac{\tau}{\mu}$ and $e_L=\chi(\frac{|z|}{L})e$. The leading tail is $r^{-2}Y_1(\omega)$ for a nonzero degree-one harmonic. We fix a smooth shell test $v_L=L^{-1}\psi(\frac{r}{L})Y_1(\omega)$ on a subannulus of the plateau of $\chi$. For a shell cutoff $\psi$ with nonzero limiting pairing we have
\begin{equation*}
 \|v_L\|_{\dot H^1}\lesssim L^{-\frac{1}{2}},\qquad
 |\langle v_L,e_L\rangle|\asymp1,\qquad \|e_L\|_{\frac{6}{5}}\lesssim L^{\frac{1}{2}}.
\end{equation*}
Testing on $v_L$ and applying Sobolev for the opposite bound, we obtain $\|e_L\|_{\dot H^{-1}}\asymp L^{\frac{1}{2}}$. The cone-supported fields
\begin{equation*}
 e_L-\frac{\|e_L\|_2^2}{\langle v_L,e_L\rangle}v_L
 \longrightarrow e\quad\text{in }\dot H^1
\end{equation*}
are exactly orthogonal to $e_L$. Their pairings with the negative eigenfunctions are exponentially small, since $e$ is orthogonal to these eigenfunctions and the modification lies at radii comparable to $L$ or larger. Subtract their pairings times the corresponding fixed columns of $D_R$ in \eqref{typeii:eq:compact-duals}. These columns are orthogonal to $e_L$ for large $L$, since $e_L=e$ on their support. The correction is exponentially small in $\dot H^1$. The limiting field $e$ has zero form and is nonzero, so the cone-supported displacement constraint has no uniform gap.
\paragraph{Norms for the physical dilation tests.}
For the mixed inverse we take $E_j$ to be the bare dilation truncated inside a fixed small fraction of the terminal cone, equal to the full dilation on the core. Then for $L=-\Delta-5P_\lambda^4-5Q_\mu^4$,
\begin{equation}\label{typeii:eq:velocity-column}
 \|LE_i\|_{\frac{6}{5}}\lesssim 1,\qquad \|LE_o\|_{\frac{6}{5}}\lesssim 1,\qquad
 \|\partial_t E_j\|_2\lesssim\tau^{-1}.
\end{equation}
The middle-potential contribution on the bare inner column is $O(\lambda^{-1}\rho^{\frac{5}{2}})$ or smaller, by splitting at $|x|=\mu$, the linear middle zero and the quadrupolar tail. The inner potential on the outer column contributes $O(\mu^{-1}\rho^{\frac{3}{2}})$. On $r\sim\tau$ the middle cutoff commutator contributes $O(\mu^{\frac{1}{2}}\tau^{-\frac{3}{2}})$ and the inner one $O(\lambda^{\frac{3}{2}}\tau^{-\frac{5}{2}})$ in $L^{\frac{6}{5}}$. A bounded $L^3$ change of potential preserves uniform boundedness of these norms by the $L^2$ bound on each column. These physical dual estimates suffice without differentiating an incident dipole in the velocity test.

\subsubsection{Solving the mixed wave equations}\label{typeii:sec:mixedinverse}
We solve the two constrained wave equations using the nested cutoff $c$ of Section~\ref{typeii:sec:terminalestimate}. With $s_i=\lambda$, $s_o=\mu$, let $E_j$ be the cone truncation of $T_j=D_{s_j}\Lambda F_j$. The Gram matrix $(\langle E_a,E_j\rangle)$ is uniformly invertible, and $\|E_j\|_2+\|\Lambda E_j\|_2\lesssim1$. By \eqref{typeii:eq:velocity-column} and \eqref{typeii:eq:corrected-potential}, the bounds $\|\dot E_j\|_2\lesssim\tau^{-1}$ and $\|LE_j\|_{\frac{6}{5}}\lesssim1$ hold also for the corrected operator. We consider
\begin{equation}\label{typeii:eq:mixed-wave}
\begin{aligned}
 h_t&=k+f_1+\sum_j\ell_jE_j, &\qquad \langle h,j_a\rangle&=0,\\
 k_t&=-Lh+f_2+\sum_j(m_jE_j+\ell_jD_j), & \langle k,E_a\rangle&=0 .
\end{aligned}
\end{equation}
Here $D_j$ may include a prescribed kinematic feedback satisfying
$\|D_j\|_2\lesssim\tau^{-1}$, in particular a bounded multiple of $(\frac{b_j}{s_j})D_{s_j}KF_j$. The two applied inputs are measured in $\dot H^1$ and $L^2$, respectively. Smooth finite-terminal equations are considered first.

\paragraph{Choice of position charts.}
By \eqref{typeii:eq:compact-duals}, the tests $j_F$ satisfy \eqref{typeii:eq:coercive-gap}, with a constant depending on the chosen cutoffs.
For $j_a=D_{s_a}j_{F_a}$ the omitted cone tails of $E_a$ tend to zero in $L^2$. Hence $\|j_a-E_a\|_2\le2\varepsilon_1$ on a late interval. The mixed Gram matrix $A_{aj}=\langle j_a,E_j\rangle$ is an $O(\varepsilon_1)$ perturbation of $(\langle E_a,E_j\rangle)$, so its inverse is bounded uniformly for small $\varepsilon_1$.
We differentiate the constraints in \eqref{typeii:eq:mixed-wave} once and obtain the triangular system
\begin{align*}
 \sum_j\langle j_a,E_j\rangle\ell_j
   &=-\langle j_a-E_a,k\rangle-\langle j_a,f_1\rangle-\langle\dot j_a,h\rangle,\\
 \sum_j\langle E_a,E_j\rangle m_j
   &=\langle LE_a,h\rangle-\langle E_a,f_2+\textstyle\sum_j\ell_jD_j\rangle-\langle\dot E_a,k\rangle .
\end{align*}
By \eqref{typeii:eq:probe-overlap}, both $A$ and $(s_j A_{aj}/s_a)_{aj}$ have bounded inverses. We use $A^{-1}$ for the term $\langle j_a-E_a,k\rangle$. For the other two terms we divide row $a$ by $s_a$ and use $\|j_a\|_{\frac{6}{5}}+\tau\|\dot j_a\|_{\frac{6}{5}}\lesssim_J s_a$. Inverting the second Gram matrix and using \eqref{typeii:eq:velocity-column}, we obtain
\begin{equation}\label{typeii:eq:two-pressures}
\begin{aligned}
 |\ell_j|&\le C\varepsilon_1\|k\|_2+C_J s_j(\tau^{-1}\|\nabla h\|_2+\|\nabla f_1\|_2),\\
 |m|&\lesssim\|f_2\|_2+\|\nabla h\|_2+\tau^{-1}(\|k\|_2+|\ell|).
\end{aligned}
\end{equation}
Both estimates involve only the applied forces and the instantaneous state.

\paragraph{Coercive energy with the two pressures.}
In the coercive form and in the virial identity we always take $L=-\Delta-\sum_j V_j$, as in Section~\ref{typeii:sec:terminalestimate}. For the physical operator $-\Delta-5U^4$, we include $(5U^4-\sum_j V_j)h$ in the second-row source and estimate it by \eqref{typeii:eq:corrected-potential}. For the continuous components $(Z,W)$ we set $N^2=\langle LZ,Z\rangle+\|W\|_2^2$. By \eqref{typeii:eq:coercive-gap}, $\|\nabla Z\|_2\lesssim_J N$, whereas $\|W\|_2\le N$ with constant one. In the three pairings of \eqref{typeii:eq:commonvirial} we use
\begin{equation*}
 \langle W,E_j\rangle=0,\quad
 |\langle c\Lambda Z,E_j\rangle|\lesssim_\chi\frac{\mu}{\tau}\|\nabla Z\|_2,\quad
 |\langle LZ,E_j\rangle|\lesssim\|\nabla Z\|_2 .
\end{equation*}
The second bound follows by local Sobolev on $r\lesssim_\chi\mu$, and the third from the column defect. Write $A_j=|p_j|/s_j+|b_j|$ for the negative-mode sizes in the decomposition below. Its columns satisfy $\|\nabla\phi_j\|_2\lesssim_J s_j^{-1}$ and $\|\psi_j\|_2\lesssim_J1$. Thus in \eqref{typeii:eq:two-pressures} the modal components of $h,k$ contribute $O_J(\sum_j A_j)$ to $|\ell|$ and $O_J(\tau^{-1}\sum_j A_j)$ to $|m|$. Since $\|c\Lambda E_j\|_2+\|D_j\|_2\lesssim\tau^{-1}$, the kinetic term obeys
\begin{equation*}
 \frac{C|\ell|}{\tau}\|W\|_2\le
 \frac{C\varepsilon_1}{\tau}\|W\|_2^2
 +C_J\frac{\mu}{\tau^2}N^2+C_J\frac{\mu}{\tau}N\|\nabla f_1\|_2+C_J\frac{N}{\tau}\sum_j A_j .
\end{equation*}
The remaining pressure terms sum to at most
\begin{equation*}
 C_{J,\chi}\Big[(1+\frac{\mu}{\tau^2})N^2+
                  N(\|\nabla f_1\|_2+\|f_2\|_2)+\frac{N}{\tau}\sum_j A_j\Big].
\end{equation*}
Upon integration, the modal term gives $\tau^{-1}\sum_j U_j$ in the first inequality of \eqref{typeii:eq:terminal-kernels}, with $U_j=A_j$. Since $A_j(v)\le D s_j(T-v)v^p$, its contribution to the weighted continuous bound is $O(\nu D)$ and is absorbed in \eqref{typeii:eq:weighted-bounds}. The coefficient of $\varepsilon_1\tau^{-1}\|W\|_2^2$ is independent of the position gap. We choose $\varepsilon_1$ first, then the virial tails and slopes as in Section~\ref{typeii:sec:terminalestimate}, and finally $\tau_0$.

\paragraph{Negative modes.}
We use the packets $\eta_j,\phi_j$ of Section~\ref{typeii:sec:pairinverse}. Let $\Pi_{\rm vel}$ be the $L^2$ projection onto the span of the velocity tests $E_a$ and set $\psi_j=(I-\Pi_{\rm vel})\eta_j$. Since $\psi_j-\phi_j=(\Pi_{\rm pos}-\Pi_{\rm vel})\eta_j$, we have
\begin{equation*}
\|\psi_j-\phi_j\|_2\lesssim_J\rho^{\frac{3}{2}},\qquad
 \|\nabla(\psi_j-\phi_j)\|_2\lesssim_J\frac{\rho^{\frac{3}{2}}}{\lambda}\lesssim_J1,\qquad
 \|\dot\psi_j\|_2\lesssim_J\tau^{-1}.
\end{equation*}
For the first two estimates, the same-frame overlaps with the physical dilation are zero before cutoff. A middle negative field has overlap $O(\rho^{\frac{3}{2}})$ with the true inner dilation by its $r^{-3}$ Kelvin tail and the linear zero of the negative field at the origin. The opposite overlap is $O(\rho^{\frac{5}{2}})$ by exponential decay and the linear middle zero. The bound for $\psi_j-\phi_j$ also implies uniform invertibility of the velocity Gram matrix. All columns have gradient norm $O(\lambda^{-1})$. For the contribution of the middle velocity column to the inner mode, we use
\begin{equation*}
 \rho^{\frac{3}{2}}\|E_o\|_{\frac{6}{5}}\lesssim\rho^{\frac{3}{2}}(\mu\tau)^{\frac{1}{2}}=o(\lambda).
\end{equation*}
Thus $\|\psi_j-\phi_j\|_2\lesssim_J s_i$ and $\|\psi_j\|_{\frac{6}{5}}\lesssim_J s_j$. The latter bound also holds after one spatial dilation derivative.
We write $h=Z+\sum p_j\phi_j$, $k=W+\sum b_j\psi_j$, with continuous components perpendicular to the respective mode columns. We set $A_j=\frac{|p_j|}{s_j}+|b_j|$, $G_j=p_j'-b_j$ and $H_j=b_j'-\frac{\kappa_j^2 p_j}{s_j^2}$. The first continuous equation has the additional source
\begin{equation*}
 -\sum G_j\phi_j+\sum b_j(\psi_j-\phi_j)-\sum p_j\dot\phi_j ,
\end{equation*}
and the second has
\begin{equation*}
 -\sum H_j\psi_j-\sum b_j\dot\psi_j-\sum p_j r_j+
   \sum (\frac{\kappa_j^2p_j}{s_j^2})(\phi_j-\psi_j),\qquad
 r_j=(L+\frac{\kappa_j^2}{s_j^2})\phi_j .
\end{equation*}
The non-multiplier pieces are bounded by $C_J\tau^{-1}\sum A_j$ in $\dot H^1\times L^2$. We differentiate $\langle Z,\phi_j\rangle=\langle W,\psi_j\rangle=0$ once. Inverting the resulting Gram matrices, we obtain
\begin{equation*}
 |G_j|\lesssim_J s_j\{\|\nabla f_1\|_2+\tau^{-1}(N+\textstyle\sum A_l)\},\qquad
 |H_j|\lesssim_J\|f_2\|_2+(1+\frac{\mu}{\tau})\|\nabla f_1\|_2+
                  C_J\tau^{-1}(N+\textstyle\sum A_l).
\end{equation*}
For $G_j$, we use $W\perp\psi_j$, $\langle E_a,\psi_j\rangle=0$, and the $L^2$ difference bound, and test $\dot\phi_j$ in $L^{\frac{6}{5}}$ against $Z$. The terms $p_l\langle\dot\phi_l,\phi_j\rangle$ with $s_l>s_j$ have the additional overlap $O_J(\rho^{\frac{5}{2}})$. The Gram off-diagonal entries across scales are also $O_J(\rho^{\frac{5}{2}})$, so inversion preserves the indicated relative weights. For $H_j$, $m$ vanishes by the velocity-mode constraint. Testing $LZ$ on $\psi_j$ reduces to $r_j$ and the bounded gradient of $\psi_j-\phi_j$. In the virial inequality, the $G_j$-term is bounded by $C_J(1+\tau^{-1})|G_j|N$ by the position-mode defect and $\|c\Lambda\phi_j\|_2\lesssim\tau^{-1}$. The $H_j$-term is bounded by $C_J(\frac{s_j}{\tau})|H_j|N$ by $W\perp\psi_j$ and integration by parts in the virial pairing against $\psi_j$. All leading $\frac{N^2}{\tau}$ terms from the pressures still have coefficient $C\varepsilon_1$. The $A_j$ terms require only a bounded coefficient.

In the variable $d\sigma_j=\frac{dt}{s_j}$, the combinations $b_j\pm\frac{\kappa_jp_j}{s_j}$ satisfy
\begin{equation}\label{typeii:eq:mode-clock}
 \left(\partial_{\sigma_j}\mp\kappa_j\right)
 \left(b_j\pm\kappa_j\frac{p_j}{s_j}\right)
 =s_jH_j\pm\kappa_j\left(G_j-\dot s_j\frac{p_j}{s_j}\right).
\end{equation}
We set the minus combination to zero at the left endpoint and the plus combination to zero at the right endpoint. Integration of \eqref{typeii:eq:mode-clock} and the continuous virial bound yield \eqref{typeii:eq:terminal-kernels} with $U_j=A_j$ and $f_c,f_j\lesssim\|\nabla f_1\|_2+\|f_2\|_2$. For zero terminal continuous data, \eqref{typeii:eq:weighted-bounds} implies
\begin{equation}\label{typeii:eq:mixed-wave-bound}
 \|\nabla h\|_2+\|k\|_2\lesssim_p M\tau^{p+1},\qquad
 A_j\lesssim_p M s_j\tau^p,\qquad p>-1 ,
\end{equation}
whenever the applied forces satisfy \eqref{typeii:eq:clock-mass} with constant $M$.
For a potential perturbation $V_{\rm pert}$ bounded in $L^\infty_tL^3_x$, Sobolev's inequality gives $\|V_{\rm pert}h\|_2\lesssim\|V_{\rm pert}\|_3\|\nabla h\|_2$. At energy weight $\tau^{p+1}$ this source has weight $\tau^p$ with a factor $O(\tau_0)$, so it is absorbed in the same estimate. This applies to the physical potential by \eqref{typeii:eq:corrected-potential}.

\paragraph{Finite-dimensional terminal selection.}
At a positive cutoff, the two constraint matrices determine $\ell,m$, and the free Duhamel formula defines the evolution in both time directions. With zero terminal continuous data, $h=\sum_jp_j\phi_j$ and $k=\sum_jb_j\psi_j$ at the right endpoint. The right mode conditions impose $b_j=-\kappa_j p_j/s_j$. It remains to invert the square homogeneous map
\begin{equation*}
 (p_j)_{\rm right}\longmapsto
       (b_j-\kappa_jp_j/s_j)_{\rm left}.
\end{equation*}
Its kernel is zero by \eqref{typeii:eq:mixed-wave-bound}. Linearity therefore enforces the left conditions for every forcing, with the same uniform estimate. This argument applies to one or two cores and to the physical columns whenever their Gram matrices are invertible.

For Lipschitz paths we use the reference rates in \eqref{typeii:eq:commonvirial}. Since $\int_0^{\tau_0}\tau^{d-1}d\tau=\frac{\tau_0^d}{d}$, the extra energy coefficient is integrable. Both constraint systems and \eqref{typeii:eq:mode-clock} use only the first derivatives of the tests.

Since $\psi_j=\eta_j-\Pi_{\rm vel}\eta_j$ is a linear combination of $\eta_j$ and the tests $E_a$, it has the same cone support. Hence for an applied acceleration $f_{\rm ext}$ supported in $r>\frac{\tau}{16}$,
\begin{equation}\label{typeii:eq:exterior-orthogonality}
 \operatorname{supp}\psi_j\ \cup\ \operatorname{supp}E_a\subset \{r\le\frac{\tau}{32}\},
 \qquad \langle f_{\rm ext},E_a\rangle=\langle f_{\rm ext},\psi_j\rangle=0 .
\end{equation}
A physical coefficient column may extend to $\frac{\tau}{8}$.
Its extra outer dilation tail on that annulus has $L^2$ norm
$O((\frac{\mu}{\tau})^{\frac{1}{2}})$ and gradient norm $O(\frac{\big(\frac{\mu}{\tau}\big)^{\frac{1}{2}}}{\tau})$.
More generally a first-row column discrepancy costing $O(1)$ in
$\dot H^1$ and $o(1)$ in $L^2$, and a second-row one costing $o(1)$
in $L^2$, contribute bounded errors or $\frac{o(1)}{\tau}$ in the energy
inequality. The corresponding mixed Gram matrices retain invertibility.
For negative modal forcing an $O(1)$ first-row $\dot H^1$ error is
tested on $\phi_j$ with dual norm $O(s_j)$.

For the application with an already selected middle wave $h_o$, the
additional potential is $5[(U+h_o)^4-U^4]$. Its action on the smaller
normal is estimated at the actual energy and cell weights in
Section~\ref{typeii:sec:transport}, including its instantaneous
pairings on the velocity tests.

\paragraph{Fully coupled true columns and the instantaneous velocity output.}
For a phase $(U,V)$ let the first-row coefficient of the acceleration
pressure be $-U_{p_j}$, in addition to the already handled position-rate
column $-U_{q_j}$. Let the second-row columns be $-V_{p_j},-V_{q_j}$.
We take the kinetic tests supported strictly inside the exterior nonlinear
cutoff as above. The two constraint systems together have matrix
\begin{equation*}
 \begin{pmatrix}
  A_{aj}=\langle j_a,-U_{q_j}\rangle&
  B_{aj}=\langle j_a,-U_{p_j}\rangle\\
  C_{aj}=\langle E_a,-V_{q_j}\rangle&
  D_{aj}=\langle E_a,-V_{p_j}\rangle
 \end{pmatrix}.
\end{equation*}
The diagonal blocks $A,D$ are perturbations of the invertible Gram matrices above. In the physical velocity row the apparent term of size $\frac{|p_j|}{s_j}$ cancels on the same core. Write $e=\Lambda F_j$, $\Lambda_1=\frac{3}{2}+x\cdot\nabla$, $E_j=\chi_v D_{s_j}e$ and $K_j^0=\zeta D_{s_j}\Lambda_1 e$, with $\zeta=1$ on the support of $\chi_v$. Integration by parts yields
\begin{equation}\label{typeii:eq:tested-dilation}
 \langle E_j,K_j^0\rangle
 =-\frac{1}{2}\int (x\cdot\nabla\chi_v)|D_{s_j}e|^2\,dx,\qquad
 |\langle E_o,K_o^0\rangle|\lesssim\frac{\mu}{\tau},\quad
 |\langle E_i,K_i^0\rangle|\lesssim\left(\frac{\lambda}{\tau}\right)^3 .
\end{equation}
The derivative of $\chi_v$ is supported at $r\asymp\tau$, where the unit-scale profiles have tails of orders $-2$ and $-3$, respectively. On opposite cores $|\langle E_a,K_j^0\rangle|\lesssim\rho^{\frac{3}{2}}$. In middle coordinates the contribution is bounded by $\rho^{\frac{3}{2}}\int |z|^{-3}\min(|z|,|z|^{-2})\,dz$. With $\frac{|p_j|}{s_j}\lesssim\tau^{-1}$, the entire tested principal block is $O(\tau^2)$.

For the physical columns the remainders in \eqref{typeii:eq:middle-physical-columns} and \eqref{typeii:eq:physical-columns} below are $O(\tau^{\frac{1}{2}})$ in the velocity row. In the position row $\|\nabla U_{p_j}\|_2\lesssim\tau^2$. Thus
\begin{equation}\label{typeii:eq:physical-gram}
 |C_{aj}|\lesssim\tau^{\frac{1}{2}},\qquad
 |B_{aj}|\lesssim_J s_a\tau^2 .
\end{equation}
The bounds for $A^{-1}$, $D^{-1}$ and the inverse of $\big(\frac{s_j A_{aj}}{s_a}\big)_{aj}$ used here follow from \eqref{typeii:eq:probe-overlap}. Indeed the position error tested against $j_a$ is $O_J(s_a\tau)$, and the velocity error in $D$ is $O(\tau^3)$. For right sides $g_1,g_2$ we substitute $\ell=A^{-1}(g_1-Bm)$ in $Dm=g_2-C\ell$. Since $\|CA^{-1}B\|\lesssim_J\mu\tau^{\frac{5}{2}}$, absorption gives
\begin{equation*}
 |m|\lesssim_J |g_2|+\tau^{\frac{1}{2}}|g_1|,\qquad
 |\ell_j-(A^{-1}g_1)_j|\lesssim_J s_j\tau^2|m| .
\end{equation*}
The same position-column bound costs $O(\tau^2|m|)$ in the energy estimate and $O_J(s_j\tau^2|m|)$ on the $j$th negative position test.

For an inner correction of size $O(\tau^q)$ on the bounded energy background $U+h_o$, \eqref{typeii:eq:cell-remainder} bounds the nonlinear pairing by $Cs_j^{-1}\tau^{2q}$. Relative to the required velocity forcing $\tau^{q-1}$, we have
\begin{equation*}
 s_j^{-1}\tau^{q+1}\lesssim\tau^{q-15}\quad(q>15),\qquad
 s_o^{-1}\tau^{7}\lesssim\tau^3\quad\text{for the middle wave at }q=6 .
\end{equation*}
The velocity tests annihilate the exterior force by support and the negative-mode accelerations by $\langle E_a,\psi_j\rangle=0$. These identities, the tested nonlinear bound and the $\dot H^{-1}$ column defect yield $|m|\lesssim\tau^{q-1}$ pointwise. On each $\lambda$-cell the exterior input is bounded by its dyadic $L^1_tL^2_x$ norm $O(\tau^{\frac{33}{2}})$, which suffices for the cell Strichartz estimate at $q\le\frac{33}{2}$. Its direct modal projections vanish by \eqref{typeii:eq:exterior-orthogonality}.

\subsection{The retarded dipole on the reference path}\label{typeii:sec:matchedcone}
We match the dipole of the elliptic core to a free wave on the light cone.
The difference in their common region determines both the matching error and the regularity needed at the origin.
Let $\mu=a\tau^4$, $r=|x|$, and write the dipole harmonic of $Q$ as
$Y(\omega)=b\cdot\omega$. For $F(s)=a^{\frac{3}{2}}(s_+)^6$, we set
\begin{equation*}
 C=-b\cdot\nabla_x\frac{F(\tau-r)}{r},\qquad
 D=\frac{F(\tau)}{r^2}Y,\qquad A=-\frac{1}{2}F''(\tau)Y .
\end{equation*}
Off the origin $\Box C=0$, since $\frac{F(\tau-r)}{r}$ is a radial free wave. It is the retarded dipole for increasing $\tau$. In physical time $t=T-\tau$, its characteristics $\tau-r=\mathrm{constant}$ travel inward.
Its support lies in $r\le\tau$. For $\mu=a\tau^4$, \eqref{typeii:eq:explicit-dipole-tail} identifies the
constant term in the elliptic tail with that of the free dipole,
\begin{equation*}
 \mu Z_\infty=-15a^{\frac{3}{2}}\tau^4Y=-\frac{F''(\tau)Y}{2}=A,
 \qquad C=\frac{Y}{r^2}[F(\tau-r)+rF'(\tau-r)] .
\end{equation*}
The derivative of $F(\tau-r)+rF'(\tau-r)$ in $r$ is $-rF''(\tau-r)$. Integrating twice, we obtain for $0<r<\tau$
\begin{equation}\label{typeii:eq:retarded-remainder}
 H=C-D-A=\frac{b\cdot x}{2}\int_0^1(1-v^2)F'''(\tau-rv)\,dv .
\end{equation}
This identity holds for every $C^3$ amplitude. Here $F'''(s)=120a^{\frac{3}{2}}s^3$ for $s>0$. Differentiating under the integral and using $\Delta A=D_{tt}$ and $\Box C=\Delta D=0$, we obtain
\begin{equation*}
 |\nabla_x^j\partial_\tau^kH|\lesssim_{j,k}\tau^{3-k}r^{1-j},
 \qquad \Box H=-A_{tt}\quad(0<r<\tau).
\end{equation*}
We take $\eta=\eta\bigl(\frac{r}{\tau}\bigr)$ equal to one up to $\frac{\tau}{16}$ and zero beyond
$\frac{\tau}{8}$, and $\chi_\mu=\chi\bigl(\frac{r}{\mu}\bigr)$ equal to one up to $\mu$ and zero
beyond $2\mu$. We define
\begin{equation*}
 V=B+(1-\eta)C+\eta\{Q_\mu+\mu Z\bigl(\frac{x}{\mu}\bigr)+(1-\chi_\mu)H\}.
\end{equation*}
Thus $V=B+Q_\mu+\mu Z\bigl(\frac{x}{\mu}\bigr)$ on $r\le\mu$ and $V=B$ on $r>\tau$.
In particular the deterministic cone modification vanishes throughout the
support $r\le2\sqrt{\lambda\mu}$ of the inner elliptic cells for small
$\tau$. We use the derivative at zero in \eqref{typeii:eq:average-slope} for the selected middle wave.

\paragraph{Structured source and the weight-six central estimate.}
We claim, for $0\le j\le3$,
\begin{equation}\label{typeii:eq:matched-middle-residual}
 \|(\mu\nabla_x)^j\mathcal R(V)(t)\|_2\lesssim_j\tau^5,\qquad
 |\mu^{\frac{5}{2}}\mathcal R(V)(t,\mu z)|\lesssim\tau^{10}|z|
 \quad(|z|\le1).
\end{equation}
On $r<\frac{\tau}{16}$ we use the elliptic cancellation
\begin{equation*}
 (Q_\mu)_{tt}-\Delta(\mu Z)-5Q_\mu^4(\nabla B(0)\cdot x+\mu Z)=0 .
\end{equation*}
The inertial terms left after adding $(1-\chi_\mu)H$ are
\begin{equation*}
 (\mu Z)_{tt}-(1-\chi_\mu)A_{tt}
                  +[\Box,1-\chi_\mu]H .
\end{equation*}
On $2\mu\le r<\frac{\tau}{16}$ we apply \eqref{typeii:eq:explicit-dipole-tail} to the first two terms,
\begin{equation*}
 |(\mu Z-A)_{tt}|\lesssim\frac{\tau^2\mu}{r},\qquad
 \|(\mu Z-A)_{tt}\|_{L^2(2\mu\le r<\frac{\tau}{16})}
 \lesssim\tau^2\mu\tau^{\frac{1}{2}}\lesssim\tau^{\frac{13}{2}} .
\end{equation*}
On $r<2\mu$ the same terms are $O(\tau^2)$ in magnitude. The
commutator is supported on $r\asymp\mu$, where each spatial derivative
of the cutoff introduces a factor of $O(\mu^{-1})$ and each time derivative a factor of
$O(\tau^{-1})$. By \eqref{typeii:eq:retarded-remainder},
\begin{equation*}
 \|(\mu\nabla_x)^j[\Box,1-\chi_\mu]H\|_2
        \lesssim_j \tau^3\mu^{\frac{1}{2}}=O(\tau^5).
\end{equation*}

We set $G=B+\mu Z+(1-\chi_\mu)H$. In the interior $|G|\lesssim r$ and $|G-B|\lesssim\mu$. By scaling $Q_\mu$ we have $\||x|^kQ_\mu^4\|_2\lesssim\mu^{k-\frac{1}{2}}$ for $k=1,2$. The two linear errors and the extreme mixed powers therefore satisfy
\begin{align*}
 \|Q_\mu^4(B-\nabla B(0)\cdot x)\|_2&\lesssim\mu^{\frac{3}{2}},&
 \|Q_\mu^4(1-\chi_\mu)H\|_2&\lesssim\tau^3\mu^{\frac{1}{2}},\\
 \|Q_\mu^3G^2\|_2&\lesssim\mu^2,&
 \|Q_\mu G^4\|_2&\lesssim\mu^{\frac{3}{2}}\tau^{\frac{7}{2}}.
\end{align*}
The last row bounds every intermediate mixed power by interpolation, and $\|G^5-B^5\|_2\lesssim\mu\tau^{\frac{11}{2}}$. On $r\ge\mu$ each spatial derivative costs $O(r^{-1})$, so these bounds persist under $(\mu\nabla_x)^j$.

For $|z|\le1$ the elliptic core is uncut. In normalized $C^j$ the first linear error is $O_j(\mu^{\frac{5}{2}})$, the cell inertia is $O_j(\mu^{\frac{5}{2}}\tau^2)$, and the higher powers are smaller. The residual vanishes at $z=0$ by symmetry. Its $C^1$ bound thus supplies the factor $|z|$ in \eqref{typeii:eq:matched-middle-residual}.

On the transition $r\asymp\tau$ the matched expression is $B+C+\eta(Q_\mu-D+\mu Z-A)$. By the dipole and elliptic tails,
\begin{equation*}
 |Q_\mu-D+\mu Z-A|\lesssim\mu^{\frac{5}{2}}\tau^{-3}+\frac{\mu^2}{\tau}
 \lesssim\tau^7 .
\end{equation*}
Two space-time derivatives cost $O(\tau^{-2})$, and the square root of the annulus volume is $O(\tau^{\frac{3}{2}})$. The resulting $L^2$ norm is therefore $O(\tau^{\frac{13}{2}})$.
The nonlinear difference
from the stationary exterior is $O(\tau^8)$, supported in a volume
$O(\tau^3)$ since the added field is $O(\tau^4)$ on this annulus and up
to $r=\tau$. The same estimates hold after $(\mu\nabla_x)^j$:
each annular derivative introduces a factor of $O_j(\tau^{-1})$ and the retarded field
vanishes to fifth order at the outer characteristic surface.

Let $h$ solve the cone-supported outer linear problem with potential $5V^4$.
Assume that its energy is $O(\tau^6)$ and that its terminal packets and pressures
satisfy the selection bounds at this weight. For each $0<\delta<\frac{1}{2}$, we claim
\begin{equation}\label{typeii:eq:central-middle-jet}
 \|\mu^{\frac{1}{2}}h(t,\mu z)\|_{C^{2,\delta}(|z|\le2)}
 +\|\mu^{\frac{3}{2}}h_t(t,\mu z)\|_{C^{1,\delta}(|z|\le2)}\lesssim_\delta \tau^6 .
\end{equation}
To apply the local wave estimate, we extend $C$ across $r=\tau$ using the polynomial
$F(s)=a^{\frac{3}{2}}s^6$. We multiply the extension by a smooth radial cutoff supported in
$\frac{\tau}{16}<r<2\tau$ and equal to one for $\frac{\tau}{8}\le r\le\tau$. Call the result $C_e$.
On the support of $h$, we can write
\begin{equation}\label{typeii:eq:split-cone-potential}
 5V^4=5(B+C_e)^4+W_c,\qquad
 W_c=\mathbf1_{\{r\le\frac{\tau}{8}\}}\big[5V^4-5(B+C_e)^4\big].
\end{equation}
The two fields match smoothly at $r=\frac{\tau}{8}$. On a slab $\tau\asymp\sigma$,
$|\nabla_{t,x}^j C_e|\lesssim_j\sigma^{4-j}$ for every fixed $j$, by the polynomial extension and the fixed cutoff in $r/\tau$. In the variables $\theta=\tau/\sigma$ and $X=x/\sigma$, the continued dipole thus has $C^M$ norm $O_M(\sigma^4)$ and $B(\sigma X)$ has $C^M(B_3)$ norm $O_M(\sigma)$ because $B(0)=0$. Consequently the macroscopic potential $5\sigma^2(B+C_e)^4$ has $C^M([1/8,1]\times B_3)$ norm $O_M(\sigma^6)$, including its time derivatives. We set $m=\mu(\sigma)$ and $L=\frac{\sigma}{m}$.
For $|z|\le L$, the normalized shrinking terms are
$O(\langle z\rangle^{-2})$, since $m^{\frac{3}{2}}\lesssim L^{-2}$.
The normalized smooth background is $O(L^{-2}\langle z\rangle)$.
We expand the fourth powers in \eqref{typeii:eq:split-cone-potential} and
use \eqref{typeii:eq:matched-middle-residual} for the source. Thus
\begin{equation*}
 |\nabla_z^j(m^2 W_c)|\lesssim_j\langle z\rangle^{-5-j},\qquad
 \|m^{\frac{5}{2}}f(v,m\cdot)\|_{H^3}\lesssim m\sigma^5\lesssim\frac{\sigma^6}{L}.
\end{equation*}
On the middle unit ball the potential agrees with that of the uncut elliptic core.
These are precisely the source and potential hypotheses of the local
estimate in Section~\ref{typeii:sec:outermoment}, with weight $\sigma^6$.
The resulting $H^4\times H^3$ bound on unit middle balls implies
\eqref{typeii:eq:central-middle-jet} by Sobolev embedding. The signed quadrupole has the
same weight by \eqref{typeii:eq:quadrupolereadout}.

\subsection{Varying scales and a fixed dipole correction}\label{typeii:sec:perturbedmatched}
We keep the dipole correction fixed as the scale and its velocity vary.
This reduces the new residual to a parameter difference and a mixed quintic term.
We set $q_*=a\tau^4$, $b_*=-4a\tau^3$, $q\asymp q_*$,
$b\asymp -\tau^3$, $A(q)=\alpha_o q^{\frac{1}{2}}$, and
$D_0=\partial_t^{\rm exp}+b\partial_q+A(q)\partial_b$.
We assume
\begin{equation*}
 \left|\frac{q}{q_*}-1\right|+\left|\frac{b}{b_*}-1\right|\le D\tau^{\frac{5}{2}}.
\end{equation*}
We use the explicit velocity correction $Z_K(z)=\frac{|z|^2\Lambda Q(z)}{4}$
from Section~\ref{typeii:sec:stationarycells}. It satisfies $L_Q Z_K=-KQ$ and
$\nabla Z_K(0)=0$. By the dipole decay of $Q$, we have
$|\nabla^j Z_K(z)|\lesssim_j\langle z\rangle^{-j}$ to each fixed order.
We write $R=q^{-\frac{1}{2}}$, $L=\frac{\tau}{q}$, $z=\frac{x}{q}$, and use the same physical cone
cutoff $\eta$ as in the matched reference. The middle displacement before dipole matching is
\begin{equation*}
 U_0=B+\eta\,q^{-\frac{1}{2}}
 \{Q(z)+\chi\bigl(\frac{z}{R}\bigr)[q^{\frac{3}{2}}Z(z)+(b^2-\beta_oq^{\frac{3}{2}})Z_K(z)]\},
 \qquad V_0=D_0U_0 .
\end{equation*}
The coefficient $b^2-\beta_o q^{\frac{3}{2}}$ is a first integral of the reference scale equation,
\begin{equation}\label{typeii:eq:middle-first-integral}
 D_0(b^2-\beta_o q^{\frac{3}{2}})
 =bq^{\frac{1}{2}}\left(2\alpha_o-\frac{3}{2}\beta_o\right)=0 .
\end{equation}
Thus $D_0$ differentiates the profile and its cutoff with this coefficient held fixed. The uncut angular coefficient is the second derivative of the dipole strength. Indeed, for $f(q)=q^{\frac{3}{2}}$, $D_0q=b$ and $D_0b=A(q)=\alpha_o q^{\frac{1}{2}}$, so
\begin{equation*}
 D_0^2 f=\frac{3}{4} q^{-\frac{1}{2}}b^2+\frac{3}{2} q^{\frac{1}{2}}A(q).
\end{equation*}
Using $Z_{K,\infty}=-\frac{3}{8}Y$ and $Z_\infty=-\frac{3}{8}(\beta_o+2\alpha_o)Y$, we obtain for every trial pair
\begin{equation*}
 qZ_\infty+q^{-\frac{1}{2}}(b^2-\beta_o q^{\frac{3}{2}})Z_{K,\infty}
       =-\frac{1}{2}(D_0^2 f)Y .
\end{equation*}

\paragraph{One residual estimate for the three spatial commutations.}
We write $\mathcal D=q_*\nabla_x$ and $R_0=D_0 V_0-\Delta U_0-U_0^5$.
The leading terms in middle units are
$-qA(q)\Lambda Q+b^2KQ-5q^{\frac{3}{2}}Q^4H_B$ and the two elliptic cells
cancel them. The errors from their cutoff, two $D_0$ derivatives, and
the cone cutoff are bounded in normalized $L^2$ by
\begin{equation*}
 q^{\frac{3}{2}}R^{-\frac{1}{2}},\qquad
 q^{\frac{3}{2}}(\frac{q}{\tau})^2 R^{\frac{3}{2}},\qquad q^{\frac{3}{2}} L^{-\frac{1}{2}}.
\end{equation*}
The Taylor and higher-power errors are $O(q^{\frac{5}{2}})$ in middle units by the interior product bounds of Section~\ref{typeii:sec:matchedcone}. These bounds also hold after
$q\partial_q$ or $b\partial_b$ at fixed $(t,x)$ and at most three
applications of $\mathcal D$. Indeed each logarithmic derivative
differentiates a coefficient or applies an Euler derivative to a
profile or cutoff. In $\frac{|x|}{q}\le2$, Taylor's integral remainder and
the smooth cells retain their powers of $q$ under normalized spatial
derivatives. On every larger annulus each such derivative introduces a factor of at most
$\frac{C_jq}{|x|}$ in the envelope, including on moving cutoff annuli.
For $q\asymp q_*$, $-b\asymp\tau^3$, the three normalized bounds above are $O(\tau^7)$. Dividing by $q$ proves the common estimate
\begin{equation}\label{typeii:eq:middle-symbol}
 \max_{j\le3}\sum_{X\in\{1,q\partial_q,b\partial_b\}}\|\mathcal D^j X R_0\|_2\lesssim\tau^3 .
\end{equation}
At the reference parameters this bounds $F_o$ in \eqref{typeii:eq:reference-residual}. Integration along the segment in $(\log q,\log(-b))$ yields $\|\mathcal D^j(R_0-R_{0,*})\|_2\lesssim D\tau^{\frac{11}{2}}$ under the stated relative error bound.

Let $V_*$ denote the matched reference position and put
$c=V_*-U_0(t,q_*,b_*)$, $(U,V)=(U_0+c,V_0+c_t)$.
We write $U_{0,*}=U_0(t,q_*,b_*)$. The fixed terms $c_{tt}-\Delta c$
cancel in the difference of residuals, and we obtain the exact identity
\begin{equation}\label{typeii:eq:frozen-dipole-difference}
 \begin{aligned}
 D_0V-\Delta U-U^5-\mathcal R(V_*)
  ={}&R_0-R_{0,*}\\
     &-\big[(U_0+c)^5-U_0^5-(U_{0,*}+c)^5+U_{0,*}^5\big].
 \end{aligned}
\end{equation}
We normalize the positions
by $q_*^{\frac{1}{2}}$, put $r=\frac{|x|}{q_*}$, $L_*=\frac{\tau}{q_*}$, and write
$\delta\widehat U_0=\widehat U_0-\widehat U_{0,*}$.
The Kelvin and cone cutoffs imply, on $r\lesssim L_*$,
\begin{align*}
 |\delta\widehat U_0|&\lesssim D\tau^{\frac{5}{2}}\langle r\rangle^{-2},&
 |\widehat U_0|+|\widehat U_{0,*}|+|\widehat c|
   &\lesssim\langle r\rangle^{-2}+\tau^6\langle r\rangle,\\
 |\widehat c|&\lesssim\min\{\tau^9\langle r\rangle^{\frac{3}{2}},\tau^6\}.
\end{align*}
For the last bound, $c$ vanishes on $r\le1$.
Up to the elliptic cutoff at $r\asymp q_*^{-\frac{1}{2}}$, \eqref{typeii:eq:retarded-remainder}
bounds $\widehat c$ by $O(\tau^9r)$. Farther out, $\widehat c=O(\tau^6)$. These envelopes also bound the first three normalized spatial derivatives, by smoothness on the core and differentiation on the annuli. For the mixed quintic in \eqref{typeii:eq:frozen-dipole-difference} we use the identity
\begin{equation*}
 (u+c)^5-u^5-(v+c)^5+v^5
 =20c(u-v)\int_0^1\!\int_0^1
        \big(v+\theta(u-v)+\eta c\big)^3\,d\theta\,d\eta .
\end{equation*}
The product rule preserves one factor from each of the envelopes for $\widehat c$ and $\delta\widehat U_0$. We use $|\widehat c|\lesssim\tau^9\langle r\rangle^{\frac{3}{2}}$ with $\langle r\rangle^{-6}$, and $|\widehat c|\lesssim\tau^6$ with $\tau^{18}\langle r\rangle^3$. For each of the four derivative orders the physical $L^2$ norm is at most
\begin{equation*}
 Cq_*^{-1}D\tau^{\frac{5}{2}}
      (\tau^9+\tau^{24}L_*^{\frac{5}{2}})=O(D\tau^{\frac{15}{2}}).
\end{equation*}
Combining \eqref{typeii:eq:frozen-dipole-difference} with \eqref{typeii:eq:matched-middle-residual}, we obtain
\begin{equation}\label{typeii:eq:varying-middle-residual}
 \max_{0\le j\le3}\|\mathcal D^j(D_0V-\Delta U-U^5)\|_2
       \lesssim \tau^5+D\tau^{\frac{11}{2}}.
\end{equation}
Thus the forcing has normalized $H^3$ norm $O(\tau^9)$.

\paragraph{Transport columns and the two outer gauges.}
We fix a trial path $(q,b)$, with $b$ the independent velocity parameter.
We use the outer path ball of Section~\ref{typeii:sec:closing}, in which the relative errors in $q,b$ and
their first logarithmic derivatives are $O(D\tau^{\frac{5}{2}})$.
In the middle projected evolution, differentiating the displacement and
momentum constraints determines two scalar outputs $a_o,c_o$. The
physical columns $-(U_q,V_q)$ and $-(U_b,V_b)$ carry these outputs.
In middle coordinates we write $e=\eta\Lambda Q$ and $\Lambda_1=\frac{3}{2}+z\cdot\nabla$, where $\eta$ is the physical cone cutoff. The normalized sources are
\begin{equation*}
 G=-q^{\frac{3}{2}}(a_oU_q+c_oU_b),\qquad
 F=-q^{\frac{5}{2}}(a_oV_q+c_oV_b).
\end{equation*}
Here the subscript $b$ denotes differentiation in the velocity parameter.
The estimates
\begin{equation}\label{typeii:eq:middle-physical-columns}
 \begin{aligned}
 \|\nabla(q^{\frac{3}{2}}U_q+e)\|_2&\lesssim\tau^5,&
 \|\nabla(q^{\frac{1}{2}}U_b)\|_2&\lesssim\tau^2,\\
 \|q^{\frac{5}{2}}V_q-b\Lambda_1e\|_2&\lesssim\tau^{\frac{9}{2}},&
 \|q^{\frac{3}{2}}V_b+e\|_2&\lesssim\tau^3
 \end{aligned}
\end{equation}
follow by differentiating only $U_0,V_0$. For the bounded cells,
$\|\nabla(\chi_RZ)\|_2\lesssim R^{\frac{1}{2}}$ and
$\|\chi_RZ\|_2\lesssim R^{\frac{3}{2}}$, including fixed logarithmic
derivatives. In the order $(q,b)$, the two rows of \eqref{typeii:eq:middle-physical-columns} have errors bounded componentwise by
\begin{equation*}
 \left(q^{\frac{3}{2}}R^{\frac{1}{2}},\,|b|R^{\frac{1}{2}}\right),\qquad
 \left(\tau^9R^{\frac{3}{2}}+\frac{q}{\tau}L^{-\frac{1}{2}},\,\tau^6R^{\frac{3}{2}}\right)
\end{equation*}
up to a fixed constant. Derivatives of $\eta$ are supported on its annulus. The leading terms of $(G,F)$ are
\begin{equation*}
 \big(a_oe,\,-ba_o\Lambda_1e+q c_oe\big).
\end{equation*}
The cutoff commutator in $L_Q e$ has dual energy norm $O(L^{-\frac{3}{2}})$ by annular rescaling in $L^{\frac{6}{5}}$. In physical units
its $L^2$ and $L^{\frac{6}{5}}$ norms are $O(q^{\frac{1}{2}}\tau^{-\frac{5}{2}})$ and
$O(q^{\frac{1}{2}}\tau^{-\frac{3}{2}})$.

We use the one-core choice of position and velocity tests from
Section~\ref{typeii:sec:mixedinverse}. The displacement probe is compact in
a fixed middle annulus, orthogonal to the negative fields and arbitrarily
close in $L^2$ to the dilation. The velocity test $e_v$ is the dilation
cut off between physical radii $\frac{\tau}{64}$ and $\frac{\tau}{32}$, hence
$\|e-e_v\|_2\lesssim L^{-\frac{1}{2}}$ in middle units.
In physical units the four column errors just estimated have norms
$O(\tau),O(\tau^2),O(\tau^{\frac{1}{2}}),O(\tau^3)$, respectively.
They are the allowed errors in the mixed inverse with only the middle
index. Its common virial bound and finite-dimensional selection
therefore control the energy on each prescribed path.

For the shell estimate we also need the two individual pressures.
We write $N=\|(\nabla X,Y)\|_2$ in middle units and
$F_*=q\|D_0V-\Delta U-U^5\|_2$, omitting further applied sources.
On the fixed displacement annulus $q^{\frac{1}{2}}U_b=O(|b|)$, so its
acceleration column has pairing $O(\tau^3)$. The differentiated
velocity test contributes $O(|q'|+\frac{q}{\tau})=O(\tau^3)$, and its spatial
operator defect contributes $O(L^{-\frac{3}{2}}+q)$ against the displacement.
Both diagonal pairings tend to nonzero constants. Solving the two
rows yields
\begin{equation}\label{typeii:eq:middle-pressures}
 |a_o|\lesssim N+\tau^3 F_*,\qquad
 |q c_o|\lesssim\tau^3N+F_* .
\end{equation}
Only the first derivatives of the trial path enter these gauges.
The defects $q'-b,b'-A(q)$ enter the reconstruction in
Section~\ref{typeii:sec:transport}.

\paragraph{Shell sources and the full-cone exterior.}
If \eqref{typeii:eq:middle-pressures} yields $|a_o|\lesssim\tau^d$ and
$|q c_o|\lesssim\tau^{3+d}$, the pressure sources on the middle slab, through any fixed spatial
order $j$, satisfy
\begin{equation*}
 |\partial_z^jG|\lesssim_j\tau^d\langle z\rangle^{-2-j},\qquad
 |\partial_z^jF|\lesssim_j\tau^{3+d}\langle z\rangle^{-2-j}.
\end{equation*}
On $r\lesssim R$ the bounded angular terms in $G,F$ have coefficients $O(\tau^{d+6})$, $O(\tau^{d+9})$, respectively. The remaining profiles and their cone cutoffs obey the differentiated dipole bounds. Thus \eqref{typeii:eq:sphere-shell} applies to both free source terms, and \eqref{typeii:eq:slab-potential-energy} controls the contribution of the exterior potential in the same Duhamel formula. The local $H^n$ bound on $G$ also controls the difference between momentum and the time derivative of displacement.

Both variable trial fields equal their monomial values for
$r>\frac{\tau}{8}$. Hence on the support of the cone solution the same smooth
exterior $5(B+C_e)^4$ in \eqref{typeii:eq:split-cone-potential} applies to every path. Its macroscopic derivative bounds are therefore uniform even for Lipschitz trial paths. In the smaller cone the
variable cells obey the same differentiated dipole envelopes as the
reference terms. Thus the concentrated potential still satisfies
$|\nabla_z^j(m^2W_c)|\lesssim_j\langle z\rangle^{-5-j}$.
By \eqref{typeii:eq:varying-middle-residual}, we retain the global $H^3_z$ source bound
$O(\frac{\sigma^6}{L})$ as on the reference path.

\subsection{The selected middle wave and its exterior force}
\label{typeii:sec:nonlinearmiddle}
We solve the middle equation with its quintic remainder cut off near the core.
We estimate the remaining exterior term in the norm of the two-core inverse.
We write $m=\mu\asymp\tau^4$, $L=\frac{\tau}{m}$, and $A(\tau)=C_0\tau^e$, with
$e>3$. Suppose the prescribed middle trial and its structured linear
residual satisfy the energy, shell and packet hypotheses of
Section~\ref{typeii:sec:outermoment} at weight $A$. The source and potential
calculations in Section~\ref{typeii:sec:perturbedmatched} verify these hypotheses at $e=6$.
We take the velocity tests supported in $|x|\le\frac{\tau}{32}$ as above. Let
$\chi$ equal one for $|x|\le\frac{\tau}{16}$ and vanish for $|x|\ge\frac{\tau}{8}$.
In the prescribed middle equation we include the source
$\chi N(V,h)$, where $N(V,h)=(V+h)^5-V^5-5V^4h$.
Then there is a selected solution at the same energy and local regularity
weights, on a sufficiently short terminal interval. The function
$(1-\chi)N(V,h)$ is left as an applied source for the two-core problem.

We use the sum of the weighted energy and local $H^4$ norms of
Section~\ref{typeii:sec:outermoment}, taking the supremum over unit balls centered
in $|x|<\frac{\tau}{4}$ in the middle frame. On a terminal slab we set
$H(z)=m^{\frac{1}{2}}h(mz)$ and $D(z)=m^{\frac{1}{2}}V(mz)$, with the time variable
suppressed. On the support of the inner cutoff all fixed spatial
derivatives of $D$ are bounded by $C\langle z\rangle^{-1}$. Indeed, the normalized regular field is $O(L^{-2}\langle z\rangle)$ and the
shrinking profile tails are $O(\langle z\rangle^{-2})$. On $|z|\lesssim L$ we have $A\langle z\rangle\lesssim LA$, where $LA\lesssim\tau^{e-3}$. We take $\tau_0$ small enough that $LA\le1$. On each unit ball meeting $\operatorname{supp}\chi$, the $H^3$ norms of $D,H,\widetilde H$ are $O(\langle z_0\rangle^{-1})$, and those of the two wave inputs are also $O(A)$. We use the exact first variation
\begin{equation*}
 \partial_H N(D,H)=20H\int_0^1(D+\theta H)^3\,d\theta .
\end{equation*}
The $H^3$ algebra estimate on the ball bounds the derivative by $CA\langle z_0\rangle^{-3}$. Integrating between the inputs and summing squares over a bounded-overlap unit cover, we obtain
\begin{equation}\label{typeii:eq:cutoff-quintic-h3}
 \|\chi[N(D,H)-N(D,\widetilde H)]\|_{H^3_z}
 \lesssim A\sup_{z_0}\|H-\widetilde H\|_{H^3(B_1(z_0))},
\end{equation}
since $\langle z\rangle^{-3}\in L^2$. Taking $\widetilde H=0$ bounds the normalized source by $CA^2$. The linear estimates require a normalized $H^3$ source of size $A/L$ and a physical $L^2$ source of size $A/\tau$. The latter norm costs $m^{-1}$ under rescaling. Since $mL=\tau$, both ratios are $O(LA)$, and \eqref{typeii:eq:cutoff-quintic-h3} supplies the same Lipschitz factor for the weighted differences. Thus the linear estimate of Section~\ref{typeii:sec:outermoment}, including its pressure and terminal terms, yields a contraction for $LA\lesssim\tau_0^{e-3}$ sufficiently small.

\paragraph{The nonlinear source near the cone.}
On a dyad $I$ with $\sigma\le\tau\le2\sigma$, we choose a fixed radial cutoff $\eta$ equal to zero for $r<\frac{\sigma}{64}$ and one for $r>\frac{\sigma}{16}$. In the physical equation for $(\eta h,\eta k)$ the only cutoff terms are $-2\nabla\eta\cdot\nabla h-(\Delta\eta)h$ in the second row. Their $L^2$ norm is $O(A/\sigma)$ by Hardy's inequality. The data at either end have energy $O(A)$ by the same inequality. By cone support and the middle dipole bound, $|V|\lesssim\sigma$ on $\operatorname{supp}(\eta h)\subset\{r\asymp\sigma\}$, so $\|\eta V^4h\|_2\lesssim\sigma^5 A$. The structured residual and $\chi N(V,h)$ already have $L^2$ norm $O(A/\sigma)$.

For the two pressure sources we use $z=x/m$, $L=\sigma/m$ and the shell bounds following \eqref{typeii:eq:middle-pressures}. In middle units $|\nabla_z^jG|\lesssim A\langle z\rangle^{-2-j}$ for $j=0,1$ and $|F|\lesssim A L^{-1}\langle z\rangle^{-2}$. On $|z|\gtrsim L$, rescaling to physical units yields
\begin{equation*}
 m^{-1}\big(\|\nabla_z[\eta(mz)G]\|_2+\|\eta(mz)F\|_2\big)
 \lesssim \frac{A}{m L^{\frac{3}{2}}}\lesssim\frac{A}{\sigma}.
\end{equation*}
We apply \eqref{p5:eq:mixed-free-wave} on $I$. For the second norm we use $\|w\|_{10}^5\le\|w\|_6\|w\|_{12}^4$ with $w=\eta h$. The energy bound therefore implies
\begin{equation}\label{typeii:eq:exterior-strichartz}
 \|\eta h\|_{L^4(I;L^{12})}+\|\eta h\|_{L^5(I;L^{10})}\lesssim A .
\end{equation}
On the support of $1-\eta$, the local regularity bound implies $|h|\lesssim A m^{-\frac{1}{2}}$. Interpolation with $\|h\|_6\lesssim A$ therefore yields $\|(1-\eta)h\|_{12}^4\lesssim A^4/m$ at each time. In particular both critical norms of $h$ are $O(A)$ on every middle cell, using interpolation for $(5,10)$.

The support of $(1-\chi)N(V,h)$ lies in $\frac{\tau}{16}\le r\le\tau$, where $\eta=1$. Taylor's formula, spatial H\"older and \eqref{typeii:eq:exterior-strichartz} imply
\begin{equation}\label{typeii:eq:exterior-dyad}
 \|(1-\chi)N(V,h)\|_{L^1(I;L^2)}
 \lesssim \sigma\,\sigma^3|B_{2\sigma}|^{\frac{1}{6}}\sup_I\|h\|_6^2
              +\|\eta h\|_{L^5(I;L^{10})}^5
 \lesssim\sigma^{\frac{9}{2}}A^2+A^5 .
\end{equation}
Since $A=C_0\sigma^e$ up to a fixed factor and $e>3$, the first term dominates for small $\sigma$. At $e=6$ the two powers are $\frac{33}{2}$ and $30$.

\paragraph{The wave forced by the exterior source.}
Suppose the exterior source has $L^1_tL^2_x$ norm $O(r^{15+d})$ on each dyad $v\asymp r$. In \eqref{typeii:eq:clock-mass} take $p=14+d$, $f_c=\|f_{\rm ext}\|_2$ and $f_j=0$ by \eqref{typeii:eq:exterior-orthogonality}. The simultaneous estimate \eqref{typeii:eq:weighted-bounds} now yields
\begin{equation*}
 n(\tau)\lesssim\tau^{15+d},\qquad U_j(\tau)\lesssim s_j(T-\tau)\tau^{14+d}.
\end{equation*}
The second bound includes the response induced by the continuous field.

\section{Coupling the scales and reaching the singular time}\label{typeii:merge:closing}
\subsection{Transport of the two phases}
\label{typeii:sec:transport}
We derive the scale equations by adding the two constrained wave equations. With $q=(\mu,\lambda)$ and independent velocities $p=(p_o,p_i)$, we use
\begin{equation*}
 A=(\alpha_o\mu^{\frac{1}{2}},\,\alpha_i\lambda^{\frac{3}{2}}\mu^{-\frac{5}{2}}),\qquad
 D_0=\partial_t^{\rm exp}+p\cdot\partial_q+A\cdot\partial_p,\qquad
 \mathcal F(u,v)=(v,\Delta u+u^5).
\end{equation*}
With $z=\frac{x}{\mu}$, $R=\mu^{-\frac{1}{2}}$ and $E_j^0=\zeta D_{q_j}\Lambda F_j$, we write the matched phase as
\begin{align*}
 U_0&=B+\zeta\mu^{-\frac{1}{2}}\{Q(z)+\chi\bigl(\frac{z}{R}\bigr)[\mu^{\frac{3}{2}}Z(z)
                 +(p_o^2-\beta_o\mu^{\frac{3}{2}})Z_K(z)]\},\\
 V_0&=D_0U_0,\qquad
 c(t,x)=U_*(t,x)-U_0(t,a\tau^4,-4a\tau^3),\\
 \Pi_b&=(U_0+c,V_0+c_t),\qquad
 \Xi=(J,-p_iE_i^0),\qquad \Pi=\Pi_b+\Xi .
\end{align*}
Here $U_*$ is the fixed retarded-dipole trial on the monomial path and $\zeta$ is the physical cone cutoff. We take the inner position of Section~\ref{typeii:sec:stationarycells} for $J$, with cutoff $\zeta$ and incident slope
\begin{equation*}
 d=\mu^{\frac{3}{2}}\partial_3\Pi_b^u(t,0)+d_h .
\end{equation*}
Thus $J$ includes both $\kappa\rho^{\frac{5}{2}}V_\times$ and $p_i^2 S_P$ in inner units. We hold $d_h$ fixed when differentiating in $q,p$. In $D_0\Xi$ we differentiate $d_h$ by \eqref{typeii:eq:slope-test}, which uses only $q_o'$ and the first middle equation.

We write $H_j=(h_j,k_j)$. With $f_{\rm ext}=(1-\chi)N(\Pi_b^u,h_o)$, the middle equation is
\begin{equation}\label{typeii:eq:prescribed-middle-wave}
 \partial_tH_o=\mathcal F(\Pi_b+H_o)-(0,f_{\rm ext})-D_0\Pi_b-
                     (\Pi_b)_{q_o}a_o-(\Pi_b)_{p_o}c_o .
\end{equation}
Once we select $H_o$, we use $f_{\rm ext}$ as a prescribed source in the inner equation,
\begin{equation}\label{typeii:eq:prescribed-inner-wave}
\begin{aligned}
 \partial_tH_i={}&\mathcal F(\Pi+H_o+H_i)-\mathcal F(\Pi_b+H_o)
                   +(0,f_{\rm ext})-D_0\Xi\\
 &{}-\Xi_{q_o}a_o-\Xi_{p_o}c_o
    -\partial_q\Pi\,a^i-\partial_p\Pi\,c^i .
\end{aligned}
\end{equation}
The vector outputs $a^i=(a_o^i,a_i^i)$ and $c^i=(c_o^i,c_i^i)$ enforce the two-core constraints. Along a trial path the chain rule reads
\begin{equation*}
 \partial_t\Pi-D_0\Pi=\partial_q\Pi\,(q'-p)+\partial_p\Pi\,(p'-A).
\end{equation*}
We add \eqref{typeii:eq:prescribed-middle-wave} and \eqref{typeii:eq:prescribed-inner-wave}. The exterior forces cancel, and for $X=\Pi+H_o+H_i$ we obtain
\begin{equation}\label{typeii:eq:sumtransport}
 \partial_tX-\mathcal F(X)
 =\partial_q\Pi\,[q'-p-(a_o,0)-a^i]
  +\partial_p\Pi\,[p'-A-(c_o,0)-c^i].
\end{equation}
We therefore impose
\begin{equation*}
 q'=p+(a_o,0)+a^i,\qquad p'=A+(c_o,0)+c^i .
\end{equation*}

\paragraph{The inherited columns.}
At fixed $d_h$, the only possible $p_o$ dependence of $J$ is through the incident slope. Since $\nabla Z_K(0)=\nabla S_Q(0)=0$, we have
\begin{equation*}
 \partial_{p_o}\Xi=0,\qquad \partial_{q_o}\Xi^v=0,\qquad
 \|\partial_{q_o}J\|_{\dot H^1}\lesssim\frac{\rho^{\frac{3}{2}}}{\mu}.
\end{equation*}
For the last estimate we differentiate $\rho^{\frac{3}{2}}d\chi W$ and its cutoff radius $\rho^{-\frac{1}{2}}$ in inner coordinates. The cutoff $r^{-2}$ and $r^{-1}$ cells and their Euler derivatives have bounded $\dot H^1$ norms. If $|a_o|\lesssim\tau^{3+\delta}$ and the middle jets have size $O(\tau^e)$, then \eqref{typeii:eq:slope-rate} and \eqref{typeii:eq:transfer-scale} imply
\begin{equation*}
 \|a_o\partial_{q_o}J\|_{\dot H^1}\lesssim\tau^{17+\delta},\qquad
 \|\dot d_h\,\partial_{d_h}J\|_{\dot H^1}
 \lesssim\rho^{\frac{3}{2}}|\dot d_h|\lesssim\tau^{14+e}.
\end{equation*}
\paragraph{The physical transport columns.}
We write $(U,V)=\Pi$ and $K_j^0=\zeta D_{q_j}KF_j$. For the columns $-\partial_{q_j}\Pi,-\partial_{p_j}\Pi$, the remainders satisfy
\begin{equation}\label{typeii:eq:physical-columns}
\begin{array}{c|cc|cc}
 &-U_{q_j}-E_j^0 & -U_{p_j}
 & -V_{q_j}+(\frac{p_j}{q_j})K_j^0 & -V_{p_j}-E_j^0
 \\[-1mm]
 j& \dot H^1&\dot H^1&L^2&L^2\\ \hline
 o&O(\tau)&O(\tau^2)&O(\tau^{\frac{1}{2}})&O(\tau^3)\\
 i&O(\tau^2)&O(\tau^{15})&0&0 .
\end{array}
\end{equation}
For $j=o$ we rescale \eqref{typeii:eq:middle-physical-columns} and add $\partial_{q_o}J$ to the position column. In particular, the velocity error and the derivative of the cone cutoff satisfy
\begin{equation*}
 \mu^{-1}\tau^{\frac{9}{2}}\lesssim\tau^{\frac{1}{2}},\qquad
 \frac{|p_o|}{\mu}\left(\frac{\tau}{\mu}\right)^{-\frac{1}{2}}
 \lesssim\tau^{\frac{1}{2}} .
\end{equation*}
The middle cutoff cells and their Euler derivatives have $\dot H^1$ norms $O(R^{\frac{1}{2}})$, whereas the inner norms are $O(1)$. The frozen $c,c_t$ have zero parameter derivatives. For $j=i$ we obtain
\begin{align*}
 \|U_{q_i}+E_i^0\|_{\dot H^1}&\lesssim\frac{\rho^{\frac{3}{2}}}{\lambda}
 \lesssim\tau^2,& \|U_{p_i}\|_{\dot H^1}&\lesssim |p_i|\lesssim\tau^{15},\\
 V_{q_i}&=\frac{p_i}{q_i}K_i^0,& V_{p_i}&=-E_i^0 .
\end{align*}
The second line follows by differentiating $-p_iE_i^0$ using \eqref{typeii:eq:unitary-dilation}. We test velocity against the cone-truncated bare dilations $E_j$. Eliminating the position rates from the velocity row, we obtain
\begin{equation*}
 |a|+\tau|c_{\rm vel}|
 \lesssim \|\nabla h\|_2+\|k\|_2+\mu\|f_1\|_{\dot H^1}+\tau\max_j|\langle f_2,E_j\rangle| .
\end{equation*}
For the column errors in the table, $\|E_j\|_2\lesssim1$ and $\mu\asymp\tau^4$. Thus
\begin{equation}\label{typeii:eq:rate-absorption}
\begin{aligned}
 \mu\|\delta f_1\|_{\dot H^1}
    +\tau\max_j|\langle\delta f_2,E_j\rangle|
 &\lesssim \mu\tau(|a|+\tau|c_{\rm vel}|)
         +\tau^{\frac{3}{2}}|a|+\tau^4|c_{\rm vel}|\\
 &\lesssim\tau^{\frac{3}{2}}(|a|+\tau|c_{\rm vel}|).
\end{aligned}
\end{equation}
This bounds the energy-source errors from the physical columns. For the tested rates, \eqref{typeii:eq:physical-gram} already proves uniform invertibility, including the principal column $-\frac{p_j}{q_j}K_j^0$.

\paragraph{Bounds for the mixed scale vector field.}
Suppose the middle energy, the uniform local norms of Section~\ref{typeii:sec:nonlinearmiddle}, and the two middle-cell norms in \eqref{typeii:eq:cell-strichartz} are $O(\tau^6)$, and that the inner energy and inner-cell Strichartz norms are $O(\tau^{16})$. We differentiate the constraints to obtain
\begin{equation}\label{typeii:eq:scale-rate-bounds}
 |a_o|\lesssim\tau^6,\qquad |c_o|\lesssim\tau^5,\qquad
 |a^i|\lesssim\tau^{16},\qquad |c^i|\lesssim\tau^{15}.
\end{equation}
For the outer rows we apply \eqref{typeii:eq:middle-pressures}. The compact position test contributes a factor $\mu$ on the first-row source. In the velocity equation we use
\begin{equation*}
 \|\dot e_o\|_2\lesssim\tau^{-1},\qquad
 \|Le_o\|_{\dot H^{-1}}\lesssim1,\qquad
 \frac{\|\nabla h_o\|_2^2}{\mu}\lesssim\frac{\tau^{12}}{\mu}\lesssim\tau^8 .
\end{equation*}
The applied acceleration is $O(\tau^5)$. For the inner wave the smooth first-row source is $O(\tau^{17})$ in $\dot H^1$. We absorb the column errors by \eqref{typeii:eq:rate-absorption}. The moving velocity tests, the velocity feedback and the quintic projection are controlled by
\begin{equation*}
 \tau^{-1}\big(\|\nabla h_i\|_2+\|k_i\|_2\big)\lesssim\tau^{15},\qquad
 q_j^{-1}\|\nabla h_i\|_2^2\lesssim\lambda^{-1}\tau^{32}\lesssim\tau^{16},
\end{equation*}
where the last estimate uses the uniform $L^6$ bound on the reference field.
By \eqref{typeii:eq:exterior-dyad} the exterior source at $e=6$ has dyadic $L^1_tL^2_x$ norm $O(\tau^{\frac{33}{2}})$ and its direct pairings vanish by \eqref{typeii:eq:exterior-orthogonality}.

On the support of the inner position $J$, $|x|\le\frac{\tau}{8}$. The uniform local middle-frame bound of Section~\ref{typeii:sec:nonlinearmiddle} gives $|h_o|\lesssim\tau^6\mu^{-\frac{1}{2}}$ throughout this support. For $|x|\le\mu$, we also use $h_o(t,0)=0$ and its first derivative bound. With $y=\frac{x}{\lambda}$ and $u_\rho$ as in \eqref{typeii:eq:mixedendpoints}, we obtain
\begin{equation}\label{typeii:eq:common-inner-envelope}
 \lambda^{\frac{1}{2}}|h_o(t,\lambda y)|\lesssim\tau^6u_\rho(|y|),\qquad
 \lambda^{\frac{1}{2}}|U_b(t,\lambda y)|\lesssim u_\rho(|y|),\qquad
 \lambda^{\frac{1}{2}}|J(t,\lambda y)|\lesssim\langle y\rangle^{-3}.
\end{equation}
For the second estimate use the signed zero and the middle profile bounds. For the last one, the corrections to $P$ are bounded by $C\rho^{\frac{3}{2}}\langle y\rangle^{-2}+C\rho^{\frac{5}{2}}\langle y\rangle^{-1}$ and vanish beyond $2\rho^{-\frac{1}{2}}$ in inner units. By \eqref{typeii:eq:quintic-taylor} the difference $N_{U_b+J}(h_o)-N_{U_b}(h_o)$ contains a factor $Jh_o^2$. It vanishes off the support of $J$. Taylor's formula bounds the remaining factor by $C[(|U_b|+|J|)^2+|h_o|^2]$. The proof of \eqref{typeii:eq:mixedendpoints} applies to the envelope $\langle y\rangle^{-3}$ itself. Thus
\begin{equation}\label{typeii:eq:inner-quintic-force}
 \begin{aligned}
 \|N_{U_b+J}(h_o)-N_{U_b}(h_o)\|_2
 &\lesssim \frac{\tau^{12}}{\lambda}
 \left(\|\langle y\rangle^{-9}u_\rho^2\|_2+
       \|\langle y\rangle^{-3}u_\rho^4\|_2\right)\\
 &\lesssim\frac{\tau^{12}}{\lambda}(\rho^3+\rho^{\frac{7}{2}})
 \lesssim\tau^{32}.
 \end{aligned}
\end{equation}
Since $\|E_j\|_2\lesssim1$, the same pointwise-in-time estimate controls each tested force.

For the inner increment $v$, the exact polynomial identity
\begin{equation*}
 N_U(h_o+v)-N_U(h_o)=5[(U+h_o)^4-U^4]v+N_{U+h_o}(v)
\end{equation*}
separates the change of potential from the inner remainder. Since $U=U_b+J$, we have $\|U_b\|_{18}^3\lesssim\mu^{-1}$ by scaling the bounded middle profile and its tails (the smooth background has bounded $L^{18}$ norm). Using \eqref{typeii:eq:common-inner-envelope} with $u_\rho\le\rho^{\frac{3}{2}}|y|$, we have
\begin{equation*}
 \||U|^3h_o\|_3\lesssim
 \|U_b\|_{18}^3\|h_o\|_6+
 \frac{\tau^6\rho^{\frac{3}{2}}}{\lambda}\||y|\langle y\rangle^{-9}\|_3
 \lesssim\tau^2 .
\end{equation*}
For any subinterval $I$ of a dyad $\tau\asymp\sigma$, take $\eta$ as in \eqref{typeii:eq:exterior-strichartz}. That estimate bounds the exterior fourth power in $L^1(I;L^3)$ by $C\sigma^{24}$. The interior estimate $\|(1-\eta)h_o\|_{12}^4\lesssim\sigma^{24}/\mu\lesssim\sigma^{20}$ holds pointwise in time. Since $|(U+h_o)^4-U^4|\lesssim |U|^3|h_o|+|h_o|^4$, Sobolev's inequality implies
\begin{equation}\label{typeii:eq:potential-multiplier}
 \|[(U+h_o)^4-U^4]v\|_{L^1(I;L^2)}
 \lesssim (|I|\sigma^2+\sigma^{24})\sup_I\|\nabla v\|_2 .
\end{equation}
For the instantaneous pairing the cubic coefficient has the same $L^3$ bound. On $\operatorname{supp}E_j\subset B_{\frac{\tau}{32}}$ we have $|h_o|\lesssim\tau^4$ and $\|v\|_{L^2(B_\tau)}\lesssim \tau\|\nabla v\|_2$. Thus the quartic term in the pairing costs at most $C\tau^{17}\|\nabla v\|_2$. These estimates apply in particular at $\|\nabla h_i\|_2\lesssim\tau^{16}$.

\paragraph{Scale boundary conditions at a positive cutoff.}
Fix $0<\epsilon<\tau_0$ and set $s=-\log\tau$, $S=-\log\tau_0$, $S_\epsilon=-\log\epsilon$. Write $Z_j=(u_j,v_j)$ with
\begin{equation*}
 q_j=c_j^0\tau^{k_j}(1+u_j),\qquad p_j=-k_j c_j^0\tau^{k_j-1}(1+v_j),
 \qquad (c_o^0,c_i^0)=(a,b),\quad(k_o,k_i)=(4,16).
\end{equation*}
Write $q_j^*=c_j^0\tau^{k_j}$ and let primes in this paragraph denote $\tau$-derivatives. We subtract from each force its derivative in its own scale at $q^*$. With $(\kappa_o,\kappa_i)=(6,360)$, the two remaining sources are
\begin{equation}\label{typeii:eq:normalized-rates}
 \begin{aligned}
 g&=(a_o+a_o^i,a_i^i),\\
 f&=(c_o+c_o^i,c_i^i)+
 A(q)-A(q^*)-\frac{1}{\tau^2}\big(\kappa_j(q_j-q_j^*)\big)_{j=o,i}.
 \end{aligned}
\end{equation}
The rates are evaluated at $T-\tau$. Since $(q_j^*)''=A_j(q^*)$, the output displacements $x_j=\widehat q_j-q_j^*$ and velocities $z_j=-\widehat p_j-(q_j^*)'$ must satisfy
\begin{equation}\label{typeii:eq:scale-linear-pair}
 x_j'=z_j-g_j,\qquad z_j'=\frac{\kappa_j}{\tau^2}x_j+f_j .
\end{equation}
At a fixed point $x_j=q_j-q_j^*$, so \eqref{typeii:eq:normalized-rates} restores the full force in both physical momentum equations.

We prescribe $x_o(\epsilon)=z_o(\epsilon)=0$, $x_i(\tau_0)=0$ and $\epsilon z_i(\epsilon)=17x_i(\epsilon)$. Let $(m_o,m_i)=(7,17)$. For each scale we put $(\kappa,m)=(\kappa_j,m_j)$, omit the index $j$ and let $r_-<r_+$ be the roots of $r(r-1)=\kappa$. By \eqref{typeii:eq:scale-linear-pair},
\begin{equation}\label{typeii:eq:scale-modes}
 P_\pm=\tau^{-m}(\tau z-r_\mp x),\qquad
 \tau P_\pm'=(r_\pm-m)P_\pm+\tau^{2-m}f+r_\mp\tau^{1-m}g .
\end{equation}
Assume $|g|+\tau|f|\lesssim\tau^{m-1}$. At the middle scale both coefficients $3-7$ and $-2-7$ are negative and $P_\pm(\epsilon)=0$. Integration in $\log\tau$ from $\epsilon$ bounds both modes uniformly.

At the inner scale, $19<r_+<20$ and $r_-=1-r_+$. The endpoint conditions become
\begin{equation*}
 P_+(\tau_0)=P_-(\tau_0),\qquad
 P_-(\epsilon)=-\frac{r_+-17}{17-r_-}P_+(\epsilon),\qquad
 0<\frac{r_+-17}{17-r_-}<\frac{3}{35}.
\end{equation*}
We integrate $P_+$ from $\tau_0$ and $P_-$ from $\epsilon$. The exponential kernels in $\log\tau$ have integrals at most $(r_+-17)^{-1}$ and $(17-r_-)^{-1}$. Thus
\begin{equation*}
 \|P_+\|_\infty\le\|P_-\|_\infty+C,\qquad
 \|P_-\|_\infty\le\frac{3}{35}\|P_+\|_\infty+C .
\end{equation*}
Absorption bounds both modes independently of $\epsilon$. With zero forcing the same inequalities have $C=0$, so the homogeneous two-point problem has trivial kernel. Its two boundary conditions therefore determine a unique solution for every forcing. At both scales we recover $x,z$ from \eqref{typeii:eq:scale-modes}, and their derivatives from \eqref{typeii:eq:scale-linear-pair}. Hence
\begin{equation}\label{typeii:eq:scale-displacement-bound}
 |x_j|+\tau(|x_j'|+|z_j|)+\tau^2|z_j'|\lesssim\tau^{m_j}.
\end{equation}

Set $(d_o,d_i)=(3,1)$ and $\|f\|_d=\sup_{S\le s\le S_\epsilon} e^{ds}|f(s)|$. We take the convex trial set
\begin{equation}\label{typeii:eq:compact-path-set}
 \begin{gathered}
 \underset{S\le s\le S_\epsilon}{\operatorname{ess\,sup}}\,
     e^{d_js}(|Z_j|+|\partial_s Z_j|)\le R\quad(j=o,i),\\
 u_i(S)=0,\qquad 16v_i(S_\epsilon)=17u_i(S_\epsilon).
 \end{gathered}
\end{equation}
For $R\ge1$ and $R^2\tau_0\le1$ the inputs obey the weaker path bounds at weights $\frac{5}{2}$ and $\frac{1}{4}$ with constant one. Taylor's formula for the two powers in $A$ yields
\begin{equation*}
 \begin{aligned}
 |f_o-c_o-c_o^i|&\lesssim\tau^2 |u_o|^2,\\
 |f_i-c_i^i|&\lesssim\tau^{14}(|u_o|+|u_i|^2).
 \end{aligned}
\end{equation*}
In particular $\tau^{14}|u_i|^2\le R^2\tau_0\tau^{15}$. By \eqref{typeii:eq:scale-rate-bounds}, $|g_j|+\tau|f_j|\lesssim\tau^{m_j-1}$ uniformly in $R$ and $\epsilon$. Dividing \eqref{typeii:eq:scale-displacement-bound} by $q_j^*$ and using $\partial_s=-\tau\partial_\tau$, we obtain
\begin{equation}\label{typeii:eq:sharp-forcings}
 \|\widehat Z_o\|_3+\|\partial_s\widehat Z_o\|_3+
 \|\widehat Z_i\|_1+\|\partial_s\widehat Z_i\|_1\lesssim1 .
\end{equation}
We choose $R$ larger than the uniform constant and then take $\tau_0$ small. The map preserves \eqref{typeii:eq:compact-path-set} for all $0<\epsilon<\tau_0$.

\subsection{Continuity on a finite interval}
\label{typeii:sec:finitecompact}
We equip \eqref{typeii:eq:compact-path-set} with the uniform topology. This convex set contains zero. On $I=[S,S_\epsilon]$ the paths are bounded in $W^{1,\infty}(I)$. By Arzel\'a--Ascoli we may extract a uniformly convergent subsequence. On $I$ we have
\begin{equation*}
 Z_{j,n}\to Z_j\ \text{uniformly},\qquad
 \partial_s Z_{j,n}\overset{*}{\rightharpoonup}\partial_s Z_j
 \quad\text{in }L^\infty(I).
\end{equation*}
The second limit follows by integration by parts. Weak-star lower semicontinuity preserves \eqref{typeii:eq:compact-path-set}, so the path set is compact.

\paragraph{Compactness of energy traces.}
On a compact physical time interval $I$ we assume
\begin{equation*}
 \sup_n\sup_{t\in I}\big(\|\nabla h_n(t)\|_2+\|k_n(t)\|_2\big)<\infty,\qquad
 \sup_n\|(\partial_t h_n,\partial_t k_n)\|_{L^\infty(I;L^2(B)\times H^{-1}(B))}<\infty
\end{equation*}
for every ball $B$. After extraction we claim
\begin{equation}\label{typeii:eq:compact-wave-traces}
 h_n\longrightarrow h\quad\text{in }C_tL^p_{x,\rm loc}\quad(2\le p<6),
 \qquad (\nabla h_n,k_n)\rightharpoonup(\nabla h,k)
       \quad\text{at every time}.
\end{equation}
For a fixed spatial cutoff $\eta$, by Rellich compactness and equicontinuity we have $(\eta h_n,\eta k_n)\to(\eta h,\eta k)$ in $C(I;L^2\times H^{-1})$ after extraction. We pass to a diagonal subsequence over all balls. Wherever $\eta h_n$ is bounded in $L_t^\infty H_x^m$, interpolation yields
\begin{equation}\label{typeii:eq:compact-interpolation}
 \sup_I\|\eta(h_n-h)\|_{H^s}
 \lesssim\left(\sup_I\|\eta(h_n-h)\|_2\right)^{1-\frac{s}{m}}\left(\sup_j\|\eta h_j\|_{L_t^\infty H_x^m}\right)^{\frac{s}{m}},
 \qquad 0\le s<m .
\end{equation}
At $m=1$, Sobolev embedding implies strong convergence in every local $L^p$, $2\le p<6$, hence convergence of the quintic products in $C_tL^1_{\rm loc}$. Density of smooth tests and the energy bound extend weak energy trace convergence uniformly to compact sets of $L^2$ tests.
If the waves also have a uniform $L^5(I;L^{10})$ bound, weak compactness and the preceding local convergence imply
\begin{equation}\label{typeii:eq:duhamel-limit}
 \|h\|_{L^5(I;L^{10})}\le\liminf_n\|h_n\|_{L^5(I;L^{10})},\qquad
 \|h^5\|_{L^1(I;L^2)}=\|h\|_{L^5(I;L^{10})}^5 .
\end{equation}
For a limiting distributional system $h_t=k+f_1$, $k_t=\Delta h+f_2$ with $(f_1,f_2)\in L^1(I;\dot H^1\times L^2)$, the Duhamel pair with the same weak trace at one time belongs to $C(I;\dot H^1\times L^2)$. Their difference is a free wave with zero trace, hence vanishes by duality against free test waves.

\paragraph{Closedness of the two prescribed wave equations.}
Suppose $(Z_{o,n},Z_{i,n})\to(Z_o,Z_i)$ uniformly on $[S,S_\epsilon]$. On $I=[T-\tau_0,T-\epsilon]$ both waves have support in a fixed ball and the scales are bounded away from zero. By the two wave equations,
\begin{equation*}
 \sup_n\big(\|\partial_t h_n\|_{L^\infty(I;L^2)}+
            \|\partial_t k_n\|_{L^\infty(I;H^{-1})}\big)<\infty .
\end{equation*}
We use the velocity and pressure bounds in the first row and $L^{\frac{6}{5}}\hookrightarrow H^{-1}$ for the quintic and exterior forces in the second. Thus \eqref{typeii:eq:compact-wave-traces} applies.
On a small fixed ball about the origin the middle waves are bounded in $L^\infty(I;H^4)$ by Section~\ref{typeii:sec:outermoment}. With $m=4$ and $s=3$, \eqref{typeii:eq:compact-interpolation} implies $C^1$ convergence on a smaller ball. Hence \eqref{typeii:eq:average-slope} and \eqref{typeii:eq:slope-rate} imply
\begin{equation}\label{typeii:eq:central-trace-limit}
 d_{h_n}\to d_h\quad\text{uniformly on }I,\qquad
 \dot d_{h_n}\overset{*}{\rightharpoonup}\dot d_h\quad\text{in }L^\infty(I).
\end{equation}
The derivative convergence follows by integration by parts and density in $L^1(I)$. The moving spatial tests range over compact subsets of $L^2$, so the uniform pairing convergence just proved applies.

The columns depend on the path values and $d_h$, and the Gram inverses are uniformly bounded. In the tested equations each path derivative and $\dot d_h$ occurs linearly, the latter only in the inner position row. For such a derivative $a_n'$ and its tested coefficient $g_n$, we use
\begin{equation}\label{typeii:eq:weak-path-rates}
a_n'g_n\rightharpoonup a'g\quad\text{in distributions},\qquad
 a_n'\overset{*}{\rightharpoonup}a'\quad\text{in }L^\infty(I),\qquad
 g_n\to g\quad\text{in }L^1(I).
\end{equation}
We first replace $g_n$ by $g$ using the uniform derivative bound, then test weak-star convergence against $g\varphi$ for $\varphi\in C_c^\infty(I^\circ)$.
The quintics converge in distributions by \eqref{typeii:eq:compact-wave-traces} at $p=5$. A finite cover by time cells supplies the hypothesis of \eqref{typeii:eq:duhamel-limit} for both waves. The remaining applied and pressure sources have the indicated $L^1$ energy bounds on $I$, so both limits satisfy their Duhamel equations. At $T-\epsilon$ the terminal fields are finite sums of the prescribed packets. Their coefficients are bounded and the packets converge in energy as the endpoint scales converge. Their weak limits retain zero terminal continuous data and the prescribed relation between position and velocity coefficients. At $T-\tau_0$ the weak traces retain the past mode condition. The finite-interval inverse in Section~\ref{typeii:sec:mixedinverse} selects the linear solution uniquely, and the same-path contractions of Sections~\ref{typeii:sec:nonlinearmiddle} and \ref{typeii:sec:closing} identify first the middle wave and then the inner wave. This identifies the limit of every subsequence.

On the fixed interval the force terms depending on $q$ in \eqref{typeii:eq:normalized-rates} converge uniformly along convergent paths. The rates are bounded in $L^\infty$ and converge in distributions by \eqref{typeii:eq:weak-path-rates}. By \eqref{typeii:eq:scale-displacement-bound} the output pairs $(x_j,z_j)$ are equicontinuous on the fixed interval. Every uniform subsequential limit solves \eqref{typeii:eq:scale-linear-pair} in distributions and obeys the same endpoint conditions. Uniqueness in \eqref{typeii:eq:scale-modes} therefore proves uniform convergence of the full sequence.

\subsection{Selection of the scales and the exact wave}\label{typeii:sec:closing}
For each path in \eqref{typeii:eq:compact-path-set} we first use the one-core inverse. The velocity tests and terminal packets lie in $r<\frac{\tau}{32}$, and the nonlinear cutoff equals one for $r\le\frac{\tau}{16}$ and vanishes for $r\ge\frac{\tau}{8}$. For the middle pressure column and velocity test of Section~\ref{typeii:sec:perturbedmatched}, the annular tail estimates in physical units are
\begin{equation*}
 \|E_o^0-E_o\|_2+\tau\|\nabla(E_o^0-E_o)\|_2
 \lesssim\left(\frac{\mu}{\tau}\right)^{\frac{1}{2}} .
\end{equation*}
Thus we may use the two cutoffs in the mixed inverse. By \eqref{typeii:eq:varying-middle-residual}, the middle residual and its first three middle-scaled derivatives are $O(\tau^5)$ in $L^2$. We apply Section~\ref{typeii:sec:nonlinearmiddle} to \eqref{typeii:eq:prescribed-middle-wave} and obtain the energy and interior bounds of order $\tau^6$. The transfer estimate \eqref{typeii:truepoltransfer} now has size
\begin{equation*}
 \frac{\rho^{\frac{3}{2}}}{\mu}\,\tau^6
 =\frac{\rho^{\frac{5}{2}}}{\lambda}\,\tau^6\lesssim\tau^{20}.
\end{equation*}
The slow inner position correction and the remaining deterministic inner force are $O(\tau^{17})$ in their respective rows.

We now solve \eqref{typeii:eq:prescribed-inner-wave}, with the middle multipliers as in \eqref{typeii:eq:sumtransport}. On a dyad $T-t\asymp r$, \eqref{typeii:eq:exterior-dyad} and the disjoint supports imply
\begin{equation*}
 \int_{T-t\asymp r}\|f_{\rm ext}(t)\|_2\,dt\lesssim r^{\frac{33}{2}},
 \qquad \langle f_{\rm ext},E_j\rangle=0 .
\end{equation*}
The direct mode pairings vanish by \eqref{typeii:eq:exterior-orthogonality}. In the difference of two inner iterates the prescribed $f_{\rm ext}$ cancels. We use the ball with weighted energy and cell $L^4_tL^{12}_x$ bounds at $q=16$, with the weighted uniform energy distance. The ball is complete since convergence in energy identifies each weak $L^4_tL^{12}_x$ limit, so the cell bounds pass to the limit. Equations \eqref{typeii:eq:mixed-wave-bound} and \eqref{typeii:eq:cell-strichartz} preserve the ball for a fixed large radius and small $\tau_0$. By \eqref{typeii:eq:cell-remainder} the wave and pressure differences gain $\tau_I^{16-15}=\tau_I$ in the energy distance alone. The potential contribution is absorbed by \eqref{typeii:eq:potential-multiplier} and its tested bound. Interpolation with $L^\infty_tL^6_x$ then controls the cell $L^5_tL^{10}_x$ norm. Thus the unique inner wave satisfies, uniformly over the admissible paths,
\begin{equation*}
 \|\nabla h_i(t)\|_2+\|k_i(t)\|_2\lesssim\tau^{16},\qquad
 \|h_i\|_{L_t^5L_x^{10}(I)}\lesssim\tau_I^{16}.
\end{equation*}

Section~\ref{typeii:sec:finitecompact} and the uniform bounds in \eqref{typeii:eq:sharp-forcings} allow us to apply Schauder's theorem on \eqref{typeii:eq:compact-path-set} for every $0<\epsilon<\tau_0$. The smallness choice for $\tau_0$ is independent of $\epsilon$. At each fixed point both brackets in \eqref{typeii:eq:sumtransport} vanish. The total pair $X_\epsilon=(u_\epsilon,v_\epsilon)$ therefore satisfies
\begin{equation*}
 \partial_tu_\epsilon=v_\epsilon,\qquad
 \partial_t v_\epsilon=\Delta u_\epsilon+u_\epsilon^5
 \quad\hbox{on }[T-\tau_0,T-\epsilon].
\end{equation*}
Finite propagation fixes $X_\epsilon=(B,0)$ on $|x|>\tau$.

The only velocity of order $\tau^3$ is $-p_oD_\mu\Lambda Q$. By \eqref{typeii:eq:compact-path-set}, $p_o=-4a\tau^3+O(\tau^6)$ and $p_i=O(\tau^{15})$. Cutting off the middle velocity at $r\asymp\tau$ costs $O(|p_o|(\mu/\tau)^{\frac{1}{2}})$ by the dipole tail. The matched annular field is $O(\tau^4)$ with first derivatives $O(\tau^3)$. Including the elliptic corrections and the two selected waves, uniformly for $\epsilon\le\tau\le\tau_0$ we obtain
\begin{equation*}
 \|\nabla(u_\epsilon-B-Q_{\mu_\epsilon}-P_{\lambda_\epsilon})\|_2+
 \|v_\epsilon-4a\tau^3 D_{\mu_\epsilon}\Lambda Q\|_2
 \lesssim \tau^{\frac{9}{2}}.
\end{equation*}
Indeed the middle position correction has size $O(\mu^{\frac{3}{2}}(\tau/\mu)^{\frac{1}{2}})$ and the inner one $O(\rho^{\frac{3}{2}})$. Their velocity corrections, the frozen cone correction and the cutoff tails are $O(\tau^{\frac{9}{2}})$ by the same moving-frame estimates.
The passage from solutions with a finite terminal cutoff to the limiting solution uses the backward compactness construction of Merle \cite{MerleBack} and Martel \cite{MartelBack}, as adapted to wave towers by Jendrej \cite[Section~3.6]{J}. Let $\epsilon\downarrow0$. On each compact $I\subset[T-\tau_0,T)$ the scale paths are equi-Lipschitz and stay away from zero. The total fields have uniformly bounded energy norm, and the exact equation bounds $\partial_t v_\epsilon$ in $L^\infty(I;H^{-1}_{\rm loc})$. On $I$ the scales have positive lower bounds, so finitely many cells control $\|u_\epsilon\|_{L^5(I;L^{10})}$ uniformly. By \eqref{typeii:eq:compact-wave-traces} and \eqref{typeii:eq:duhamel-limit}, a diagonal subsequence has locally uniformly convergent scales and a limiting strong energy solution with
\begin{equation*}
 u_\epsilon\to u\quad\text{in } C(I;L^5_{\rm loc}),\qquad
 u_t=v,\quad v_t=\Delta u+u^5 .
\end{equation*}
By weak lower semicontinuity and strong continuity of dilation at positive scales,
\begin{equation}\label{typeii:eq:terminalmodulus}
 \|\nabla(u-B-Q_\mu-P_\lambda)\|_2+\|u_t-4a\tau^3 D_\mu\Lambda Q\|_2\lesssim\tau^{\frac{9}{2}},
 \qquad \vec u=(B,0)\quad(|x|>\tau),
\end{equation}
where $\vec u=(u,u_t)$ and $\mathcal H=\dot H^1\times L^2$.

\subsection{The singular point and its two energy regions}
\label{typeii:sec:terminalcutoffs}
We apply \eqref{typeii:eq:terminalmodulus} on two shrinking balls to identify the concentrating energies and locate the singularity.
The bound \eqref{typeii:eq:compact-path-set} holds with the same $R$ for every finite-cutoff fixed point. Passing to the locally uniform limit at the two weights $3$ and $1$, we obtain
\begin{equation}\label{typeii:eq:sharp-scales}
 \mu=a\tau^4(1+O(\tau^3)),\qquad \lambda=b\tau^{16}(1+O(\tau)).
\end{equation}
Put $\vec z=(u-B,u_t)$, $w=u-B-P_\lambda-Q_\mu$, $r_1=\sqrt{\lambda\mu}$, $r_2=\sqrt{\mu\tau}$ and $R=r_2/\mu$. The differentiated tails of $P$ and $Q$ satisfy
\begin{equation*}
 N_{|x|>r_1}(P_\lambda,0)\lesssim\rho^{\frac{5}{4}},\qquad
 N_{|x|<r_1}(Q_\mu,0)\lesssim\rho^{\frac{3}{4}} .
\end{equation*}
By \eqref{p8:eq:density-continuity} we discard the other core on either region at cost $O(\rho^{\frac{3}{4}})$. We then apply \eqref{p8:eq:local-stationary} with $k=0$ to $F+w$, taking $F=P_\lambda$ on the ball and $F=Q_\mu$ outside. Since $\rho\asymp\tau^{12}$, the two static energies equal $E(P,0)$ and $E(Q,0)$ up to $O(\tau^9)$. By \eqref{typeii:eq:terminalmodulus}, the total kinetic energy is $8M_o a^2\tau^6+O(\tau^{\frac{15}{2}})$. Its part on $|x|<r_1$ is $O(\tau^9)$, since $\|\Lambda Q\|_{L^2(|y|<\sqrt\rho)}\lesssim\rho^{\frac{5}{4}}$.

Kelvin inversion of the wall jet yields $Q(r\omega)=r^{-2}Y_1+r^{-3}Y_2+O(r^{-4})$, with the differentiated expansion and $Y_1\perp Y_2$. Testing $-\Delta Q=Q^5$ against $Q$ on $|y|>R$ and using \eqref{typeii:eq:two-forces}, we obtain
\begin{equation*}
 \frac{1}{2}\int_{|y|>R}|\nabla Q|^2\,dy
 =-\frac{R^2}{2}\int_{\mathbb S^2}Q\partial_rQ\,d\omega
   +\frac{1}{2}\int_{|y|>R}Q^6\,dy
 =\frac{C_1}{3R^3}+O(R^{-5}).
\end{equation*}
On $D=\{|x|>r_2\}$ the same tail bounds imply
\begin{equation*}
 N_D(Q_\mu,0)\lesssim R^{-\frac{3}{2}},\quad
 N_D(w+P_\lambda,0)\lesssim\tau^{\frac{9}{2}}+\left(\frac{\lambda}{r_2}\right)^{\frac{5}{2}},\quad
 \|u_t\|_{L^2(D)}\lesssim\tau^3 R^{-\frac{1}{2}}+\tau^{\frac{9}{2}}.
\end{equation*}
Cauchy--Schwarz bounds the $L^1$ difference of $e(\vec z)$ and $\frac{1}{2}|\nabla Q_\mu|^2$ there by $O(\tau^{\frac{9}{2}}R^{-\frac{3}{2}}+\tau^{\frac{15}{2}})=O(\tau^{\frac{27}{4}})$. The same comparison applies to $|e(\vec z)|$. Since $R^{-5}=O(\tau^{\frac{15}{2}})$, we conclude that
\begin{equation}\label{typeii:eq:singular-energy-regions}
 \begin{aligned}
 \int_{|x|<r_1}e(\vec z)\,dx&=E(P,0)+O(\tau^9),\\
 \int_{r_1<|x|<r_2}e(\vec z)\,dx&=E(Q,0)-\frac{C_1}{3R^3}+8M_o a^2\tau^6+O(\tau^{\frac{27}{4}}),\\
 \int_{|x|>r_2}|e(\vec z)|\,dx&=\frac{C_1}{3R^3}+O(\tau^{\frac{27}{4}}).
 \end{aligned}
\end{equation}
By \eqref{typeii:eq:terminalmodulus}, $(\nabla(u-B),u_t)$ is bounded in $L^2$ and supported on $B_\tau$, hence converges weakly to zero. The two dilation frames are orthogonal since $\lambda/\mu\to0$. Thus the same estimate yields
\begin{equation*}
 \int_{B_{r_2}} |\nabla u|^2\,dx
 \longrightarrow \|\nabla P\|_2^2+\|\nabla Q\|_2^2>0 .
\end{equation*}
Indeed $r_2/\mu\to\infty$ captures both profiles, and $\|\nabla B\|_{L^2(B_{r_2})}\to0$. By the stationary equations, $E(F,0)=\frac{1}{3}\|\nabla F\|_2^2$ for $F=P,Q$. The positive gradient mass on shrinking balls excludes a strong energy trace. Since $u=B$ outside $B_\tau$ and the energy norm is bounded, only the origin is singular and the blow-up is of Type~II.
\space
\appendix
\section{The variational Hopf state and its linearized kernel}\label{cert:start}
Let $D=\{(x,y,z)\in\Sph^2:x,y>0\}$, $a=\frac{3}{4}$ and
$P=-4\Delta_{\Sph^2}+a$. We normalize the area of the lune to one.
Reflection then identifies its $L^p$ norms with the corresponding norms
for normalized area on $\Sph^2$.

\begin{theorem}[Nondegeneracy of the Hopf state]\label{cert:thm:executed}\label{thm:stationary-certificate}
The minimization in \eqref{ap:min} determines a positive Dirichlet solution
$h$ of $Ph=h^5$ on $D$. Let $H$ be its odd reflection across both sides,
$q=H\circ\pi$, $A=-\Delta_{\Sph^3}+a-5q^4$ and $\tau f(Z)=f(\bar Z)$.
Then
\begin{equation}\label{ap:kernel}
 \ker A=T_q(\operatorname{Conf}(\Sph^3)\cdot q),\qquad
 \dim_{\R}\ker A=9,\qquad \dim_{\R}(\ker A\cap\ker(\tau+I))=4 .
\end{equation}
\end{theorem}

The cone minimization \eqref{ap:min} selects the profile, and the proof below
determines its full and odd kernels. For a different family, the polygonal states of
del Pino, Musso, Pacard and Pistoia \cite{DPMPP}, Musso and Wei
\cite[Theorem~1.1 and the discussion of $n=3$]{MW} state that the bounded
kernel is generated by the conformal fields, giving dimension nine when $n=3$.
Clapp, Pistoia and Weth \cite[Theorem~1.1]{EllCPW} constructed signed-equivariant
variational solutions on the round sphere, including $\mathbb S^3$.
For a higher-dimensional gluing construction involving a Clifford torus,
see Medina, Musso and Wei \cite{EllMMW}.

We first approximate $5h^4$ by an explicit rational potential in
Section~\ref{ap:potential}. We compare the fiber charges in
Section~\ref{ap:disk}, locate the remaining eigenvalues by their moments
in Section~\ref{ap:tests}, and count the zero modes in Section~\ref{ap:count}.

\subsection{Construction and pointwise bounds on the lune}\label{ap:profile}
We set $T=x^2+y^2$, $w=2xy$ and $\sigma=\frac{w}{T}$. Let $\mathcal C$ be
the closed convex cone in $H^1_0(D)$ generated by $e_{Nn}=T^{N-n}w^n$,
where $N\ge n\ge1$ and $n$ is odd. By compactness in $L^6(D)$ we may choose a minimizer $v$ in
\begin{equation}\label{ap:min}
 S_0=\min_{\substack{v\in\mathcal C\\ \|v\|_6=1}}\langle Pv,v\rangle,
 \qquad h=S_0^{\frac{1}{4}}v .
\end{equation}
For $G=P_D^{-1}$ we have $Gv^5\in\mathcal C$. Indeed, products of five
generators are again generators, and the identity
\begin{equation}\label{ap:triangular}
 Pe_{Nn}=(16N^2+8N+a)e_{Nn}
 -16(N^2-n^2)e_{N-1,n}-16n(n-1)e_{N-1,n-2}
\end{equation}
shows by induction that $Ge_{Nn}\in\mathcal C$. Terms with zero
coefficient are omitted. The assertion for $v$ follows by continuity.
We put $f=Gv^5$. By \eqref{ap:min}, H\"older and Cauchy--Schwarz,
\[
 S_0\|f\|_6^2\le\langle Pf,f\rangle
 =\langle v^5,f\rangle\le\|f\|_6,\qquad
 1=\langle Pf,v\rangle^2\le S_0\langle Pf,f\rangle\le1.
\]
Equality implies $v=S_0 f$. Thus $Ph=h^5$ and $\int_D h^6=S_0^{\frac{3}{2}}$.
By odd reflection and elliptic regularity, $H$ is analytic on the sphere.
By the strong maximum principle, $h>0$ in $D$.

Let $M$ denote the value of $h$ at the peak $x=y=\frac{1}{\sqrt{2}}$,
$z=0$. The cone expansion takes the form $h=M\mathbb E(T^N\sigma^n)$.
For a cone polynomial $v_j=\sum c_{Nn}^{(j)}e_{Nn}$ with $v_j\to h$ in $H^1$,
\[
 \sum_{N,n}c_{Nn}^{(j)}(r')^N
 \lesssim_{r'}\int_{\{w>r'\}}v_j,\qquad 0<r'<1.
\]
A diagonal subsequence converges uniformly on $T\le r<r'$. Along $\sigma=1$ the
limit as $T\uparrow1$ identifies the sum of the coefficients with $M$.
On the meridian $T=1-z^2$, and on the equator $\sigma=\cos(2\phi)$,
where $\phi$ is longitude measured from the peak. Since
$\frac{1-t^k}{1-t}\uparrow k$ as $t\uparrow1$, the second derivatives
there are $-2M\mathbb E N$ and $-4M\mathbb E n$, respectively. Hence
\[
 8\mathbb E(N+2n)=M^4-a,\qquad
 h\ge M T^{\mathbb E N}\sigma^{\mathbb E n}.
\]
We used Jensen's inequality in the second bound. By the beta integral and
$\Gamma(t+\frac{1}{2})^2\le\Gamma(t)\Gamma(t+1)$, for $p,q\ge0$,
\[
 \int_D T^p\sigma^q
 =\frac{\Gamma(p+1)\Gamma(\frac{q+1}{2})}
 {2\Gamma(p+\frac{3}{2})\Gamma(\frac{q}{2}+1)}
 \ge\frac{1}{\sqrt{2(p+1)(q+1)}}\ge\frac{2}{p+2q+3}.
\]
Taking $(p,q)=6(\mathbb E N,\mathbb E n)$, we obtain
$\int_D h^6\ge\frac{32M^6}{12M^4+39}$.

For the upper bound we use $\psi=w+w^5$. We write $I_k=\int_D w^k$.
We integrate $-\Delta w=6w$ against $w^{i+j-1}$. For $i,j\ge1$,
\[
 \langle Pw^i,w^j\rangle=
 \left(a+\frac{24ij}{i+j-1}\right)I_{i+j},\qquad
 I_0=1,\quad I_{2k+2}=\frac{4(2k+1)^2}{(4k+3)(4k+5)}I_{2k}.
\]
The recurrence follows from the preceding beta integral with $p=q=2k$.
Substituting $(i,j)\in\{1,5\}^2$ in the energy and
$\int_D\psi^6=\sum_{j=0}^6\binom{6}{j} I_{6+4j}$ in \eqref{ap:min}, we obtain
\begin{equation}\label{ap:peak}
 \frac{1024M^{12}}{(12M^4+39)^2}\le S_0^3
 \le\frac{\langle P\psi,\psi\rangle^3}{\int_D\psi^6}
 <1660<\frac{1024\cdot240^3}{2919^2},\qquad M^4<240 .
\end{equation}

To bound $h$ away from the peak, set
\[
 c=\frac{x+y}{\sqrt{2}},\quad d=\frac{x-y}{\sqrt{2}},\quad
 X=1-c^2,\qquad v(r)=M\mathbb E(1-r^2)^N.
\]
We first prove $h\le v(\sqrt X)$. On $x>y>0$, $z>0$, the rotation
fields $F=zH_y$ and $K=z(H_y-H_x)$ solve $(P-5H^4)F=(P-5H^4)K=0$.
Here we differentiate the extension $H(T,w)$ independent of $z$.
The positive series implies $F>0$ in the triangle. On the open $y=0$
edge we have $F>0$ by the boundary point lemma. On its three edges
$K\ge0$ by reflection. For the negative part $K_-$ the ground state
identity and its equation imply
\[
 4\int F^2\left|\nabla\frac{K_-}{F}\right|^2
 \le \langle (P-5H^4)K_-,K_-\rangle=0 .
\]
The inequality follows by density from the identity for compactly
supported tests. Thus $K_-$ is a multiple of $F$, and its zero trace
on $y=0$ forces $K_-=0$. At fixed $c$, the angular derivative away
from $d=0$ equals $-\frac{K}{\sqrt{2}}\le0$. Reflection proves the claim.

On $d=0$ the equation for $v$ reads
\[
 -4(1-r^2)\left(v''+\frac{2}{r}v'\right)
 =v^5-av+\frac{16M}{1-r^2}\mathbb E[(N-n)(1-r^2)^N].
\]
We take $R^2=X$, $V=v(R)^4>a$ and $\delta=1-\frac{a}{V}$.
Then $-\Delta_{\R^3}v\ge\frac{\delta}{4}v^5$ on the ball of radius $R$.
Iteration between the constant $v(R)$ and $v$ produces a radial solution
with boundary value $v(R)$ and central value $t^{\frac{1}{4}}\le M$.
The explicit radial solution and $(1+y)^2\ge4y$ imply
\[
 V=\frac{t}{(1+\frac{\delta tX}{12})^2},\qquad
 t<240,\qquad V-a=\delta V\le\frac{3}{X}.
\]
We write $g=\frac{79}{4}$. If $gX\le1$ and $V\ge\frac{240}{(1+gX)^2}$,
then $V\ge60$ and $\delta\ge\frac{g}{20}$. On $0\le t\le240$ the
function $t(1+\frac{gXt}{240})^{-2}$ increases, hence
$V\le\frac{240}{(1+gX)^2}$. Finally $\frac{v(\sqrt X)}{1-X}$
decreases in $X$ and $v(\frac{1}{\sqrt{5}})^4\le\frac{63}{4}$.
Thus $h^4\le H_b(X)$, where $H_b$ is the minimum of the applicable bounds
\begin{equation}\label{ap:env}
 \frac{240}{(1+gX)^2}\quad(gX\le1),\qquad
 a+\frac{3}{X}\quad(X>0),\qquad
 \frac{63}{4}\left[\frac{5}{4}(1-X)\right]^4
 \quad\left(X\ge\frac{1}{5}\right).
\end{equation}

\subsection{A rational approximation to the potential}\label{ap:potential}
We estimate $5h^4$ by comparing $h$ with one explicit trial. In the variables
$u=z^2$, $v=(x^2-y^2)^2$, we have $w^2=(1-u)^2-v$. We set
\[
 L_0(t)=\frac{3}{34}+\frac{7}{6}t,\qquad
 p=\frac{N}{D_0},\quad h_*=\frac{w}{\sqrt p}.
\]
The two arrays list $10^9[u^iv^j]N$ and $10^9[u^iv^j]D_0$ by row $i$
and column $j$. Blank entries are zero.
\begin{center}\small
\setlength{\tabcolsep}{4pt}\begin{tabular}{crrrrrr}\toprule
$N$ & $j=0$ & $1$ & $2$ & $3$ & $4$ & $5$\\\midrule
0 & 88371781 & 446158659 & 573423893 & 185804138 & 1490123 & -765416\\
1 & 1952652646 & 5305706904 & 2917108940 & 146507823 & -13495811 & \\
2 & 12182019655 & 15022036590 & 2062685788 & -83482297 &  & \\
3 & 25274895106 & 9909090144 & -61372815 &  &  & \\
4 & 15528842903 & 940327368 &  &  &  & \\
5 & 2056793093 &  &  &  &  & \\
\bottomrule\end{tabular}\end{center}
\begin{center}\small
\setlength{\tabcolsep}{4pt}\begin{tabular}{crrrrrr}\toprule
$D_0$ & $j=0$ & $1$ & $2$ & $3$ & $4$ & $5$\\\midrule
0 & 1000000000 & 1905128713 & 827668942 & 65000067 & -1557564 & -115944\\
1 & 8845491633 & 8063874369 & 955345940 & -37839843 & -1373898 & \\
2 & 20286776957 & 5701656829 & -229364295 & -2093478 &  & \\
3 & 11479683948 & -403565381 & 10840947 &  &  & \\
4 & 117647467 & 34786629 &  &  &  & \\
5 & 77868796 &  &  &  &  & \\
\bottomrule\end{tabular}\end{center}
We use the potential
\begin{equation}\label{ap:surrogate}
 K=1+\frac{119}{9}u+\frac{13}{4}v,\quad
 S=w^4\left(\frac{2561}{4}+\frac{543}{4}v-36u-228uv-105u^2\right),
 \quad V_c=\frac{S}{K^2}.
\end{equation}

\begin{lemma}\label{ap:reslem}
With the area normalization above,
\begin{equation}\label{ap:res}
 L_0(u)+\frac{1}{4}v\le p<2,\quad
 \|Ph_*-h_*^5\|_2<0.8\cdot10^{-6},\quad
 \|5h_*^4-V_c\|_\infty\le\frac{17}{400}.
\end{equation}
\end{lemma}
\begin{proof}
We set $l=L_0(u)+\frac{1}{4}v$, $L=1+\frac{12}{5}u+\frac{1}{2}v$ and
\[
 J=K^2N^2,\qquad Z_1=5w^4K^2D_0^2-SN^2 .
\]
On the triangle $u,v\ge0$, $u+v\le1$, the Bernstein coefficients
of degree $16$ are
\[
 \beta_{ij}(u^av^b)=\frac{\binom{i}{a}\binom{j}{b}}
                    {\binom{16}{a+b}\binom{a+b}{a}},\qquad i+j\le16.
\]
We take out-of-range binomial coefficients to be zero. The coefficient minima are
\begingroup\setlength{\arraycolsep}{3pt}
\[
\begin{array}{c|rrrrr}
 F&2D_0-N-1&N-lD_0&D_0-L^3&\frac{17J}{400}+Z_1&\frac{17J}{400}-Z_1\\
 \lfloor10^6\min\beta_{ij}(F)\rfloor&911628&136&0&253&72 .
\end{array}
\]
\endgroup
Positivity of the basis implies $D_0\ge L^3$, $l\le p<2$ and
$\frac{|Z_1|}{J}\le\frac{17}{400}$.

For the residual let $\mathcal Lp=\sum g_{ij}p_{ij}+\ell\cdot\nabla p$,
where the indices refer to $(u,v)$ and
\[
 (g_{ij})=\begin{pmatrix}4u(1-u)&-8uv\\-8uv&16v(1-u-v)\end{pmatrix},
 \qquad \ell=(2-14u,\,8(1-u)-36v).
\]
By the product rule,
\begin{equation}\label{ap:resformula}
 R=Ph_*-h_*^5=\frac{wZ}{p^{\frac{5}{2}}D_0^4},\qquad
 Z=D_0^4\left[(24+a)p^2-w^4+2p\mathcal Lp-3\sum g_{ij}p_ip_j\right].
\end{equation}
Here $Z$ is a polynomial. Since $0\le L-1\le\frac{12}{5}$, we have
$l\ge L_0(0)L^2$. Hence $R^2\le L_0(0)^{-5}w^2 Z^2 L^{-34}$.

We use one polynomial inequality for this reciprocal and the spectral
moments below. We choose $(A,B)=(2,3)$ for $x=L$ and $(4,13)$ for $x=K$.
We write
\begin{align*}
 W(x)&=\left[\left(1-\frac{x}{A}\right)\left(1-\frac{x}{B}\right)\right]^{32},
 &h_k(x)&=\sum_{j=0}^{k-1}[t^j]\frac{1}{W(t)}\,x^j,\\
 P_k(x)&=\frac{1-W(x)h_k(x)}{x^k},
 &E_k(x)&=W(x)h_k(1).
\end{align*}
The Taylor remainder makes $P_k$ a polynomial. For $x\ge1$ the
coefficients of $h_k$ are nonnegative, hence
$0\le x^{-k}-P_k(x)=\frac{W(x)h_k(x)}{x^k}\le E_k(x)$.
For $B_*,C_*\ge0$ it follows that
\begin{equation}\label{ap:recip}
 -E_k(x) C_*\le\frac{B_*-C_*}{x^k}-P_k(x)(B_*-C_*)\le E_k(x) B_*.
\end{equation}
The area average $\mathcal M$ factors into a uniform integral in $|z|$
and an angular average of $\sin^{2j}(2\phi)$. Evaluating the beta integrals, we find
\begin{equation}\label{ap:mom}
 \mathcal M(u^iv^j)
 =\frac{(2i-1)!!(2j)!\binom{2j}{j}}{(2i+4j+1)!!},\qquad (-1)!!=1 .
\end{equation}
Applying \eqref{ap:recip} with $k=34$ and expanding by \eqref{ap:mom},
we obtain
\[
 \|R\|_2^2\le L_0(0)^{-5}\mathcal M[w^2 Z^2(P_{34}+E_{34})(L)]
 \in 10^{-18}(627508,627509).
\]
\end{proof}

\begin{lemma}\label{ap:secant}
On Dirichlet functions even under $z\mapsto-z$ and $x\leftrightarrow y$,
the secant potential $S_5=\sum_{i=0}^4 h^i h_*^{4-i}$ satisfies
\begin{equation}\label{ap:secgap}
 \lambda_1(P-S_5)<-4,\qquad \lambda_2(P-S_5)>\frac{2}{5}.
\end{equation}
\end{lemma}
\begin{proof}
By \eqref{ap:peak}, $h\le Mw<6h_*$ in $D$. The positive Dirichlet
ground state $w$ of $P$ has eigenvalue $24+a$. We test with $h$ to obtain
\[
 \langle(P-S_5)h,h\rangle<-\frac{1}{6}\int_D h^6
 \le-\frac{24+a}{6}\|h\|_2^2 .
\]

For the upper bound on $S_5$, we fix $c$ and put $X=1-c^2$. Then
\[
 0\le u\le X,\qquad v=4(1-X)(X-u),\qquad w=1-2X+u .
\]
At fixed $X$, the derivative of $\frac{w^2}{L_0(u)+\frac{v}{4}}$ has the
sign of $w[\frac{1}{102}+\frac{4}{3}X+(\frac{1}{6}+X)u]$. By
\eqref{ap:res} we may therefore take $u=X$ to obtain
\[
 h_*^4\le V_b(X)=\frac{(1-X)^4}{L_0(X)^2},\qquad
 S_5\le B=\sum_{i=0}^4 H_b^{\frac{i}{4}} V_b^{1-\frac{i}{4}} .
\]
We extend the Dirichlet trial by zero to $c>0$. Its angular frequencies
about the $c$-axis are even. By weighted arithmetic--geometric mean,
\[
 XB\le\frac{5}{2}(XH_b+XV_b)
 \le\frac{5}{2}\left(\frac{15}{4}+\frac{17}{7}\right)<16,\qquad
 B\le\frac{499}{256}H_b+\frac{94}{27}V_b .
\]
For the second inequality we apply the same inequality after scaling
$H_b$ by $(\frac{3}{4})^4$. Thus the angular energy
$\frac{4k^2}{X}$ dominates $B$ for $k\ne0$. On $c\le\frac{7}{10}$
each latitude intersects $D$ in two arcs of angular length
$2\arcsin(\frac{c}{\sqrt X})\le\frac{\pi c}{\sqrt X}$.
There the Dirichlet Poincar\'e inequality bounds the angular energy
below by $\frac{4}{c^2}$ times the mass. We use this bound for the whole
trial on those latitudes. By the last branch of \eqref{ap:env},
$B\le84c^8<\frac{4}{c^2}$ there. Hence every nonzero frequency in the
resulting radial form lies above $\frac{2}{5}$.

For both radial comparisons we use the following identity on an interval $I$.
If $e>0$ and the real function $Y$ has no interior zero, then
\begin{equation}\label{ap:groundtransform}
 \int_I(e|f'|^2-W|f|^2)
 =\int_I eY^2\left|\left(\frac{f}{Y}\right)'\right|^2
 +\int_I\left[-\frac{(eY')'}{Y}-W\right]|f|^2
 +\left[e\frac{Y'}{Y}|f|^2\right]_{\partial I}.
\end{equation}
At a simple interior zero of $Y$ we impose $f=0$ and apply
\eqref{ap:groundtransform} on each side, first to functions vanishing
near the cut and then by density.

For the angular mean we put $s=\frac{1-c}{1+c}$. On the intervals with endpoints
$0,\frac{1}{78},\frac{7}{125},\frac{3}{17}$, we take $H$ to be the
first, second and third bounds in \eqref{ap:env}, respectively. Indeed
$X(\frac{1}{78})\le\frac{4}{79}$ and $X(\frac{7}{125})\ge\frac{1}{5}$.
The resulting potential bound is
\[
 M_1=\begin{cases}
 \frac{499}{256}H+\frac{94}{27}V_b,&0\le s<\frac{3}{17},\\[2pt]
 84c^8-\frac{4}{c^2},&\frac{3}{17}\le s<1 .
 \end{cases}
\]
We set $z_0=\frac{1}{130}$ and impose $f(z_0)=0$ on the angular mean.
We take
\[
 Y=(1-130s)c^{\frac{3}{2}}(1+bs)^{-\frac{4}{3}},\qquad
 \frac{4b}{3(130+b)}=\frac{1}{2}-\frac{390}{130^2-1}.
\]
This choice of $b$ implies $(sY')'(z_0)=0$. On each of the four intervals $[l,h]$ we write
\begin{equation}\label{ap:signrule}
 \frac{P_1}{Q}=(1-s)^2\left[-\frac{(sY')'}{Y}
       -\frac{M_1-\frac{7}{20}}{4(1+s)^2}\right],\qquad
 a_i=[t^i](1+t)^{17} P_1\!\left(\frac{l+ht}{1+t}\right),
\end{equation}
in lowest terms with $Q>0$ in the interior. The pole at $z_0$ is
removable. In each interval $a_0>0$, and direct expansion determines the integer part
\[
 \left\lfloor\min_{1\le i\le17}\frac{a_i}{a_0}\right\rfloor
 =17,\ 14,\ 8,\ 29
\]
in order of increasing $s$. Thus the bracket in \eqref{ap:signrule}
is strictly positive. The radial Laplacian acts by $\Delta f=(1+s)^2(sf')'$.
We apply \eqref{ap:groundtransform} with $e=s$ and
$W=\frac{M_1-\frac{7}{20}}{4(1+s)^2}$. At $s=0$ we have
$|f(s)|\lesssim_f 1+\sqrt{|\log s|}$ by Cauchy--Schwarz, so the flux
vanishes. To justify the endpoint at $s=1$ ($c=0$), first use
a trial in $C_c^\infty(D)$, whose zero extension and angular mean
vanish for small $c$. These trials are dense in $H^1_0(D)$.
Averaging is continuous into radial $H^1$ near $c=0$, where the radial
energy density $1-c^2$ is bounded below. The mean therefore has zero
trace there. The term $\int |f|^2/c^2$ passes to the limit by the
one-dimensional Hardy inequality at this endpoint. Since $a-\frac{2}{5}=\frac{7}{20}$,
the radial form exceeds $\frac{2}{5}\|f\|_2^2$ under this single trace
condition. The min--max principle proves \eqref{ap:secgap}.
\end{proof}

Now $(P-S_5)(h-h_*)=-R$, hence $y=h-h_*$ satisfies
$\|y\|_2\le\frac{5}{2}\|R\|_2$. Its odd reflection has only even spherical
degrees $\ell\ge2$. For $C=1200+a$ set
$x_\ell=4\ell(\ell-1)+C-4$. The midpoint of $x_\ell,x_{\ell+2}$ is
$4\ell(\ell+1)+C$ and $x_{\ell+2}-x_\ell=8(2\ell+1)$. Thus
\[
 \sum_{\ell=2,4,\dots}\frac{2\ell+1}{[4\ell(\ell+1)+C]^2}
 \le\sum_{\ell=2,4,\dots}\frac{x_{\ell+2}-x_\ell}{8x_\ell x_{\ell+2}}
 =\frac{1}{8(C+4)} .
\]
By the spherical addition theorem, Cauchy--Schwarz and $0\le S_5\le1200$,
\begin{equation}\label{ap:sup}
 \|y\|_\infty\le\frac{\|(P+1200)y\|_2}{\sqrt{8(C+4)}}
 \le\frac{6001\|R\|_2}{\sqrt{8(C+4)}}<0.00004891 .
\end{equation}
Since $h^4,h_*^4\le(L_0(0)^{-\frac{1}{2}}+0.00004891)^4<129$,
\begin{equation}\label{ap:verror}
 \|5h^4-5h_*^4\|_\infty\le20\cdot129^{\frac{3}{4}}\|y\|_\infty<0.0375,\qquad
 \|5h^4-V_c\|_\infty<\varepsilon=\frac{2}{25}.
\end{equation}
For two scaled minimizers $h_1,h_2$ the secant potential between them
differs from the one between $h_*,h_2$ by at most
$10\cdot129^{\frac{3}{4}}\cdot0.00004891<0.025$. The first secant operator
is therefore invertible by \eqref{ap:secgap}. It annihilates $h_1-h_2$,
so the scaled minimizer is unique.
\space
\subsection{The free lune forms and comparison in the fiber charge}\label{ap:disk}
Let $A_{m,j}$ be the restriction of $A$ to
\[
 f(e^{it}Z)=e^{imt}f(Z),\qquad f(iZ_1,Z_2)=i^j f(Z),
 \qquad m\in\mathbb Z,\quad j\in\mathbb Z/4\mathbb Z .
\]
The charge-$m$ functions are sections of the monopole line bundle over
the base sphere, as in Wu and Yang \cite{EllWY}. For the relation
between the round three-sphere Laplacian and the magnetic Laplacian on
the unit base sphere, see \cite[Section~5.5.2]{EllCdVT}.
We relax the matching conditions on the sides of $D$ and denote the
resulting operator by $B_m$. The horizontal dilation of the Hopf map is
two, so its form and the min--max comparison are
\begin{equation}\label{ap:free}
 Q_m(f)=\int_D\big[4|(d+im\alpha)f|^2+(m^2+a-V)|f|^2\big],
 \qquad V=5h^4,\qquad \lambda_i(A_{m,j})\ge\lambda_i(B_m).
\end{equation}
We count eigenvalues with multiplicity, starting at one.
In the radial gauge about the peak, $\alpha=\frac{s}{1+s}\,d\omega$
and $4|\alpha|^2=s$. Rotation through $\pi$ about the peak axis
splits $B_m=B_m^+\oplus B_m^-$. We claim
\begin{equation}\label{ap:diskgap}
 d_m=\frac{3(1+m^2)}{4},\qquad
 \lambda_2(B_m^+)>d_m,\qquad \lambda_3(B_m^-)>d_m+\frac{16}{3}.
\end{equation}

We map $D$ to the disk by
\[
 Z=e^{-\frac{i\pi}{4}}\frac{x+iy}{1+z},\qquad
 \zeta=\frac{1-Z^2}{1+Z^2}=\rho e^{i\theta},\qquad
 r=\rho^2,\qquad D_1=\frac{1+r+|1-\zeta^2|}{2}.
\]
In the metric $J|d\zeta|^2$ we have
\[
 J=\frac{1}{D_1|1-\zeta^2|}\le\frac{1}{1-r},\quad JT=D_1^{-2},\quad
 X=\frac{r}{D_1},\quad D_1=\frac{1+r}{1+u}\ge1,\quad w=\frac{1-r}{D_1}.
\]
By \eqref{ap:res}, $p\ge L_0(u)+\frac{v}{4}$.
Since $v=(1-u)^2-w^2$ and $u\le r$, we obtain
\[
 D_1^2\left(L_0(u)+\frac{v}{4}\right)-L_0(r)
 =\frac{(r-u)[17ru+9(r+u)+1]}{102(1+u)^2}\ge0 .
\]
Thus $h_*\le(1-r)L_0(r)^{-\frac{1}{2}}$. Also $h\le Mw$ and
$h_*\le L_0(0)^{-\frac{1}{2}}w$, so $\max(h,h_*)\le3.4(1-r)$ by \eqref{ap:sup}.
Since $4\cdot3.4^3\|h-h_*\|_\infty<0.01$, we have
\[
 Jh^4\le C(r)=\frac{(1-r)^3}{L_0(r)^2}+\frac{1}{100}.
\]

Both charge comparisons use the square completion
\begin{equation}\label{ap:magnetic-square}
 |\xi+\eta|^2\ge b|\xi|^2-\frac{b}{1-b}|\eta|^2,\qquad 0<b<1 .
\end{equation}
With $\xi=dg$, $\eta=im\alpha g$ and $b=(1+4s)^{-1}$ we obtain
\[
 4|(d+im\alpha)g|^2+\frac{m^2}{4}|g|^2\ge4e|dg|^2,\qquad
 e=\frac{(2-r)^2}{4+r^2}\le\frac{1}{1+4s}.
\]
Indeed $s=\frac{X}{(1+c)^2}\le\frac{r}{(2-r)^2}$.
To prove \eqref{ap:diskgap}, we write the trial as $Fg$ and choose
$(F,\delta)=(1,0)$ in the even space and
$(F,\delta)=(1+u-\frac{u^2}{5},\frac{16}{3})$ in the odd space.
Both multipliers have zero normal derivative. In the odd space
\[
 \delta F+4\Delta F=T\left(\frac{40}{3}-\frac{224u}{15}\right),
 \qquad (1+u)^2F\left(\frac{40}{3}-\frac{224u}{15}\right)\le21 ,
\]
where we maximize $(1+u)^3(1-\frac{28u}{25})$ on its positive
set. With $t=\log\rho$ we integrate by parts in $F$ to obtain
\begin{align*}
 \frac{\pi}{4}\big[Q_m(Fg)-(d_m+\delta)\|Fg\|^2\big]
 &\ge \int\big[e(|g_t|^2+|g_\theta|^2)-p_j|g|^2\big]\,dt\,d\theta,\\
 p_0&=\frac{5rC}{4},\qquad p_1=\frac{r}{4}\left[5C(1+r)^2+\frac{21}{(1+r)^2}\right].
\end{align*}
Here the indices $0,1$ correspond to even and odd half-turn parity. Using
$e\ge\max(1-r,\frac{1}{5})$, $\frac{r}{L_0^2}\le\frac{17}{7}$,
$1-r^2\le\frac{3}{5}(2-r)$ and $(1+r)(2-r)\ge2$, we obtain
\[
 \frac{p_0}{e}\le\frac{85}{28}+\frac{1}{16}<4,\qquad
 \frac{p_1}{e}\le\frac{153(1-r)}{28}+\frac{105r}{16}+\frac{1}{4}<7 .
\]
It remains to treat frequency $0$ in the even space and $\pm1$ in the
odd space. We set $D_r=r\partial_r$ and choose
\[
 (W,h,\nu,b)=(p_0,0,25,12)\quad\text{or}\quad
 \left(p_1-e,\frac{1}{2},5,11\right),\qquad
 Y=r^h\frac{(1-\nu r)(1+dr)}{(1+br)^2}.
\]
We choose $d>0$ to cancel the pole of $U=-\frac{4D_r(eD_rY)}{Y}$ at
$r=\frac{1}{\nu}$. Explicitly,
\[
 \frac{d}{\nu+d}=\frac{2b}{\nu+b}-\frac{1}{2}-h
                 -\frac{1}{2}(D_r\log e)(\frac{1}{\nu}).
\]
On $[0,1]$ the coefficient formula \eqref{ap:signrule} for
$\frac{U-W}{r}$ yields $a_0>0$ with the respective lower ratios $44,63$
after rounding down. In both cases
\[
 U>W\quad(0<r<1),\qquad
 \frac{D_rY(1)}{Y(1)}
 =h+\frac{\nu}{\nu-1}+\frac{d}{1+d}-\frac{2b}{1+b}>\frac{1}{50}.
\]
We impose zero trace at $r=\frac{1}{\nu}$ on each remaining Fourier
coefficient. Identity \eqref{ap:groundtransform} with this $Y$ on either
side of the cut proves strict positivity. The outer flux is nonnegative.
At $t=-\infty$ the even coefficient satisfies $|g(t)|\lesssim_g 1+\sqrt{|t|}$
and $\left|\frac{Y_t}{Y}\right|\lesssim r$. The odd coefficient belongs to
$H^1(-\infty,0)$ and tends to zero. Approximation away from the cut
uses its zero trace. The codimensions are one and two, respectively,
which proves \eqref{ap:diskgap}. The reflections $z\mapsto-z$ and
$x\leftrightarrow y$ act by $\zeta\mapsto-\bar\zeta$ and
$\zeta\mapsto\bar\zeta$. For scalar functions odd under the first
and even under the second, the two modes $\pm1$ have one independent
coefficient and require only one trace condition.

For consecutive charges, on the full free space or either half-turn
space we have
\begin{equation}\label{ap:step}
 m\ge0,\quad \lambda_i(B_m)\ge L,\quad L\le m^2+a
 \quad\Longrightarrow\quad \lambda_i(B_{m+1})>L+2m+1-\frac{14}{3}.
\end{equation}
We take $F=\frac{1}{2}\sqrt{3+c}$. We apply \eqref{ap:magnetic-square} with
$\xi=(d+im\alpha)f$, $\eta=i\alpha f$, $b=F^2$, and integrate by parts
in $Q_m(Ff)$. Since $\partial_\nu c=-\frac{1}{\sqrt{2}}$ on the open sides of $D$,
the boundary term $-4F\partial_\nu F$ is nonnegative. The two kinetic
terms reduce to
\[
 \frac{4F^2|\alpha|^2}{1-F^2}-4F\Delta F
 =\frac{3+c}{1+c}+c+\frac{1-c^2}{4(3+c)}
 \le 3+\frac{1}{12}.
\]
Here $\Delta c=-2c$ and $|dc|^2=1-c^2$. The fixed-$c$ bound for $h_*$
and \eqref{ap:verror} imply $V\le\frac{5c^8}{L_0(X)^2}+0.05$. Using
$\frac{(1-F^2)c^8}{X}=\frac{c^8}{4(1+c)}\le\frac{1}{8}$ and
$0\le1-F^2\le\frac{1}{4}$, we estimate
\[
 (1-F^2)V\le\frac{85}{56}+\frac{1}{80}<\frac{19}{12}.
\]
Thus for $f\ne0$,
\[
 Q_{m+1}(f)-Q_m(Ff)>
 \left(2m+1-\frac{14}{3}\right)\|f\|^2+
 (m^2+a)\int_D(1-F^2)|f|^2 .
\]
We require $Ff$ to be orthogonal to the first $i-1$ eigenfunctions
of $B_m$. Since $LF^2+(m^2+a)(1-F^2)\ge L$, the min--max principle
proves \eqref{ap:step}.

\subsection{Polynomial moments for the remaining spectral tests}\label{ap:tests}
In the three scalar rows below, $B=B_m^-$ is the free lune operator in
odd half-turn parity. In the tangent-field row, $B$ is the charge-two
operator on the whole sphere in the real parity class $101$, identified
with tangent fields in Section~\ref{ap:count}. In each row we write
$\widetilde B$ for the same operator with $V_c$ in place of $V$.
By the spectral theorem and \eqref{ap:verror},
\[
 \operatorname{dist}(b,\operatorname{spec} B)\,\|f\|
 \le\|(B-b)f\|\le\|(\widetilde B-b)f\|+\varepsilon\|f\|.
\]
Thus, for $r>\varepsilon$,
\begin{equation}\label{ap:inclusion}
 I=\|(\widetilde B-b)f\|^2-(r-\varepsilon)^2\|f\|^2<0
 \quad\Longrightarrow\quad \operatorname{spec} B\cap(b-r,b+r)\ne\varnothing .
\end{equation}
In the folded measure below we obtain
\begin{equation}\label{ap:lowtable}
\begin{array}{clcr}
 m&\text{space and trial}&(b-r,b+r)&\text{upper bound for }10^6 I\\ \hline
 1&B_1^-,\,f_{12}&(0,\frac{27}{4})&-15205\\
 2&B_2^-,\,f_{12}&(4,9)&-10539\\
 2&\text{type }101\text{ on }\mathbb S^2,\,Y&(0,\frac{15}{4})&-198\\
 3&B_3^-,\,f_3&(0,4)&-173
\end{array}
\end{equation}

For the scalar trials we set $\phi=\arg(x+iy)-\frac{\pi}{4}$ and
$\xi=T\sin(2\phi)$, so $v=\xi^2$. In the gauge with phase
$\exp(im[\arg Z_2+\frac{1+z}{2}\phi])$ the connection is
$\frac{\phi}{2}\,dz$. We use
\[
 f=\frac{zP_1+i\xi P_2}{K^{\frac{3}{2}}},\qquad
 \widetilde B_m f=-4\left(\partial_z+\frac{im\phi}{2}\right)
 \left[T\left(\partial_z+\frac{im\phi}{2}\right)f\right]
 -\frac{4}{T}\partial_\phi^2 f+(m^2+a-V_c)f .
\]
The gauge change is even under the half-turn since both connections are
invariant and the half-turn fixes the peak. The amplitude is odd under
$(z,\phi)\mapsto(-z,-\phi)$ and $\partial_\phi f=0$ on the sides.
At geodesic distance $\varrho$ from either pole, $T=O(\varrho^2)$.
In a smooth local bundle gauge the nonconstant angular terms of the
trial and the phase correction are $O(\varrho^2)$, with first derivatives
$O(\varrho)$ and second derivatives $O(1)$ on the lune. Indeed
$\phi$ ranges over a fixed angular interval, and at each pole the
angle-dependent part of the gauge phase is a smooth multiple of
$(1-|z|)\phi$. The corner contributes zero boundary flux, and the
trial belongs to the free operator domain. For the tangent trial we set $\xi=x^2-y^2$,
$p=(x,y,z)$ and $\Pi=I-pp^{\mathsf T}$. Differentiating the projection, we obtain
\[
 C_1=\frac{(zP_1,0,xP_2)}{K^{\frac{3}{2}}},\qquad
 Y=\Pi C_1,\qquad
 \widetilde B Y=\Pi(-4\Delta+a-V_c)C_1+8\nabla(p\cdot C_1).
\]
We take the following polynomials.
\begin{center}\small
\begin{tabular}{crll}\toprule
 & $k$ & $kP_1$ & $kP_2$\\\midrule
$f_{12}$&5&$5+3u-2v-4uv$&$-5-8u-6u^2$\\
$f_3$&40&$40+22u-3v+v^2$&$14+8u-17u^2+v^2$\\
$Y$&100&$100+48\xi+5\xi^2+8\xi^3+4\xi^4$&$2-28\xi+4\xi^2-5\xi^3-2\xi^5$\\
&&$\quad+u(49+35\xi)-24u^2$&$\quad+u(63-25\xi)+29u^2$\\
\bottomrule\end{tabular}\end{center}

Let $\Phi$ denote integration over $[0,1]\times[0,\frac{\pi}{4}]$ in
$(z,\phi)$. For each row of \eqref{ap:lowtable} we apply \eqref{ap:recip} in the form
\[
 I\le\Phi[P_7(K)(B-C)+E_7(K)B],\qquad
 (B,C)=K^7\big(|(\widetilde B-b)f|^2,(r-\varepsilon)^2|f|^2\big).
\]
We use $f=Y$ in the tangent row. In each case $B,C$ are polynomials.
To evaluate all angular powers at once, we introduce $d=\frac{\pi}{2}$ and
$F_j(t)=\int_0^d e^{t\theta}\sin^j\theta\,d\theta$. We integrate twice by parts to obtain
\begin{equation}\label{ap:fold}
 \begin{aligned}
 (t^2+j^2)F_j(t)&=j(j-1)F_{j-2}(t)+te^{dt}+\mathbf1_{\{j=1\}}
 \quad(j\ge1),\\
 \Phi(u^i\xi^j\phi^k)&=
 \frac{2^{j-k-1}j!(2i-1)!!}{(2i+2j+1)!!}\,
 k![t^k]F_j(t),\qquad F_0(t)=\frac{e^{dt}-1}{t}.
 \end{aligned}
\end{equation}
We omit the $F_{j-2}$ term at $j=1$. The second identity follows by
integrating $z^{2i}(1-z^2)^j$ in $z$ and taking $\theta=2\phi$.
For $Y$ we write the longitude as $2\phi$, so $\xi=T\cos(4\phi)$.
Only even powers of $\xi$ survive, and they have the same $k=0$ moments.
Since $k\le4$, the upper bound is a polynomial $\sum_{n=0}^5c_nd^n$.
The alternating series in $d=8\arctan(\frac{1}{5})-2\arctan(\frac{1}{239})$
bounds $d$ between $d_-=\frac{314159}{200000}$ and $d_+=\frac{314160}{200000}$.
The four integers in \eqref{ap:lowtable} are
$\lceil 10^6\sum_n\max(c_nd_-^n,c_nd_+^n)\rceil$.

\paragraph{The cap trial.}
We set $p=\frac{s}{1+s}$, $e=1-p$ and $\xi=1-7p$. We take
\[
 G=\frac{e}{14}\big(\xi^2+4\xi^6+7\xi^{15}+2\xi^{38}\big)
 \quad (0\le p\le\frac{1}{7}),\qquad G=0\quad(p>\frac{1}{7}).
\]
On the support, $s\le\frac{1}{6}<\tan^2(\frac{\pi}{8})$, so the cap lies inside $D$. Both $G$ and $G'$ vanish on the cut, so $G$ belongs to the
operator domain and extends in every character. We put $t=\sin^2\omega$.
In these coordinates
\[
 X=4pe,\quad u=Xt,\quad v=4X(1-X)(1-t),\quad w=1-2X+u,\qquad
 \mathbb E(\xi^it^j)=\frac{\binom{2j}{j}}{(i+1)4^j},
\]
where $\mathbb E$ denotes the cap average. The variable $\xi$ is uniform
on $[0,1]$. For $M=K^2G$ the radial operator satisfies
\[
 D=K^2\widetilde B_{15}G
   =K^2\left[-4\partial_p(pe\,\partial_pG)
    +\left(\frac{225}{e}+a\right)G\right]-SG .
\]
We compute the first two spectral moments by \eqref{ap:recip} with $k=4$ and
$4DM=(D+M)^2-(D-M)^2$. Inserting the cap moments, we obtain
\[
 (\widetilde\mu,\widetilde\nu)
 =\frac{1}{\mathbb E(G^2)}\mathbb E\!\left[\frac{(DM,D^2)}{K^4}\right],
 \qquad
 \frac{673795}{10^6}\le\widetilde\mu\le\frac{677832}{10^6},
 \quad \widetilde\nu\le\frac{83766160}{10^6}.
\]
It follows from \eqref{ap:verror}, with $g=\frac{G}{\|G\|}$, that
\begin{equation}\label{ap:cap}
 0.59<\mu=\langle B_{15}g,g\rangle<0.76,\qquad
 \nu=\|B_{15}g\|^2<(\sqrt{84}+\varepsilon)^2<86.
\end{equation}
Since $G$ is real and radial, we also have
\[
 Q_m(G)=Q_{15}(G)-(225-m^2)\int\frac{G^2}{e}<0,\qquad 0\le m\le14 .
\]

\paragraph{A negative plane in charge one.}
In character $(1,0)$ we take the real amplitudes
$f_0=|H|$ and $f_1=\sqrt T\,A$, $A=z(1+u)K^{-\frac{3}{2}}$, in
the gauge $e^{i\arg Z_2}$. Here the magnetic and fiber energies on a
real amplitude have density $\frac{2}{1-z}|f|^2$. Reflection cancels
the imaginary cross term. At fixed longitude $H^2$ is a nonnegative
sum of powers $T^k$, $k\ge2$. Since
$\int T^{k-1}\,dz=(1+\frac{1}{2k})\int T^k\,dz$ and
$\int H^6\ge(24+a)\int H^2$, the form entries satisfy
\[
 Q_{00}=-4\int H^6+2\int\frac{H^2}{T}<-96\|H\|^2,\qquad
 Q_{01}=2\int\frac{z|H|f_1}{T},\qquad
 \frac{|Q_{01}|^2}{-Q_{00}}\le\frac{5}{96}\int uA^2 .
\]
We estimate the Schur complement directly in the quadratic form.
The even part of $\frac{2}{1-z}$ is $\frac{2}{T}$. Hence
\begin{align*}
 Q_{11}+\frac{|Q_{01}|^2}{-Q_{00}}
 &\le\mathcal M\!\left(\frac{B-C}{K^5}\right),\qquad C=K^5 V_c f_1^2,\\
 B&=K^5\left[4|\nabla f_1|^2+
 \left(\frac{2}{T}+\frac{5u}{96T}+a+\varepsilon\right)f_1^2\right].
\end{align*}
Both $B,C$ are nonnegative polynomials in $u,v$, since
\[
 |\nabla f_1|^2=T|\partial_z f_1|^2+\frac{1}{T}|\partial_\phi f_1|^2,
 \qquad \partial_z v=-\frac{4zv}{T},\qquad (\partial_\phi v)^2=16vw^2 .
\] By \eqref{ap:recip} and \eqref{ap:mom},
\[
 -324<10^6\mathcal M[P_5(K)(B-C)+E_5(K)B]<-323 .
\]
Both Schur pivots are negative, so $f_0,f_1$ span a negative plane.
\space
\subsection{The zero modes in the charge decomposition}\label{ap:count}
We prove Theorem~\ref{cert:thm:executed}. Differentiation of the conformal
action embeds the orbit tangent space in $\ker A$. Its dimension is nine
by Proposition~\ref{geo:prop:stabilizer}. It thus suffices to bound the
dimension of the kernel from above.

\paragraph{Charges $|m|\ge3$.}
The second row of \eqref{ap:lowtable} locates an eigenvalue of $B_2^-$ in $(4,9)$.
By \eqref{ap:diskgap} and \eqref{ap:step},
\[
 \lambda_3(B_2^-)>\frac{109}{12}>9,\qquad \lambda_2(B_2^-)>4,\qquad
 \lambda_2(B_3^-)>4+5-\frac{14}{3}>4.
\]
The eigenvalue of $B_3^-$ in $(0,4)$ is thus its first.
By \eqref{ap:diskgap} the even block has at most one nonpositive direction,
so $\lambda_2(B_3)>0$. We iterate
\eqref{ap:step} to obtain
\[
 \lambda_2(B_m)>m^2-9-\frac{14}{3}(m-3)\quad(m\ge3),\qquad
 \lambda_2(B_{15})>160 .
\]
For the cap trial in \eqref{ap:cap},
\[
 \nu-160\mu<86-160\cdot0.59<0 .
\]
The spectral measure therefore meets $(0,160)$. As $\lambda_2(B_{15})>160$,
we have $\lambda_1(B_{15})>0$.
For $m\ge15$ positivity follows by \eqref{ap:step}. For $3\le m\le14$,
the cap trial and \eqref{ap:free} imply $\lambda_1(A_{m,j})<0<
\lambda_2(A_{m,j})$ for every $j$. Complex conjugation identifies the charges $m$ and $-m$.

\paragraph{Charge zero.}
We write $L=P-5H^4$. The subscripts of $L$ record odd ($1$) or even ($0$)
parity in $x,y,z$. When the block is invariant under $x\leftrightarrow y$,
we also use its exchange sign $\pm$. The three rotations
\[
 zH_x\ (011),\qquad zH_y\ (101),\qquad
 J=xH_y-yH_x=2(x^2-y^2)H_w\ (000,-)
\]
solve $Lf=0$. By Section~\ref{ap:profile}, each is positive on its
fundamental region, the northern octant for the first two and its half
$x>y$ for the third. Each is therefore a simple ground state.

Let $n_-,n_0$ count negative and zero eigenvalues, with multiplicity.
The remaining scalar information used below is
\[
\begin{array}{c|l}
(n_-,n_0)&\text{scalar blocks}\\ \hline
(1,0)&000,+\quad 001,+\quad 110,+\quad 100\quad 010\\
(0,0)&001,-\quad 110,-\quad 111 .
\end{array}
\]
Indeed, the trial $|H|$ has energy $-4\int H^6$ in $000,+$, whose
second eigenvalue is positive by \eqref{ap:diskgap}. Thus
$\lambda_2(L_{000})=0$. Changing a Neumann side to a Dirichlet side
strictly raises every eigenvalue, since equality in min--max would
impose both Cauchy data on an open side of a nonzero eigenfunction.
Adding sides to $000$ thus raises the second eigenvalue in $100,010,110$,
while $h$ still tests negative in each. Adding sides to $000,-$ proves
positivity of $001,-$ and $110,-$, and adding one to $011$ proves
positivity of $111$. Removing the $y$ side condition from $011$ yields
$\lambda_1(L_{001})<0$. The one-trace bound following \eqref{ap:diskgap}
ensures $\lambda_2(L_{001,+})>0$. This proves the table and
$\dim_{\R}\ker L=3$.

\paragraph{Charge one.}
Since $d_1+\frac{16}{3}>\frac{27}{4}$, the first row of
\eqref{ap:lowtable} and \eqref{ap:diskgap} imply
$\lambda_2(B_1^\pm)>0$, hence $\lambda_3(A_{1,j})>0$.
In Section~\ref{ap:tests} we constructed a negative plane for $j=0$
and a negative trial for every $j$. Exchanging $Z_1,Z_2$ sends $j$ to $1-j$.
Consequently
\[
 \ker A_{1,0}=\ker A_{1,1}=0,\qquad
 \dim_{\C}\ker A_{1,j}\le1\quad(j=2,3).
\]

\paragraph{Charge two.}
The real charge pair identifies with tangent fields $Y\perp p=(x,y,z)$
on the base. Indeed the horizontal frame $(-\bar Z_2,\bar Z_1)$ and
its $i$ multiple turn through $-2t$ relative to transport by
$Z\mapsto e^{it}Z$. The induced connection and the Levi--Civita
connection have the same curvature $2\,d\alpha$ up to orientation,
and $\Sph^2$ is simply connected. The extrinsic derivative identity
then identifies the form,
\[
 \sum_i |dY_i|^2=|\nabla^T Y|^2+|Y|^2,\qquad
 \mathfrak q(Y)=\sum_i\langle LY_i,Y_i\rangle,\qquad p\cdot Y=0 .
\]
For vector type $b$, component $i$ has type $b\oplus e_i$, with
addition modulo two. The parallel isometry $*Y=p\times Y$ preserves
the form. The quarter turn $R(x,y,z)=(-y,x,z)$ is interchange followed
by $x$-reflection. Thus
\[
 *:\ (b,\epsilon)\longmapsto(b\oplus111,-\epsilon),\qquad
 R|_{(b,\epsilon)}=(-1)^{b_1}\epsilon I\quad(b_1=b_2).
\]

We parallel transport $e_z=(0,0,1)$ from the peak along meridians to
obtain the radial frame vector $E$. The cap $GE$ has real amplitude
$G$ and energy $Q_2(G)<0$ on $D$. Equatorial reflection reverses $E$
and interchange fixes it. Since $\operatorname{supp}G\Subset D$,
we may reflect $GE$ across $x=0,y=0$ with either sign on each side.
We take the signs of $111,+$, $001,+$ and $101$ to obtain a negative trial in each.

Restricting to the lune, we claim
\[
 \min\big\{\lambda_3^{\C}(R=-I),\,\lambda_3^{\R}(101)\big\}
 \ge\lambda_3(B_2)>\frac{15}{4}.
\]
The last inequality follows from \eqref{ap:diskgap} and
$\lambda_2(B_2^-)>4$. In both spaces restriction to $D$ preserves
Rayleigh quotients by symmetry. On $R=-I$ it is complex linear for
the complex structure $*$. In type $101$ the equatorial reflection is
antiunitary, since it anticommutes with $*$. The Hermitian norm and
form are real on its fixed-sign real subspace, so complexification
preserves the dimension and upper Rayleigh bound of every trial space.
The negative $111,+$ cap lies in $R=-I$. That whole sector therefore
contributes at most one complex dimension. In $101$ the cap tests
negative and \eqref{ap:lowtable} places an eigenvalue in
$(0,\frac{15}{4})$. Hence this block and its images $011,010,100$
under interchange and $*$ have zero nullity.

For $R=I$ it suffices by $*$ to consider $111,-$ and $001,+$. Before
imposing $p\cdot Y=0$, the direct sum $L^\oplus$ of component forms
satisfies
\[
\begin{array}{c|c|c}
 (b,\epsilon)&n_-(L^\oplus)&Y\in\ker L^\oplus\quad(A\in\R)\\ \hline
 111,- &0&(AzH_x,-AzH_y,0)\\
 001,+ &1&(AzH_y, AzH_x,0).
\end{array}
\]
Tangency in the second row forces $A=0$ by positivity. In the first
it requires $A(x^2-y^2)H_T=0$. If $H=F(w)$, then
$\Delta H=4(T-w^2)F''-6wF'$. Varying $T$ at fixed $w$ forces
$F''=0$, contrary to $PH=H^5$. Thus $\ker L^\oplus$ contains no
nonzero tangent field in either row. Projection onto the negative
space of $L^\oplus$ is injective on each nonpositive tangent subspace,
since a vector in the kernel of the projection must have zero component
energy. The $001,+$ cap already accounts for the sole negative
direction allowed by the table. Neither row contributes to the kernel.

Combining charge zero, the two conjugate charge-one sectors and the
real charge-two pair, we obtain
\[
 9=\dim T_q(\operatorname{Conf}(\Sph^3)\cdot q)
 \le\dim_{\R}\ker A\le3+2\cdot2+2=9 .
\]
The kernel is thus the orbit tangent space. Formula \eqref{common:equivariant}
identifies its $\tau$-odd part with the image of the conformal generators
commuting with conjugation. The real and imaginary coordinate planes
have dimension two, and conjugation reverses the Hopf generator.
Thus the dimension is $\binom{2}{2}+\binom{2}{2}+2=4$, proving \eqref{ap:kernel}.\qed
\space
\space
\section{Notation for the oscillating solution}\label{notation:start}
In Sections~\ref{main:start}--\ref{end:start}, $E(f,g)$ denotes the
conserved energy of a pair, whereas $E(f)$ denotes the gradient norm
of a wave position, with the time supremum in \eqref{lin:eq:observations}.
The stationary profile $h$ on $\Sph^2$ and the dipole response $h(t)$
appear in \eqref{common:lune-pde} and \eqref{non:eq:halo}, respectively.

\small
\renewcommand{\arraystretch}{1.2}
\begin{longtable}{@{}L{2.7cm}L{8.5cm}L{3.5cm}@{}}
\toprule
Notation & Definition or role & Location\\
\midrule
\endfirsthead
\toprule Notation & Definition or role & Location\\\midrule
\endhead
\bottomrule\endlastfoot
$\mathcal H,\mathrm{Sim}$ & Energy space $\dot H^1\times L^2$ and spatial similarities acting isometrically on energy pairs. & \eqref{intro:wave}, \eqref{main:eq:similarity}\\
$q_\xi,Q,K,R$ & $q_\xi=Q+K\xi+R(\xi)$, $K^*K=I$, $K^*R=0$. The columns of $K$ span the odd energy kernel. & \eqref{geo:eq:chart}\\
$D,G_0,L,V$ & $D=\sqrt{-\Delta}$, $G_0=D^{-2}$, $L=D^2-V$, $V=5Q^4$. & Section~\ref{lin:green-section}\\
$\rho,X_s,Y_s$ & $\rho=\langle x\rangle$, $X_s=L^2(\rho^{-2s}dx)$, $Y_s=L^2(\rho^{2s}dx)$. & \eqref{lin:eq:kernelbounds}, \eqref{lin:eq:mixed-newton}\\
$\Phi,P,P_0$ & Schwartz tests of $L$-spectrum below one, $\Phi K=I$, $\Phi\Psi_l=0$, $P=I-K\Phi$, $P_0=I-KK^*$. & \eqref{lin:eq:spectralgauge}\\
$E,S,\LE,\OO$ & Gradient energy, critical Strichartz, weighted position and Newton velocity norms. & \eqref{lin:eq:observations}\\
$\mathcal R(I)$ & Combined energy, Strichartz and weighted observation norm on $I$. & \eqref{lin:eq:green-norm}\\
$\varepsilon,\eta_j,T_j$ & $\varepsilon=d^3$, $\eta_j=\varepsilon2^{-j}$, $T_j=4^jd^{-5}$, $T_j\eta_j^2=d$. & \eqref{main:scales}\\
$a_j,a_j^*,v_a,\xi$ & Exact nodal impulses and normalized slope jumps, radiation train with $\Phi v_a=0$, and stationary chart path. & \eqref{main:eq:ideal-controls}, \eqref{cat:eq:train}\\
$\widetilde v,p,M_\xi$ & Corrected velocity $\bar v+P\Phi^*p$, feedback $p=-DR^*\widetilde v$, and $M_\xi=I+\Phi DR$. & Section~\ref{non:sec:field}\\
$H,F(v_a)$ & Modulation force and its leading quadratic radiation term. & \eqref{cat:eq:forcecomparison}\\
$h(t),G$ & Explicit dipole response and regular remainder $G=w-h$, with bounded position mass. & \eqref{non:eq:halo}, \eqref{non:eq:dipole-kernel}\\
$\gamma(\xi),r$ & Chart coordinate $\gamma=\Phi(q_\xi-Q)$ and displacement $(-1)^j d r$ on $I_j$. & \eqref{end:ep:coordinate}, \eqref{end:eq:control-cube}\\
$Z$ & Tangent velocity observation $Z=K^*\bar v-\Phi\Phi^*p$. & \eqref{non:eq:tangent-observation}, Lemma~\ref{end:ep:k}\\
$W,\Gamma_{\rm coord}$ & Projected position $W=w+PR(\xi)$ and coordinate $\Gamma_{\rm coord}=\gamma(\xi)$. & \eqref{non:eq:projected}\\
$(w,u_t),v_+$ & Physical remainder and its unique asymptotic free wave. & Section~\ref{end:global-section}\\
\end{longtable}
\normalsize
\renewcommand{\arraystretch}{1}
\space
\section*{Statement on use of AI.}
Our contribution is the use of the signed Hopf profiles in the three-dimensional wave constructions described here. We directed ChatGPT Astra 6.2 to assist in deriving and writing out specific technical estimates, and to carry out the computer-assisted numerical verification in its entirety.
\space


\begin{thebibliography}{99}
\bibitem{Agmon}
S.~Agmon, \emph{Spectral properties of Schr\"odinger operators and scattering theory},
Ann. Scuola Norm. Sup. Pisa Cl. Sci. (4) \textbf{2} (1975), 151--218.
\url{https://www.numdam.org/item/ASNSP_1975_4_2_2_151_0/}.

\bibitem{Aubin}
T.~Aubin, \emph{Probl\`emes isop\'erim\'etriques et espaces de Sobolev},
J. Differential Geom. \textbf{11} (1976), 573--598.
\url{https://doi.org/10.4310/jdg/1214433725}.

\bibitem{BG}
H.~Bahouri and P.~G\'erard,
\emph{High frequency approximation of solutions to critical nonlinear
wave equations}, Amer. J. Math. \textbf{121} (1999), 131--175.
\url{https://doi.org/10.1353/ajm.1999.0001}.

\bibitem{BGP}
C.~B\"ar, N.~Ginoux and F.~Pf\"affle, \emph{Wave Equations on Lorentzian
Manifolds and Quantization}, ESI Lectures in Mathematics and Physics,
European Mathematical Society, 2007, Chapter~2.

\bibitem{Bec}
M.~Beceanu,
\emph{A centre-stable manifold for the energy-critical wave equation in
$\R^3$ in the symmetric setting},
J. Hyperbolic Differ. Equ. \textbf{11} (2014), 437--476.
\url{https://arxiv.org/abs/1212.2285}.

\bibitem{BecGold}
M.~Beceanu and M.~Goldberg, \emph{Strichartz estimates and maximal operators for the wave equation in $\mathbb R^3$}, J. Funct. Anal. \textbf{266} (2014), 1476--1510.
\url{https://arxiv.org/abs/1212.0568}.

\bibitem{CGS}
L.~A.~Caffarelli, B.~Gidas and J.~Spruck,
\emph{Asymptotic symmetry and local behavior of semilinear elliptic equations with critical Sobolev growth},
Comm. Pure Appl. Math. \textbf{42} (1989), 271--297.
\url{https://doi.org/10.1002/cpa.3160420304}.

\bibitem{CHS3DR6}
X.~Chen, J.~Holmer and S.~Shen,
\emph{A counterexample to fixed-profile soliton resolution for the three-dimensional quintic wave equation},
manuscript in preparation, November 2026 (balanced single-node construction).

\bibitem{CHS5D}
X.~Chen, J.~Holmer and S.~Shen,
\emph{Radiation and persistent shape dynamics for the five-dimensional
energy-critical nonlinear wave equation}, arXiv:2611.17243v1 (2026),
Section~10. \url{https://arxiv.org/abs/2611.17243v1}.

\bibitem{CK}
M.~Christ and A.~Kiselev, \emph{Maximal functions associated to filtrations},
J. Funct. Anal. \textbf{179} (2001), no.~2, 409--425.

\bibitem{CP}
M.~Clapp and A.~Pistoia,
\emph{Symmetries, Hopf fibrations and supercritical elliptic problems},
in \emph{Mathematical Congress of the Americas},
Contemp. Math. \textbf{656}, Amer. Math. Soc., 2016, 1--12.
\url{https://arxiv.org/abs/1501.03349}.

\bibitem{EllCPW}
M.~Clapp, A.~Pistoia and T.~Weth,
\emph{An upper bound for the least energy of a nodal solution to the
Yamabe equation on the sphere}, Minimax Theory Appl. \textbf{7} (2022),
no.~2, 173--184. \url{https://arxiv.org/abs/1910.05797}.

\bibitem{EllCdVT}
Y.~Colin de Verdi\`ere and F.~Truc,
\emph{Confining quantum particles with a purely magnetic field},
Ann. Inst. Fourier (Grenoble) \textbf{60} (2010), no.~7, 2333--2356.
\url{https://doi.org/10.5802/aif.2609}.

\bibitem{CollotAnnouncement}
C.~Collot,
\emph{Singularities of the critical wave equation and Keller--Segel system formed by the collision of collapsing solitons at a point},
Princeton Analysis Seminar, November 30, 2026. Seminar announcement describing ongoing work with Li on the four-dimensional cubic wave equation.
\href{https://www.math.princeton.edu/events/singularities-critical-wave-equation-and-keller-segel-system-formed-collision-collapsing}{Princeton seminar abstract}.

\bibitem{CDKM}
C.~Collot, T.~Duyckaerts, C.~Kenig and F.~Merle,
\emph{Non-radial soliton resolution for the energy critical wave equation
for initial data with compact support},
Discrete Contin. Dyn. Syst. \textbf{56} (2026), 501--556.
\url{https://doi.org/10.3934/dcds.2026115}.

\bibitem{DPMPP}
M.~del Pino, M.~Musso, F.~Pacard and A.~Pistoia,
\emph{Large energy entire solutions for the Yamabe equation},
J. Differential Equations \textbf{251} (2011), 2568--2597.
\url{https://doi.org/10.1016/j.jde.2011.03.008}.

\bibitem{HeatTower}
M.~del Pino, M.~Musso and J.~Wei, \emph{Existence and stability of infinite time bubble towers in the energy critical heat equation}, Anal. PDE \textbf{14} (2021), 1557--1598.
\url{https://doi.org/10.2140/apde.2021.14.1557}.

\bibitem{Ding}
W.~Y.~Ding, \emph{On a conformally invariant elliptic equation on $\R^n$},
Comm. Math. Phys. \textbf{107} (1986), 331--335.
\url{https://doi.org/10.1007/BF01209398}.

\bibitem{DHKS}
R.~Donninger, M.~Huang, J.~Krieger and W.~Schlag,
\emph{Exotic blowup solutions for the $u^5$ focusing wave equation in $\mathbb R^3$},
Michigan Math. J. \textbf{63} (2014), 451--501.
\url{https://doi.org/10.1307/mmj/1409932630}.

\bibitem{DK}
R.~Donninger and J.~Krieger,
\emph{Nonscattering solutions and blowup at infinity for the critical
wave equation}, Math. Ann. \textbf{357} (2013), 89--163.
\url{https://arxiv.org/abs/1201.3258}.

\bibitem{DJKM}
T.~Duyckaerts, H.~Jia, C.~Kenig and F.~Merle,
\emph{Soliton resolution along a sequence of times for the focusing
energy critical wave equation}, Geom. Funct. Anal. \textbf{27} (2017),
no.~4, 798--862. \url{https://arxiv.org/abs/1601.01871}.

\bibitem{DKMradial}
T.~Duyckaerts, C.~Kenig and F.~Merle,
\emph{Classification of radial solutions of the focusing, energy-critical
wave equation}, Camb. J. Math. \textbf{1} (2013), no.~1, 75--144.
\url{https://arxiv.org/abs/1204.0031}.

\bibitem{DKMprofiles}
T.~Duyckaerts, C.~Kenig and F.~Merle,
\emph{Profiles for bounded solutions of dispersive equations, with applications
to energy-critical wave and Schr\"odinger equations},
Commun. Pure Appl. Anal. \textbf{14} (2015), no.~4, 1275--1326.
\url{https://arxiv.org/abs/1311.0665}.

\bibitem{DKMcompact}
T.~Duyckaerts, C.~Kenig and F.~Merle,
\emph{Solutions of the focusing nonradial critical wave equation with the
compactness property}, Ann. Sc. Norm. Super. Pisa Cl. Sci. (5)
\textbf{15} (2016), 731--808.
\url{https://arxiv.org/abs/1402.0365}.

\bibitem{DKM}
T.~Duyckaerts, C.~Kenig and F.~Merle,
\emph{Concentration-compactness and universal profiles for the non-radial
energy critical wave equation},
Nonlinear Anal. \textbf{138} (2016), 44--82.
Conjecture numbering refers to \url{https://arxiv.org/abs/1510.01750}.

\bibitem{DKModd}
T.~Duyckaerts, C.~Kenig and F.~Merle,
\emph{Soliton resolution for the radial critical wave equation in all odd
space dimensions}, Acta Math. \textbf{230} (2023), 1--92.
\url{https://arxiv.org/abs/1912.07664}.

\bibitem{EM}
M.~Engelstein and D.~Mendelson,
\emph{Non-uniqueness of bubbling for wave maps},
Ars Inven. Anal. (2022), Paper No.~5, 30 pp.
\url{https://arxiv.org/abs/2005.14128}.

\bibitem{FPsingle}
J.~C.~Fern\'andez and J.~Petean,
\emph{Low energy nodal solutions to the Yamabe equation},
J. Differential Equations \textbf{268} (2020), 6576--6597.
\url{https://doi.org/10.1016/j.jde.2019.11.043}.

\bibitem{FPP}
J.~C.~Fern\'andez, O.~Palmas and J.~Petean,
\emph{Supercritical elliptic problems on the round sphere and nodal
solutions to the Yamabe problem in projective spaces},
Discrete Contin. Dyn. Syst. \textbf{40} (2020), 2495--2514.
\url{https://doi.org/10.3934/dcds.2020123}.

\bibitem{Fuglede}
B.~Fuglede, \emph{An integral formula}, Math. Scand. \textbf{6} (1958), 207--212.
\url{https://www.mscand.dk/article/view/10545}.

\bibitem{GS}
M.~Goldberg and W.~Schlag,
\emph{A limiting absorption principle for the three-dimensional
Schr\"odinger equation with $L^p$ potentials},
Int. Math. Res. Not. \textbf{2004} (2004), no.~75, 4049--4071.
\url{https://doi.org/10.1155/S1073792804140324}.

\bibitem{GNT}
S.~Gustafson, K.~Nakanishi and T.-P.~Tsai,
\emph{Asymptotic stability, concentration, and oscillation in harmonic map
heat-flow, Landau--Lifshitz, and Schr\"odinger maps on $\R^2$},
Comm. Math. Phys. \textbf{300} (2010), 205--242.
\url{https://arxiv.org/abs/0904.0461}.

\bibitem{Helgason}
S.~Helgason, \emph{The Radon Transform}, second edition, Birkh\"auser, Boston, 1999.
\url{https://doi.org/10.1007/978-1-4757-1463-0}.

\bibitem{IN}
G.~Ivanov and M.~Nasz\'odi,
\emph{Quantitative Steinitz theorem: a polynomial bound},
Bull. Lond. Math. Soc. \textbf{56} (2024), no.~2, 796--802.
\url{https://doi.org/10.1112/blms.12965}.

\bibitem{J5}
J.~Jendrej,
\emph{Construction of type II blow-up solutions for the energy-critical
wave equation in dimension 5},
J. Funct. Anal. \textbf{272} (2017), 866--917.
\url{https://arxiv.org/abs/1503.05024}.

\bibitem{JOpposite}
J.~Jendrej, \emph{Nonexistence of two-bubbles with opposite signs for the radial energy-critical wave equation},
Ann. Sc. Norm. Super. Pisa Cl. Sci. (5) \textbf{18} (2018), 735--778.
\url{https://arxiv.org/abs/1510.03965}.

\bibitem{J}
J.~Jendrej, \emph{Construction of two-bubble solutions for energy-critical wave equations},
Amer. J. Math. \textbf{141} (2019), 55--118.
\url{https://arxiv.org/abs/1602.06524}.

\bibitem{FTJendrejKrieger}
J. Jendrej and J. Krieger, \emph{Concentric bubbles concentrating in finite time for the energy critical Wave Maps equation}, \href{https://arxiv.org/abs/2501.08396}{arXiv:2501.08396}.

\bibitem{JL}
J.~Jendrej and A.~Lawrie,
\emph{Soliton resolution for the energy-critical nonlinear wave equation
in the radial case}, Ann. PDE \textbf{9} (2023), article 18.
\url{https://arxiv.org/abs/2203.09614}.

\bibitem{JLR}
J.~Jendrej, A.~Lawrie and C.~Rodriguez,
\emph{Dynamics of bubbling wave maps with prescribed radiation},
Ann. Sci. \'Ec. Norm. Sup\'er. (4) \textbf{55} (2022), 1135--1198.
\url{https://arxiv.org/abs/1908.08512}.

\bibitem{JLS}
J.~Jendrej, A.~Lawrie and W.~Schlag,
\emph{Continuous in time bubble decomposition for the harmonic map heat flow},
Forum Math. Pi \textbf{13} (2025), e4.
\url{https://doi.org/10.1017/fmp.2024.15}.

\bibitem{JendrejMartel}
J.~Jendrej and Y.~Martel,
\emph{Construction of multi-bubble solutions for the energy-critical wave equation in dimension 5},
J. Math. Pures Appl. \textbf{139} (2020), 317--355.
\url{https://doi.org/10.1016/j.matpur.2020.02.007}.

\bibitem{JYZ}
J.~Jendrej, Y.~Yin and L.~Zhao, \emph{Nonexistence of blow-up solutions with smooth radiation for energy-critical equivariant wave maps}, arXiv:2607.04952 (2026). \url{https://arxiv.org/abs/2607.04952}.

\bibitem{JZZ}
J.~Jendrej, C.~Zhang and L.~Zhao,
\emph{Rigidity of the multi-bubble solutions to the energy critical wave equation in dimension five},
arXiv:2605.27994 (2026). \url{https://arxiv.org/abs/2605.27994}.

\bibitem{JK}
A.~Jensen and T.~Kato,
\emph{Spectral properties of Schr\"odinger operators and time-decay of the
wave functions}, Duke Math. J. \textbf{46} (1979), 583--611.
See also the historical discussion in \cite{JN}.

\bibitem{JN}
A.~Jensen and G.~Nenciu,
\emph{A unified approach to resolvent expansions at thresholds},
Rev. Math. Phys. \textbf{13} (2001), 717--754.
\url{https://doi.org/10.1142/S0129055X01000843}.

\bibitem{JLSX}
H.~Jia, B.~Liu, W.~Schlag and G.~Xu, \emph{Global center stable manifold for the defocusing energy critical wave equation with potential}, Amer. J. Math. \textbf{142} (2020), 1497--1557.
\url{https://doi.org/10.1353/ajm.2020.0038}.

\bibitem{KadarGlobal}
I.~Kadar, \emph{Construction of multi-soliton solutions for the energy critical wave equation in dimension 3},
arXiv:2409.05267 (2024). \url{https://arxiv.org/abs/2409.05267}.

\bibitem{FTKadar}
I. Kadar, \emph{Smooth finite time singularity formation without quantization}, \href{https://arxiv.org/abs/2603.14985}{arXiv:2603.14985}.

\bibitem{KT}
M.~Keel and T.~Tao, \emph{Endpoint Strichartz estimates},
Amer. J. Math. \textbf{120} (1998), no.~5, 955--980.

\bibitem{KM}
C.~E.~Kenig and F.~Merle,
\emph{Global well-posedness, scattering and blow-up for the energy-critical
focusing non-linear wave equation},
Acta Math. \textbf{201} (2008), no.~2, 147--212.
\url{https://arxiv.org/abs/math/0610801}.

\bibitem{KNS}
J.~Krieger, K.~Nakanishi and W.~Schlag,
\emph{Center-stable manifold of the ground state in the energy space for
the critical wave equation}, Math. Ann. \textbf{361} (2015), 1--50.
\url{https://arxiv.org/abs/1303.2400}.

\bibitem{FTKriegerPalacios}
J. Krieger and J. M. Palacios, \emph{Long finite time bubble trees for two co-rotational wave maps}, \href{https://arxiv.org/abs/2602.22825}{arXiv:2602.22825}.

\bibitem{KS}
J.~Krieger and W.~Schlag,
\emph{On the focusing critical semi-linear wave equation},
Amer. J. Math. \textbf{129} (2007), 843--913.
\url{https://arxiv.org/abs/math/0508574}.

\bibitem{KSFull}
J.~Krieger and W.~Schlag, \emph{Full range of blow up exponents for the quintic wave equation in three dimensions},
J. Math. Pures Appl. \textbf{101} (2014), 873--900.
\url{https://arxiv.org/abs/1212.3795}.

\bibitem{KST}
J.~Krieger, W.~Schlag and D.~Tataru,
\emph{Slow blow-up solutions for the $H^1(\R^3)$ critical focusing
semilinear wave equation},
Duke Math. J. \textbf{147} (2009), 1--53.
\url{https://arxiv.org/abs/math/0702033}.

\bibitem{MartelBack}
Y.~Martel, \emph{Asymptotic $N$-soliton-like solutions of the subcritical and critical generalized Korteweg--de Vries equations},
Amer. J. Math. \textbf{127} (2005), 1103--1140.
\url{https://doi.org/10.1353/ajm.2005.0033}.

\bibitem{MartelRaphael}
Y.~Martel and P.~Rapha\"el, \emph{Strongly interacting blow up bubbles for the mass critical nonlinear Schr\"odinger equation},
Ann. Sci. \'Ec. Norm. Sup\'er. (4) \textbf{51} (2018), 701--737.
\url{https://doi.org/10.24033/asens.2364}.

\bibitem{EllMMW}
M.~Medina, M.~Musso and J.~Wei,
\emph{Desingularization of Clifford torus and nonradial solutions to
the Yamabe problem with maximal rank}, J. Funct. Anal. \textbf{276}
(2019), no.~8, 2470--2523. \url{https://doi.org/10.1016/j.jfa.2019.02.001}.

\bibitem{MerleBack}
F.~Merle, \emph{Construction of solutions with exactly $k$ blow-up points for the Schr\"odinger equation with critical nonlinearity},
Comm. Math. Phys. \textbf{129} (1990), 223--240.
\url{https://doi.org/10.1007/BF02096981}.

\bibitem{MW}
M.~Musso and J.~Wei,
\emph{Nondegeneracy of nodal solutions to the critical Yamabe problem},
Comm. Math. Phys. \textbf{340} (2015), 1049--1107.
\url{https://doi.org/10.1007/s00220-015-2462-1}. \href{https://purehost.bath.ac.uk/ws/files/173510250/MussoWeiCMPRevised.pdf}{Accepted author manuscript}.

\bibitem{NPsingle}
S.~L.~N.~Neves and A.~Pistoia,
\emph{Nondegeneracy of low-energy nodal solutions to the Yamabe equation},
arXiv:2610.40960 (2026). \url{https://arxiv.org/abs/2610.40960}.

\bibitem{PillaiQuintic}
M.~Pillai, \emph{Global, non-scattering solutions to the quintic, focusing semilinear wave equation on $\mathbb R^{1+3}$},
arXiv:2108.12024 (2021). \url{https://arxiv.org/abs/2108.12024}.

\bibitem{RaphaelSzeftel}
P.~Rapha\"el and J.~Szeftel, \emph{Existence and uniqueness of minimal blow-up solutions to an inhomogeneous mass critical NLS},
J. Amer. Math. Soc. \textbf{24} (2011), 471--546.
\url{https://doi.org/10.1090/S0894-0347-2010-00688-1}.

\bibitem{ShenPure}
R.~Shen, \emph{Nonexistence of pure bubble type II blow-up solutions to energy critical wave equation in the 3D radial case}, arXiv:2602.04205 (2026). \url{https://arxiv.org/abs/2602.04205}.

\bibitem{RS}
R.~Shen, \emph{Nonexistence of multi-bubble radial solutions to the 3D energy critical wave equation},
arXiv:2603.21602 (2026). \url{https://arxiv.org/abs/2603.21602}.

\bibitem{Syl}
J.~Sylvester, \emph{An estimate for the free Helmholtz equation that scales},
Inverse Probl. Imaging \textbf{3} (2009), no.~2, 333--351.
\url{https://doi.org/10.3934/ipi.2009.3.333}.
The weighted estimate cited in the text is also displayed as (1.5) in the
author's June 2, 2008 manuscript,
\url{https://sites.math.washington.edu/~sylvest/papers/agest.pdf}.

\bibitem{Talenti}
G.~Talenti, \emph{Best constant in Sobolev inequality},
Ann. Mat. Pura Appl. (4) \textbf{110} (1976), 353--372.
\url{https://doi.org/10.1007/BF02418013}.

\bibitem{Top97}
P.~Topping, \emph{Rigidity in the harmonic map heat flow},
J. Differential Geom. \textbf{45} (1997), 593--610.
\url{https://doi.org/10.4310/jdg/1214459844}.

\bibitem{Top04}
P.~Topping, \emph{Winding behaviour of finite-time singularities of the
harmonic map heat flow}, Math. Z. \textbf{247} (2004), 279--302.
\url{https://doi.org/10.1007/s00209-003-0582-3}.

\bibitem{EllWY}
T.~T.~Wu and C.~N.~Yang,
\emph{Dirac monopole without strings: Monopole harmonics},
Nuclear Phys. B \textbf{107} (1976), no.~3, 365--380.
\url{https://doi.org/10.1016/0550-3213(76)90143-7}.

\bibitem{YuanFive}
X.~Yuan, \emph{Construction of excited multi-solitons for the 5D energy-critical wave equation},
J. Hyperbolic Differ. Equ. \textbf{18} (2021), no.~2, 397--434. \url{https://doi.org/10.1142/S0219891621500120}.

\bibitem{YuanFour}
X.~Yuan, \emph{Construction of excited multi-solitons for the focusing 4D cubic wave equation},
J. Funct. Anal. \textbf{282} (2022), Art.~109336.
\url{https://arxiv.org/abs/2103.04168}.
\space
\end{thebibliography}
\end{document}